\documentclass[twoside]{article}
\usepackage{amsfonts}

\newcommand{\ssection}[1]{\section{\sf{\textbf{#1}}}}
\newcommand{\ssubsection}[1]{\subsection[\sf{#1}]{\sf{\textbf{#1}}}}
\newtheorem{theorem}{Theorem}[section]
\newtheorem{proposition}[theorem]{Proposition} 
\newtheorem{lemma}[theorem]{Lemma} 
\newtheorem{corollary}[theorem]{Corollary} 
\newtheorem{definition}[theorem]{Definition}

\newcommand{\CC}{{\mathbb C}}

\newcommand{\GG}{{\mathbb G}}

\newcommand{\PP}{{\mathbb P}}
\newcommand{\QQ}{{\mathbb Q}}
\newcommand{\RR}{{\mathbb R}}
\newcommand{\TT}{{\mathbb T}}

\newcommand{\ZZ}{{\mathbb Z}}

\newcommand{\cA}{{\mathcal A}}
\newcommand{\cB}{{\mathcal B}}

\newcommand{\cG}{{\mathcal G}}
\newcommand{\cH}{{\mathcal H}}

\newcommand{\cJ}{{\mathcal J}}

\newcommand{\cL}{{\mathcal L}}
\newcommand{\cM}{{\mathcal M}}
\newcommand{\cN}{{\mathcal N}}
\newcommand{\cO}{{\mathcal O}}
\newcommand{\cP}{{\mathcal P}}

\newcommand{\cU}{{\mathcal U}}
\newcommand{\cV}{{\mathcal V}}

\newcommand{\SU}{\mathrm{SU}}
\newcommand{\SO}{\mathrm{SO}}
\newcommand{\GL}{\mathrm{GL}}
\newcommand{\SL}{\mathrm{SL}}
\newcommand{\PSL}{\mathrm{PSL}}
\newcommand{\Sp}{\mathrm{Sp}}
\newcommand{\UU}{\mathrm{U}}
\newcommand{\OO}{\mathrm{O}}

\newcommand{\fa}{{\mathfrak a}}
\newcommand{\fg}{{\mathfrak g}}
\newcommand{\fh}{{\mathfrak h}}
\newcommand{\fn}{{\mathfrak n}}
\newcommand{\ft}{{\mathfrak t}}
\newcommand{\fu}{{\mathfrak u}}

\newcommand \id {\mathrm{id}}
\newcommand \Id {\mathrm{Id}}

\newcommand \curv {\mathrm{curv}}

\newcommand \red {\mathrm{red}}

\newcommand {\ip}[2] {\langle #1, #2 \rangle}

\newcommand{\del}{\partial}
\newcommand{\delbar}{\bar{\partial}}

\begin{document}

\title{\sf{Symplectic Geometry}\\
\vspace*{2ex}
{\normalsize overview written for the
{\em Handbook of Differential Geometry}, vol.~2}\\
{\normalsize (F.J.E.\ Dillen and L.C.A.\ Verstraelen, eds.)}\\
\vspace*{5ex}}

\author{\sf{Ana Cannas da Silva}\thanks{E-mail:
{\tt acannas@math.ist.utl.pt} or
{\tt acannas@math.princeton.edu}}}
\date{\sf{September 2004}}

\maketitle

\def\contentsname{\sf{\textbf{Contents}}}
\addtocontents{toc}{\protect\vspace{2ex}}
\tableofcontents

\newpage

\pagestyle{headings}

\sffamily{

%%%%%%%%%%%%%%%%%%%%%%%%%%%%%%%%%%%%%%%%%%%%%%%%%%%%%%%%%%%%%%%%%%%%%%%%%%%%%
%%%%%%%%%%%%%%%%%%%%%%%%%%%%%%%%%%%%%%%%%%%%%%%%%%%%%%%%%%%%%%%%%%%%%%%%%%%%%
% --> Introduction
%%%%%%%%%%%%%%%%%%%%%%%%%%%%%%%%%%%%%%%%%%%%%%%%%%%%%%%%%%%%%%%%%%%%%%%%%%%%%
%%%%%%%%%%%%%%%%%%%%%%%%%%%%%%%%%%%%%%%%%%%%%%%%%%%%%%%%%%%%%%%%%%%%%%%%%%%%%

\section*{\sf{\textbf{Introduction}}}
\label{introduction}

\markboth{\sf{INTRODUCTION}}{\sf{INTRODUCTION}}

\addcontentsline{toc}{section}{\sf{\textbf{Introduction}}}

\thispagestyle{empty}

This is an overview of symplectic geometry\footnote{The
word {\em symplectic} in mathematics
was coined in the late 1930's by Weyl~\cite[p.165]{we:classical}
who substituted the Latin root in {\em complex} by
the corresponding Greek root in order to label the symplectic group
(first studied be Abel).
An English dictionary is likely to list {\em symplectic}
as the name for a bone in a fish's head.}
-- the geometry of {\em symplectic manifolds}.
From a language for classical mechanics in the XVIII century,
symplectic geometry has matured since the 1960's
to a rich and central branch of differential geometry and topology.
A current survey can thus only aspire
to give a partial flavor on this exciting field.
The following six topics have been chosen for this handbook:

\vspace*{2ex}

\textbf{1. Symplectic manifolds}
are manifolds equipped with {\em symplectic forms}.
A symplectic form is a closed nondegenerate 2-form.
The algebraic condition (nondegeneracy)
says that the top exterior power of a symplectic
form is a volume form, therefore symplectic manifolds
are necessarily even-dimensional and orientable.
The analytical condition (closedness) is a natural differential equation
that forces all symplectic manifolds to being locally indistinguishable:
they all locally look like an even-dimensional euclidean space
equipped with the $\sum dx_i \wedge dy_i$ symplectic form.
All cotangent bundles admit canonical symplectic forms,
a fact relevant for analysis of differential operators,
dynamical systems, classical mechanics, etc.
Basic properties, major classical examples, equivalence notions,
local normal forms of symplectic manifolds and symplectic submanifolds
are discussed in Chapter~\ref{section1}.

\vspace*{2ex}

\textbf{2. Lagrangian submanifolds}\footnote{The
name {\em lagrangian manifold} was introduced by
Maslov~\cite{ma:perturbation} in the 1960's,
followed by {\em lagrangian plane}, etc., introduced
by Arnold~\cite{ar:mathematical}.}
are submanifolds of symplectic manifolds
of half dimension and where the restriction of
the symplectic form vanishes identically.
By the {\em lagrangian creed}~\cite{we:lectures},
everything is a lagrangian submanifold, starting
with closed 1-forms, real functions modulo constants
and symplectomorphisms (diffeomorphisms that respect
the symplectic forms).
Chapter~\ref{section2} also describes normal neighborhoods
of lagrangian submanifolds with applications.

\vspace*{2ex}

\textbf{3. Complex structures} or almost complex structures
abound in symplectic geometry:
any symplectic manifold possesses almost complex
structures, and even so in a {\em compatible} sense.
%establishes a link from symplectic geometry to complex geometry, and
This is the point of departure for the modern technique
of studying pseudoholomorphic curves,
as first proposed by Gromov~\cite{gr:pseudo}\index{Gromov !
pseudo-holomorphic curve}\index{pseudo-holomorphic curve}.
K\"ahler geometry lies at the intersection of
complex, riemannian and symplectic geometries,
and plays a central role in these three fields.
Chapter~\ref{section3} includes the local normal form
for K\"ahler manifolds
and a summary of Hodge theory for K\"ahler manifolds.

\vspace*{2ex}

\textbf{4. Symplectic geography} is concerned with
existence and uniqueness of symplectic forms on a given manifold.
Important results from K\"ahler geometry remain true in the
more general symplectic category, as shown using pseudoholomorphic methods.
This viewpoint was more recently continued
with work on the existence of certain symplectic submanifolds,
in the context of Seiberg-Witten invariants,
and with topological descriptions in terms of Lefschetz pencils.
Both of these directions are particularly relevant to 4-dimensional
topology and to mathematical physics, where symplectic manifolds
occur as building blocks or as key examples.
Chapter~\ref{section4} treats constructions of symplectic manifolds
and invariants to distinguish them. 

\vspace*{2ex}

\textbf{5. Hamiltonian geometry} is the geometry of
symplectic manifolds equipped with a {\em moment map},
that is, with a collection of quantities conserved by symmetries.
With roots in hamiltonian mechanics, moment maps became a
consequential tool in geometry and topology.
The notion of a moment map arises from the fact that,
to any real function on a symplectic manifold,
is associated a vector field whose flow preserves the
symplectic form and the given function;
this is called the {\em hamiltonian vector field}
of that (hamiltonian) function.
The Arnold conjecture in the 60's regarding hamiltonian dynamics was a
major driving force up to the establishment of Floer homology in the 80's.
Chapter~\ref{section5} deals mostly with the geometry of moment maps,
including the classical Legendre transform, integrable systems and convexity.

\vspace*{2ex}

\textbf{6. Symplectic reduction} 
is at the heart of many symplectic arguments.
There are infinite-dimensional analogues with
amazing consequences for differential geometry,
as illustrated in a symplectic approach to Yang-Mills theory.
Symplectic toric manifolds provide examples of extremely
symmetric symplectic manifolds that arise from
symplectic reduction using just the data of a polytope.
All properties of a symplectic toric manifold
may be read from the corresponding polytope.
There are interesting interactions with algebraic geometry,
representation theory and geometric combinatorics.
The variation of reduced spaces is also addressed in Chapter~\ref{section6}.

%%%%%%%%%%%%%%%%%%%%%%%%%%%%%%%%%%%%%%%%%%%%%%%%%%%%%%%%%%%%%%%%%%%%%%%%%%%%%
%%%%%%%%%%%%%%%%%%%%%%%%%%%%%%%%%%%%%%%%%%%%%%%%%%%%%%%%%%%%%%%%%%%%%%%%%%%%%
% --> Section 1
%%%%%%%%%%%%%%%%%%%%%%%%%%%%%%%%%%%%%%%%%%%%%%%%%%%%%%%%%%%%%%%%%%%%%%%%%%%%%
%%%%%%%%%%%%%%%%%%%%%%%%%%%%%%%%%%%%%%%%%%%%%%%%%%%%%%%%%%%%%%%%%%%%%%%%%%%%%

\newpage

\ssection{Symplectic Manifolds}
\label{section1}

%%%%%%%%%%%%%%%%%%%%%%%%%%%%%%%%%%%%%%%%%%%%%%%%%%%%%%%%%%%%%%%%%%%%%%%%%%%%%
%%%%%%%%%%%%%%%%%%%%%%%%%%%%%%%%%%%%%%%%%%%%%%%%%%%%%%%%%%%%%%%%%%%%%%%%%%%%%

\ssubsection{Symplectic Linear Algebra}
\label{symplectic_linear_algebra}

Let $V$ be a vector space over $\RR$,
and let $\Omega: V \times V \to \RR$ be a skew-symmetric bilinear map.
By a skew-symmetric version of the Gram-Schmidt process,\footnote{Let
$u_1,\dots,u_k$ be a basis of
$U := \{u \in V \mid \Omega(u,v) = 0, \mbox{ for all } v \in V \}$,
and $W$ a complementary subspace such that $V = U \oplus W$.
Take any nonzero $e_1 \in W$.
There is $f_1 \in W$ with $\Omega(e_1,f_1) = 1$.
Let $W_1$ be the span of $e_1,f_1$ and
$W_1^\Omega := \{ v \in V \, | \, \Omega (v,u) = 0 \;
\forall u \in W_1 \}$.
Then $W = W_1 \oplus W_1^\Omega$.
Take any nonzero $e_2 \in W_1^\Omega$.
There is $f_2 \in W_1^\Omega$ for which $\Omega(e_2,f_2) = 1$.
Let $W_2$ be the span of $e_2,f_2$, and so on.}
there is a basis $u_1,\dots,u_k$, $e_1,\dots,e_n$, $f_1,\dots,f_n$
of $V$ for which
$\Omega(u_i,v) = \Omega(e_i,e_j) = \Omega(f_i,f_j) = 0$ and
$\Omega(e_i,f_j) = \delta_{ij}$ for all $i,j$ and all $v \in V$.
Although such a basis is not unique, it is commonly referred to
as a \textbf{canonical basis}.
The dimension $k$ of the subspace
$U = \{u \in V \mid \Omega(u,v) = 0, \mbox{ for all } v \in V \}$
is an invariant of the pair $(V,\Omega)$.
Since $k + 2n = \dim V$, the even number $2n$ is also an invariant
of $(V,\Omega)$, called the \textbf{rank} of
$\Omega$.\index{rank}\index{skew-symmetric bilinear map ! rank}
We denote by ${\widetilde \Omega}: V \to V^*$ the linear map
defined by ${\widetilde \Omega}(v)(u) := \Omega(v,u)$.
We say that $\Omega$ is
\textbf{symplectic}\index{skew-symmetric bilinear map !
symplectic}\index{symplectic !
bilinear map}\index{symplectic ! linear symplectic structure}
(or \textbf{nondegenerate}\index{skew-symmetric bilinear map !
nondegenerate}\index{nondegenerate ! bilinear map}) if
the associated ${\widetilde \Omega}$ is bijective (i.e.,
the kernel $U$ of ${\widetilde \Omega}$ is the trivial space $\{0\}$).
In that case, the map $\Omega$ is called a \textbf{linear symplectic
structure} on $V$, and the pair $(V,\Omega)$ is called a
\textbf{symplectic vector space}.\index{symplectic !
vector space}\index{vector space ! symplectic}
A linear symplectic structure
$\Omega$\index{symplectic ! properties of linear symplectic structures}
expresses a {\em duality} by the bijection ${\widetilde \Omega}:
V \stackrel{\simeq}{\longrightarrow} V^*$,\index{symplectic ! duality}
similar to the (symmetric) case of an inner product.
By considering a canonical basis, we see that
the dimension of a symplectic vector space $(V,\Omega)$
{\em must be even}, $\dim V = 2n$,
and that $V$ admits a basis $e_1,\dots,e_n,f_1,\dots,f_n$ satisfying
$\Omega(e_i,f_j) = \delta_{ij}$ and
$\Omega(e_i,e_j) = 0 = \Omega(f_i,f_j)$.
Such a basis is then called a \textbf{symplectic basis} of
$(V,\Omega)$,\index{symplectic ! basis}
and, in terms of exterior algebra,
$\Omega = e_1^* \wedge f_1^* + \ldots + e_n^* \wedge f_n^*$,
where $e_1^*, \ldots, e_n^*, f_1^*, \ldots,  f_n^*$
is the dual basis.
With respect to a symplectic basis,
the map $\Omega$ is represented by the matrix
\[
        \left[ \begin{array}{cc}
        0 & \mbox{Id} \\
        -\mbox{Id} & 0
        \end{array} \right] \ .
\]

%%%%%%%%%%%%%%%%%%%%%%%%%%%%%%%%%%%%%%%%%%%%%%%%%%%%%%%%%%%%%%%%%%%%%%%%%%%%%

\begin{examples}
\begin{enumerate}
\item
The \textbf{prototype of a symplectic vector space} is
$(\RR^{2n},\Omega_0)$ with $\Omega_0$ such that the canonical basis
$e_1=(1,0,\ldots,0), \ldots, e_n, f_1, \ldots, f_n=(0,\ldots,0,1)$
is a symplectic basis.
Bilinearity then determines $\Omega_0$ on other vectors.

\item
For any real vector space $E$, the direct sum $V = E \oplus E^*$
has a \textbf{canonical symplectic structure} determined by the formula
$\Omega_0 (u \oplus \alpha, v \oplus \beta) = \beta (u) - \alpha (v)$.
If $e_1, \ldots ,e_n$ is a basis of $E$, and
$f_1, \ldots ,f_n$ is the dual basis,
then $e_1 \oplus 0, \ldots ,e_n \oplus 0,
0 \oplus f_1, \ldots ,0 \oplus f_n$ is a symplectic basis for $V$.
\end{enumerate}
\end{examples}

%%%%%%%%%%%%%%%%%%%%%%%%%%%%%%%%%%%%%%%%%%%%%%%%%%%%%%%%%%%%%%%%%%%%%%%%%%%%%

Given a linear subspace $W$ of a symplectic vector space
$(V, \Omega)$, its \textbf{symplectic orthogonal}\index{symplectic !
orthogonal} is the subspace
$W^\Omega := \{ v \in V \, | \, \Omega (v,u) = 0 \;
\mbox{for all } u \in W \}$.
By nondegeneracy, we have $\dim W + \dim W^\Omega = \dim V$
and $(W^\Omega)^\Omega = W$.
For subspaces $W$ and $Y$, we have
$(W \cap Y)^\Omega = W^\Omega + Y^\Omega$, and
if $W \subseteq Y$ then $Y^\Omega \subseteq W^\Omega$.

There are special types of linear subspaces
of a symplectic vector space $(V,\Omega)$.
A subspace $W$ is a \textbf{symplectic subspace}
if the restriction $\Omega|_W$ is nondegenerate,
that is, $W \cap W^\Omega = \{ 0 \}$, or equivalently
$V = W \oplus W^\Omega$.\index{subspace ! symplectic}
A subspace $W$ is an \textbf{isotropic subspace}
if $\Omega|_W \equiv 0$, that is,
$W \subseteq W^\Omega$.\index{subspace ! isotropic}
A subspace $W$ is a \textbf{coisotropic subspace}
if $W^\Omega \subseteq W$.\index{coisotropic !
subspace}\index{subspace ! coisotropic}
A subspace $W$ is a \textbf{lagrangian subspace}
if it is both isotropic and coisotropic, or equivalently,
if it is an isotropic subspace with
$\dim W = {1 \over 2} \dim V$.\index{lagrangian
subspace}\index{subspace ! lagrangian}
A basis $e_1, \ldots ,e_n$ of a lagrangian subspace
can be extended to a symplectic basis:
choose $f_1$ in the symplectic orthogonal to the linear span
of $\{e_2, \ldots ,e_n \}$, etc.

\begin{examples}
\begin{enumerate}
\item
For a symplectic basis as above,
the span of $e_1,f_1$ is symplectic,
that of $e_1,e_2$ isotropic,
that of $e_1,\dots,e_n,f_1$ coisotropic,
and that of $e_1,\dots,e_n$ lagrangian.

\item
The graph of a linear map $A : E \to E^*$ is a lagrangian subspace
of $E \oplus E^*$ with the canonical symplectic structure
if and only if $A$ is symmetric (i.e., $(Au)v = (Av)u$).
Therefore, the grassmannian of all lagrangian subspaces
in a $2n$-dimensional symplectic vector space has dimension
$\frac{n(n+1)}{2}$.
\end{enumerate}
\end{examples}

%%%%%%%%%%%%%%%%%%%%%%%%%%%%%%%%%%%%%%%%%%%%%%%%%%%%%%%%%%%%%%%%%%%%%%%%%%%%%

A \textbf{symplectomorphism}\index{symplectomorphism ! linear}
$\varphi$ between symplectic vector spaces
$(V,\Omega)$ and $(V',\Omega')$ is a linear isomorphism
$\varphi: V \stackrel{\simeq}{\longrightarrow} V'$ such that
$\varphi^*\Omega' = \Omega$.\footnote{By definition,
$(\varphi^*\Omega')(u,v) = \Omega'(\varphi(u),\varphi(v))$.}
If a symplectomorphism exists,
$(V,\Omega)$ and $(V',\Omega')$ are said to be
\textbf{symplectomorphic}\index{symplectomorphic}.
Being symplectomorphic is clearly
an equivalence relation
in the set of all even-dimensional vector spaces.
The existence of canonical bases shows that
every $2n$-dimensional symplectic vector space
$(V,\Omega)$ is symplectomorphic to the prototype $(\RR^{2n},\Omega_0)$;
a choice of a symplectic basis for $(V,\Omega)$ yields
a symplectomorphism to $(\RR^{2n},\Omega_0)$.
Hence, nonnegative even integers classify
equivalence classes for the relation of being symplectomorphic.

Let $\Omega(V)$ be the space of all linear symplectic structures
on the vector space $V$.
Take a $\Omega \in \Omega(V)$, and let $\Sp (V,\Omega)$
be the \textbf{group of symplectomorphisms}\index{group
of symplectomorphisms} of $(V,\Omega)$.
The group $\GL(V)$ of all isomorphisms of $V$
acts {\em transitively} on $\Omega(V)$ by pullback
(i.e., all symplectic structures are related by a linear isomorphism),
and $\Sp (V,\Omega)$ is the stabilizer of the given $\Omega$.
Hence, $\Omega(V) \simeq \GL (V) / \Sp (V,\Omega)$.

%%%%%%%%%%%%%%%%%%%%%%%%%%%%%%%%%%%%%%%%%%%%%%%%%%%%%%%%%%%%%%%%%%%%%%%%%%%%%
%%%%%%%%%%%%%%%%%%%%%%%%%%%%%%%%%%%%%%%%%%%%%%%%%%%%%%%%%%%%%%%%%%%%%%%%%%%%%

\ssubsection{Symplectic Forms}
\label{symplectic_forms}
\index{symplectic ! form}

Let $\omega$ be a de Rham 2-form on a manifold\footnote{Unless
otherwise indicated, all vector spaces are real and finite-dimensional,
all maps are smooth (i.e., $C^\infty$) and all manifolds are smooth,
Hausdorff and second countable.} $M$.
For each point $p \in M$, the map $\omega_p:T_pM \times T_pM \rightarrow \RR$
is skew-symmetric and bilinear on the tangent space to $M$ at $p$,
and $\omega_p$ varies smoothly in $p$.

\begin{definition}
The 2-form $\omega$ is \textbf{symplectic}\index{symplectic !
form}\index{form ! symplectic} if
$\omega$ is closed (i.e., its exterior derivative $d \omega$ is zero)
and $\omega_p$ is symplectic for all $p \in M$.
A \textbf{symplectic manifold}\index{symplectic !
manifold}\index{manifold ! symplectic} is a pair $(M, \omega)$
where $M$ is a manifold and $\omega$ is a symplectic form.
\end{definition}

%%%%%%%%%%%%%%%%%%%%%%%%%%%%%%%%%%%%%%%%%%%%%%%%%%%%%%%%%%%%%%%%%%%%%%%%%%%%%

Symplectic manifolds must be {\em even-dimensional}.
Moreover, the $n$th exterior power $\omega^n$ of a
symplectic form $\omega$ on a $2n$-dimensional manifold is a
{\em volume form}\index{volume}.\footnote{A \textbf{volume form}
is a nonvanishing form of top degree.
If $\Omega$ is a symplectic structure
on a vector space $V$ of dimension $2n$, its $n$th exterior power
$\Omega^n = \Omega \wedge \ldots \wedge \Omega$ does not vanish.
Actually, a skew-symmetric bilinear map $\Omega$ is symplectic
if and only if $\Omega^n \neq 0$.}
Hence, any symplectic manifold $(M,\omega)$ is
{\em canonically oriented}.
The form $\frac{\omega^n}{n!}$ is called the
\textbf{symplectic volume}\index{symplectic ! volume}\index{volume}
or \textbf{Liouville volume}\index{Liouville ! volume}\index{volume !
Liouville} of $(M,\omega)$.
%Non-orientable manifolds can never be symplectic.
When $(M,\omega)$ is a {\em compact} $2n$-dimensional symplectic manifold,
the de Rham cohomology\index{de
Rham cohomology}\index{cohomology ! de Rham} class
$[\omega ^n] \in H^{2n} (M;\RR)$ must be non-zero
by Stokes theorem\index{Stokes theorem}\index{theorem ! Stokes}.
Therefore, the class $[\omega]$ must be non-zero,
as well as its powers $[\omega]^k = [\omega^k] \neq 0$.
{\em Exact symplectic forms} can only exist on noncompact manifolds.
Compact manifolds with a trivial even cohomology group
$H^{2k} (M;\RR)$, $k = 0,1,\ldots,n$,
such as spheres $S^{2n}$ with $n > 1$, can thus never be symplectic.
On a manifold of dimension greater than 2, a function multiple $f \omega$
of a symplectic form $\omega$ is symplectic if and only if
$f$ is a nonzero locally constant function
(this follows from the existence of a symplectic basis).

%%%%%%%%%%%%%%%%%%%%%%%%%%%%%%%%%%%%%%%%%%%%%%%%%%%%%%%%%%%%%%%%%%%%%%%%%%%%%

\begin{examples}
\begin{enumerate}
\index{example ! of symplectic manifold}
\item
Let $M = \RR^{2n}$ with linear coordinates $x_1,\dots,x_n,y_1,\dots,y_n$.
The form
\[
   \omega_0 = \sum \limits_{i=1}^n dx_i \wedge dy_i
\]
is symplectic, and the vectors
$\left( \frac {\partial}{\partial x_1} \right)_p,\dots,\left( \frac
{\partial}{\partial x_n} \right)_p, \left( \frac {\partial}{\partial y_1}
\right)_p,\dots,\left( \frac {\partial}{\partial y_n} \right)_p$
constitute a symplectic basis of $T_pM$.

\vspace*{-1ex}

\item
Let $M = \CC^{n}$ with coordinates $z_1,\dots,z_n$.
The form $\omega_0 = \frac i2 \sum dz_k \wedge d\bar z_k$ is symplectic.
In fact, this form coincides with that of the previous example
under the identification $\CC^{n} \simeq \RR^{2n}$,
$z_k = x_k + iy_k$.

\vspace*{-1ex}

\item
The 2-sphere $S^2$, regarded as the set of unit vectors in $\RR^3$,
has tangent vectors at $p$ identified with vectors orthogonal to $p$.
The standard symplectic form on $S^2$ is induced by the
standard inner (dot) and exterior (vector) products:
$\omega_p (u,v) := \langle p, u \times v \rangle$,
for $u,v \in T_p S^2 = \{ p \} ^\perp$.
%This form is closed because it is of top degree;
%it is nondegenerate because $\langle p, u \times v \rangle \neq 0$
%when $u \neq 0$ and we take, for instance, $v = u \times p$.
This is the standard area form on $S^2$ with total area $4\pi$.
In terms of cylindrical polar coordinates
$0 \leq \theta < 2\pi$ and $-1 \leq z \leq 1$ away from the poles,
it is written $\omega = d \theta \wedge dz$.

\vspace*{-1ex}

\item
On any Riemann surface,
regarded as a 2-dimensional oriented manifold,
any area form, that is, any never vanishing 2-form,
is a symplectic form.

\vspace*{-1ex}

\item
Products of symplectic manifolds are naturally symplectic
by taking the sum of the pullbacks of the symplectic forms from the factors.

\vspace*{-1ex}

\item
If a $(2n+1)$-dimensional manifold $X$
admits a \textbf{contact form}\index{contact form},
that is, a 1-form $\alpha$ such that $\alpha \wedge (d \alpha)^n$
is never vanishing, then the 2-form $d(e^t \alpha)$ is
symplectic on $X \times \RR$,
and the symplectic manifold $(X \times \RR, d(e^t \alpha))$
is called the \textbf{symplectization}\index{symplectization}
of the {\em contact manifold} $(X , \alpha)$.
For more on {\em contact geometry}\index{contact geometry},
see for instance the corresponding contribution in this volume.

\end{enumerate}
\end{examples}

\vspace*{-2ex}

%%%%%%%%%%%%%%%%%%%%%%%%%%%%%%%%%%%%%%%%%%%%%%%%%%%%%%%%%%%%%%%%%%%%%%%%%%%%%

\begin{definition}
Let $(M_1,\omega_1)$ and $(M_2,\omega_2)$
be symplectic manifolds.
A (smooth) map $\psi:M_1\to M_2$ is
\textbf{symplectic}\index{symplectic map} if
$\psi^{*}\omega_2=\omega_1$.\footnote{By definition
of \textbf{pullback}\index{pullback}, we have
$(\psi^{*}\omega_2)_p (u,v) =
(\omega_2)_{\psi(p)} (d\psi_p (u),d\psi_p (v))$,
at tangent vectors $u,v \in T_p M_1$.}
A symplectic diffeomorphism $\varphi:M_1\to M_2$
is a \textbf{symplectomorphism}.\index{symplectomorphism
! definition}
$(M_1,\omega_1)$ and $(M_2,\omega_2)$
are said to be \textbf{symplectomorphic}\index{symplectomorphic}
when there exists a symplectomorphism between them.
\end{definition}

The classification of symplectic manifolds
up to symplectomorphism
is an open problem in symplectic geometry.
However, the local classification is taken care of
by the {\em Darboux theorem} (Theorem~\ref{thm:darboux})\index{theorem !
Darboux}\index{Darboux ! theorem}:
the dimension is the only local invariant of symplectic
manifolds up to symplectomorphisms.
That is, just as any $n$-dimensional manifold
is locally diffeomorphic to $\RR^n$,
any symplectic manifold $(M^{2n},\omega)$
is locally symplectomorphic to $({\RR}^{2n},\omega_{0})$.
As a consequence, if we prove for $(\RR^{2n},\omega_{0})$
a local assertion that is invariant under symplectomorphisms,
then that assertion holds for any symplectic manifold.
We will hence refer to ${\RR}^{2n}$, with linear coordinates
$(x_1,\ldots,x_n,y_1,\ldots, y_n)$, and with symplectic form
$\omega_0=\sum_{i=1}^n dx_i\wedge dy_i$,
as the \textbf{prototype of a local piece of a $2n$-dimensional
symplectic manifold}.

%%%%%%%%%%%%%%%%%%%%%%%%%%%%%%%%%%%%%%%%%%%%%%%%%%%%%%%%%%%%%%%%%%%%%%%%%%%%%
%%%%%%%%%%%%%%%%%%%%%%%%%%%%%%%%%%%%%%%%%%%%%%%%%%%%%%%%%%%%%%%%%%%%%%%%%%%%%

\ssubsection{Cotangent Bundles}
\label{cotangent_bundles}

Cotangent bundles are major examples of symplectic manifolds.
Let $(\cU,x_1,\ldots,x_n)$ be a coordinate chart for a manifold
$X$, with associated cotangent coordinates
$(T^* \cU,x_1,\ldots, x_n,\xi_1,\ldots,\xi_n)$.\footnote{If
an $n$-dimensional manifold $X$ is described by
coordinate charts $(\cU,x_1,\ldots,x_n)$ with $x_i: \cU \to \RR$,
then, at any $x \in \cU$, the differentials
$(dx_i)_x$ form a basis of $T_x^*X$, inducing a map
\[
\begin{array}{rcl}
        T^* \cU & \longrightarrow & \RR^{2n} \\
        (x, \xi) & \longmapsto &
        (x_1, \ldots , x_n, \xi_1, \ldots, \xi_n)\ ,
\end{array}
\]
where $\xi_1, \ldots, \xi_n \in \RR$ are the corresponding
coordinates of $\xi \in T_x^*X$: $\xi = \sum_{i=1}^n \xi_i (dx_i)_x$.
Then $(T^* \cU,x_1,\ldots,x_n, \xi_1, \ldots, \xi_n)$
is a coordinate chart for the cotangent bundle $T^*X$;
the coordinates $x_1,\ldots,x_n, \xi_1, \ldots, \xi_n$ are called the
\textbf{cotangent coordinates}\index{cotangent bundle ! coordinates}
associated to the coordinates $x_1,\ldots,x_n$ on $\cU$.
One verifies that the transition functions on the overlaps are smooth,
so $T^* X$ is a $2n$-dimensional manifold.}
Define a symplectic form on $T^* \cU$ by
\index{form ! canonical}
\index{canonical form on $T^*X$ ! coordinate definition}
\[
        \omega = \sum \limits_{i=1}^n dx_i \wedge d\xi_i \ .
\]
One can check that this $\omega$ is intrinsically defined
by considering the 1-form on $T^* \cU$
\index{form ! tautological}\index{form !
canonical}\index{canonical form on $T^*X$ ! intrinsic definition}
\index{tautological form on $T^*X$ ! coordinate definition}
\[
        \alpha=\sum\limits_{i=1}^n \xi_i \ dx_i
\]
which satisfies $\omega=-d\alpha$ and is coordinate-independent:
in terms of the natural projection
$\pi: M \to X$, $p=(x,\xi) \mapsto x$,
the form $\alpha$\index{form ! tautological}
\index{tautological form on $T^*X$ ! intrinsic definition}
may be equivalently defined pointwise without coordinates by
\[
        \alpha_{p}=(d\pi_{p})^* \xi \quad \in T_p ^* M\ ,
\]
where $(d\pi_{p})^*: T_x^*X \to T_p ^*M$ is the transpose of $d\pi_{p}$,
that is, $\alpha_{p}(v) = \xi ( (d\pi_{p})v )$ for $v\in T_{p}M$.
Or yet, the form $\alpha$ is
uniquely characterized by the property that $\mu^* \alpha = \mu$
for every 1-form $\mu: X \to T^* X$\index{tautological
form on $T^*X$ ! property} (see Proposition~\ref{prop:tautological_property}).
The 1-form $\alpha$ is the \textbf{tautological form}
(or the \textbf{Liouville 1-form}\index{Liouville 1-form}) and the
2-form $\omega$ is the \textbf{canonical symplectic form}\index{form !
tautological}\index{tautological form on $T^*X$ !
coordinate definition}\index{form ! canonical}\index{canonical form
on $T^*X$ ! coordinate definition} on $T^*X$.
When referring to a cotangent bundle as a symplectic manifold,
the symplectic structure is meant to be given by this canonical
$\omega$.\index{example ! of symplectic manifold}

%%%%%%%%%%%%%%%%%%%%%%%%%%%%%%%%%%%%%%%%%%%%%%%%%%%%%%%%%%%%%%%%%%%%%%%%%%%%%

Let $X_1$ and $X_2$ be $n$-dimensional manifolds
with cotangent bundles $M_1=T^*X_1$ and $M_2=T^*X_2$,
and tautological 1-forms $\alpha_1$ and $\alpha_2$.
Suppose that $f:X_1\to X_2$ is a diffeomorphism.
Then there is a natural diffeomorphism $f_{\sharp}:M_1\to M_2$
which \textbf{lifts}\index{lift ! of a diffeomorphism} $f$; namely,
for $p_1=(x_1,\xi_1) \in M_1$ we define
\[
        f_{\sharp}(p_1) =p_2 =(x_2,\xi_2) \ ,
        \quad \mbox{ with } \left\{ \begin{array}{l}
        x_2=f(x_1) \in X_2 \quad \mbox{ and } \\
        \xi_1=(df_{x_1})^* \xi_2 \in T_{x_1}^* X_1 \ ,
        \end{array} \right.
\]
where $(df_{x_1})^* :
T_{x_2}^* X_2 \stackrel{\simeq}{\longrightarrow} T_{x_1}^* X_1$,
so $f_{\sharp}|_{T_{x_1}^* }$ is the inverse map of $(df_{x_1})^*$.

\begin{proposition}
\label{prop:diffeo_lift}
The lift $f_{\sharp}$ of a diffeomorphism $f: X_1 \rightarrow X_2$
pulls the tautological form on $T^* X_2$ back to the
tautological form on $T^* X_1$, i.e., $(f_{\sharp})^* \alpha_2=\alpha_1$.
\end{proposition}

\vspace*{-2ex}

\begin{proof}
At $p_1=(x_1,\xi_1)\in M_1$, the claimed identity says
$\left( d f_{\sharp} \right)^*_{p_1} (\alpha_2)_{p_2} = (\alpha_1)_{p_1}$,
where $p_2=f_{\sharp}(p_1)$, that is,
$p_2 = (x_2,\xi_2)$ where $x_2 = f(x_1)$ and $(df_{x_1})^* \xi_2 = \xi_1$.
This can be proved as follows:
\[
\begin{array}{rclcl}
        (df_{\sharp})^* _{p_1}(\alpha_2)_{p_2}
        & = & (df_{\sharp})^* _{p_1}(d\pi_2)^* _{p_2}\xi_2
        & & \mbox{ by definition of $\alpha_2$} \\
        & = & \left( d(\pi_2 \circ f_{\sharp}) \right) ^* _{p_1}\xi_2
        & & \mbox{ by the chain rule} \\
        & = & \left( d(f \circ \pi_1) \right) ^* _{p_1}\xi_2
        & & \mbox{ because $\pi_2 \circ f_{\sharp} = f \circ \pi_1$} \\
        & = & (d\pi_1)^* _{p_1}(df)^* _{x_1} \xi_2
        & & \mbox{ by the chain rule} \\
        & = & (d\pi_1)^* _{p_1} \xi_1
        & & \mbox{ by definition of $f_{\sharp}$} \\
        & = & (\alpha_1)_{p_1}
        & & \mbox{ by definition of $\alpha_1$} \ .
\end{array}
\]
\end{proof}

As a consequence of this naturality for the tautological
form,\index{canonical form on $T^*X$ ! naturality}\index{tautological form
on $T^*X$ ! naturality}\index{cotangent bundle ! canonical symplectomorphism}
a diffeomorphism of manifolds induces a canonical
symplectomorphism\index{canonical !
symplectomorphism}\index{symplectomorphism !
canonical}\index{cotangent bundle ! canonical symplectomorphism}
of cotangent bundles:

\begin{corollary}
The lift $f_{\sharp}: T^*X_1 \to T^*X_2$ of a diffeomorphism
$f: X_1 \rightarrow X_2$ is a symplectomorphism
for the canonical symplectic forms, i.e.,
$(f_{\sharp})^* \omega_2=\omega_1$.
\end{corollary}

%%%%%%%%%%%%%%%%%%%%%%%%%%%%%%%%%%%%%%%%%%%%%%%%%%%%%%%%%%%%%%%%%%%%%%%%%%%%%

In terms of the group (under composition)
of diffeomorphisms $\mathrm{Diff}(X)$
of a manifold $X$, and the \textbf{group of symplectomorphisms}
$\mathrm{Sympl} (T^*X,\omega)$\index{symplectomorphism ! group of
symplectomorphisms}\index{group ! of symplectomorphisms}
of its cotangent bundle, we see that the injection
$\mathrm{Diff}(X) \to \mathrm{Sympl} (T^*X,\omega)$,
$f \mapsto f_{\sharp}$ is a group homomorphism.
Clearly this is not surjective: for instance,
consider the symplectomorphism $T^*X \to T^*X$ given by
translation along cotangent fibers.

\begin{example}
Let $X_1=X_2=S^1$.
Then $T^* S^1$ is a cylinder $S^1 \times \RR$.
The canonical form is the area form $\omega=d\theta\wedge d\xi$.
If $f:S^1\rightarrow S^1$ is any diffeomorphism, then
$f_{\sharp}: S^1 \times \RR \rightarrow S^1 \times \RR$
is a symplectomorphism, i.e., is an area-preserving
diffeomorphism of the cylinder.
Translation along the $\RR$ direction is area-preserving
but is not induced by a diffeomorphism of the base manifold $S^1$.
\end{example}

%%%%%%%%%%%%%%%%%%%%%%%%%%%%%%%%%%%%%%%%%%%%%%%%%%%%%%%%%%%%%%%%%%%%%%%%%%%%%%%

There is a criterion for which cotangent symplectomorphisms
arise as lifts of diffeomorphisms in terms of the tautological form.
First note the following feature of symplectic manifolds
with \textbf{exact symplectic forms}.
Let $\alpha$ be a 1-form on a manifold $M$
such that $\omega = - d\alpha$ is symplectic.
There exists a unique vector field $v$ whose
interior product with $\omega$ is $\alpha$,
i.e., $\imath_v \omega = - \alpha$.
If $g : M \to M$ is a symplectomorphism that preserves $\alpha$
(that is, $g^* \alpha = \alpha$), then $g$ commutes
with the flow\footnote{For
$p \in M$, $(\exp tv) (p) $ is the unique
curve in $M$ solving the initial value problem
\[
\left\{ \begin{array}{l}
        {d \over dt} (\exp tv (p)) = v (\exp tv (p)) \\
        (\exp tv) (p) |_{t=0} = p
\end{array} \right.
\]
for $t$ in some neighborhood of $0$.
The one-parameter group of diffeomorphisms $\exp tv$
is called the \textbf{flow} of the vector field $v$.}
of $v$, i.e., $(\exp tv) \circ g = g \circ (\exp tv)$.
When $M = T^* X$ is the cotangent bundle of an
arbitrary $n$-dimensional manifold $X$, and
$\alpha$ is the tautological 1-form on $M$,
the vector field $v$
is just $\sum \xi_i \, {\partial \over \partial \xi_i}$
with respect to a cotangent coordinate chart
$(T^* \cU,x_1,\ldots,x_n, \xi_1, \ldots, \xi_n)$.
The flow $\exp tv$, $-\infty < t < \infty$, satisfies
$(\exp tv) (x, \xi) = (x, e^t \xi)$, for every $(x, \xi)$ in $M$.

\begin{theorem}
A symplectomorphism $g: T^*X \to T^*X$ is a lift of a
diffeomorphism $f: X \rightarrow X$ if and only if
it preserves the tautological form: $g^* \alpha = \alpha$.
\end{theorem}

\vspace*{-2ex}

\begin{proof}
By Proposition~\ref{prop:diffeo_lift},
a lift $f_{\sharp} : T^*X \to T^*X$ of a
diffeomorphism $f: X \rightarrow X$
preserves the tautological form.
Conversely, if $g$ is a symplectomorphism of $M$ that preserves
$\alpha$, then $g$ preserves the cotangent fibration:
by the observation above,
$g (x, \xi) = (y,\eta) \Rightarrow g (x, \lambda \xi) = (y, \lambda \eta)$
for all $(x, \xi) \in M$ and $\lambda > 0$,
and this must hold also for $\lambda \leq 0$ by the
differentiability of $g$ at $(x,0)$.
Therefore, there exists a diffeomorphism $f: X \to X$
such that $\pi \circ g = f \circ \pi$, where $\pi : M \to X$ is the
projection map $\pi (x, \xi) = x$, and $g = f_\#$.
\end{proof}

%%%%%%%%%%%%%%%%%%%%%%%%%%%%%%%%%%%%%%%%%%%%%%%%%%%%%%%%%%%%%%%%%%%%%%%%%%%%%

The canonical form is natural also in the following way.
Given a smooth function $h: X \to \RR$,
the diffeomorphism $\tau_h$ of $M = T^*X$ defined by
$\tau_h (x, \xi) = (x, \xi + dh_x)$
turns out to be always a symplectomorphism.
Indeed, if $\pi : M \to X$, $\pi (x, \xi) = x$,
is the projection, we have
$\tau_h^* \alpha = \alpha + \pi^* dh$,
so that $\tau_h^* \omega = \omega$.

%%%%%%%%%%%%%%%%%%%%%%%%%%%%%%%%%%%%%%%%%%%%%%%%%%%%%%%%%%%%%%%%%%%%%%%%%%%%%
%%%%%%%%%%%%%%%%%%%%%%%%%%%%%%%%%%%%%%%%%%%%%%%%%%%%%%%%%%%%%%%%%%%%%%%%%%%%%

\ssubsection{Moser's Trick}
\label{sec:trick}
\index{symplectic ! equivalence}\index{Moser ! trick}

There are other relevant notions of equivalence
for symplectic manifolds\footnote{Understanding these notions
and the normal forms requires tools,
such as isotopies (by \textbf{isotopy}\index{isotopy}
we mean a smooth one-parameter family of diffeomorphisms
starting at the identity, like the
flow of a vector field), Lie derivative,
tubular neighborhoods and the homotopy formula
in de Rham theory, covered in differential geometry
or differential topology texts.}
besides being symplectomorphic.\index{symplectomorphic}
Let $M$ be a manifold with
two symplectic forms $\omega_0 , \omega_1$.

\begin{definition}
The symplectic manifolds $(M,\omega_0)$ and $(M,\omega_1)$ are
\textbf{strongly isotopic}\index{strong isotopy}\index{symplectic !
strong isotopy} if there is an isotopy $\rho_t:M\to M$ such that
$\rho^*_1 \omega_1=\omega_0$.
$(M,\omega_0)$ and $(M,\omega_1)$ are
\textbf{deformation-equivalent}\index{deformation
equivalence}\index{symplectic ! deformation equivalence}
if there is a smooth family $\omega_t$ of symplectic
forms joining $\omega_0$ to $\omega_1$.
$(M,\omega_0)$ and $(M,\omega_1)$ are
\textbf{isotopic}\index{isotopy ! symplectic}\index{symplectic ! isotopy}
if they are deformation-equivalent and the de Rham cohomology
class\index{de Rham cohomology class} $[\omega_t]$
is independent of $t$.
\end{definition}

Hence, being strongly isotopic implies being symplectomorphic,
and being isotopic implies being deformation-equivalent.
We also have that being strongly isotopic implies being isotopic,
because, if $\rho_t : M \to M$ is an isotopy such that
$\rho_1 ^* \omega_1 = \omega_0$, then
$\omega_t := \rho_t ^* \omega_1$
is a smooth family of symplectic forms joining $\omega_1$
to $\omega_0$ and $[\omega_t]=[\omega_1]$, $\forall t$,
by the homotopy invariance of de Rham cohomology.

Moser~\cite{mo:volume} proved that, on a compact manifold,
being isotopic implies being strongly isotopic (Theorem~\ref{thm:moser}).
McDuff\index{McDuff counterexample}\index{example ! McDuff}
showed that deformation-equivalence is indeed a necessary
hypothesis: even if $[\omega_0]=[\omega_1]\in H^2(M;\RR)$,
there are compact examples where $(M,\omega_0)$ and $(M,\omega_1)$
are not strongly isotopic; see Example~7.23 in~\cite{mc-sa:introduction}.
In other words, fix $c\in H^2(M)$ and define
$S_c$ as the set of symplectic forms $\omega$ in $M$ with $[\omega]=c$.
On a compact manifold,
all symplectic forms in the same path-connected component
of $S_c$ are symplectomorphic according to the Moser theorem,
though there might be symplectic forms in different components
of $S_c$ that are not symplectomorphic.

%%%%%%%%%%%%%%%%%%%%%%%%%%%%%%%%%%%%%%%%%%%%%%%%%%%%%%%%%%%%%%%%%%%%%%%%%%%%%

\begin{theorem}
\label{thm:moser}\index{Moser ! theorem}\index{theorem !
Moser}
\textbf{(Moser)} $\;$
Let $M$ be a compact manifold with symplectic forms
$\omega_0$ and $\omega_1$.
Suppose that $\omega_t$, $0\leq t\leq 1$, is a smooth
family of symplectic forms joining $\omega_0$ to $\omega_1$
with cohomology class $[\omega_t]$ independent of $t$.
Then there exists an isotopy $\rho:M \times \RR \to M$ such that
$\rho^*_t\omega_t=\omega_0$, $0\leq t\leq 1$.
\end{theorem}

Moser applied an extremely useful argument, known as
\textbf{Moser's trick}\index{Moser ! trick},
starting with the following observation.
If there existed an isotopy $\rho:M \times \RR \to M$
such that $\rho^*_t\omega_t=\omega_0$, $0\leq t\leq 1$,
in terms of the associated time-dependent vector field
\[
        v_t := \frac{d\rho_t}{dt}\circ\rho^{-1}_t \ , \qquad t\in\RR \ ,
\]
we would then have for all $0\leq t\leq 1$ that
\[
        \displaystyle{0=\frac{d}{dt} (\rho^*_t\omega_t) =\rho^*_t
        \big(\cL_{v_t}\omega_t+\frac{d\omega_t}{dt}\big)}
        \iff \cL_{v_t}\omega_t +
        \displaystyle{\frac{d\omega_t}{dt}} =0\ .
\]
Conversely, the existence of a smooth time-dependent vector
field $v_t$, $t \in \RR$, satisfying the last equation
is enough to produce by integration (since $M$ is compact)
the desired isotopy $\rho:M \times \RR \to M$ satisfying
$\rho^*_t\omega_t= \rho^*_0 \omega_0=\omega_0$, for all $t$.
So everything boils down to solving the equation
$\cL_{v_t}\omega_t + \frac{d\omega_t}{dt} =0$ for $v_t$.

\begin{proof}
By the cohomology assumption that $\big[ \frac{d}{dt} \omega_t\big] = 0$,
there exists a {\em smooth} family of 1-forms $\mu_t$ such that
\[
        \displaystyle{\frac{d\omega_t}{dt} = d\mu_t}\ ,
        \quad 0 \leq t \leq 1\ .
\]
The argument involves the Poincar\'e lemma for compactly-supported forms,
together with the Mayer-Vietoris sequence in order to use
induction on the number of charts in a good cover of $M$;
for a sketch, see page~95 in~\cite{mc-sa:introduction}.
In the simplest case where
$\omega_t=(1-t)\omega_0 + t\omega_1$ with $[\omega_0]=[\omega_1]$,
we have that $\frac{d\omega_t}{dt} = \omega_1-\omega_0 =d\mu$ is exact.

The nondegeneracy assumption on $\omega_t$,
guarantees that we can pointwise solve the equation,
known as \textbf{Moser's equation},\index{Moser ! equation}
\[
        \imath_{v_t} \omega_t + \mu_t=0
\]
to obtain a unique smooth family of vector fields $v_t$, $0 \leq t \leq 1$.
Extend $v_t$ to all $t \in \RR$.
Thanks to the compactness of $M$, the vector fields $v_t$
generate an isotopy $\rho$ satisfying
$\frac{d\rho_t}{dt} = v_t \circ \rho_t$.
Then we indeed have
\[
        \displaystyle{\frac{d}{dt} (\rho_t^*\omega_t) =
        \rho^*_t(\cL_{v_t}\omega_t + \frac{d\omega_t}{dt}) =
        \rho^*_t(d \imath_{v_t} \omega_t + d\mu_t)  =
        \rho^*_t d (\imath_{v_t} \omega_t + \mu_t) = 0}\ ,
\]
where we used Cartan's magic formula\index{Cartan ! magic formula}
in $\cL_{v_t}\omega_t = d \imath_{v_t}\omega_t + \imath_{v_t} d\omega_t$.
\end{proof}

\vspace*{-1ex} 

\begin{example}
On a compact oriented 2-dimensional manifold $M$,
a symplectic form is just an area form.\index{form ! area}
Let $\omega_0$ and $\omega_1$ be two area forms on $M$.
If $[\omega_0] = [\omega_1]$,
i.e., $\omega_0$ and $\omega_1$ give the same total area,
then any convex combination of them is symplectic
(because they induce the same orientation),
and there is an isotopy $\varphi_t : M \to M$, $t \in [0,1]$,
such that $\varphi_1^* \omega_0 = \omega_1$.
Therefore, up to strong isotopy, there is a unique symplectic representative
in each non-zero 2-cohomology class of $M$.
\end{example}

On a {\em noncompact} manifold, given $v_t$, we would need to check
the existence for $0\leq t\leq 1$ of an isotopy $\rho_t$
solving the differential equation $\frac{d\rho_t}{dt} = v_t \circ \rho_t$.

%%%%%%%%%%%%%%%%%%%%%%%%%%%%%%%%%%%%%%%%%%%%%%%%%%%%%%%%%%%%%%%%%%%%%%%%%%%%%
%%%%%%%%%%%%%%%%%%%%%%%%%%%%%%%%%%%%%%%%%%%%%%%%%%%%%%%%%%%%%%%%%%%%%%%%%%%%%

\ssubsection{Darboux and Moser Theorems}
\label{moser_relative_theorem}

By a \textbf{submanifold} of a manifold $M$ we mean either
a manifold $X$ with a {\em closed embedding}\footnote{A \textbf{closed
embedding} is a {\em proper} injective immersion.
A map is \textbf{proper}\index{proper map} when
its preimage of a compact set is always compact.}
$i: X \hookrightarrow M$,
or an \textbf{open submanifold} (i.e., an open subset of $M$).
%Regarding the embedding $i: X \hookrightarrow M$ as an inclusion,
%we identify points $p = i(p)$ and tangent vectors $v = di_p (v)$.

Given a $2n$-dimensional manifold $M$,
a $k$-dimensional submanifold $X$,
neighborhoods $\cU_0,\cU_1$ of $X$,
and symplectic forms $\omega_0,\omega_1$ on $\cU_0,\cU_1$,
we would like to know whether there exists a
\textbf{local symplectomorphism preserving $X$}, i.e.,
a diffeomorphism $\varphi:\cU_0\to \cU_1$
with $\varphi^*\omega_1=\omega_0$ and $\varphi(X)=X$.
Moser's Theorem~\ref{thm:moser}
addresses the case where $X=M$.
At the other extreme, when $X$ is just one point,
there is the classical Darboux theorem (Theorem~\ref{thm:darboux}).
In general, we have:

\begin{theorem}
\label{thm:moser_relative}\index{Moser !
theorem -- relative version}\index{theorem !
Moser -- relative version}
\textbf{(Moser Theorem -- Relative Version)} $\;$
Let $\omega_0$ and $\omega_1$ be symplectic forms on a manifold $M$,
and $X$ a compact submanifold of $M$.
Suppose that the forms coincide, $\omega_0|_p = \omega_1|_p$,
at all points $p\in X$.
Then there exist neighborhoods $\cU_0$ and $\cU_1$ of $X$ in $M$,
and a diffeomorphism $\varphi : \cU_0 \to \cU_1$ such that
$\varphi^*\omega_1=\omega_0$ and
$\varphi$ restricted to $X$ is the identity map.
\end{theorem}

\vspace*{-2ex}

\begin{proof}
Pick a tubular neighborhood $\cU_0$ of $X$.
The 2-form $\omega_1-\omega_0$ is closed on $\cU_0$, and
satisfies $(\omega_1-\omega_0)_p=0$ at all $p\in X$.
By the homotopy formula on the tubular neighborhood,
there exists a 1-form $\mu$ on
$\cU_0$ such that $\omega_1-\omega_0=d\mu$ and $\mu_p=0$ at all $p\in X$.
Consider the family
$\omega_t=(1-t)\omega_0+t\omega_1=\omega_0+ td\mu$
of closed 2-forms on $\cU_0$.
Shrinking $\cU_0$ if necessary, we can assume that $\omega_t$ is symplectic
for $t \in [0,1]$, as nondegeneracy is an open property.
Solve Moser's equation, $\imath_{v_t} \omega_t = - \mu$, for $v_t$
By integration,
shrinking $\cU_0$ again if necessary, there exists
a local isotopy $\rho: \cU_0 \times [0,1] \to M$ with
$\rho^*_t\omega_t=\omega_0$, for all $t \in [0,1]$.
Since $v_t |_X =0$, we have $\rho_t |_X =\id_X$.
Set $\varphi=\rho_1$, $\cU_1=\rho_1(\cU_0)$.
\end{proof}

\vspace*{-1ex}

\begin{theorem}
\label{thm:darboux}\index{theorem ! Darboux}\index{Darboux ! theorem}
\textbf{(Darboux)} $\;$
Let $(M,\omega)$ be a symplectic manifold, and
let $p$ be any point in $M$.
Then we can find a chart
$(\cU,x_1, \ldots ,x_n,y_1,\ldots y_n)$
centered at $p$ where
\[
        \omega = \displaystyle{\sum_{i=1}^n dx_i\wedge dy_i}\ .
\]
\end{theorem}

Such a coordinate chart $(\cU,x_1,\dots,x_n,y_1,\dots,y_n)$
is called a \textbf{Darboux chart},\index{theorem !
Darboux}\index{Darboux ! chart}\index{chart ! Darboux}
and the corresponding coordinates are called \textbf{Darboux
coordinates}.\index{Darboux ! coordinates}\index{coordinates ! Darboux}
%The vector space version of Darboux coordinates is the
%symplectic basis discussed in Section~\ref{symplectic_linear_algebra}.

The classical proof of Darboux's theorem is by induction
on the dimension of the manifold~\cite{ar:mathematical},
in the spirit of the argument for a symplectic basis
(Section~\ref{symplectic_linear_algebra}).
The proof below, using Moser's theorem,
was first provided by Weinstein~\cite{we:lagrangian}.

\begin{proof}
Apply Moser's relative theorem to $X=\{p\}$.
More precisely, use any symplectic basis for $(T_pM, \omega_p)$
to construct coordinates $(x'_1,\ldots ,x'_n,$ $y'_1,\ldots y'_n)$
centered at $p$ and valid on some neighborhood $\cU'$, so that
$\omega_p=  \left. \sum dx'_i\wedge dy'_i \right|_p$.
There are two symplectic forms on $\cU'$: the given $\omega_0=\omega$
and $\omega_1=\sum dx'_i\wedge dy'_i$.
By Theorem~\ref{thm:moser_relative}, there are neighborhoods
$\cU_0$ and $\cU_1$ of $p$, and a diffeomorphism
$\varphi : \cU_0 \to \cU_1$ such that $\varphi (p)=p$ and
$\varphi ^*(\sum dx'_i\wedge dy'_i)=\omega$.
Since $\varphi ^*(\sum dx'_i \wedge dy'_i)
= \sum d(x'_i\circ \varphi) \wedge d(y'_i\circ \varphi)$,
we simply set new coordinates $x_i=x'_i\circ \varphi$,
$y_i = y'_i\circ \varphi$.
\end{proof}

Darboux's theorem is easy in the 2-dimensional case.\index{Darboux ! theorem
in dimension two}\index{theorem ! Darboux}
Being closed $\omega$ is locally exact,
$\omega = d \alpha$.
Every nonvanishing 1-form on a surface can be written
locally as $\alpha = g \, dh$ for suitable functions $g,h$,
where $h$ is a coordinate on the local leaf space
of the kernel foliation of $\alpha$.
The form $\omega = dg \wedge dh$ is nondegenerate
if and only if $(g,h)$ is a local diffeomorphism.
By the way, transversality shows that the normal
form for a {\em generic}\footnote{\textbf{Generic} here means
that the subset of those 2-forms having this behavior
is open, dense and invariant under diffeomorphisms of the manifold.}
2-form is $x dx \wedge dy$ near a point where it is degenerate.

%%%%%%%%%%%%%%%%%%%%%%%%%%%%%%%%%%%%%%%%%%%%%%%%%%%%%%%%%%%%%%%%%%%%%%%%%%%%%
%%%%%%%%%%%%%%%%%%%%%%%%%%%%%%%%%%%%%%%%%%%%%%%%%%%%%%%%%%%%%%%%%%%%%%%%%%%%%

\ssubsection{Symplectic Submanifolds}
\label{symplectic_submanifolds}

Moser's argument permeates many other proofs,
including those of the next two results
regarding {\em symplectic submanifolds}.
Let $(M,\omega)$ be a symplectic manifold.

\begin{definition}
A \textbf{symplectic submanifold}\index{symplectic submanifold !
definition} of $(M,\omega)$
is a submanifold $X$ of $M$ where, at each $p \in X$,
the space $T_pX$ is a symplectic subspace of $(T_pM,\omega_p)$.
\end{definition}

If $i: X \hookrightarrow M$ is the inclusion of a
symplectic submanifold $X$,
then the restriction of $\omega$ to $X$ is a symplectic form,
so that $(X,i^*\omega)$ is itself a symplectic manifold.

%%%%%%%%%%%%%%%%%%%%%%%%%%%%%%%%%%%%%%%%%%%%%%%%%%%%%%%%%%%%%%%%%%%%%%%%%%%%%

Let $X$ be a symplectic submanifold of $(M,\omega)$.
At each $p \in X$, we have $T_pM = T_p X \oplus (T_p X)^{\omega_p}$
(Section~\ref{symplectic_linear_algebra}),
so the map $(T_p X)^{\omega_p} \to T_p M/T_p X$ is an isomorphism.
This canonical identification of the
\textbf{normal space}\index{normal ! space}\index{space ! normal}
of $X$ at $p$, $N_p X := T_p M/T_p X$, with the
symplectic orthogonal $(T_p X)^{\omega_p}$,
yields a \-canonical identification
of the \textbf{normal bundle} $NX$ with the {\em symplectic vector bundle}
$(TX)^\omega$.
A \textbf{symplectic vector bundle}\index{symplectic ! vector bundle}
is a vector bundle $E \to X$ equipped with a
smooth\footnote{Smoothness means that, for any pair of (smooth) sections
$u$ and $v$ of $E$, the real-valued function $\Omega (u,v) : X \to \RR$
given by evaluation at each point is smooth.} field
$\Omega$ of fiberwise nondegenerate
skew-symmetric bilinear maps $\Omega_p : E_p \times E_p \to \RR$.
The \textbf{symplectic normal bundle} is the normal
bundle of a symplectic submanifold, with the symplectic structure
induced by orthogonals.
The next theorem, due to Weinstein~\cite{we:lagrangian},
states that a neighborhood of a symplectic submanifold $X$
is determined by $X$ and (the isomorphism class of)
its symplectic normal bundle.

\begin{theorem}
\label{thm:weinstein_symplectic}\index{Weinstein !
symplectic neighborhood theorem}\index{theorem ! Weinstein
symplectic neighborhood}\index{neighborhood ! Weinstein
symplectic neighborhood}
\textbf{(Symplectic Neighborhood Theorem)} $\;$
Let $(M_0,\omega_0)$, $(M_1,\omega_1)$ be symplectic manifolds
with diffeomorphic compact symplectic submanifolds $X_0$, $X_1$.
Let $i_0 : X_0 \hookrightarrow M_0$, $i_1 : X_1 \hookrightarrow M_1$
be their inclusions.
Suppose there is an isomorphism $\widetilde \phi : NX_0 \to NX_1$
of the corresponding symplectic normal bundles covering a
symplectomorphism $\phi : (X_0, i_0^* \omega_0) \to (X_1, i_1^* \omega_1)$.
Then there exist neighborhoods $\cU_0 \subset M_0$, $\cU_1 \subset M_1$
of $X_0$, $X_1$ and a symplectomorphism
$\varphi: \cU_0 \to \cU_1$ extending $\phi$ such that the restriction
of $d \varphi$ to the normal bundle $NX_0$ is $\widetilde \phi$.
\end{theorem}

As first noted by Thurston~\cite{th:examples},
the form $\Omega + \pi^* \omega_X$ is symplectic
in some neighborhood of the zero section in $NX$,
where $\pi : NX \to X$ is the bundle projection and
$\omega_X$ is the restriction of $\omega$ to $X$.
Therefore, {\em a compact symplectic submanifold $X$ always admits
a tubular neighborhood in the ambient $(M,\omega)$
symplectomorphic to a tubular neighborhood of the zero section
in the symplectic normal bundle $NX$}.

\begin{proof}
By the Whitney extension theorem\footnote{
\textbf{Whitney Extension Theorem: }\index{Whitney
extension theorem}\index{theorem ! Whitney extension}
{\em Let $M$ be a manifold and $X$ a submanifold of $M$.
Suppose that at each $p\in X$ we are given a linear isomorphism
$L_p:T_pM\stackrel{\simeq}{\longrightarrow} T_pM$ such that
$L_p|_{T_pX}= \Id_{T_pX}$ and $L_p$ depends smoothly on $p$.
Then there exists an embedding $h : \cN \to M$ of some
neighborhood $\cN$ of $X$ in $M$ such that $h|_X= \id_X$
and $dh_p=L_p$ for all $p\in X$.}
A proof relies on a tubular neighborhood model.}\index{theorem ! Whitney
extension}\index{Whitney extension theorem}
there exist neighborhoods $\cU_0 \subset M_0$ and $\cU_1 \subset M_1$
of $X_0$ and $X_1$, and a diffeomorphism
$h: \cU_0 \to \cU_1$ such that $h \circ i_0 = i_1 \circ \phi$
and the restriction of
$d h$ to the normal bundle $NX_0$ is the given $\widetilde \phi$.
Hence $\omega_0$ and $h^* \omega_1$ are two symplectic
forms on $\cU_0$ which coincide at all points $p\in X_0$.
The result now follows from Moser's relative theorem
(Theorem~\ref{thm:moser_relative}).
\end{proof}

%%%%%%%%%%%%%%%%%%%%%%%%%%%%%%%%%%%%%%%%%%%%%%%%%%%%%%%%%%%%%%%%%%%%%%%%%%%%%

Carefully combining Moser's argument
with the existence of an ambient isotopy
that produces a given deformation of
a compact submanifold, we can show:

\begin{theorem}
Let $X_t$, $t \in [0,1]$, be a (smooth) family of compact
symplectic submanifolds of a compact symplectic manifold $(M,\omega)$.
Then there exists an isotopy $\rho:M \times \RR \to M$ such that
for all $t \in [0,1]$ we have $\rho^*_t\omega=\omega$
and $\rho_t (X_0) = X_t$.
\end{theorem}

%%%%%%%%%%%%%%%%%%%%%%%%%%%%%%%%%%%%%%%%%%%%%%%%%%%%%%%%%%%%%%%%%%%%%%%%%%%%%

Inspired by complex geometry, Donaldson~\cite{do:almost}
proved the following theorem
on the existence of symplectic submanifolds.
A major consequence is the characterization
of symplectic manifolds in terms of {\em Lefschetz pencils};
see Section~\ref{sec:pencils}.

\begin{theorem}
\label{thm:donaldson_submanifolds}
\textbf{(Donaldson)} $\;$
Let $(M,\omega)$ be a compact symplectic manifold.
Assume that the cohomology class $[\omega]$ is integral,
i.e., lies in $H^2 (M;\ZZ)$.
Then, for every sufficiently large integer $k$,
there exists a connected codimension-2 symplectic
submanifold $X$ representing the Poincar\'e dual of
the integral cohomology class $k[\omega]$.
\end{theorem}

Under the same hypotheses, Auroux extended this result
to show that given $\alpha \in H_{2m} (M;\ZZ)$ there exist
positive $k,\ell \in \ZZ$ such that
$k \mathrm{PD} [\omega^{n-m}] + \ell \alpha$
is realized by a $2m$-dimensional symplectic submanifold.

%%%%%%%%%%%%%%%%%%%%%%%%%%%%%%%%%%%%%%%%%%%%%%%%%%%%%%%%%%%%%%%%%%%%%%%%%%%%%
%%%%%%%%%%%%%%%%%%%%%%%%%%%%%%%%%%%%%%%%%%%%%%%%%%%%%%%%%%%%%%%%%%%%%%%%%%%%%
% --> Section 2
%%%%%%%%%%%%%%%%%%%%%%%%%%%%%%%%%%%%%%%%%%%%%%%%%%%%%%%%%%%%%%%%%%%%%%%%%%%%%
%%%%%%%%%%%%%%%%%%%%%%%%%%%%%%%%%%%%%%%%%%%%%%%%%%%%%%%%%%%%%%%%%%%%%%%%%%%%%

\newpage

\ssection{Lagrangian Submanifolds}
\index{lagrangian submanifold ! definition}
\label{section2}

%%%%%%%%%%%%%%%%%%%%%%%%%%%%%%%%%%%%%%%%%%%%%%%%%%%%%%%%%%%%%%%%%%%%%%%%%%%%%
%%%%%%%%%%%%%%%%%%%%%%%%%%%%%%%%%%%%%%%%%%%%%%%%%%%%%%%%%%%%%%%%%%%%%%%%%%%%%

\ssubsection{First Lagrangian Submanifolds}
\label{lagrangian_submanifolds}

Let $(M,\omega)$ be a symplectic manifold.

\begin{definition}
A submanifold $X$ of $(M,\omega)$ is
\textbf{lagrangian}\index{lagrangian submanifold ! definition}
(respectively, \textbf{isotropic}\index{isotropic submanifold}
and \textbf{coisotropic}\index{coisotropic submanifold})
if, at each $p \in X$, the space $T_pX$ is a lagrangian
(respectively, isotropic and coisotropic) subspace of $(T_pM,\omega_p)$.
\end{definition}

If $i: X \hookrightarrow M$ is the inclusion map,
then $X$ is a \textbf{lagrangian submanifold}
if and only if $i^*\omega = 0$ and $\dim X = \frac {1}{2} \dim M$.

The problem of embedding\footnote{An \textbf{embedding}\index{embedding}
is an immersion that is a homeomorphism onto its image.}
a compact manifold as a lagrangian submanifold of a given
symplectic manifold is often global.
For instance, Gromov~\cite{gr:pseudo} proved that
there can be no lagrangian spheres in $(\CC^n,\omega_0)$,
except for the circle in $\CC^2$, and more generally
no compact \textbf{exact lagrangian} submanifolds, in the sense that
$\alpha_0 =\sum y_j \ dx_j$ restricts to an exact 1-form.
The argument uses {\em pseudoholomorphic curves}
(Section~\ref{sec:pseudoholomorphic}).
Yet there are {\em immersed} lagrangian spheres
(Section~\ref{sec:special_lagrangians}).
More recently were found topological and geometrical constraints
on manifolds that admit lagrangian embeddings into {\em compact}
symplectic manifolds;
see for instance~\cite{bi:intersections,bi-ci:subcritical,se:graded}.

%%%%%%%%%%%%%%%%%%%%%%%%%%%%%%%%%%%%%%%%%%%%%%%%%%%%%%%%%%%%%%%%%%%%%%%%%%%%%

\begin{examples}
\begin{enumerate}
\item
Any 1-dimensional submanifold of a symplectic surface is lagrangian
(because a 1-dimensional subspace of a symplectic vector space
is always isotropic).

Therefore, any product of $n$ embedded curves arises as a
lagrangian submanifold of (a neighborhood of zero in)
the prototype $( \RR ^{2n},\omega_0)$.
In particular, a \textbf{torus}\index{torus}
$\TT^n = S^1 \times \ldots \times S^1$
can be embedded as a lagrangian submanifold of any $2n$-dimensional
symplectic manifold, by Darboux's theorem (Theorem~\ref{thm:darboux}).

\item
Let $M = T^*X$ be the cotangent bundle of a manifold $X$.
With respect to a cotangent coordinate chart
$(T^*U, x_1,\dots,x_n,\xi_1,\dots,\xi_n)$, the tautological
form is $\alpha = \sum \xi_idx_i$ and
the canonical form is $\omega = -d\alpha = \sum dx_i \wedge d\xi_i$.

The \textbf{zero section}\index{lagrangian submanifold !
zero section}\index{cotangent bundle !
zero section}\index{example ! of lagrangian submanifold}
$X_0 := \{(x,\xi) \in T^*X \mid \xi = 0 \mbox{ in } T_x^*X\}$
is an $n$-dimensional submanifold of $T^*X$ whose intersection with $T^*U$
is given by the equations $\xi_1 = \dots = \xi_n = 0$.
Clearly $\alpha$ vanishes on $X_0 \cap T^*U$.
Hence, if $i_0: X_0 \hookrightarrow T^*X$ is the inclusion map,
we have $i_0^*\omega = i_0^*d\alpha = 0$, and so $X_0$ is lagrangian.

A \textbf{cotangent fiber}\index{lagrangian submanifold !
zero section}\index{cotangent bundle !
zero section}\index{example ! of lagrangian submanifold}
$T_{x_0}^*X$ is an $n$-dimensional submanifold of $T^*X$
given by the equations $x_i = (x_0)_i$, $i = 1, \dots , n$, on $T^*U$.
Since the $x_i$'s are constant, the form
$\alpha$ vanishes identically,
% on $T_{x_0}^*X$,
and $T_{x_0}^*X$ is a lagrangian submanifold.
\end{enumerate}
\end{examples}

%%%%%%%%%%%%%%%%%%%%%%%%%%%%%%%%%%%%%%%%%%%%%%%%%%%%%%%%%%%%%%%%%%%%%%%%%%%%%

Let $X_\mu $ be (the image of) an arbitrary section, that is,
an $n$-dimensional submanifold of $T^*X$ of the form
$X_\mu  = \{(x,\mu_x) \mid x \in X,\ \mu_x \in T_x^*X\}$,
where the covector $\mu_x$ depends smoothly on $x$,
so $\mu: X \rightarrow T^*X$ is a de Rham 1-form.
We will investigate when such an $X_\mu$ is lagrangian.
Relative to the inclusion $i: X_\mu  \hookrightarrow T^*X$ and the
cotangent projection $\pi: T^*X \rightarrow X$,
these $X_\mu$'s are exactly the submanifolds for which
$\pi \circ i: X_\mu  \rightarrow X$ is a diffeomorphism.

\begin{proposition}
\label{prop:tautological_property}
The tautological 1-form $\alpha$ on $T^*X$ satisfies
$\mu^* \alpha = \mu$, for any 1-form $\mu: X \to T^*X$.
\end{proposition}

\vspace*{-2ex}

\begin{proof}
Denote by $s_{\mu}: X \rightarrow T^*X$, $x \mapsto (x,\mu_x)$,
the 1-form $\mu$ regarded exclusively as a map.
From the definition, $\alpha_p = (d\pi_p)^*\xi$ at $p = (x,\xi) \in M$.
For $p = s_{\mu}(x) = (x,\mu_x)$, we have $\alpha_p = (d\pi_p)^*\mu_x$.
Then, since $\pi \circ s_{\mu} = \mathrm{id}_X$, we have
\[
        (s_{\mu}^*\alpha)_x = (ds_{\mu})_x^* \alpha_p
        = (ds_{\mu})_x^*(d\pi_p)^*\mu_x
        = (d (\pi \circ s_{\mu}) )_x^* \mu_x = \mu_x \ .
\]
\end{proof}

%%%%%%%%%%%%%%%%%%%%%%%%%%%%%%%%%%%%%%%%%%%%%%%%%%%%%%%%%%%%%%%%%%%%%%%%%%%%%

The map $s_{\mu}: X \rightarrow T^*X$, $s_{\mu} (x) = (x,\mu_x)$
is an embedding with image the section $X_\mu$.
The diffeomorphism $\tau : X \to X_\mu$,
$\tau (x) := (x, \mu_x)$, satisfies $i \circ \tau = s_\mu$.

\begin{proposition}
\label{prop:closed_1_forms}
The sections of $T^*X$ that are lagrangian
are those corresponding to closed 1-forms on $X$.\index{lagrangian
submanifold ! closed 1-form}\index{cotangent bundle !
lagrangian submanifold}\index{lagrangian submanifold ! of TX@of $T^*X$}
\end{proposition}

\vspace*{-2ex}

\begin{proof}
Using the previous notation, the
condition of $X_\mu$ being lagrangian becomes:
$i^*d\alpha = 0 \Leftrightarrow \tau^* i^* d\alpha = 0
%\Leftrightarrow (i \circ \tau)^* d\alpha = 0
\Leftrightarrow s_{\mu}^* d\alpha = 0
\Leftrightarrow d(s_{\mu}^*\alpha) = 0
\Leftrightarrow d\mu = 0$.
\end{proof}

%%%%%%%%%%%%%%%%%%%%%%%%%%%%%%%%%%%%%%%%%%%%%%%%%%%%%%%%%%%%%%%%%%%%%%%%%%%%%

When $\mu = dh$ for some $h \in C^{\infty}(X)$,
such a primitive $h$ is called a
\textbf{generating function}\index{generating function}
for the lagrangian submanifold\index{lagrangian submanifold !
generating function} $X_\mu $.
Two functions generate the same lagrangian submanifold
if and only if they differ by a locally constant function.
When $X$ is simply connected,
or at least $H_{\mathrm{de Rham}}^1(X) = 0$,
every lagrangian $X_\mu $ admits a generating function.

%%%%%%%%%%%%%%%%%%%%%%%%%%%%%%%%%%%%%%%%%%%%%%%%%%%%%%%%%%%%%%%%%%%%%%%%%%%%%

Besides the cotangent fibers,
there are lots of lagrangian submanifolds of
$T^*X$\index{cotangent bundle ! lagrangian submanifold}
not covered by the description in terms of closed 1-forms.
Let $S$ be any submanifold of an
$n$-dimensional manifold $X$.
The \textbf{conormal space}\index{conormal ! space} of $S$ at $x \in S$ is
\[
        N_x^*S = \{\xi \in T_x^*X \mid \xi(v) = 0 \ ,
        \mbox{ for all } v \in T_xS\} \ .
\]
The \textbf{conormal bundle}\index{conormal ! bundle} of $S$ is
$N^*S = \{(x,\xi) \in T^*X \mid x \in S,\ \xi \in N_x^*S\}$.
This is an $n$-dimensional submanifold of $T^*X$.
In particular, taking $S = \{x\}$ to be one point,
the conormal bundle is the corresponding cotangent fiber $T_x^*X$.
Taking $S = X$, the conormal bundle is the zero section $X_0$ of $T^*X$.

\begin{proposition}
If $i: N^*S \hookrightarrow T^*X$ is the inclusion
of the conormal bundle of a submanifold $S \subset X$,
and $\alpha$ is the tautological 1-form on $T^*X$, then $i^*\alpha = 0$.
\end{proposition}

\vspace*{-2ex}

\begin{proof}
Let $(\cU,x_1,\dots,x_n)$ be a coordinate chart on $X$
adapted to $S$, so that $\cU \cap S$ is
described by $x_{k+1} = \dots = x_n = 0$.
Let $(T^*\cU,x_1,\dots,x_n,\xi_1,\dots,\xi_n)$ be the associated
cotangent coordinate chart.
The submanifold $N^*S \cap T^*\cU$ is described by
$x_{k+1} = \dots = x_n = 0$ and $\xi_1 = \dots = \xi_k = 0$.
Since $\alpha = \sum \xi_i dx_i$ on $T^*\cU$,
we conclude that, at $p \in N^*S$,
\[
        (i^*\alpha)_p = \alpha_p|_{T_p(N^*S)} =
        \left. \sum \limits_{i>k} \xi_i dx_i
        \right|_{\mathrm{span} \{ \frac{\partial}{\partial x_i},
        i \leq k \}} = 0 \ .
\]
\end{proof}

\vspace*{-2ex}

\begin{corollary}
\index{lagrangian submanifold ! conormal bundle}\index{cotangent
bundle ! lagrangian submanifold}\index{cotangent bundle !
conormal bundle}
For any submanifold $S$ of $X$, the conormal bundle $N^*S$
is a lagrangian submanifold of $T^*X$.
\end{corollary}

%%%%%%%%%%%%%%%%%%%%%%%%%%%%%%%%%%%%%%%%%%%%%%%%%%%%%%%%%%%%%%%%%%%%%%%%%%%%%
%%%%%%%%%%%%%%%%%%%%%%%%%%%%%%%%%%%%%%%%%%%%%%%%%%%%%%%%%%%%%%%%%%%%%%%%%%%%%

\ssubsection{Lagrangian Neighborhood Theorem}
\label{weinstein_lagrangian_theorem}

Weinstein~\cite{we:lagrangian} proved that,
if a compact submanifold $X$ is lagrangian
with respect to two symplectic forms $\omega_0$ and $\omega_1$,
then the conclusion of the Moser relative theorem
(Theorem~\ref{thm:moser_relative}) still holds.
We need some algebra for the Weinstein theorem\index{Weinstein !
lagrangian neighborhood theorem}\index{theorem ! Weinstein
lagrangian neighborhood}\index{neighborhood ! Weinstein
lagrangian neighborhood}.

%%%%%%%%%%%%%%%%%%%%%%%%%%%%%%%%%%%%%%%%%%%%%%%%%%%%%%%%%%%%%%%%%%%%%%%%%%%%%

Suppose that $U,W$ are $n$-dimensional vector spaces,
and $\Omega: U \times W \to \RR$ is a bilinear pairing;
the map $\Omega$ gives rise to a linear map $\widetilde\Omega:U\to W^*$,
$\widetilde\Omega(u) = \Omega(u,\cdot)$.
Then $\Omega$ is nondegenerate if and only if
$\widetilde\Omega$ is bijective.

\begin{proposition}
\label{prop:complement}
Let $(V,\Omega)$ be a symplectic vector space,
$U$ a lagrangian subspace\index{lagrangian
subspace}\index{subspace ! lagrangian} of $(V,\Omega)$,
and $W$ any vector space complement to $U$, not necessarily lagrangian.
Then from $W$ we can canonically build a
lagrangian complement\index{lagrangian complement} to $U$.
\end{proposition}

\vspace*{-2ex}

\begin{proof}
From $\Omega$ we get a nondegenerate pairing
$\Omega' : U \times W \to \RR$,
so $\widetilde\Omega':U \to W^*$ is bijective.
We look for a lagrangian complement to $U$ of the form
$W'=\{w+Aw\mid w\in W\}$ for some linear map $A:W \to U$.
For $W'$ to be lagrangian we need that $\Omega( w_1,w_2) =
\widetilde\Omega'(Aw_2)(w_1) -\widetilde\Omega'(Aw_1)(w_2)$.
Let $A'=\widetilde\Omega'\circ A$, and look for $A'$ such that
$\Omega( w_1,w_2)=A'(w_2)(w_1)-A'(w_1)(w_2)$ for all $w_1,w_2\in W$.
The canonical choice is $A'(w)=-\frac 12\Omega(w,\cdot)$.
Set $A=(\widetilde\Omega')^{-1}\circ A'$.
\end{proof}

\vspace*{-2ex}

\begin{proposition}
\label{prop:canonical_iso}
Let $V$ be a vector space,
let $\Omega_0$ and $\Omega_1$ be symplectic forms on $V$,
let $U$ be a subspace of $V$ lagrangian for $\Omega_0$ and $\Omega_1$,
and let $W$ be any complement to $U$ in $V$.
Then from $W$ we can canonically construct a linear isomorphism
$L:V\stackrel{\simeq} \longrightarrow V$ such that $L|_U= \Id_U$ and
$L^*\Omega_1=\Omega_0$.
\end{proposition}

\vspace*{-2ex}

\begin{proof}
By Proposition~\ref{prop:complement},
from $W$ we canonically obtain complements $W_0$ and $W_1$ to
$U$ in $V$ such that $W_0$ is lagrangian for $\Omega_0$ and $W_1$ is
lagrangian for $\Omega_1$.
The nondegenerate bilinear pairings
$\Omega_i : W_i \times U \to \RR$, $i=0,1$, give isomorphisms
$\widetilde\Omega_i:W_i \stackrel{\simeq}\longrightarrow U^*$, $i=0,1$,
respectively.
Let $B : W_0 \to W_1$ be the linear map satisfying
$\widetilde\Omega_1 \circ B = \widetilde\Omega_0$,
i.e., $\Omega_0(w_0,u)=\Omega_1(Bw_0,u)$,
$\forall w_0 \in W_0$, $\forall u \in U$.
Let $L := \Id_U\oplus B: U\oplus W_0 \to U\oplus W_1$
be the extension of $B$ to the rest of $V$ by setting
it to be the identity on $U$.
It satisfies:
\[
\begin{array}{rcl}
        (L^*\Omega_1)(u\oplus w_0,u'\oplus w'_0) & = &
        \Omega_1(u\oplus B w_0,u'\oplus B w'_0) \\
        & = & \Omega_1(u,B w'_0) + \Omega_1(B w_0,u') \\
        & = & \Omega_0(u, w'_0) + \Omega_0(w_0,u') \\
        & = & \Omega_0(u\oplus w_0, u'\oplus w'_0)\ .
\end{array}
\]
\end{proof}

\vspace*{-2ex}

%%%%%%%%%%%%%%%%%%%%%%%%%%%%%%%%%%%%%%%%%%%%%%%%%%%%%%%%%%%%%%%%%%%%%%%%%%%%%

\begin{theorem}
\label{thm:weinstein_lagrangian}\index{Weinstein !
lagrangian neighborhood theorem}\index{theorem ! Weinstein
lagrangian neighborhood}\index{neighborhood ! Weinstein
lagrangian neighborhood}
\textbf{(Weinstein Lagrangian Neighborhood Theorem)} $\;$
Let $M$ be a $2n$-dimensional manifold,
$X$ a compact $n$-dimensional submanifold,
$i: X \hookrightarrow M$ the inclusion map, and
$\omega_0$ and $\omega_1$ symplectic forms on $M$
such that $i^*\omega_0=i^*\omega_1=0$, i.e., $X$ is a lagrangian
submanifold of both $(M,\omega_0)$ and $(M,\omega_1)$.
Then there exist neighborhoods $\cU_0$ and $\cU_1$ of $X$ in $M$
and a diffeomorphism $\varphi: \cU_0 \to \cU_1$ such that
$\varphi^*\omega_1=\omega_0$ and $\varphi$ is the identity on $X$,
i.e., $\varphi (p) = p$, $\forall p\in X$.
\end{theorem}

\vspace*{-2ex}

\begin{proof}
Put a riemannian metric $g$ on $M$.
Fix $p\in X$, and let $V=T_pM$, $U=T_pX$ and $W=U^{\perp}$,
the orthocomplement of $U$ in $V$ relative to
the inner product $g_p(\cdot,\cdot)$.
Since $i^*\omega_0=i^*\omega_1=0$, the subspace $U$ is lagrangian
for both $(V,\omega_0|_p)$ and $(V,\omega_1|_p)$.
By Proposition~\ref{prop:canonical_iso},
we canonically get from $U^{\perp}$
a linear isomorphism $L_p:T_pM\to T_pM$ depending smoothly on $p$,
such that $L_p|_{T_pX}= \Id_{T_pX}$ and $L^*_p\omega_1|_p=\omega_0|_p$.
By the Whitney extension theorem (Section~\ref{moser_relative_theorem}),
there exist a neighborhood $\cN$ of $X$
and an embedding $h: \cN \hookrightarrow M$ with $h|_X= \id_X$ and
$dh_p=L_p$ for $p\in X$.
Hence, at any $p\in X$, we have $(h^*\omega_1)_p = (dh_p)^* \omega_1|_p
= L^*_p \omega_1|_p = \omega_0|_p$.
Applying the Moser relative theorem (Theorem~\ref{thm:moser_relative})
to $\omega_0$ and $h^*\omega_1$,
we find a neighborhood $\cU_0$ of $X$ and an embedding
$f : \cU_0 \to \cN$ such that $f|_X= \id _X$ and
$f^*(h^*\omega_1)=\omega_0$ on $\cU_o$.
Set $\varphi = h \circ f$ and $\cU_1 = \varphi (\cU_0)$.
\end{proof}

Theorem~\ref{thm:weinstein_lagrangian} has
the following generalization.
For a proof see, for instance, either
of~\cite{go:coisotropic,gu-st:techniques,we:isotropic}.

\begin{theorem}
\label{thm:coisotropic}
\index{theorem ! coisotropic embedding}\index{coisotropic !
embedding}\index{embedding ! coisotropic}
\textbf{(Coisotropic Embedding Theorem)} $\;$
Let $M$ be a manifold of dimension $2n$,
$X$ a submanifold of dimension $k \geq n$,
$i: X \hookrightarrow M$ the inclusion, and
$\omega_0$ and $\omega_1$ symplectic forms on $M$,
such that $i^*\omega_0=i^*\omega_1$ and $X$ is
coisotropic for both $(M,\omega_0)$ and $(M,\omega_1)$.
Then there exist neighborhoods $\cU_0$ and $\cU_1$ of $X$ in $M$
and a diffeomorphism $\varphi: \cU_0 \to \cU_1$ such that
$\varphi^*\omega_1=\omega_0$ and $\varphi |_X = \id_X$.
\end{theorem}

%%%%%%%%%%%%%%%%%%%%%%%%%%%%%%%%%%%%%%%%%%%%%%%%%%%%%%%%%%%%%%%%%%%%%%%%%%%%%
%%%%%%%%%%%%%%%%%%%%%%%%%%%%%%%%%%%%%%%%%%%%%%%%%%%%%%%%%%%%%%%%%%%%%%%%%%%%%

\ssubsection{Weinstein Tubular Neighborhood Theorem}
\label{weinstein_tubular_theorem}
\index{Weinstein ! tubular neighborhood theorem}\index{theorem !
Weinstein tubular neighborhood}\index{neighborhood !
Weinstein tubular neighborhood}

Let $(V,\Omega)$ be a symplectic linear space,
and let $U$ be a lagrangian subspace.\index{symplectic ! linear algebra}
Then there is a canonical nondegenerate bilinear pairing
$\Omega': V/U \times U \to \RR$ defined by
$\Omega'([v],u) = \Omega(v,u)$ where $[v]$ is the
equivalence class of $v$ in $V/U$.
Consequently, we get a canonical isomorphism
${\widetilde \Omega}': V/U \to U^*$,
${\widetilde \Omega}'([v]) = \Omega'([v],\cdot)$.

In particular, if $(M,\omega)$ is a symplectic manifold,
and $X$ is a lagrangian submanifold, then
$T_p X$ is a lagrangian subspace of $(T_p M,\omega_p )$ for each $p \in X$
and there is a canonical identification of the
\textbf{normal space}\index{normal ! space}\index{space ! normal}
of $X$ at $p$, $N_p X := T_p M/T_p X$, with the cotangent fiber $T_p ^*X$.
Consequently the normal bundle $NX$ and the cotangent bundle $T^*X$
are canonically identified.

%%%%%%%%%%%%%%%%%%%%%%%%%%%%%%%%%%%%%%%%%%%%%%%%%%%%%%%%%%%%%%%%%%%%%%%%%%%%%

\begin{theorem}
\index{Weinstein ! tubular neighborhood theorem}\index{theorem !
Weinstein tubular neighborhood}\index{neighborhood !
Weinstein tubular neighborhood}\index{tubular neighborhood !
Weinstein theorem}\index{Weinstein ! lagrangian
embedding}\index{embedding ! lagrangian}\label{thm:weinstein_tubular}
\textbf{(Weinstein Tubular Neighborhood Theorem)} $\;$
Let $(M,\omega)$ be a symplectic manifold,
$X$ a compact lagrangian submanifold,
$\omega_0$ the canonical symplectic form on $T^*X$,
$i_0: X \hookrightarrow T^*X$ the lagrangian embedding
as the zero section, and
$i: X \hookrightarrow M$ the lagrangian embedding
given by inclusion.
Then there are neighborhoods $\cU_0$ of $X$ in $T^*X$,
$\cU$ of $X$ in $M$, and a diffeomorphism
$\varphi: \cU_0 \to \cU$ such that
$\varphi^*\omega=\omega_0$ and $\varphi \circ i_0 = i$.
\end{theorem}

\vspace*{-2ex}

\begin{proof}
By the standard tubular neighborhood\index{tubular
neighborhood ! theorem}\index{theorem !
tubular neighborhood}\index{neighborhood ! tubular neighborhood}
theorem\footnote{\textbf{Tubular Neighborhood Theorem:}
{\em Let $M$ be a manifold, $X$ a submanifold,
$NX$ the normal bundle of $X$ in $M$,
$i_0: X \hookrightarrow NX$ the zero section, and
$i: X \hookrightarrow M$ the inclusion.
Then there are neighborhoods $\cU_0$ of $X$ in $NX$,
$\cU$ of $X$ in $M$ and a diffeomorphism
$\psi: \cU_0 \to \cU$ such that $\psi \circ i_0 = i$.}
This theorem can be proved with the exponential map
using a riemannian metric;
see for instance~\cite{spivak:comprehensive}.}
and since $NX \simeq T^*X$ are canonically identified,
we can find a neighborhood
$\cN_0$ of $X$ in $T^*X$, a neighborhood $\cN$ of $X$ in $M$,
and a diffeomorphism $\psi: \cN_0 \to \cN$ such that $\psi \circ i_0 = i$.
Let $\omega_0$ be the canonical form on $T^*X$ and
$\omega_1 = \psi^*\omega$.
The submanifold $X$ is lagrangian for both of these
symplectic forms on $\cN_0$.
By the Weinstein lagrangian neighborhood theorem
(Theorem~\ref{thm:weinstein_lagrangian}),
there exist neighborhoods $\cU_0$ and $\cU_1$ of $X$ in $\cN_0$
and a diffeomorphism $\theta: \cU_0 \to \cU_1$ such that
$\theta^*\omega_1 = \omega_0$ and $\theta \circ i_0 = i_0$.
Take $\varphi = \psi \circ \theta$ and $\cU = \varphi(\cU_0)$.
Then $\varphi^*\omega = \theta^* \psi^*\omega = \theta^* \omega_1 = \omega_0$.
\end{proof}

Theorem~\ref{thm:weinstein_tubular} classifies
compact lagrangian embeddings:
up to local symplectomorphism, the set of lagrangian embeddings
is the set of embeddings of manifolds into their cotangent
bundles as zero sections.

The classification of compact {\em isotropic} embeddings is also
due to Weinstein in~\cite{we:lectures,we:isotropic}\index{Weinstein !
isotropic embedding}\index{embedding ! isotropic}\index{isotropic ! embedding}.
An \textbf{isotropic embedding} of a manifold $X$
into a symplectic manifold $(M,\omega)$ is a closed embedding
$i: X \hookrightarrow M$ such that $i^*\omega = 0$.
Weinstein showed that neighborhood equivalence of isotropic
embeddings is in one-to-one correspondence with isomorphism
classes of symplectic vector bundles.

The classification of compact {\em coisotropic} embeddings is due to
Gotay~\cite{go:coisotropic}.\index{embedding !
coisotropic}\index{coisotropic ! embedding}\index{embedding !
coisotropic}\index{Gotay ! coisotropic embedding}
A \textbf{coiso\-tro\-pic embedding} of a manifold $X$ carrying a closed 2-form
$\alpha$ of constant rank into a symplectic manifold $(M,\omega)$
is an embedding $i: X \hookrightarrow M$ such that $i^*\omega = \alpha$
and $i(X)$ is coisotropic as a submanifold of $M$.
Let $E$ be the \textbf{characteristic distribution}\index{characteristic
distribution} of a closed form $\alpha$ of constant rank on $X$,
i.e., $E_p$ is the kernel of $\alpha_p$ at $p \in X$.
Gotay showed that then the total space $E^*$ carries a symplectic
structure in a neighborhood of the zero section, such that
$X$ embeds  coisotropically onto this zero section and, moreover,
every coisotropic embedding is equivalent to this
in some neighborhood of the zero section.

%%%%%%%%%%%%%%%%%%%%%%%%%%%%%%%%%%%%%%%%%%%%%%%%%%%%%%%%%%%%%%%%%%%%%%%%%%%%%
%%%%%%%%%%%%%%%%%%%%%%%%%%%%%%%%%%%%%%%%%%%%%%%%%%%%%%%%%%%%%%%%%%%%%%%%%%%%%

\ssubsection{Application to Symplectomorphisms}
\label{application_symplectomorphisms}
\index{symplectomorphism ! group of symplectomorphisms}\index{group !
of symplectomorphisms}\index{symplectomorphism ! vs.\ lagrangian submanifold}

Let $(M_1,\omega_1)$ and $(M_2,\omega_2)$ be two
$2n$-dimensional symplectic manifolds.
Given a diffeomorphism $f: M_1 \stackrel{\simeq}{\longrightarrow}
M_2$, there is a way to express the condition of
$f$ being a symplectomorphism in terms of a certain
submanifold being lagrangian.
Consider the two projection maps ${\mathrm{pr}}_i : M_1 \times M_2 \to M_i$,
$(p_1,p_2) \mapsto p_i$, $i=1,2$.
The \textbf{twisted product form}\index{twisted product form}
on $M_1 \times M_2$ is the symplectic\footnote{More generally,
$\lambda_1({\mathrm{pr}}_1)^*\omega_1 +
\lambda_2({\mathrm{pr}}_2)^*\omega_2$ is
symplectic for all $\lambda_1,\lambda_2 \in \RR{\backslash}\{0\}$.} form
\[
        {\widetilde \omega} =
        ({\mathrm{pr}}_1)^*\omega_1 - ({\mathrm{pr}}_2)^*\omega_2 \ .
\]

\begin{proposition}
\index{symplectomorphism ! vs.\ lagrangian submanifold}\index{lagrangian
submanifold ! vs.\ symplectomorphism}\index{theorem ! symplectomorphism
vs.\ lagrangian submanifold}
A diffeomorphism $f: M_1 \stackrel{\simeq}{\longrightarrow} M_2$
is a symplectomorphism
if and only if the graph of $f$ is a lagrangian submanifold of
$(M_1 \times M_2,{\widetilde \omega})$.
\end{proposition}

\vspace*{-2ex}

\begin{proof}
The graph of $f$ is the $2n$-dimensional submanifold
${\mathrm{Graph}} \, f
= \{(p,f(p)) \mid p \in M_1\} \subseteq M_1 \times M_2$,
which is the image of the embedding
$\gamma: M_1 \to M_1 \times M_2$, $p \mapsto (p,f(p))$.
We have $\gamma^*{\widetilde \omega}
= \gamma^* {\mathrm{pr}}_1^* \ \omega_1 -
\gamma^* {\mathrm{pr}}_2^* \ \omega_2
= ({\mathrm{pr}}_1 \circ \gamma)^*\omega_1 -
({\mathrm{pr}}_2 \circ \gamma)^*\omega_2$,
and ${\mathrm{pr}}_1 \circ \gamma$ is the identity map on $M_1$
whereas ${\mathrm{pr}}_2 \circ \gamma = f$.
So ${\mathrm{Graph}} \, f$ is lagrangian, i.e.,
$\gamma^*{\widetilde \omega} = 0$, if and only if
$f^*\omega_2 = \omega_1$, i.e., $f$ is a symplectomorphism.
\end{proof}

Lagrangian submanifolds of $(M_1 \times M_2,{\widetilde \omega})$
are called \textbf{canonical relations}, when viewed as
morphisms between $(M_1,\omega_1)$ and $(M_2,\omega_2)$,
even if $\dim M_1 \neq \dim M_2$.
Under a reasonable assumption, there is a notion of
composition~\cite{we:lectures}.

%%%%%%%%%%%%%%%%%%%%%%%%%%%%%%%%%%%%%%%%%%%%%%%%%%%%%%%%%%%%%%%%%%%%%%%%%%%%%

Take $M_1 = M_2 = M$ and
suppose that $(M,\omega)$ is a {\em compact} symplectic manifold
and $f \in \mathrm{Sympl}(M,\omega)$.
The graphs $\mathrm{Graph} \, f$ and $\Delta$, of $f$ and of the
identity map $\id : M \to M$, are lagrangian submanifolds of
$M \times M$ with ${\widetilde \omega} =
\mathrm{pr}_1^*\omega - \mathrm{pr}_2^*\omega$.
By the Weinstein tubular neighborhood theorem, there exist
a neighborhood $\cU$ of $\Delta$ in $(M \times M,{\widetilde \omega})$
and a neighborhood $\cU_0$ of $M$ in $(T^*M,\omega_0)$
with a symplectomorphism $\varphi: \cU \to \cU_0$
satisfying $\varphi (p,p) = (p,0)$, $\forall p \in M$.

Suppose that $f$ is sufficiently
\textbf{\boldmath{$C^1$}-close}\index{C-topology@$C^1$-topology}\footnote{Let
$X$ and $Y$ be manifolds.
A sequence of maps $f_i: X \to Y$ \textbf{converges in the
\boldmath{$C^0$}-topology} (a.k.a.\ the \textbf{compact-open topology})
to $f: X \to Y$ if and only if $f_i$ converges uniformly on compact sets.
A sequence of $C^1$ maps $f_i: X \to Y$
\textbf{converges in the
\boldmath{$C^1$}-topology}\index{C-topology@$C^1$-topology}
to $f: X \to Y$ if and only if it and
the sequence of derivatives $df_i: TX \to TY$ converge uniformly on
compact sets.} to $\id$,
i.e., $f$ is in some sufficiently small neighborhood of
the identity $\id$ in the $C^1$-topology\index{C-topology@$C^1$-topology}.
Hence we can assume that $\mathrm{Graph} \, f \subseteq \cU$.
Let $j: M \hookrightarrow \cU$, $j(p)=(p,f(p))$,
be the embedding as $\mathrm{Graph} \, f$,
and $i: M \hookrightarrow \cU$, $i(p)=(p,p)$,
be the embedding as $\Delta = \mathrm{Graph} \, \id$.
The map $j$ is sufficiently $C^1$-close to $i$.
These maps induce embeddings
$\varphi \circ j = j_0: M \hookrightarrow \cU_0$
and $\varphi \circ i = i_0: M \hookrightarrow \cU_0$ as 0-section,
respectively.
Since the map $j_0$ is sufficiently $C^1$-close to $i_0$,
the image set $j_0 (M)$ intersects each fiber $T_p^*M$
at one point $\mu_p$ depending smoothly on $p$.
Therefore, the image of $j_0$ is the image of a smooth section
$\mu: M \to T^*M$, that is, a 1-form $\mu = j_0 \circ (\pi \circ j_0)^{-1}$.
We conclude that $\mathrm{Graph} \, f \simeq
\{(p,\mu_p) \ | \ p \in M ,\ \mu_p \in T_p^*M \}$.
Conversely, if $\mu$ is a 1-form sufficiently $C^1$-close
to the zero 1-form, then
$\{(p,\mu_p) \ | \ p \in M,\ \mu_p \in T^*_pM \} \simeq \mathrm{Graph} \, f$,
for some diffeomorphism $f: M \to M$.

By Proposition~\ref{prop:closed_1_forms},
$\mathrm{Graph} \, f$ is lagrangian if and only if $\mu$ is closed.
A small $C^1$-neighborhood of $\id$ in $\mathrm{Sympl}(M,\omega)$
is thus homeomorphic to a $C^1$-neighborhood of zero
in the vector space of closed 1-forms on $M$.
So we obtain the model:
\[
        T_\id (\mathrm{Sympl}(M,\omega)) \simeq
        \{\mu \in \Omega^1(M) \ | \ d\mu = 0\}\ .
\]
In particular, $T_\id (\mathrm{Sympl}(M,\omega))$ contains
the space of exact 1-forms that correspond to generating functions,
%$\{\mu = dh \ | \ h \in C^\infty (M) \} \ \simeq \ 
$C^\infty (M) / \{\mbox{locally constant functions}\}$.

%%%%%%%%%%%%%%%%%%%%%%%%%%%%%%%%%%%%%%%%%%%%%%%%%%%%%%%%%%%%%%%%%%%%%%%%%%%%%

\begin{theorem}
\index{symplectomorphism ! fixed point}\index{fixed point}
Let $(M,\omega)$ be a compact symplectic manifold
(and not just one point) with $H^1_\mathrm{deRham}(M) = 0$.
Then any symplectomorphism of $M$ that is sufficiently $C^1$-close
to the identity has at least two fixed points.
\end{theorem}

\vspace*{-2ex}

\begin{proof}
If $f \in \mathrm{Sympl}(M,\omega)$ is sufficiently
$C^1$-close to $\id$, then its graph
corresponds to a closed 1-form $\mu$ on $M$.
As $H^1_{\mathrm{de Rham}}(M) = 0$,
we have that $\mu = dh$ for some $h \in C^{\infty}(M)$.
But $h$ must have at least two critical points
because $M$ is compact.
A point $p$ where $\mu_p = dh_p = 0$ corresponds to
a point in the intersection of the graph of $f$ with the diagonal,
%$\mathrm{Graph} \, f \cap \Delta$,
that is, a fixed point of $f$.
\end{proof}

This result has the following analogue in terms
of \textbf{lagrangian intersections}\index{lagrangian
submanifold ! intersection problem}\index{intersection of lagrangian
submanifolds}:
{\em if $X$ is a compact lagrangian submanifold of
a symplectic manifold $(M,\omega)$ with
$H^1_\mathrm{deRham}(X) = 0$,
then every lagrangian submanifold of $M$ that is
$C^1$-close\footnote{We say that a submanifold $Y$ of $M$ is
\textbf{\boldmath{$C^1$}-close} to another submanifold $X$
when there is a diffeomorphism $X \to Y$ that is, as a map into $M$,
$C^1$-close to the inclusion $X \hookrightarrow M$.}
to $X$ intersects $X$ in at least two points.}

%%%%%%%%%%%%%%%%%%%%%%%%%%%%%%%%%%%%%%%%%%%%%%%%%%%%%%%%%%%%%%%%%%%%%%%%%%%%%
%%%%%%%%%%%%%%%%%%%%%%%%%%%%%%%%%%%%%%%%%%%%%%%%%%%%%%%%%%%%%%%%%%%%%%%%%%%%%

\ssubsection{Generating Functions}
\label{generating_function}
\index{generating function}

We focus on symplectomorphisms between
the cotangent bundles $M_1 = T^*X_1$,
$M_2 = T^*X_2$\index{symplectomorphism ! recipe}\index{recipe !
for symplectomorphisms}\index{example ! of symplectomorphism}
of two $n$-dimensional manifolds $X_1$, $X_2$.
Let $\alpha_1,\alpha_2$ and $\omega_1,\omega_2$ be the
corresponding tautological and canonical forms.
Under the natural identification
\[
        M_1 \times M_2 = T^*X_1 \times T^*X_2 \simeq T^*(X_1 \times X_2) \ ,
\]
the tautological 1-form on $T^*(X_1 \times X_2)$ is
$\alpha = {\mathrm{pr}}_1^* \alpha_1 + {\mathrm{pr}}_2^*\alpha_2$,
the canonical 2-form on $T^*(X_1 \times X_2)$ is $\omega = -d\alpha =
{\mathrm{pr}}_1^*\omega_1 + {\mathrm{pr}}_2^*\omega_2$,
and the twisted product form\index{twisted product form}
is ${\widetilde \omega} =
{\mathrm{pr}}_1^*\omega_1 - {\mathrm{pr}}_2^*\omega_2$.
We define the involution $\sigma_2: M_2 \to M_2$,
$(x_2,\xi_2) \mapsto (x_2,-\xi_2)$,
which yields $\sigma_2^*\alpha_2 = -\alpha_2$.
Let $\sigma = {\mathrm{id}}_{M_1}
\times \sigma_2: M_1 \times M_2 \rightarrow M_1 \times M_2$.
Then $\sigma^*{\widetilde \omega} =
{\mathrm{pr}}_1^*\omega_1 + {\mathrm{pr}}_2^*\omega_2 = \omega$.
If $L$ is a lagrangian submanifold of $(M_1 \times M_2, \omega)$,
then its \textbf{twist} $L^{\sigma} := \sigma(L)$
is a lagrangian submanifold of
$(M_1 \times M_2,{\widetilde \omega})$.

For producing a symplectomorphism
$M_1 = T^*X_1 \rightarrow M_2 = T^*X_2$\index{symplectomorphism ! recipe}
we can start with a lagrangian submanifold $L$ of
$(M_1 \times M_2,\omega)$,
twist it to obtain a lagrangian submanifold
$L^{\sigma}$ of $(M_1 \times M_2,{\widetilde \omega})$,
and, if $L^{\sigma}$ happens to be the graph of some diffeomorphism
$\varphi: M_1 \rightarrow M_2$, then $\varphi$ is a symplectomorphism.

A method to obtain lagrangian submanifolds of
$M_1 \times M_2 \simeq T^*(X_1 \times X_2)$ relies
on generating functions.\index{generating function}
For any $f \in C^{\infty}(X_1 \times X_2)$,
$df$ is a closed 1-form on $X_1 \times X_2$.
The \textbf{lagrangian submanifold generated by $f$}\index{lagrangian
submanifold ! generating function} is
$L_{f} := \{((x,y),(df)_{(x,y)}) \mid (x,y) \in X_1 \times X_2\}$
(cf.\ Section~\ref{lagrangian_submanifolds}).
We adopt the loose notation
\[
\begin{array}{rcccl}
        d_x f & := & d_x f (x,y) & := &
        (df)_{(x,y)} \mbox{ projected to } T_x^*X_1 \times \{0\}, \\
        d_y f & := & d_y f (x,y) & := &
        (df)_{(x,y)} \mbox{ projected to } \{0\} \times T_y^*X_2 \ ,
\end{array}
\]
which enables us to write
$L_{f} = \{(x,y,d_x f,d_y f) \mid (x,y) \in X_1 \times X_2\}$ and
\[
        L_{f}^{\sigma} =
        \{(x,y,d_x f,-d_y f) \mid (x,y) \in X_1 \times X_2\} \ .
\]
When $L_{f}^{\sigma}$ is in fact the graph of a diffeomorphism
$\varphi: M_1=T^*X_1 \rightarrow M_2=T^*X_2$,
we call $\varphi$ the \textbf{symplectomorphism generated
by $f$},\index{symplectomorphism ! generating function}
and call $f$ the \textbf{generating function}\index{generating
function} of $\varphi$.
The issue now is to determine whether a given $L_{f}^{\sigma}$
is the graph of a diffeomorphism $\varphi: M_1 \rightarrow M_2$.
Let $(\cU_1,x_1,\dots,x_n),(\cU_2,y_1,\dots,y_n)$ be
coordinate charts for $X_1,X_2$, with
associated charts $(T^*\cU_1,x_1,\dots,x_n,\xi_1,\dots,\xi_n)$,
$(T^*\cU_2,y_1,\dots,y_n,\eta_1,\dots,\eta_n)$ for $M_1,M_2$.
The set $L_{f}^{\sigma}$
is the graph of $\varphi: M_1 \rightarrow M_2$ exactly when,
for any $(x,\xi) \in M_1$ and $(y,\eta) \in M_2$, we have
$\varphi(x,\xi) = (y,\eta) \Leftrightarrow
\xi = d_x f \mbox{ and } \eta = -d_y f$.
Therefore, given a point $(x,\xi) \in M_1$, to find its image
$(y,\eta) = \varphi(x,\xi)$ we must solve the
{\em Hamilton look-alike equations}\index{Hamilton equations}
\[
\left\{
\begin{array}{rll}
        \xi_i & = & \displaystyle{\phantom{-}
        \frac {\partial f}{\partial x_i} (x,y)} \\
        \eta_i &= &\displaystyle{-\frac {\partial f}{\partial y_i} (x,y)} \ .
\end{array} \right.
\]
If there is a solution $y = \varphi_1(x,\xi)$ of the first equation,
we may feed it to the second thus obtaining
$\eta = \varphi_2(x,\xi)$,
so that $\varphi(x,\xi) = (\varphi_1(x,\xi),\varphi_2(x,\xi))$. 
By the implicit function theorem\index{theorem !
implicit function}, in order to solve the first equation
locally and smoothly for $y$ in terms of $x$ and $\xi$, we need the condition
\[
        \det\left[ \frac {\partial}{\partial y_j} \left( \frac
        {\partial f}{\partial x_i}
        \right)\right]^n_{i,j=1} \ne 0 \ .
\]
This is a necessary condition for $f$ to generate a
symplectomorphism $\varphi$.
Locally this is also sufficient, but globally there
is the usual bijectivity issue.

\begin{example}
Let $X_1 = X_2 = \RR^n$,
and $f (x,y) = -\frac{|x-y|^2}{2}$, the square of euclidean
distance\index{euclidean ! distance} up to a constant.
In this case, the Hamilton equations\index{Hamilton equations} are
\[
\left\{
\begin{array}{rllll}
        \xi_i & = & \displaystyle{\phantom{-}
        \frac {\partial f}{\partial x_i}}
        & = & y_i - x_i \\
        \eta_i & = &\displaystyle{-\frac {\partial f}{\partial y_i}}
        & = & y_i - x_i
\end{array} \right. \qquad \Longleftrightarrow \qquad \left\{
\begin{array}{rll}
        y_i & = & x_i + \xi_i \\
        \\
        \eta_i & = & \xi_i \ .
\end{array} \right.
\]
The symplectomorphism generated by $f$ is
$\varphi (x,\xi) = (x+\xi , \xi)$.
If we use the euclidean inner product\index{euclidean !
inner product} to identify $T^* \RR^n$ with $T\RR^n$,
and hence regard $\varphi$ as $\widetilde \varphi : T\RR^n \to T\RR^n$
and interpret $\xi$ as the velocity vector, then the symplectomorphism
$\varphi$ corresponds to free translational motion in
euclidean space\index{euclidean ! space}.
\end{example}

%%%%%%%%%%%%%%%%%%%%%%%%%%%%%%%%%%%%%%%%%%%%%%%%%%%%%%%%%%%%%%%%%%%%%%%%%%%%%

The previous example can be generalized to the
{\em geodesic flow on a riemannian manifold}.\footnote{A
\textbf{riemannian metric}\index{riemannian ! metric}\index{metric}
on a manifold $X$ is a smooth function $g$ that assigns
to each point $x \in X$ an {\em inner product}\index{positive !
inner product} $g_x$ on $T_x X$, that is, a symmetric
positive-definite bilinear map $g_x: T_x X \times T_x X \to \RR$.
Smoothness means that
for every (smooth) vector field $v: X \to TX$ the real-valued
function $x \mapsto g_x (v_x, v_x)$ is smooth.
A \textbf{riemannian manifold}\index{riemannian ! manifold}
is a pair $(X,g)$ where $g$ is
a riemannian metric on the manifold $X$.
The \textbf{arc-length}\index{arc-length} of a piecewise
smooth curve $\gamma: [a,b] \to X$ on a riemannian $(X,g)$ is
$\int_a^b \left| \frac{d \gamma}{dt} \right| \, dt$,
where $\frac{d \gamma}{dt} (t) = d \gamma_t (1) \in T_{\gamma (t)} X$
and $\left| \frac{d \gamma}{dt} \right| =
\sqrt{g_{\gamma(t)} ( \frac{d \gamma}{dt}, \frac{d \gamma}{dt} )}$
is the \textbf{velocity} of $\gamma$.
A \textbf{reparametrization} of a curve $\gamma : [a,b] \to X$
is a curve of the form $\gamma \circ \tau : [c,d] \to X$
for some $\tau : [c,d] \to [a,b]$.
By the change of variable formula for the integral, we see that
the arc-length of $\gamma$ is invariant by reparametrization.
The \textbf{riemannian distance}\index{riemannian ! distance}
between two points $x$ and $y$ of a connected riemannian manifold $(X,g)$
is the infimum $d(x,y)$ of the set of all arc-lengths for piecewise smooth
curves joining $x$ to $y$.
A \textbf{geodesic}\index{geodesic ! curve}
is a curve that locally minimizes
distance and whose velocity is constant.
Given any curve $\gamma : [a,b] \to X$ with
${d \gamma \over dt}$ never vanishing, there is
a reparametrization $\gamma \circ \tau : [a,b] \to X$ of constant velocity.
A \textbf{minimizing geodesic}\index{geodesic ! minimizing}
from $x$ to $y$ is a geodesic joining $x$ to $y$ whose
arc-length is the
riemannian distance\index{riemannian ! distance} $d(x,y)$.
A riemannian manifold $(X,g)$ is
\textbf{geodesically convex}\index{geodesic ! geodesically convex}
if every point $x$ is joined to every other point $y$
by a unique (up to reparametrization) minimizing geodesic.
For instance,
$(\RR^n , \langle \cdot , \cdot \rangle)$ is a geodesically
convex riemannian manifold (where $g_x (v,w) = \langle v,w \rangle$
is the euclidean inner product on $T\RR^n \simeq \RR^n \times \RR^n$),
for which the riemannian distance is the
usual euclidean distance\index{euclidean ! distance} $d(x,y) = |x-y|$.}
Let $(X,g)$ be a geodesically convex riemannian manifold,
where $d(x,y)$ is the riemannian distance between points $x$ and $y$.
Consider the function
\[
        f: X \times X \longrightarrow \RR\ ,
        \qquad f (x,y) = - \frac{d(x,y)^2}{2} \ .
\]
We want to investigate if $f$ generates a symplectomorphism
$\varphi: T^*X \to T^*X$.
Using the identification
$\widetilde g _x : T_xX \stackrel{\simeq}{\longrightarrow} T_x^*X$,
$v \mapsto g_x(v,\cdot)$, induced by the metric,
we translate $\varphi$ into a map $\widetilde \varphi : TX \to TX$.
We need to solve 
\begin{eqnarray}
\label{eqn:geodesic}
\left\{
\begin{array}{lllll}
        \widetilde g_x(v) & = & \xi  & = & \phantom{-} d_x f (x,y) \\
        \widetilde g_y(w) & = & \eta & = & -d_y f (x,y)
\end{array} \right.
\end{eqnarray}
for $(y,\eta)$ in terms of $(x,\xi)$ in order to find $\varphi$, or,
equivalently,
for $(y,w)$ in terms $(x,v)$ in order to find $\widetilde \varphi$.
Assume that $(X,g)$ is \textbf{geodesically complete},
that is, every geodesic can be extended indefinitely.

\begin{proposition}
\label{prop:geodesic}
Under the identification $T_xX \simeq T_x^*X$ given by the metric,
the symplectomorphism generated by $f$ corresponds to the map
\[
\begin{array}{rrcl}
        \widetilde \varphi: & TX & \longrightarrow & TX \\
        & (x,v) & \longmapsto & (\gamma(1), \frac{d\gamma}{dt} (1)) \ ,
\end{array}
\]
where $\gamma$ is the geodesic with initial conditions
$\gamma(0) = x$ and $\frac{d\gamma}{dt} (0) = v$.
\end{proposition}

This map $\widetilde \varphi$
is called the \textbf{geodesic flow}\index{geodesic ! flow} on $(X,g)$.

\begin{proof}
Given $(x,v) \in TX$, let $\exp (x,v): \RR \to X$ be the unique
geodesic with initial conditions
$\exp (x,v) (0) = x$ and ${d \exp (x,v) \over dt} (0) = v$.
In this notation, we need to show that
the unique solution of the system of equations~(\ref{eqn:geodesic}) is
$\widetilde \varphi (x,v) = (\exp (x,v) (1) , d {\exp (x,v) \over dt} (1))$.

The Gauss lemma\index{Gauss lemma} in riemannian geometry
(see, for instance,~\cite{spivak:comprehensive})
asserts that geodesics are orthogonal to the level sets of the
distance function.
To solve the first equation for $y = \exp (x,u) (1)$ for some $u \in T_xX$,
evaluate both sides at $v$ and at vectors $v' \in T_x X$ orthogonal to $v$
\[
   |v|^2 = {d \over dt} \left[
   {-d(\exp (x,v)(t),y)^ 2 \over 2} \right]_{t=0}
   \quad \mbox{ and } \quad
   0 = {d \over dt} \left[
   {-d(\exp (x,v')(t),y)^ 2 \over 2} \right]_{t=0}
\]
to conclude that $u=v$, and thus $y = \exp (x,v) (1)$.

We have $-d_y f (x,y) (w') =0$ at vectors $w' \in T_y X$
orthogonal to $W:={d \exp (x,v) \over dt} (1)$,
because $f(x,y)$ is essentially the arc-length of a minimizing geodesic.
Hence $w=kW$ must be proportional to $W$,
and $k=1$ since
\[
   k |v|^2 = g_y (kW, W) = - {d \over dt} \left[
   {- d(x,\exp (x,v)(1-t))^ 2 \over 2} \right]_{t=0} = |v|^ 2 \ .
\]
\end{proof}

%%%%%%%%%%%%%%%%%%%%%%%%%%%%%%%%%%%%%%%%%%%%%%%%%%%%%%%%%%%%%%%%%%%%%%%%%%%%
%%%%%%%%%%%%%%%%%%%%%%%%%%%%%%%%%%%%%%%%%%%%%%%%%%%%%%%%%%%%%%%%%%%%%%%%%%%%

\ssubsection{Fixed Points}
\index{fixed_points}

Let $X$ be an $n$-dimensional manifold,
and $M = T^*X$ its cotangent bundle equipped with
the canonical symplectic form $\omega$.
Let $f: X \times X \to \RR$ be a smooth function
generating a symplectomorphism $\varphi: M \to M$,
$\varphi(x,d_x f) = (y, -d_y f)$,
with the notation of Section~\ref{generating_function}.
To describe the fixed points of $\varphi$\index{fixed point},
we introduce the function $\psi : X \to \RR$, $\psi(x) = f(x,x)$.

\begin{proposition}
\label{prop:fixed_vs_critical}
There is a one-to-one correspondence between the fixed points
of the symplectomorphism $\varphi$ and the critical points of $\psi$.
\end{proposition}

\vspace*{-2ex}

\begin{proof}
At $x_0 \in X$, $d_{x_0} \psi = ( d_x f + d_y f)|_{(x,y)=(x_0,x_0)}$.
Let $\xi = d_x f |_{(x,y)=(x_0,x_0)}$.
Recalling that $L_f^\sigma$
is the graph of $\varphi$, we have that
$x_0$ is a critical point of $\psi$, i.e.,
$d_{x_0} \psi =0$, if and only if $d_y f |_{(x,y)=(x_0,x_0)} = -\xi$,
which happens if and only if the point in $L_f^\sigma$
corresponding to $(x,y) = (x_0,x_0)$ is $(x_0,x_0,\xi,\xi)$,
i.e., $\varphi(x_0,\xi) = (x_0,\xi)$ is a fixed point.
\end{proof}

Consider the iterates
$\varphi^N = \varphi \circ \varphi \circ \ldots \circ \varphi$,
$N=1,2,\ldots$, given by $N$ successive applications of $\varphi$.
According to the previous proposition, if the symplectomorphism
$\varphi^N: M \to M$ is generated by some function $f^{(N)}$,
then there is a one-to-one correspondence
between the set of fixed points of $\varphi^N$
and the set of critical points of
$\psi^{(N)} : X \to \RR\ , \ \psi^{(N)} (x) = f^{(N)} (x,x)$.
It remains to know whether $\varphi^N$ admits a generating
function.\index{function ! generating}\index{generating function}
We will see that to a certain extent it does.

For each pair $x,y \in X$, define a map
$X \to \RR$, $z \mapsto f(x,z) + f (z,y)$.
Suppose that this map has a unique critical point $z_0$
and that $z_0$ is nondegenerate.
As $z_0$ is determined for each $(x,y)$ implicitly
by the equation $d_y f (x,z_0) + d_x f (z_0,y) =0$,
by nondegeneracy, the implicit function theorem
assures that $z_0 = z_0 (x,y)$ is a smooth function.
Hence, the function
\[
   f^{(2)}: X \times X \longrightarrow \RR\ , \quad
   f^{(2)} (x,y) := f(x,z_0) + f (z_0,y)
\]
is smooth.
Since $\varphi$ is generated by $f$, and $z_0$ is critical, we have
\[
\begin{array}{crclcl}
        & \varphi^2 (x,d_x f^{(2)} (x,y) )
        & = & \varphi ( \varphi (x, d_xf (x,z_0))
        & = & \varphi ( z_0, -d_y f (x,z_0)) \\
        = & \varphi ( z_0, d_x f (z_0,y) )
        & = & (y, -d_y f (z_0,y) )
        & = & (y, -d_y f^{(2)} (x,y) )\ .
\end{array}
\]
We conclude that the function $f^{(2)}$ is
a generating function for $\varphi^2$,
as long as, for each $\xi \in T^*_x X$, there is a unique
$y \in X$ for which $d_x f^{(2)} (x,y)$ equals $\xi$.

There are similar partial recipes for generating functions
of higher iterates.
In the case of $\varphi^3$, suppose that the function
$X \times X \to \RR$, $(z,u) \mapsto f(x,z) + f(z,u) + f (u,y)$,
has a unique critical point $(z_0,u_0)$ and that
it is a nondegenerate critical point.
A generating function would be
$f^{(3)} (x,y) = f(x,z_0) + f(z_0,u_0) + f (u_0,y)$.

When the generating functions $f$, $f^{(2)}$, $f^{(3)}$, \ldots , $f^{(N)}$
exist given by these formulas, the \textbf{$N$-periodic points}
of $\varphi$, i.e., the fixed points of $\varphi^N$,
are in one-to-one correspondence with the critical points of
\[
        (x_1, \ldots , x_N) \longmapsto
        f(x_1,x_2) + f(x_2,x_3) + \ldots +
        f (x_{N-1},x_N) + f (x_N,x_1) \ .
\]

%%%%%%%%%%%%%%%%%%%%%%%%%%%%%%%%%%%%%%%%%%%%%%%%%%%%%%%%%%%%%%%%%%%%%%%%%%%%%

\begin{example}
Let $\chi : \RR \to \RR^2$ be a smooth plane curve
that is 1-periodic, i.e., $\chi (s+1) = \chi (s)$, and parametrized
by arc-length, i.e., $\left| \frac{d \chi}{ds} \right| = 1$.
Assume that the region $Y$ enclosed by the image of $\chi$ is
\textbf{convex}, i.e., for any $s \in \RR$, the tangent line
$\{ \chi(s) + t \frac{d \chi}{ds} \mid t \in \RR \}$ intersects
the image $X := \partial Y$ of $\chi$ only at the point $\chi (s)$.

Suppose that a ball\index{billiards} is thrown into a billiard table
of shape $Y$ rolling with constant velocity and bouncing off
the boundary subject to the usual law of reflection.
The map describing successive points on the orbit of the ball is
\[
\begin{array}{rrcl}
        \varphi: & \RR / \ZZ \times (-1,1) & \longrightarrow &
        \RR / \ZZ \times (-1,1) \\
        & (x,v) & \longmapsto & (y,w) \ ,
\end{array}
\]
saying that when the ball bounces off $\chi(x)$
with angle $\theta = \arccos v$,
it will next collide with $\chi(y)$
and bounce off with angle $\nu = \arccos w$.
Then the function $f :  \RR / \ZZ \times  \RR / \ZZ \to \RR$
defined by $f (x,y) = -|\chi(x)-\chi(y)|$ is smooth off the diagonal,
and for $\varphi (x,v) = (y,w)$ satisfies
\[
\left\{ \begin{array}{lclcccc}
        \displaystyle{\frac{\partial f}{\partial x} (x,y)} & = &
        \displaystyle{\left. \frac{\chi(y)-\chi(x)}{|\chi(x)-\chi(y)|}
        \cdot \frac{d\chi}{ds} \right|_{s=x}}
        & = & \cos \theta & = & v \\
        \\
        \displaystyle{\frac{\partial f}{\partial y} (x,y)} & = &
        \displaystyle{\left. \frac{\chi(x)-\chi(y)}{|\chi(x)-\chi(y)|}
        \cdot \frac{d\chi}{ds} \right|_{s=y}}
        & = & - \cos \nu & = & -w \ .
\end{array} \right.
\]
We conclude that $f$ is a generating function for $\varphi$.
Similar approaches work for higher-dimensional billiard problems.
Periodic points are obtained by finding critical points of
real functions of $N$ variables in $X$,
\[
        (x_1, \ldots , x_N) \longmapsto
        |\chi(x_1)-\chi(x_2)| + \ldots +
        |\chi(x_{N-1})-\chi(x_N)| + |\chi(x_N)-\chi(x_1)| \ ,
\]
that is, by finding the $N$-sided (generalized) polygons
inscribed in $X$ of critical perimeter.
Notice that $\RR / \ZZ \times (-1,1) \simeq
\{ (x,v) \mid x \in X, v \in T_xX, |v|<1 \}$
is the open unit tangent ball bundle of a circle $X$,
which is an open annulus $A$, and the
map $\varphi: A \to A$ is area-preserving,
as in the next two theorems.
\end{example}

%%%%%%%%%%%%%%%%%%%%%%%%%%%%%%%%%%%%%%%%%%%%%%%%%%%%%%%%%%%%%%%%%%%%%%%%%%%%%

While studying {\em Poincar\'e return maps} in dynamical systems,
Poincar\'e arrived at the following results.

\begin{theorem}\label{thm:poincare_recurrence}\index{theorem !
Poincar\'e recurrence}\index{recurrence}\index{Poincar\'e !
recurrence theorem}
\textbf{(Poincar\'e Recurrence Theorem)} $\;$
Let $\varphi:A \to A$ be a volume-preserving diffeomorphism
of a finite-volume manifold $A$,
and $\cU$ a nonempty open set in $A$.
Then there is $q \in \cU$ and a positive
integer $N$ such that $\varphi^N (q) \in \cU$.
\end{theorem}

Hence, under iteration, a mechanical system governed by $\varphi$
will eventually return arbitrarily close to the initial state.

\begin{proof}
Let $\cU_0 = \cU, \cU_1=\varphi(\cU), \cU_2 =
\varphi^2 (\cU), \ldots$.
If all of these sets were disjoint, then,
since $\mbox{Volume } (\cU_i) = \mbox{ Volume } (\cU) > 0$ for all $i$,
the volume of $A$ would be greater or equal to 
$\sum_i \mbox{ Volume } (\cU_i) = \infty$.
To avoid this contradiction we must have
$\varphi^k (\cU) \cap \varphi^\ell (\cU) \ne \emptyset$ for some $k > \ell$,
which implies $\varphi^{k-\ell} (\cU) \cap \cU \ne \emptyset$.
\end{proof}

\vspace*{-2ex}

\begin{theorem}\index{Poincar\'e ! last geometric
theorem}\index{theorem ! Poincar\'e's last geometric
theorem}\label{thm:poincare_birkhoff}
\textbf{(Poincar\'e's Last Geometric Theorem)} $\;$
Suppose that $\varphi:A \to A$ is an area-preserving diffeomorphism
of the closed annulus $A =\RR / \ZZ \times [-1,1]$ that preserves
the two components of the boundary and twists them in opposite directions.
Then $\varphi$ has at least two fixed points.
\end{theorem}

This theorem was proved in 1913 by Birkhoff\index{Birkhoff !
Poincar\'e-Birkhoff theorem}~\cite{bi:dynamical},
and hence is also called the
\textbf{Poincar\'e-Birkhoff theorem}\index{Poincar\'e !
Poincar\'e-Birkhoff theorem}\index{theorem ! Poincar\'e-Birkhoff}.
It has important
applications to dynamical systems\index{dynamical system}
and celestial mechanics\index{mechanics ! celestial}.
The {\em Arnold conjecture}\index{Arnold ! conjecture}\index{conjecture !
Arnold} on the existence of fixed points for
symplectomorphisms\index{fixed point}\index{symplectomorphism !
fixed point}\index{symplectomorphism ! Arnold conjecture}
of compact manifolds (see Section~\ref{sec:arnold_floer})
may be regarded as a generalization of the Poincar\'e-Birkhoff theorem.
This conjecture has motivated a significant amount of
research involving a more general notion of generating function;
see, for instance,~\cite{el-gr:lagrangian,gi:periodic}.

%%%%%%%%%%%%%%%%%%%%%%%%%%%%%%%%%%%%%%%%%%%%%%%%%%%%%%%%%%%%%%%%%%%%%%%%%%%%%
%%%%%%%%%%%%%%%%%%%%%%%%%%%%%%%%%%%%%%%%%%%%%%%%%%%%%%%%%%%%%%%%%%%%%%%%%%%%%

\ssubsection{Lagrangians and Special Lagrangians in $\CC^n$}
\label{sec:special_lagrangians}

The standard \textbf{hermitian inner product}
$h(\cdot,\cdot)$ on $\CC^n$
has real and imaginary parts given by the euclidean inner product
$\langle \cdot , \cdot \rangle$ and (minus) the symplectic form
$\omega_0$, respectively: for $v=(x_1 + i y_1, \ldots , x_n + i y_n),
u=(a_1 + i b_1, \ldots , a_n + i b_n) \in \CC^n$,
\[
   \begin{array}{rcl}
   h (v,u) & = & \textstyle{\sum \limits_{k=1}^n
   (x_k + i y_k) (a_k - i b_k)} \\
   & = & \textstyle{\sum \limits_{k=1}^n (x_k a_k + y_k b_k) -
   i \sum \limits_{k=1}^n (x_k b_k - y_k a_k)} \\
   & = & \langle v, u \rangle - i \omega_0 (v,u) \ .
   \phantom{\sum \limits_{k=1}^n}
   \end{array}
\]

\begin{lemma}
\label{lem:lagrangian}
Let $W$ be a subspace of $(\CC^n, \omega_0)$
and $e_1, \ldots , e_n$ vectors in $\CC^n$.
Then:
\begin{itemize}
\item[(a)]
$W$ is lagrangian if and only if $W^\perp = i W$;
\item[(b)]
$(e_1, \ldots , e_n)$ is an orthonormal basis of a lagrangian subspace
if and only if $(e_1, \ldots , e_n)$ is a unitary basis of $\CC^n$.
\end{itemize}
\end{lemma}

\vspace*{-2ex}

\begin{proof}
\begin{itemize}
\item[(a)]
We always have
$\omega_0 (v,u) = - \mathrm{Im}\, h (v,u)
= \mathrm{Re}\, h (iv,u) = \langle iv, u \rangle$.
It follows that, if $W$ is lagrangian,
so that $\omega_0 (v,u)=0$ for all $v,u \in W$, then $i W \subseteq W^\perp$.
These spaces must be equal because they have the same dimension.
Reciprocally, when $\langle iv, u \rangle = 0$ for all $v,u \in W$,
the equality above shows that $W$ must be isotropic.
Since $\dim W = \dim iW = \dim W^\perp = 2n - \dim W$,
the dimension of $W$ must be $n$.
\item[(b)]
If $(e_1, \ldots , e_n)$ is an orthonormal basis of a
lagrangian subspace $W$, then, by the previous part,
$(e_1, \ldots , e_n, ie_1, \ldots , ie_n)$ is an
orthonormal basis of $\CC^n$ as a real vector space.
Hence $(e_1, \ldots , e_n)$ must be a complex basis of $\CC^n$
and it is unitary because
$h (e_j,e_k) = \langle e_j, e_k \rangle - i \omega_0 (e_j,e_k)
= \delta _{jk}$.
Conversely, if $(e_1, \ldots , e_n)$ is a unitary basis of $\CC^n$,
then the real span of these vectors is lagrangian
($\omega_0 (e_j,e_k) = - \mathrm{Im}\, h (e_j,e_k) = 0$)
and they are orthonormal
($\langle e_j, e_k \rangle = \mathrm{Re}\, h (e_j,e_k) = \delta _{jk}$).
\end{itemize}
\end{proof}

The \textbf{lagrangian grassmannian} $\Lambda_n$ is the set
of all lagrangian subspaces of $\CC^n$.
It follows from part~(b) of Lemma~\ref{lem:lagrangian} that
$\Lambda_n$ is the set of all subspaces of $\CC^n$
admitting an orthonormal basis that is a unitary basis of $\CC^n$.
Therefore, we have
\[
   \Lambda_n \simeq \UU(n) / \OO (n) \ .
\]
Indeed $\UU (n)$ acts transitively on $\Lambda_n$:
given $W,W' \in \Lambda_n$ with orthonormal bases
$(e_1, \ldots , e_n)$, $(e_1', \ldots , e_n')$ respectively,
there is a unitary transformation of $\CC^n$ that maps
$(e_1, \ldots , e_n)$ to $(e_1', \ldots , e_n')$ as
unitary bases of $\CC^n$.
And the stabilizer of $\RR^n \in \Lambda_n$
is the subgroup of those unitary transformations that
preserve this lagrangian subspace, namely $\OO(n)$.
It follows that $\Lambda_n$ is a compact connected manifold
%\footnote{The fundamental group of $\Lambda_n$
%is isomorphic to $\ZZ$.
%The \textbf{Maslov class} is the generator
%of $H^1 (\Lambda_n ; \ZZ)$.}
of dimension $\frac{n(n+1)}{2}$;
cf.\ the last example of Section~\ref{symplectic_linear_algebra}.

The lagrangian grassmannian comes with a \textbf{tautological
vector bundle}
\[
   \tau_n := \{ (W,v) \in \Lambda_n \times \CC^n \mid v \in W \} \ ,
\]
whose fiber over $W \in \Lambda_n$ is the $n$-dimensional real space $W$.
It is a consequence of part~(a) of Lemma~\ref{lem:lagrangian}
that the following map gives a well-defined global isomorphism
of the complexification $\tau_n \otimes_\RR \CC$ with the
trivial bundle $\underline{\CC^n}$ over $\Lambda_n$
(i.e., {\em a global trivialization}):
$(W,v \otimes c) \mapsto (W, cv)$, for $W \in \Lambda_n, v \in W, c \in \CC$.

\begin{definition}
A \textbf{lagrangian immersion} of a manifold $X$
is an immersion $f : X \to \CC^n$ such that $df_p (T_pX)$
is a lagrangian subspace of $(\CC^n, \omega_0)$, for every $p \in X$.
\end{definition}

\vspace*{-2ex}

\begin{example}
The graph of a map $h: \RR^n \to i \RR^n$
is an embedded $n$-dimensional submanifold $X$ of $\CC^n$.
Its tangent space at $(p,h(p))$
is $\{ v + dh_p (v) \mid v \in \RR^n \}$.
Let $e_1, \ldots , e_n$ be the standard basis of $\RR^n$.
Since $\omega_0 (e_k + dh_p (e_k) , e_j + dh_p (e_j)) =
\langle e_k , -i \, dh_p (e_j) \rangle
+ \langle e_j , i \, dh_p (e_k) \rangle$,
we see that $X$ is lagrangian if and only if
$\frac{\partial h_k}{\partial x_j} = \frac{\partial h_j}{\partial x_k}$,
$\forall j,k$, which in $\RR^n$ is if and only if
$h$ is the gradient of some $H : \RR^n \to i \RR$.
\end{example}

If $f : X \to \CC^n$ is a lagrangian immersion,
we can define a \textbf{Gauss map}
\[
\begin{array}{rrcl}
   \lambda_f : & X & \longrightarrow & \Lambda_n \\
   & p & \longmapsto & df_p (T_p X) \ .
\end{array}
\]
Since $\lambda_f^* \tau_n = TX$ and
$\tau_n \otimes \CC \simeq \underline{\CC^n}$,
we see that a necessary condition for an immersion $X \to \CC^n$
to exist is that the complexification of $TX$ be trivializable.
Using the h-principle (Section~\ref{sec:compatible_almost}),
Gromov~\cite{gr:partial} showed that this is also sufficient:
{\em an $n$-dimensional manifold $X$ admits a lagrangian immersion
into $\CC^n$ if and only if the complexification of its tangent bundle
is trivializable.}

\begin{example}
For the unit sphere
$S^n = \{ (t,x) \in \RR \times \RR^n \, : \, t^2 + |x|^2 = 1 \}$,
the \textbf{Whitney sphere immersion} is the map
\[
\begin{array}{rrcl}
   f : & S^n & \longrightarrow & \CC^n \\
   & (t,x) & \longmapsto & x + itx \ .
\end{array}
\]
The only self-intersection is at the origin
where $f (-1,0,\ldots,0) = f (1,0,\ldots,0)$.
Since $T_{(t,x)} S^n = (t,x)^\perp$, the differential
$df_{(t,x)} : (u,v) \mapsto v + i (tv + ux)$ is always injective:
$v + i (tv + ux)= 0 \Leftrightarrow v=0 \mbox{ and } ux=0$,
but when $x = 0$ it is $t=\pm 1$ and
$T_{(\pm 1,0)} S^n = \{0\} \times \RR^n$, so it must be $u=0$.
We conclude that $f$ is an immersion.
By computing $\omega_0$ at two vectors of the form $v + i (tv + ux)$,
we find that the image $df_p (T_p S^n)$
is an $n$-dimensional isotropic subspace of $\CC^n$.
Therefore, $f$ is a lagrangian immersion of $S^n$,
% of a sphere,
and the complexification $TS^n \otimes \CC$ must be always
trivializable, though the tangent bundle $TS^n$ is only
trivializable in dimensions $n=0,1,3,7$.
\end{example}

The \textbf{special lagrangian grassmannian} $S\Lambda_n$ is the set
of all {\em oriented} subspaces of $\CC^n$
admitting a {\em positive} orthonormal basis $(e_1, \ldots , e_n)$
that is a {\em special} unitary basis of $\CC^n$.
By the characterization of lagrangian in the
part~(b) of Lemma~\ref{lem:lagrangian},
it follows that the elements of $S\Lambda_n$
are indeed lagrangian submanifolds.
Similarly to the case of the lagrangian grassmannian, we have that
\[
   S \Lambda_n \simeq \SU(n) / \SO (n)
\]
is a compact connected manifold of dimension $\frac{n(n+1)}{2} - 1$.

We can single out the {\em special} lagrangian subspaces by
expressing the condition on the determinant in terms of the real $n$-form
in $\CC^n$
\[
   \beta : = \mathrm{Im}\, \Omega \  , \quad
   \mbox{ where } \quad \Omega : = dz_1 \wedge \ldots \wedge dz_n \ .
\]
Since for $A \in \SO(n)$, we have $\det A = 1$ and
$\Omega (e_1, \ldots , e_n) = \Omega (Ae_1, \ldots , Ae_n)$,
we see that, for an oriented real $n$-dimensional subspace $W \subset \CC^n$,
the number $\Omega (e_1, \ldots , e_n)$ does not depend
on the choice of a positive orthonormal basis $(e_1, \ldots , e_n)$ of $W$,
thus can be denoted $\Omega (W)$ and its imaginary part $\beta (W)$.

\begin{proposition}
\label{prop:special_lagrangian}
A subspace $W$ of $(\CC^n, \omega_0)$ has an orientation
for which it is a special lagrangian
if and only if $W$ is lagrangian and $\beta (W) = 0$.
\end{proposition}

\vspace*{-2ex}

\begin{proof}
Any orthonormal basis $(e_1, \ldots , e_n)$ of
a lagrangian subspace $W \subset \CC^n$
is the image of the canonical basis
of $\CC^n$ by some $A \in \UU (n)$, and $\Omega (W) = \det A \in S^1$.
Therefore, $W$ admits an orientation for which such a {\em positive}
$(e_1, \ldots , e_n)$ is a {\em special} unitary basis of $\CC^n$
if and only if $\det A = \pm 1$, i.e., $\beta (W) = 0$.
\end{proof}

\vspace*{-1ex}

\begin{definition}
A \textbf{special lagrangian immersion} of an oriented manifold $X$
is a lagrangian immersion $f : X \to \CC^n$
such that, at each $p \in X$, the space $df_p (T_pX)$
is a special lagrangian subspace of $(\CC^n, \omega_0)$.
\end{definition}

For a special lagrangian immersion $f$,
the Gauss map $\lambda_f$ takes values in $S\Lambda_n$.

By Proposition~\ref{prop:special_lagrangian},
the immersion $f$ of an $n$-dimensional
manifold $X$ in $(\CC^n,\omega_0)$ is {\em special lagrangian}
if and only if $f^* \omega_0 = 0$ and $f^* \beta = 0$

\begin{example}
In $\CC^2$, writing $z_k = x_k + i y_k$,
we have $\beta = dx_1 \wedge dy_2 + dy_1 \wedge dx_2$.
We have seen that the graph of the gradient $i \nabla H$
is lagrangian, for any function $H : \RR^2 \to \RR$.
So $f(x_1,x_2) = (x_1,x_2,i \frac{\partial H}{\partial x_1},
i \frac{\partial H}{\partial x_2})$ is a lagrangian immersion.
For $f$ to be a {\em special} lagrangian immersion, we need
the vanish of
\[
   f^* \beta = dx_1 \wedge d \left( \frac{\partial H}{\partial x_2} \right)
   + d \left( \frac{\partial H}{\partial x_1} \right) \wedge dx_2
   = \left( \frac{\partial^2 H}{\partial x_1^2} 
   + \frac{\partial^2 H}{\partial x_2^2} \right) dx_1 \wedge dx_2 \ .
\]
Hence the graph of $\nabla H$ is special lagrangian
if and only if $H$ is {\em harmonic}.
\end{example}

If $f : X \to \CC^n$ is a special lagrangian immersion,
then $f^* \Omega$ is an exact (real) volume form:
$f^* \Omega = d \mathrm{Re}\, (z_1 dz_2 \wedge \ldots \wedge dz_n)$.
We conclude, by Stokes theorem, that there can be no
special lagrangian immersion of a compact manifold in $\CC^n$.
{\em Calabi-Yau manifolds}\footnote{\textbf{Calabi-Yau manifolds}
are compact {\em K\"ahler manifolds} (Section~\ref{sec:kahler})
with vanishing first Chern class.}
are more general manifolds where a definition of special lagrangian
submanifold makes sense and where the space of special lagrangian
embeddings of a compact manifold is interesting.
Special lagrangian geometry was introduced by
Harvey and Lawson~\cite{ha-la:calibrated}.
For a treatment of lagrangian and special lagrangian submanifolds
with many examples, see for instance~\cite{au:barcelona}.

%%%%%%%%%%%%%%%%%%%%%%%%%%%%%%%%%%%%%%%%%%%%%%%%%%%%%%%%%%%%%%%%%%%%%%%%%%%%%
%%%%%%%%%%%%%%%%%%%%%%%%%%%%%%%%%%%%%%%%%%%%%%%%%%%%%%%%%%%%%%%%%%%%%%%%%%%%%
% --> Section 3
%%%%%%%%%%%%%%%%%%%%%%%%%%%%%%%%%%%%%%%%%%%%%%%%%%%%%%%%%%%%%%%%%%%%%%%%%%%%%
%%%%%%%%%%%%%%%%%%%%%%%%%%%%%%%%%%%%%%%%%%%%%%%%%%%%%%%%%%%%%%%%%%%%%%%%%%%%%

\newpage

\ssection{Complex Structures}
\index{complex structure}
\label{section3}

%%%%%%%%%%%%%%%%%%%%%%%%%%%%%%%%%%%%%%%%%%%%%%%%%%%%%%%%%%%%%%%%%%%%%%%%%%%%%
%%%%%%%%%%%%%%%%%%%%%%%%%%%%%%%%%%%%%%%%%%%%%%%%%%%%%%%%%%%%%%%%%%%%%%%%%%%%%

\ssubsection{Compatible Linear Structures}
\label{compatible_linear_structures}

A \textbf{complex structure}\index{complex structure !
on a vector space} on a vector space $V$ is a linear map
$J: V \to V$ such that $J^2 = -\Id$.
The pair $(V,J)$ is then called a
\textbf{complex vector space}\index{complex vector space}\index{vector
space ! complex}.
A complex structure $J$ on $V$ is equivalent to a structure of
vector space over $\CC$, the map $J$ corresponding to multiplication by $i$.
If $(V, \Omega)$ is a symplectic vector space,
a complex structure $J$ on $V$ is said to be
\textbf{compatible}\index{compatible ! complex structure}\index{complex
structure ! compatibility} (with $\Omega$, or $\Omega$-compatible) if
the bilinear map $G_{_J} : V \times V \to \RR$
defined by $G_{_J}(u,v) = \Omega(u,Jv)$ is an inner product on $V$.
This condition comprises $J$ being a symplectomorphism
(i.e., $\Omega(Ju,Jv) = \Omega(u,v)$ $\forall u,v$) and the so-called
\textbf{taming}\index{taming}: $\Omega(u, Ju) > 0$, $\forall u \neq 0$.

\begin{example}
For the symplectic vector space $(\RR^{2n}, \Omega_0)$
with symplectic basis
$e_1=(1,0,\ldots,0), \ldots, e_n, f_1, \ldots, f_n=(0,\ldots,0,1)$,
there is a standard compatible complex structure $J_0$ determined
by $J_0(e_j) = f_j$ and $J_0(f_j) = -e_j$ for all $j=1,\ldots,n$.
This corresponds to a standard identification
of $\RR^{2n}$ with $\CC^n$,
and $\Omega_0 (u,J_0v) = \langle u,v \rangle$ is
the standard euclidean inner product.
With respect to the symplectic basis $e_1, \ldots, e_n, f_1, \ldots, f_n$,
the map $J_0$ is represented by the matrix
\[
        \left[ \begin{array}{cc}
        0 & -\mbox{Id} \\
        \mbox{Id} & 0
        \end{array} \right] \ .
\]
The \textbf{symplectic linear group}\index{symplectic !
linear group}, $\Sp (2n) := \{ A \in \GL (2n;\RR) \, | \,
\Omega_0 (Au , Av) = \Omega_0 (u , v)$
$\mbox{for all } u,v \in \RR^{2n} \}$,
is the group of all linear
transformations of $\RR^{2n}$ that preserve the
standard symplectic structure.
The \textbf{orthogonal group} $\OO (2n)$ is the group
formed by the linear transformations $A$ that preserve the euclidean
inner product, $\langle Au , Av \rangle = \langle u , v \rangle$,
for all $u,v \in \RR^{2n}$.
The \textbf{general complex group} $\GL (n;\CC)$
is the group of linear transformations $A: \RR^{2n} \to \RR^{2n}$
commuting with $J_0$, $A(J_0v) = J_0 (Av)$,
for all $v \in \RR^{2n}$.\footnote{Identify the complex
$n \times n$ matrix $X + i Y$ with the real $2n \times 2n$ matrix
$\left[ \begin{array}{cc} X & -Y \\ Y & X \end{array} \right]$.}
The compatibility between the structures
$\Omega_0$, $\langle \cdot , \cdot \rangle$ and $J_0$
implies that the intersection of
{\em any two} of these subgroups of $\GL (2n;\RR)$ is the same group,
namely the \textbf{unitary group}\index{unitary group} $\UU (n)$.
\end{example}

As $(\RR^{2n}, \Omega_0)$ is the prototype of a $2n$-dimensional
symplectic vector space, the preceding example shows that
compatible complex structures always exist on
symplectic vector spaces.\footnote{Conversely,
given $(V,J)$, there is a symplectic $\Omega$
with which $J$ is compatible:
take $\Omega (u,v) = G(Ju,v)$ for an inner product $G$
such that $J^t = -J$.}
There is yet a way to produce a {\em canonical} compatible
complex structure $J$ after the choice of an inner product $G$
on $(V, \Omega)$, though the starting $G(u,v)$
is usually different from $G_{_J}(u,v) := \Omega(u,Jv)$.

\begin{proposition}
\label{polar_decomposition}
Let $(V, \Omega)$ be a symplectic vector space,
with an inner product $G$.
Then there is a canonical compatible complex structure $J$ on $V$.
\end{proposition}

\vspace*{-2ex}

\begin{proof}
By nondegeneracy of $\Omega$ and $G$,
the maps $u \mapsto \Omega(u, \cdot)$ and
$w \mapsto G(w, \cdot)$ are both isomorphisms between $V$ and $V^*$.
Hence, $\Omega(u,v) = G(Au,v)$ for some linear $A: V \to V$.
The map $A$ is skew-symmetric, and the
product $AA^t$ is symmetric\footnote{A map $B:V \to V$
is \textbf{symmetric}, respectively
\textbf{skew-symmetric}, when $B^t = B$, resp.\ $B^t = -B$,
where the transpose $B^t : V \to V$ is determined by
$G(B^tu, v) = G(u,Bv)$.}
and positive: $G(AA^tu,u) = G(A^tu, A^tu) > 0$, for $u \neq 0$.
By the spectral theorem, these properties imply that $AA^t$ diagonalizes
with positive eigenvalues $\lambda_i$, say
$AA^t = B \ \mbox{diag} \,
( \lambda_1, \ldots, \lambda_{2n} ) \ B^{-1}$.
We may hence define an arbitrary real power of $AA^t$ by
rescaling the eigenspaces, in particular,
\[
   \sqrt{AA^t} := B \, \mbox{diag} \,
   ( \sqrt{\lambda_1}, \ldots, \sqrt{\lambda_{2n}} ) \ B^{-1} \ .
\]
The linear transformation $\sqrt{AA^t}$ is symmetric, positive-definite
and does not depend on the
choice of $B$ nor of the ordering of the eigenvalues.
It is completely determined by its effect on each eigenspace of $AA^t$:
on the eigenspace corresponding to the eigenvalue $\lambda_k$,
the map $\sqrt{AA^t}$ is defined to be multiplication by $\sqrt{\lambda_k}$.

Let $J := (\sqrt{AA^t})^{-1}A$.
Since $A$ and $\sqrt{AA^t}$ commute,
$J$ is orthogonal ($JJ^t = \Id$),
as well as skew-symmetric ($J^t = -J$).
It follows that $J$ is a complex structure on $V$.
Compatibility is easily checked:
\[
\begin{array}{c}
        \Omega(Ju, Jv) = G(AJu,Jv) = G(JAu,Jv) = G(Au,v)
        = \Omega(u,v) \mbox{ and }\\
        \Omega(u,Ju) = G(Au,Ju) = G(-JAu, u)
        = G(\sqrt{AA^t} \, u,u) > 0\ , \mbox{ for }u \neq 0 \ .
\end{array}
\]
\end{proof}

The factorization $A = \sqrt{AA^t} \, J$ is called the
\textbf{polar decomposition}\index{polar decomposition}\index{complex
structure ! polar decomposition} of $A$.

\begin{remark}
Being {\em canonical}, this construction may be {\em smoothly} performed:
when $(V_t, \Omega_t)$ is a family of symplectic vector spaces
with a family $G_t$ of inner products,
all depending smoothly on a parameter $t$,
an adaptation of the previous proof
shows that there is a smooth family $J_t$ of
compatible complex structures on $(V_t, \Omega_t)$.
\end{remark}

Let $(V, \Omega)$ be a symplectic vector space of dimension $2n$,
and let $J$ be a complex structure on $V$.
If $J$ is $\Omega$-compatible and $L$ is a
lagrangian subspace of $(V, \Omega)$, then $JL$ is also lagrangian
and $JL = L^\perp$, where $\perp$ indicates orthogonality with respect
to the inner product $G_{_J} (u,v) = \Omega (u, Jv)$.
Therefore, a complex structure $J$ is $\Omega$-compatible
{\em if and only if}
there exists a symplectic basis for $V$ of the form
\[
        e_1, e_2, \ldots , e_n, f_1=Je_1, f_2=Je_2, \ldots, f_n=Je_n \ .
\]

Let ${\cal J} (V,\Omega)$ be the \textbf{set of all
compatible complex structures
in a symplectic vector space} $(V,\Omega)$.

\begin{proposition}
\label{prop:linear_contractible}
The set ${\cal J} (V,\Omega)$ is
contractible.\footnote{\textbf{Contractibility} of ${\cal J} (V,\Omega)$
means that there exists a homotopy $h_t: \cJ(V,\Omega) \to \cJ(V, \Omega)$,
$0 \leq t \leq 1$, starting at the identity $h_0 = \Id$,
finishing at a trivial map $h_1: \cJ(V,\Omega) \to \{J_0\}$,
and fixing $J_0$ (i.e., $h_t(J_0) = J_0$, $\forall t$)
for some $J_0  \in \cJ(V,\Omega)$.}
\end{proposition}

\vspace*{-2ex}

\begin{proof}
Pick a lagrangian subspace $L_0$ of $(V,\Omega)$.
Let ${\cal L} (V,\Omega,L_0)$ be the space of all lagrangian subspaces
of $(V,\Omega)$ that intersect $L_0$ transversally.
Let ${\cal G} (L_0)$ be the space of all inner products on $L_0$.
The map
\[
\begin{array}{rrcl}
        \Psi : & {\cal J} (V,\Omega) & \longrightarrow &
        {\cal L} (V,\Omega,L_0)\times{\cal G}(L_0) \\
        & J & \longmapsto & (J L_0, G_{_J}|_{L_0})
\end{array}
\]
is a homeomorphism, with inverse as follows.
Take $(L,G) \in {\cal L} (V,\Omega,L_0)\times{\cal G}(L_0)$.
For $v\in L_0$, $v^\perp = \{ u \in L_0 \, | \, G(u,v) = 0 \}$
is a $(n-1)$-dimensional space of $L_0$; its symplectic orthogonal
$(v^\perp)^\Omega$ is $(n+1)$-dimensional.
Then $(v^\perp)^\Omega \cap L$ is $1$-dimensional.
Let $Jv$ be the unique vector in this line such that $\Omega (v,Jv) =1$.
If we take $v$'s in some $G$-orthonormal basis of $L_0$,
this defines an element $J \in {\cal J} (V,\Omega)$.

The set ${\cal L} (V,\Omega,L_0)$ can be identified with
the vector space of all symmetric $n \times n$ matrices.
In fact, any $n$-dimensional subspace $L$ of $V$ that is transverse
to $L_0$ is the graph of a linear map $JL_0 \to L_0$,
and the lagrangian ones correspond to symmetric maps
(cf.\ Section~\ref{symplectic_linear_algebra}).
Hence, ${\cal L} (V,\Omega,L_0)$ is contractible.
Since ${\cal G} (L_0)$ is contractible (it is even convex),
we conclude that ${\cal J} (V,\Omega)$ is contractible.
\end{proof}

%%%%%%%%%%%%%%%%%%%%%%%%%%%%%%%%%%%%%%%%%%%%%%%%%%%%%%%%%%%%%%%%%%%%%%%%%%%%%
%%%%%%%%%%%%%%%%%%%%%%%%%%%%%%%%%%%%%%%%%%%%%%%%%%%%%%%%%%%%%%%%%%%%%%%%%%%%%

\ssubsection{Compatible Almost Complex Structures}
\label{sec:compatible_almost}
\index{almost complex structure !
definition}\index{almost complex structure ! compatible}

An \textbf{almost complex structure}\index{almost complex structure !
definition} on a manifold $M$ is a
smooth\footnote{\textbf{Smoothness} means that for any vector field
$v$, the image $Jv$ is a (smooth) vector field.}
field of complex structures on the tangent spaces,
$J_p: T_pM \to T_pM$, $p \in M$.
The pair $(M,J)$ is then called an
\textbf{almost complex manifold}\index{almost complex manifold}.

\begin{definition}
An almost complex structure $J$ on
a symplectic manifold $(M,\omega)$
is \textbf{compatible}\index{almost complex structure !
compatibility}\index{compatible ! almost complex structure}
(with $\omega$ or $\omega$-compatible) if the map that
assigns to each point $p \in M$ the bilinear pairing
$g_p: T_pM \times T_pM \to \RR$, $g_p(u,v) := \omega_p(u, J_pv)$
is a riemannian metric on $M$.\index{riemannian ! metric}\index{metric}
A triple $(\omega, g, J)$ of a symplectic form,
a riemannian metric and an almost complex structure
on a manifold $M$ is a
\textbf{compatible triple}\index{compatible ! triple}
when $g(\cdot, \cdot) = \omega(\cdot, J\cdot)$.
\end{definition}

If $(\omega, J, g)$ is a
compatible triple,\index{compatible ! triple}
each of $\omega$, $J$ or $g$ can be written in terms of the other two.

\begin{examples}
\begin{enumerate}
\item
If we identify $\RR^{2n}$ with $\CC^n$ using
coordinates $z_j = x_j + i y_j$, multiplication
by $i$ induces a constant linear map $J_0$
on the tangent spaces such that $J_0^2 = - \Id$, known as the
\textbf{standard almost complex structure} on $\RR^{2n}$:
\[
   J_0 \left (\frac{\partial}{\partial x_j} \right)
   = \frac{\partial}{\partial y_j}\ ,
   \qquad
   J_0 \left( \frac{\partial}{\partial y_j} \right)
   = -\frac{\partial}{\partial x_j}\ .
\]
For the standard symplectic form $\omega_0 = \sum dx_j \wedge dy_j$
and the euclidean inner product $g_0 = \langle \cdot, \cdot \rangle$,
the compatibility relation holds: $\omega_0(u,v) = g_0(J_0(u),v)$.

\item
   Any oriented hypersurface $\Sigma \subset \RR^3$
carries a natural symplectic form and a natural
compatible almost complex structure induced by the
standard inner (or dot) and exterior (or vector) products.
They are given by the formulas
$\omega_p (u,v) := \langle \nu_p , u \times v \rangle$
and $J_p (v) = \nu_p \times v$ for $v \in T_p\Sigma$,
where $\nu_p$ is the outward-pointing unit normal vector at $p \in \Sigma$
(in other words, $\nu : \Sigma \to S^2$ is the {\em Gauss map}).
Cf.\ Example~3 of Section~\ref{symplectic_forms}.
The corresponding riemannian metric is the restriction to $\Sigma$
of the standard euclidean metric $\langle \cdot , \cdot \rangle$.

\item
   The previous example generalizes to the
oriented hypersurfaces $M \subset \RR^7$.
Regarding $u,v \in \RR^7$ as imaginary {\em octonions}
(or {\em Cayley numbers}),
the natural vector product $u \times v$ is the imaginary part of the
product of $u$ and $v$ as octonions.
This induces a natural almost complex structure on $M$ given by
$J_p (v) = \nu_p \times v$, where $\nu_p$
is the outward-pointing unit normal vector at $p \in M$.
In the case of $S^6$, at least, this $J$ is not compatible
with any symplectic form, as $S^6$ cannot be a symplectic manifold.
\end{enumerate}
\end{examples}

As a consequence of the remark in
Section~\ref{compatible_linear_structures}, we have:

\begin{proposition}
On any symplectic manifold $(M,\omega)$ with
a riemannian metric $g$,
there is a canonical compatible almost complex structure $J$.
\end{proposition}

Since riemannian metrics always exist, we conclude that
{\em any symplectic manifold has compatible almost complex structures}.
The metric $g_{_J}(\cdot,\cdot) := \omega(\cdot, J\cdot)$
tends to be different from the given $g(\cdot, \cdot)$.

%%%%%%%%%%%%%%%%%%%%%%%%%%%%%%%%%%%%%%%%%%%%%%%%%%%%%%%%%%%%%%%%%%%%%%%%%%%%%

\begin{proposition}
Let $(M,J)$ be an almost complex manifold
where $J$ is compatible with two
symplectic forms $\omega_0, \omega_1$
Then $\omega_0$ and $\omega_1$ are deformation-equivalent.\index{deformation
equivalent}
\end{proposition}

\vspace*{-2ex}

\begin{proof}
Simply take the convex combinations
$\omega_t = (1-t)\omega_0 + t\omega_1$, $0 \leq t \leq 1$.
\end{proof}

A counterexample to the converse of this proposition
is provided by the family
$\omega_t = \cos \pi t~ dx_1 \wedge dy_1 + \sin \pi t ~dx_1 \wedge dy_2
+ \sin \pi t ~dy_1 \wedge dx_2 + \cos \pi t ~dx_2 \wedge dy_2$
for $0 \leq t \leq 1$.
There is no $J$ in $\RR^4$ compatible with both $\omega_0$
and $\omega_1 = - \omega_0$.

A submanifold $X$ of an almost complex manifold $(M,J)$
is an \textbf{almost complex submanifold}\index{submanifold !
almost complex}\index{almost complex submanifold} when
$J(TX) \subseteq TX$, i.e., we have $J_p v \in T_p X$,
$\forall p \in X, v \in T_p X$.

\begin{proposition}
Let $(M,\omega)$ be a symplectic manifold
equipped with a compatible almost complex structure $J$.
Then any almost complex submanifold $X$ of $(M,J)$ is a symplectic
submanifold of $(M, \omega)$.
\end{proposition}

\vspace*{-2ex}

\begin{proof}
Let $i: X \hookrightarrow M$ be the inclusion.
Then $i^*\omega$ is a closed 2-form on $X$.
Since $\omega_p(u,v) = g_p(J_pu,v)$, $\forall p \in X$,
$\forall u, v \in T_pX$, and since $g_p|_{T_pX}$ is nondegenerate,
so is $\omega_p|_{T_pX}$, and $i^*\omega$ is nondegenerate.
\end{proof}

%%%%%%%%%%%%%%%%%%%%%%%%%%%%%%%%%%%%%%%%%%%%%%%%%%%%%%%%%%%%%%%%%%%%%%%%%%%%%

It is easy to see that the
\textbf{set $\cJ(M, \omega)$ of all compatible almost complex structures
on a symplectic manifold} $(M,\omega)$ is path-connected.
From two almost complex structures $J_0, J_1$ compatible with $\omega$,
we get two riemannian metrics
$g_0 (\cdot, \cdot) = \omega (\cdot, J_0 \cdot)$,
$g_1 (\cdot, \cdot) = \omega (\cdot, J_1 \cdot)$.
Their convex combinations
\[
        g_t(\cdot, \cdot) = (1-t)g_0(\cdot, \cdot) + tg_1(\cdot, \cdot)\ ,
        \qquad 0 \leq t \leq 1\ ,
\]
form a smooth family of riemannian metrics.
Applying the polar decomposition to the family
$(\omega, g_t)$, we obtain a smooth path of compatible
almost complex structures $J_t$ joining $J_0$ to $J_1$.
The set $\cJ(M, \omega)$ is even
{\em contractible} (this is important for defining invariants).
The first ingredient is the contractibility
of the set of compatible complex structures on a vector space
(Proposition~\ref{prop:linear_contractible}).
Consider the fiber bundle $\cJ \to M$ with fiber over $p \in M$
being the space $\cJ_p := \cJ(T_pM, \omega_p)$
of compatible complex structures on the tangent space at $p$.
A compatible almost complex structure on $(M,\omega)$ is
a section of $\cJ$.
The space of sections of $\cJ$ is contractible
because the fibers are contractible.\footnote{The base being a
(second countable and Hausdorff) manifold,
a contraction can be produced using a countable cover
by trivializing neighborhoods
whose closures are compact subsets of larger trivializing neighborhoods,
and such that each $p \in M$ belongs to only a finite number
of such neighborhoods.}

The \textbf{first Chern class} $c_1(M,\omega)$ of a
symplectic manifold $(M,\omega)$ is the first
Chern class of $(TM,J)$ for any compatible $J$.
The class $c_1(M,\omega) \in H^2 (M;\ZZ)$
is invariant under deformations of $\omega$.

We never used the closedness of $\omega$ to
obtain compatible almost complex structures.
The construction holds for an
\textbf{almost symplectic manifold}\index{manifold !
almost symplectic}\index{almost symplectic manifold}\index{symplectic !
almost symplectic manifold}
$(M,\omega)$, that is, a pair of a manifold $M$ and
a nondegenerate 2-form $\omega$, not necessarily closed.
We could further work with a
\textbf{symplectic vector bundle},\index{symplectic ! vector bundle}
that is, a vector bundle $E \to M$
equipped with a smooth field $\omega$ of fiberwise nondegenerate
skew-symmetric bilinear maps (Section~\ref{symplectic_submanifolds}).
The existence of such a field $\omega$ is equivalent to
being able to reduce the structure group of the bundle
from the general linear group to the linear symplectic group.
As both $\Sp (2n)$ and $\GL (n;\CC)$ retract to their common
maximal compact subgroup $\UU (n)$, a symplectic vector bundle
can be always endowed with a structure of complex vector bundle,
and vice-versa.

Gromov showed in his thesis~\cite{gr:stable} that
any {\em open}\footnote{A manifold is \textbf{open} if it has
no closed connected components, where \textbf{closed} means
compact and without boundary.} almost complex manifold
admits a symplectic form.
The books~\cite[\S 10.2]{el-mi:principle}
and~\cite[\S 7.3]{mc-sa:introduction}
contain proofs of this statement using different techniques.

\begin{theorem}
\textbf{(Gromov)} $\;$ For an open manifold
the existence of an almost complex structure $J$ implies
that of a symplectic form $\omega$ in any
given 2-cohomology class and such that $J$ is homotopic
to an almost complex structure compatible with $\omega$.
\end{theorem}

From an almost complex structure $J$ and a metric $g$,
one builds a nondegenerate 2-form
$\omega (u,v) = g (Ju,v)$, which will not be closed in general.
Closedness is a {\em differential relation}, i.e., a condition
imposed on the partial derivatives, encoded as a subset of {\em jet space}.
One says that a differential relation satisfies the
\textbf{h-principle}\footnote{There are in fact different h-principles
depending on the different possible coincidences of homotopy groups
for the spaces of formal solutions and of holonomic solutions.}
if any \textbf{formal solution}
(i.e., a solution for the associated algebraic problem,
in the present case a nondegenerate 2-form) is homotopic
to a \textbf{holonomic solution} (i.e, a genuine solution,
in the present case a closed nondegenerate 2-form).
Therefore, when the h-principle holds,
one may concentrate on a purely topological question
(such as the existence of an almost complex structure)
in order to prove the existence of a differential solution.
Gromov showed that, for an open differential relation
on an open manifold, when the relation is invariant under
the group of diffeomorphisms of the underlying manifold,
the inclusion of the space of holonomic solutions into
the space of formal solutions is a weak homotopy equivalence,
i.e., induces isomorphisms of all homotopy groups.
The previous theorem fits here as an application.

For {\em closed} manifolds there is no such theorem:
as discussed in Section~\ref{symplectic_forms},
the existence of a 2-cohomology class whose top power
is nonzero is also necessary for the existence of a symplectic form
and there are further restrictions coming from
{\em Gromov-Witten theory} (see Section~\ref{sec:invariants}).

%%%%%%%%%%%%%%%%%%%%%%%%%%%%%%%%%%%%%%%%%%%%%%%%%%%%%%%%%%%%%%%%%%%%%%%%%%%%%
%%%%%%%%%%%%%%%%%%%%%%%%%%%%%%%%%%%%%%%%%%%%%%%%%%%%%%%%%%%%%%%%%%%%%%%%%%%%%

\ssubsection{Integrability}
\label{sec:integrability}
\index{Dolbeault decompositions}\index{integrability}

Any {\em complex manifold}\index{complex manifold}\footnote{A
\textbf{complex manifold}\index{complex manifolds}\index{manifold !
complex}\index{complex ! manifold} of (complex) dimension $n$
is a set $M$ with a complete complex atlas
$\left\{(\cU_\alpha, \cV_\alpha, \varphi_\alpha)\ ,
\alpha \in \mbox{ index set } I \right\}$
where $M = \cup_\alpha \cU_\alpha$, the $\cV_\alpha$'s are open
subsets of $\CC^n$,
and the maps $\varphi_\alpha: \cU_\alpha \to \cV_\alpha$ are
bijections such that the transition maps
$\psi_{\alpha\beta} = \varphi_\beta \circ \varphi_\alpha^{-1}:
\cV_{\alpha \beta} \to \cV_{\beta \alpha}$
are {\em biholomorphic}\index{biholomorphic map}
(i.e., bijective, holomorphic and with holomorphic inverse)
as maps on open subsets of $\CC^n$,
$\cV_{\alpha \beta} = \varphi_\alpha (\cU_\alpha \cap \cU_\beta)$.}
has a canonical almost complex structure $J$.
It is defined locally over the domain $\cU$ of a complex chart
$\varphi: \cU \rightarrow \cV \subseteq \CC^n$, by
$J_p \left( \left. \frac {\partial}{\partial x_j} \right|_p \right) =
\left. \frac {\partial}{\partial y_j} \right|_p$ and
$J_p\left( \left. \frac {\partial}{\partial y_j} \right|_p \right) =
\left. -\frac {\partial}{\partial x_j} \right|_p$,
where these are the tangent vectors induced by the
real and imaginary parts of the coordinates of $\varphi = (z_1,\ldots,z_n)$,
$z_j = x_j + iy_j$.
This yields a globally well-defined $J$, thanks to
the {\em Cauchy-Riemann equations}\index{Cauchy-Riemann
equations}\index{Riemann ! Cauchy-Riemann equations}
satisfied by the components of the transition maps.

An almost complex structure $J$ on a manifold $M$
is called \textbf{integrable}\index{integrable !
almost complex structure}\index{almost complex structure ! integrability}
when $J$ is induced by some underlying structure of
complex manifold on $M$ as above.
The question arises whether some compatible almost complex structure
$J$ on a symplectic manifold $(M,\omega)$ is integrable.\index{integrability}
To understand what is involved,
we review Dolbeault theory and the Newlander-Nirenberg theorem.

%%%%%%%%%%%%%%%%%%%%%%%%%%%%%%%%%%%%%%%%%%%%%%%%%%%%%%%%%%%%%%%%%%%%%%%%%%%%%

Let $(M,J)$ be a $2n$-dimensional almost complex manifold.
The fibers of the complexified tangent bundle, $TM \otimes \CC$,
are $2n$-dimensional vector spaces over $\CC$.
We may extend $J$ linearly to $TM \otimes \CC$ by
$J(v\otimes c) = Jv\otimes c$, $v \in TM$, $c \in \CC$.
Since $J^2 = -\Id$, on the complex vector space $(TM \otimes \CC)_p$
the linear map $J_p$ has eigenvalues $\pm i$.
The $(\pm i)$-eigenspaces of $J$ are denoted $T_{1,0}$ and $T_{0,1}$,
respectively, and called the spaces of
\textbf{$J$-holomorphic}\index{J-holomorphic
tangent@($J$-)holomorphic tangent
vectors}\index{holomorphic tangent@($J$-)holomorphic tangent vectors}
and of \textbf{$J$-anti-holomorphic
tangent vectors}\index{J-anti-holomorphic tangent@($J$-)anti-holomorphic
tangent vectors}\index{anti-holomorphic tangent@($J$-)anti-holomorphic
tangent vectors}.
We have an isomorphism
\[
\begin{array}{rrcl}
        (\pi_{1,0}, \pi_{0,1}) : &
        TM \otimes \CC & \stackrel{\simeq}\longrightarrow &
        T_{1,0}\oplus T_{0,1} \\
        & v & \longmapsto &
        \frac{1}{2} (v - iJv , v + iJv)
\end{array}
\]
where the maps to each summand satisfy
$\pi_{1,0} \circ J = i\pi_{1,0}$ and 
$\pi_{0,1} \circ J = - i\pi_{0,1}$.
Restricting $\pi_{1,0}$ to $TM$, we see that
$(TM, J) \simeq T_{1,0} \simeq \overline{T_{0,1}}$,
as complex vector bundles, where the multiplication by $i$
is given by $J$ in $(TM, J)$ and where $\overline{T_{0,1}}$
denotes the complex conjugate bundle of $T_{0,1}$.

Similarly, $J^*$ defined on $T^*M \otimes \CC$ by
$J^* \xi = \xi \circ J$ has ($\pm i$)-eigenspaces
$T^{1,0} = (T_{1,0})^*$ and $T^{0,1} = (T_{0,1})^*$, respectively,
called the spaces of \textbf{complex-linear}\index{complex-linear
cotangent vectors}
and of \textbf{complex-antilinear cotangent
vectors}\index{complex-antilinear cotangent vectors}.
Under the two natural projections $\pi^{1,0}, \pi^{0,1}$ ,
the complexified cotangent bundle splits as
\[
\begin{array}{rrcl}
        (\pi^{1,0}, \pi^{0,1}) : &
        T^*M \otimes \CC & \stackrel{\simeq}\longrightarrow &
        T^{1,0} \oplus T^{0,1} \\
        & \xi & \longmapsto &
        \frac{1}{2} (\xi - iJ^* \xi, \xi + i J^* \xi) \ .
\end{array}
\]

%%%%%%%%%%%%%%%%%%%%%%%%%%%%%%%%%%%%%%%%%%%%%%%%%%%%%%%%%%%%%%%%%%%%%%%%%%%%%

Let
\[
   \Lambda^k (T^*M\otimes\CC) :=
   \Lambda^k(T^{1,0}\oplus T^{0,1}) =
   \oplus_{\ell+m=k} \Lambda^{\ell,m} \ ,
\]
where $\Lambda^{\ell,m} := (\Lambda^\ell T^{1,0}) \wedge (\Lambda^mT^{0,1})$,
and let $\Omega^k(M;\CC)$ be the space of sections of
$\Lambda^k (T^*M \otimes \CC)$, called
\textbf{complex-valued $k$-forms on $M$}\index{form !
complex-valued}\index{complex-valued form}.
%In particular, $\Lambda^{1,0} = T^{1,0}$ and $\Lambda^{0,1} = T^{0,1}$.
The \textbf{differential forms of type $(\ell,m)$}\index{form ! type}
\index{form ! complex-valued}on $(M,J)$ are the sections of $\Lambda^{\ell,m}$,
and the space of these differential forms is denoted $\Omega^{\ell,m}$.
The decomposition of forms by Dolbeault type is
$\Omega^k(M;\CC) = \oplus_{\ell+m = k} \Omega^{\ell,m}$.
Let $\pi^{\ell,m}: \Lambda^k(T^*M \otimes \CC) \to \Lambda^{\ell,m}$
be the projection map, where $\ell + m = k$.
The usual exterior derivative $d$
(extended linearly to smooth complex-valued forms)
composed with two of these projections
induces the \textbf{del} and \textbf{del-bar}
differential operators, $\del$ and $\delbar$,
on forms of type $(\ell, m)$:
\[
\begin{array}{rcl}
        \del & := & \pi^{\ell+1,m}\circ d :
        \Omega^{\ell,m} \longrightarrow \Omega^{\ell+1,m}
\; \quad \mbox{ and } \\
        \delbar & := & \pi^{\ell,m+1}\circ d :
        \Omega^{\ell,m} \longrightarrow \Omega^{\ell,m+1} \ .
\end{array}
\]
If $\beta \in \Omega^{\ell,m} (M)$, with $k = \ell+m$,
then $d\beta \in \Omega^{k+1}(M ; \CC)$:
\[
        d\beta = \displaystyle{\sum_{r+s = k+1}} \pi^{r,s}d\beta =
        \pi^{k+1,0}d\beta + \cdots + \del\beta + \delbar\beta+ \cdots + 
        \pi^{0,k+1}d\beta\ .
\]

In particular, on complex-valued functions we have
$df = d (\mathrm{Re} f) +i \, d(\mathrm{Im} f)$
and $d = \del + \delbar$, where
$\del = \pi^{1,0}\circ d$ and $\delbar = \pi^{0,1}\circ d$.
A function $f: M \to \CC$ is
\textbf{$J$-holomorphic at} $p \in M$\index{J-holomorphic
function@$J$-holomorphic function}
if $df_p$ is complex linear, i.e., $df_p\circ J_p = i \, df_p$
(or $df_p \in T_p^{1,0}$).
A function $f$ is \textbf{$J$-holomorphic}
if it is holomorphic at all $p \in M$.
A function $f: M \to \CC$ is
\textbf{$J$-anti-holomorphic at} $p \in M$\index{J-anti-holomorphic
function@$J$-anti-holomorphic function}
if $df_p$ is complex antilinear, i.e., $df_p\circ J_p = -i \, df_p$
(or $df_p \in T_p^{0,1}$), that is, when the conjugate function
$\bar{f}$ is holomorphic at $p \in M$.
In terms of $\del$ and $\delbar$,
a function $f$ is $J$-holomorphic if and only if $\delbar f = 0$,
and $f$ is $J$-anti-holomorphic if and only if $\del f = 0$.

%%%%%%%%%%%%%%%%%%%%%%%%%%%%%%%%%%%%%%%%%%%%%%%%%%%%%%%%%%%%%%%%%%%%%%%%%%%%%

When $M$ is a {\em complex manifold} and $J$ is its
canonical almost complex structure, the splitting
$\Omega^k(M;\CC) = \oplus_{\ell+m=k} \Omega^{\ell,m}$\index{form !
on a complex manifold}
is particularly interesting.
Let $\cU \subseteq M$ be the domain of a complex coordinate chart
$\varphi = (z_1,\ldots,z_n)$, where the corresponding real coordinates
$x_1,y_1,\ldots,x_n,y_n$ satisfy $z_j = x_j + iy_j$.
In terms of
\[
        \displaystyle{
        \frac {\partial}{\partial z_j} :=
        \frac {1}{2} \left( \frac {\partial}{\partial x_j}
        - i\frac {\partial}{\partial y_j} \right)
        \quad \mbox{ and } \quad
        \frac {\partial}{\partial {\bar z}_j} :=
        \frac {1}{2} \left( \frac {\partial}{\partial x_j} +
        i\frac {\partial}{\partial y_j} \right)\ ,}
\]
the $(\pm i)$-eigenspaces of $J_p$ ($p \in \cU$) can be written
\[
   (T_{1,0})_p =
\CC \mbox{-span}\left\{ \left. \frac {\partial}{\partial z_j}
\right|_p: j = 1,\ldots,n\right\} \quad \mbox{ and } \quad
   (T_{0,1})_p =
\CC \mbox{-span}\left\{ \left. \frac {\partial}{\partial {\bar z}_j}
\right|_p \right\}\ .
\]
Similarly, putting $dz_j = dx_j + idy_j$ and
$d{\bar z}_j = dx_j - idy_j$, we obtain
simple formulas for the differentials of a $b \in C^{\infty}({\cU};\CC)$,
$\partial b = \sum \frac {\partial b}{\partial z_j} dz_j$ and
${\bar \partial}b =
\sum \frac {\partial b}{\partial {\bar z}_j} d{\bar z}_j$, and
we have $T^{1,0} = \CC \mbox{-span}\{dz_j: j = 1,\ldots,n\}$ and
$T^{0,1} = \CC \mbox{-span}\{d{\bar z}_j: j = 1,\ldots,n\}$.
If we use multi-index notation $J = (j_1,\ldots,j_\ell)$
where $1 \leq j_1 < \ldots < j_\ell \leq n$, $|J| = \ell$ and
$dz_{_J} = dz_{j_1} \wedge dz_{j_2} \wedge \ldots \wedge dz_{j_\ell}$,
then the set of $(\ell,m)$-forms on $\cU$ is
\[
        \Omega^{\ell,m} =
        \left\{ \displaystyle{\sum_{|J| = \ell,|K| = m}}
        b_{_{J,K}} dz_{_J} \wedge d{\bar z}_{_K}
        \mid b_{_{J,K}} \in C^{\infty}({\cU};\CC) \right\}\ .
\]
A form $\beta \in \Omega^k(M;\CC)$ may be written over $\cU$ as
\[
        \beta = \displaystyle{\sum_{\ell+m=k}}
        \left( \displaystyle{\sum_{|J| = \ell, |K| = m}}
        b_{_{J,K}}dz_{_J} \wedge d{\bar z}_{_K}\right) \ .
\]
Since $d = \partial + {\bar \partial}$ on functions, we get
\[
\begin{array}{rl}
        d\beta = & \displaystyle{\sum_{\ell+m=k}}
        \left( \displaystyle{\sum_{|J| = \ell, |K| = m}}
        db_{_{J,K}} \wedge dz_{_J} \wedge d{\bar z}_{_K}\right) \\
        \\
        = & \displaystyle{\sum_{\ell+m=k}}
        \underbrace{\left( \displaystyle{\sum_{|J| = \ell, |K| = m}}
        \partial b_{_{J,K}} \wedge dz_{_J} \wedge
        d{\bar z}_{_K} \right.}_{\in \Omega^{\ell+1,m}} +
        \underbrace{\left. \displaystyle{\sum_{|J| = \ell, |K| = m}}
        {\bar \partial}b_{_{J,K}} \wedge dz_{_J}
        \wedge d{\bar z}_{_K} \right)}_{\in \Omega^{\ell,m+1}} \\
        \\
        = & \partial \beta + {\bar \partial}\beta\ ,
\end{array}
\]
and conclude that, {\em on a complex manifold,
$d = \partial + {\bar \partial}$
on forms of any degree}.\index{complex ! differentials}
This cannot be proved for an almost complex manifold, because there are
no coordinate functions $z_j$ to give a suitable basis of 1-forms.

When $d = \del + \delbar$, for any form $\beta \in \Omega^{\ell,m}$, we have
\[
        0 = d^2\beta =
        \underbrace{\del^2 \beta}_{\in \Omega^{\ell+2,m}} +
        \underbrace{\del \delbar \beta + \delbar\del \beta}_
        {\in \Omega^{\ell+1,m+1}} +
        \underbrace{\delbar^2 \beta}_{\in \Omega^{\ell,m+2}}
\quad \Longrightarrow \quad
\left\{ \begin{array}{l}
        \delbar^2 = 0 \\
        \del \delbar + \delbar\del = 0 \\
        \del^2 = 0
\end{array} \right.
\]
Since $\delbar^2 = 0$, the chain
$0 \longrightarrow \Omega^{\ell,0} \stackrel{\delbar}{\longrightarrow}
\Omega^{\ell,1} \stackrel{\delbar}{\longrightarrow}
\Omega^{\ell,2} \stackrel{\delbar}{\longrightarrow} \cdots$
is a differential complex.
Its cohomology groups
\[
        H^{\ell,m}_{\mathrm{Dolbeault}}(M) :=
        \frac{\ker ~ \delbar: \Omega^{\ell,m} \to \Omega^{\ell,m+1}}
        {\mathrm{im} ~\delbar:\Omega^{\ell,m-1} \to \Omega^{\ell,m}}
\]
are called the \textbf{Dolbeault cohomology} groups.\index{Dolbeault !
cohomology}\index{cohomology ! Dolbeault}
The Dolbeault theorem\index{theorem !
Dolbeault}\index{Dolbeault ! theorem}
states that for complex manifolds
$H_{\mathrm{Dolbeault}}^{\ell,m} (M) \simeq H^m(M; \cO(\Omega^{(\ell,0)}))$,
where $\cO(\Omega^{(\ell,0)})$ is the sheaf of forms of type
$(\ell,0)$ over $M$.

It is natural to ask whether the identity $d = \partial + {\bar \partial}$
could hold for manifolds other than complex manifolds.
Newlander and Nirenberg\index{Newlander-Nirenberg theorem}~\cite{ne-ni:complex}
showed that the answer is no:
for an almost complex manifold $(M,J)$, the following are equivalent
\[
        \mbox{$M$ is a complex manifold }
        \iff \cN \equiv 0
        \iff d = \partial + {\bar \partial}
        \iff {\bar \partial}^2 = 0 \ ,
\]
where $\cN$ is the \textbf{Nijenhuis tensor}\index{Nijenhuis tensor}:
\[
        {\cal N} (X,Y) := [JX, JY] - J[JX,Y] - J[X,JY] - [X,Y]\ ,
\]
for vector fields $X$ and $Y$ on $M$, $[ \cdot , \cdot ]$ being
the usual bracket.\footnote{The \textbf{bracket} of vector fields
$X$ and $Y$ is the vector field $[X,Y]$
characterized by the property that
$\cL_{[X,Y]} f := \cL_X (\cL_Y f) - \cL_Y (\cL_X f)$,
for $f \in C^\infty (M)$, where $\cL_X f = df (X)$.}
The Nijenhuis tensor can be thought of as a measure
of the existence of $J$-holomorphic functions:
if there exist $n$ $J$-holomorphic functions, $f_1, \ldots ,f_n$,
on $\RR^{2n}$, that are independent at some point $p$,
i.e., the real and imaginary parts of $(df_1)_p, \ldots ,(df_n)_p$
form a basis of $T^*_p \RR^{2n}$, then ${\cal N}$ vanishes identically at $p$.
More material related to Dolbeault theory or to the
Newlander-Nirenberg theorem can be found
in~\cite{ch:potential,du:heat,gr-ha:principles,ho:several,we:complex}.

\begin{example}\index{example !
of almost complex manifold}\index{example !
of non-almost-complex manifold}
Out of all spheres, only $S^2$ and $S^6$ admit almost
complex structures~\cite[\S41.20]{st:fibre_bundles}.
As a complex manifold, $S^2$ if referred to as the
{\em Riemann sphere} $\CC \PP^1$.
The almost complex structure on $S^6$
from Example~3 of Section~\ref{sec:compatible_almost}
is not integrable, but it is not yet known whether
$S^6$ admits a structure of complex manifold.
\end{example}

In the (real) 2-dimensional case ${\cal N}$ always vanishes
simply because ${\cal N}$ is a tensor, i.e.,
${\cal N} (fX,gY) = fg {\cal N} (X,Y)$ for any $f,g \in C^\infty (M)$,
and ${\cal N} (X,JX) =0$ for any vector field $X$.
Combining this with the fact that any orientable surface is symplectic,
we conclude that any orientable surface is a complex manifold,
a result already known to Gauss.
However, most almost complex structures on higher dimensional manifolds
are not integrable.
In Section~\ref{sec:hodge}
we see that the existence of a complex structure compatible
with a symplectic structure on a compact manifold imposes
significant topological constraints.

%%%%%%%%%%%%%%%%%%%%%%%%%%%%%%%%%%%%%%%%%%%%%%%%%%%%%%%%%%%%%%%%%%%%%%%%%%%%%
%%%%%%%%%%%%%%%%%%%%%%%%%%%%%%%%%%%%%%%%%%%%%%%%%%%%%%%%%%%%%%%%%%%%%%%%%%%%%

\ssubsection{K\"ahler Manifolds}
\label{sec:kahler}
\index{K\"ahler ! manifold}\index{manifold ! K\"ahler}

\begin{definition}
A \textbf{K\"ahler manifold}\index{K\"ahler !
manifold}\index{manifold ! K\"ahler} is a symplectic manifold $(M,\omega)$
equipped with an integrable compatible almost complex structure $J$.
The symplectic form $\omega$ is then called a
\textbf{K\"ahler form}.\index{K\"ahler ! form}\index{form ! K\"ahler}
\end{definition}

As a complex manifold, a K\"ahler manifold $(M,\omega,J)$
has Dolbeault cohomology.
As it is also a symplectic manifold, it is interesting to
understand where the symplectic form $\omega$ sits
with respect to the Dolbeault type decomposition.

\begin{proposition}
A K\"ahler form $\omega$ is a $\partial$- and ${\bar \partial}$-closed
$(1,1)$-form that is given on a local complex chart $(\cU,z_1,\ldots,z_n)$ by
\[
   \omega = \frac i2 \sum_{j,k=1}^n h_{jk}\ dz_j \wedge d{\bar z}_k
\]
where, at every point $p \in \cU$, $(h_{jk}(p))$ is a positive-definite
hermitian matrix.
\end{proposition}

In particular, $\omega$ defines a Dolbeault $(1,1)$-cohomology class,
$[\omega] \in H_{\mathrm{Dolbeault}}^{1,1}(M)$.

\begin{proof}
Being a form in $\Omega^2(M;\CC) =
\Omega^{2,0} \oplus \Omega^{1,1} \oplus \Omega^{0,2}$,
with respect to a local complex chart, $\omega$ can be written
\[
        \omega = \sum a_{jk} \ dz_j \wedge dz_k + \sum b_{jk} \ dz_j \wedge
        d{\bar z}_k + \sum c_{jk} \ d{\bar z}_j \wedge d{\bar z}_k
\]
for some $a_{jk},b_{jk},c_{jk} \in C^{\infty}(\cU;\CC)$.
By the compatibility of $\omega$ with the complex structure,
$J$ is a symplectomorphism, that is, $J^*\omega = \omega$
where $(J^* \omega ) (u,v) := \omega (Ju,Jv)$.
Since $J^*dz_j = dz_j \circ J = idz_j$ and
$J^*d{\bar z}_j = d{\bar z}_j \circ J = -id{\bar z}_j$,
we have $J^*\omega = \omega$ if and only if the coefficients
$a_{jk}$ and $c_{jk}$ all vanish identically, that is,
if and only if $\omega \in \Omega^{1,1}$.
Since $\omega$ is closed, of type $(1,1)$ and
$d\omega = \partial \omega + {\bar \partial}\omega$,
we must have $\partial\omega = 0$ and ${\bar \partial}\omega = 0$.
Set $b_{jk} = \frac i2 h_{jk}$.
As $\omega$ is real-valued, i.e.,
$\omega = \frac i2 \sum h_{jk}\ dz_j \wedge d{\bar z}_k$ and
${\overline \omega} = -\frac {i}{2} \sum
\overline{h_{jk}} \ d{\bar z}_j \wedge dz_k$ coincide,
we must have $h_{jk} = \overline{h_{kj}}$ for all $j$ and $k$.
In other words, at every point $p \in \cU$,
the $n \times n$ matrix $(h_{jk}(p))$ is hermitian.
The nondegeneracy amounts to the nonvanishing of
\[
        \omega^n = n!\left( \frac {i}{2} \right)^n
        \textstyle{\det} (h_{jk})
        \, dz_1 \wedge d{\bar z}_1 \wedge \ldots
        \wedge dz_n \wedge d{\bar z}_n \ .
\]
Therefore, at every $p \in M$, the matrix $(h_{jk}(p))$ must be nonsingular.
Finally, the positivity condition $\omega(v,Jv) > 0$, $\forall v \neq 0$,
from compatibility, implies that, at each $p \in \cU$,
the matrix $(h_{jk}(p))$ is positive-definite.
\end{proof}

%%%%%%%%%%%%%%%%%%%%%%%%%%%%%%%%%%%%%%%%%%%%%%%%%%%%%%%%%%%%%%%%%%%%%%%%%%%%%

Consequently,
%there is the result due to Banyaga
%\index{theorem ! Banyaga}\index{Banyaga theorem}
if $\omega_0$ and $\omega_1$ are both K\"ahler forms
on a compact manifold $M$
with $[\omega_0] = [\omega_1] \in H_{\mathrm{deRham}}^2(M)$,
then $(M,\omega_0)$ and $(M,\omega_1)$ are strongly isotopic
by Moser's Theorem~\ref{thm:moser}.
Indeed $\omega_t = (1-t)\omega_0 + t\omega_1$ is
symplectic for $t \in [0,1]$, as convex combinations
of positive-definite matrices are still positive-definite.

%%%%%%%%%%%%%%%%%%%%%%%%%%%%%%%%%%%%%%%%%%%%%%%%%%%%%%%%%%%%%%%%%%%%%%%%%%%%%

Another consequence is the following recipe
for K\"ahler forms.\index{recipe
! for K\"ahler forms}\index{K\"ahler ! recipe}
A smooth real function $\rho$ on a complex manifold $M$
is \textbf{strictly plurisubharmonic}\index{strictly
plurisubharmonic}\index{potential ! strictly plurisubharmonic}
(\textbf{s.p.s.h.})\index{s.p.s.h.}
if, on each local complex chart $(\cU,z_1,\ldots,z_n)$, the matrix
$\left( \frac {\partial^2\rho}{\partial z_j \partial {\bar z}_k} (p)\right)$
is positive-definite at all $p \in \cU$.
If $\rho \in C^{\infty}(M;\RR)$ is s.p.s.h., then the form
\[
        \omega = \frac i2 \partial {\bar \partial} \rho
\]
is K\"ahler.
The function $\rho$ is then called a (global)
\textbf{K\"ahler potential}.\index{K\"ahler !
potential}\index{potential ! K\"ahler}

\begin{example}
Let $M = \CC^n \simeq \RR^{2n}$, with complex coordinates
$(z_1,\ldots,z_n)$ and corresponding real coordinates
$(x_1,y_1,\ldots,x_n,y_n)$ via $z_j = x_j + iy_j$.
The function
\[
        \rho(x_1,y_1,\ldots,x_n,y_n) = \sum_{j=1}^n (x_j^2 + y_j^2)
        = \sum |z_j|^2 = \sum z_j{\bar z}_j
\]
is s.p.s.h.\ and is a K\"ahler potential for the standard K\"ahler form:
\[
   \frac i2 \partial{\bar \partial}\rho
   = \frac i2 \sum \limits_{j,k} \delta_{jk} \ dz_j \wedge d{\bar z}_k
   = \frac i2 \sum \limits_j dz_j \wedge d{\bar z}_j
   = \sum \limits_j dx_j \wedge dy_j
   = \omega_0 \ .
\]
\end{example}

%%%%%%%%%%%%%%%%%%%%%%%%%%%%%%%%%%%%%%%%%%%%%%%%%%%%%%%%%%%%%%%%%%%%%%%%%%%%%

There is a local converse to the previous construction of K\"ahler forms.

\begin{proposition}
\index{local form}\index{K\"ahler ! local form}
Let $\omega$ be a closed real-valued $(1,1)$-form
on a complex manifold $M$ and let $p \in M$.
Then on a neighborhood $\cU$ of $p$ we have
$\omega = \frac i2 \partial{\bar \partial}\rho$
for some $\rho \in C^{\infty}(\cU;\RR)$.
\end{proposition}

The proof of this theorem requires holomorphic versions of Poincar\'e's lemma,
namely, the local triviality of Dolbeault groups
(the fact that any point in a complex manifold admits a
neighborhood $\cU$ such that
$H^{\ell,m}_{\mathrm{Dolbeault}} (\cU) = 0$ for all $m >0$)
and the local triviality of the holomorphic de Rham groups;
see~\cite{gr-ha:principles}.

For a K\"ahler $\omega$, such a local function $\rho$ is called a
\textbf{local K\"ahler potential}.\index{K\"ahler !
potential}\index{potential ! K\"ahler}

\begin{proposition}
Let $M$ be a complex manifold,
$\rho \in C^{\infty}(M;\RR)$ s.p.s.h., $X$ a complex submanifold, and
$i: X \hookrightarrow M$ the inclusion map.
Then $i^*\rho$ is s.p.s.h..
\end{proposition}

\begin{proof}
It suffices to verify this locally by considering
a complex chart $(z_1,\ldots,z_n)$ for $M$ adapted to $X$
so that $X$ is given there by the equations $z_1 = \ldots = z_m =0$.
Being a principal minor of the positive-definite matrix
$\left( \frac {\partial^2}{\partial z_j\partial {\bar z}_k}
(0,\ldots,0,z_{m+1},\ldots,z_n)\right)$
the matrix
$\left( \frac
        {\partial^2\rho}{\partial z_{m+j}\partial {\bar z}_{m+k}}
        (0,\ldots,0,z_{m+1},\ldots,z_n)\right)$
is also positive-definite.
\end{proof}

\vspace*{-1ex}

\begin{corollary}\label{thm:subkahler}
Any complex submanifold of a K\"ahler manifold is also
K\"ahler.\index{example ! complex submanifold of a K\"ahler manifold}
\end{corollary}

\vspace*{-1ex}

\begin{definition}
Let $(M,\omega)$ be a K\"ahler manifold,
$X$ a complex submanifold, and $i: X \hookrightarrow M$ the inclusion.
Then $(X, i^* \omega)$ is called a
\textbf{K\"ahler submanifold}.\index{K\"ahler !
submanifold}\index{submanifold ! K\"ahler}
\end{definition}

\vspace*{-1ex}

%%%%%%%%%%%%%%%%%%%%%%%%%%%%%%%%%%%%%%%%%%%%%%%%%%%%%%%%%%%%%%%%%%%%%%%%%%%%%

\begin{examples}\index{example ! of K\"ahler submanifold}
\begin{enumerate}
\item
Complex vector space $(\CC^n,\omega_0)$ where
$\omega_0 = \frac i2 \sum dz_j \wedge d{\bar z}_j$ is K\"ahler.
According to Corollary~\ref{thm:subkahler},
every complex submanifold of $\CC^n$ is K\"ahler.

\item
In particular, {\em Stein manifolds}\index{example !
Stein manifold}\index{Stein manifold} are K\"ahler.
\textbf{Stein manifolds} are the properly embedded
complex submanifolds of $\CC^n$.
They can be alternatively characterized as being
the K\"ahler manifolds $(M,\omega)$
that admit a global proper K\"ahler potential,
i.e., $\omega = \frac {i}{2} \partial{\bar \partial} \rho$ for some
proper function $\rho: M \to \RR$.

\item
The function $z \mapsto \log (|z|^2 + 1)$ on $\CC^n$
is strictly plurisubharmonic.
Therefore the 2-form
\[
        \omega_{_{\mathrm{FS}}} = \textstyle{\frac{i}{2}}
        \partial \bar{\partial} \log (|z|^2 + 1)
\]
is another K\"ahler form on $\CC^n$
This is called the \textbf{Fubini-Study form} on $\CC^n$.

\item
Let $\{ ({\cal U}_j , \CC^n , \varphi_j), j=0, \ldots ,n\}$
be the usual complex atlas for \textbf{complex projective
space}.\index{complex ! projective space}\index{example !
complex projective space}\footnote{The \textbf{complex projective space}
$\CC \PP ^n$ is the complex $n$-dimensional
manifold\index{example ! of complex manifold}
given by the space of complex lines in $\CC^{n+1}$.
It can be obtained from
$\CC^{n+1} \setminus \{ 0 \}$ by making the identifications
$(z_0, \ldots, z_n) \sim (\lambda z_0, \ldots, \lambda z_n)$
for all $\lambda \in \CC \setminus \{ 0 \}$.
One denotes by $[z_0, \ldots, z_n]$ the equivalence class of
$(z_0, \ldots, z_n)$, and calls $z_0, \ldots, z_n$ the
\textbf{homogeneous coordinates} of the point $p = [z_0, \ldots, z_n]$.
(Homogeneous coordinates are, of course, only determined up to
multiplication by a non-zero complex number $\lambda$.)
Let ${\cal U}_j$ be the subset of $\CC \PP ^n$ consisting
of all points $p = [z_0, \ldots, z_n]$ for which $z_j \neq 0$.  Let
$\varphi_j : {\cal U}_j \to \CC^n$ be the map defined by
\[
        \varphi_j ([z_0, \ldots, z_n]) = \displaystyle{
        \left( \textstyle{\frac{z_0}{z_j}} , \ldots ,
        \textstyle{\frac{z_{j-1}}{z_j}} ,
        \textstyle{\frac{z_{j+1}}{z_j}} , \ldots ,
        \textstyle{\frac{z_n}{z_j}} \right)} \ .
\]
The collection $\{ ({\cal U}_j , \CC^n , \varphi_j), j=0, \ldots ,n\}$
is the \textbf{usual complex atlas}\index{complex ! atlas} for $\CC \PP ^n$.
For instance, the transition map from
$({\cal U}_0 , \CC^n , \varphi_0)$ to $({\cal U}_1 , \CC^n , \varphi_1)$
is $\varphi_{0,1} (z_1 , \ldots, z_n) =
(\textstyle{\frac{1}{z_1}} , \textstyle{\frac{z_2}{z_1}},
\ldots, \textstyle{\frac{z_n}{z_1}})$
defined from the set
$\{ (z_1 , \ldots, z_n) \in \CC^n \, | \, z_1 \neq 0\}$ to itself.}
The form $\omega_{_{\mathrm{FS}}}$ is preserved by
the transition maps, hence $\varphi_j ^* \omega_{_{\mathrm{FS}}}$
and $\varphi_k ^* \omega_{_{\mathrm{FS}}}$ agree on the
overlap ${\cal U}_j \cap {\cal U}_k$.
The \textbf{Fubini-Study form} on $\CC \PP ^n$
is the K\"ahler form obtained by gluing together
the $\varphi_j ^* \omega_{_{\mathrm{FS}}}$, $j=0, \ldots, n$.

\item
Consequently, all \textbf{non-singular projective
varieties}\index{non-singular projective variety}\index{example !
non-singular projective variety} are K\"ahler submanifolds.
Here by non-singular we mean smooth,
and by projective variety we mean the zero locus of some
collection of homogeneous polynomials.

\item
All \textbf{Riemann surfaces}\index{example !
Riemann surface}\index{Riemann ! surface} are K\"ahler,
since any compatible almost complex structure is integrable
for dimension reasons (Section~\ref{sec:integrability}).

\item
%$\CC \PP ^1$ is diffeomorphic to $S^2$ as a real 2-dimensional manifold.
The Fubini-Study form
on the chart ${\cal U}_0 = \{ [z_0, z_1] \in \CC \PP ^1 \, | z_0 \neq 0 \}$
of the \textbf{Riemann sphere}\index{Riemann sphere}
$\CC \PP ^1$ is given by the formula
\[
        \omega_{_{\mathrm{FS}}} = \frac{ dx \wedge dy }{(x^2+y^2+1)^2}
\]
where $\frac{z_1}{z_0} = z = x+iy$ is the usual coordinate on $\CC$.
The standard area form $\omega_{_{\mathrm{std}}} = d \theta \wedge dh$
is induced by regarding $\CC \PP ^1$ as the unit sphere $S^2$ in $\RR^3$
(Example~3 of Section~\ref{symplectic_forms}).
Stereographic projection\index{stereographic projection}
shows that $\omega_{_{\mathrm{FS}}} = \frac{1}{4} \omega_{_{\mathrm{std}}}$.

\item
        \textbf{Complex tori}\index{example !
complex torus}\index{complex torus} are K\"ahler.
Complex tori look like quotients $\CC^n / \ZZ^n$ where
$\ZZ^n$ is a lattice in $\CC^n$.
The form $\omega = \sum dz_j \wedge d \bar z_j$ induced by the
euclidean structure is K\"ahler.

\item
        Just like products of symplectic manifolds are symplectic,
also products of K\"ahler manifolds\index{example !
product of K\"ahler manifolds} are K\"ahler.

\end{enumerate}
\end{examples}

%%%%%%%%%%%%%%%%%%%%%%%%%%%%%%%%%%%%%%%%%%%%%%%%%%%%%%%%%%%%%%%%%%%%%%%%%%%%%
%%%%%%%%%%%%%%%%%%%%%%%%%%%%%%%%%%%%%%%%%%%%%%%%%%%%%%%%%%%%%%%%%%%%%%%%%%%%%

\ssubsection{Hodge Theory}
\label{sec:hodge}

Hodge~\cite{ho:theory} identified the spaces of cohomology classes of
forms with spaces of actual forms, by picking {\em the} representative
from each class that solves a certain differential equation, namely
the {\em harmonic} representative.\index{harmonic form}\index{form !
harmonic}
We give a sketch of Hodge's idea.
The first part makes up ordinary Hodge theory,
which works for any compact oriented riemannian manifold $(M,g)$,
not necessarily K\"ahler.

At a point $p \in M$, let $e_1,\ldots,e_n$ be a positively
oriented orthonormal basis of the cotangent space $T^*_pM$,
with respect to the induced inner product and orientation.
The \textbf{Hodge star operator}\index{Hodge !
star operator@$\ast$-operator}
is the linear operator on the exterior algebra of $T^*_pM$ defined by
\[
\begin{array}{rcl}
        \ast (1) & = & e_1\wedge \ldots \wedge e_n \\
        \ast (e_1\wedge \ldots \wedge e_n) & = & 1 \\
        \ast (e_1 \wedge \ldots \wedge e_k) & = &
        e_{k+1} \wedge \ldots \wedge e_n \ .
\end{array}
\]
We see that $\ast: \Lambda^k(T^*_pM) \to \Lambda^{n-k}(T^*_pM)$ and
satisfies $\ast\ast =(-1)^{k(n-k)}$.
The \textbf{codifferential}\index{codifferential}
and the \textbf{laplacian}\index{laplacian}
are the operators defined by
\[
\begin{array}{cclcl}
        \delta & = & (-1)^{n(k+1)+1} \ast d \ast & : &
        \Omega^k(M) \to \Omega^{k-1}(M) \  , \\
        \Delta & = & d\delta + \delta d & : &
        \Omega^k(M) \to \Omega^k(M)\ .
\end{array}
\]
The operator $\Delta$ is also called the
\textbf{Laplace-Beltrami operator} and satisfies
$\Delta \ast = \ast\Delta$.\index{Laplace-Beltrami
operator}\index{Beltrami ! Laplace-Beltrami
operator}\index{operator ! Laplace-Beltrami}
On $\Omega^0(\RR^n) = C^{\infty}(\RR^n)$, it is simply the usual laplacian
$\Delta = -\sum_{i=1}^n \frac {\partial^2}{\partial x_i^2}$.
The \textbf{inner product on forms} of any degree,
\[
        \langle \cdot , \cdot \rangle:
        \Omega^k (M) \times \Omega^k (M) \longrightarrow \RR \ ,
        \qquad
        \langle \alpha,\beta \rangle := \int_M \alpha \wedge \ast\beta \ ,
\]
satisfies
$\langle d\alpha,\beta \rangle = \langle \alpha,\delta\beta\rangle$,
so the codifferential $\delta$ is often denoted by
$d^*$ and called the {\em adjoint\footnote{When $M$ is not compact,
we still have a {\em formal adjoint} of $d$
with respect to the nondegenerate bilinear pairing
$\langle \cdot , \cdot \rangle: \Omega^k (M) \times \Omega^k_c (M) \to \RR$
defined by a similar formula, where $\Omega^k_c (M)$ is the space of
compactly supported $k$-forms.} of $d$}.
Also, $\Delta$ is self-adjoint (i.e.,
$\langle \Delta\alpha,\beta\rangle = \langle \alpha,\Delta\beta\rangle$),
and
$\langle \Delta\alpha,\alpha\rangle = |d\alpha|^2 +|\delta\alpha|^2 \geq 0$,
where $| \cdot |$ is the norm with respect to this inner product.
The \textbf{harmonic $k$-forms}\index{form !
harmonic}\index{harmonic form} are the elements of
${\cH}^k := \{\alpha \in \Omega^k \mid \Delta\alpha = 0\}$.
Note that $\Delta\alpha = 0$ if and only if $d\alpha = \delta\alpha = 0$.
Since a harmonic form is $d$-closed, it defines a de Rham cohomology class.

\begin{theorem}\index{Hodge ! theorem}\index{theorem ! Hodge}
\textbf{(Hodge)} $\;$
        Every de Rham cohomology class on a compact oriented riemannian
manifold $(M,g)$ possesses a unique harmonic representative, i.e.,
there is an isomorphism $\cH^k \simeq H_{\mathrm{deRham}}^k(M;\RR)$.
In particular, the spaces $\cH^k$ are finite-dimensional.
We also have the following orthogonal decomposition
with respect to the inner product on forms:
$\Omega^k \simeq \cH^k \oplus \Delta(\Omega^k(M))
\simeq \cH^k \oplus d\Omega^{k-1} \oplus \delta\Omega^{k+1}$.
\end{theorem}

This decomposition is called the \textbf{Hodge decomposition on
forms}.\index{Hodge ! decomposition}
The proof of this and the next theorem involves functional analysis,
elliptic differential operators, pseudodifferential operators and
Fourier analysis;
see for instance~\cite{gr-ha:principles,ko:harmonic,we:complex}.

Here is where \textbf{complex Hodge theory}\index{Hodge ! complex
Hodge theory}\index{complex ! Hodge theory} begins.
When $M$ is K\"ahler, the laplacian satisfies
$\Delta = 2( \bar \partial \bar \partial ^* + \bar \partial ^* \bar \partial)$
(see, for example,~\cite{gr-ha:principles})
and preserves the decomposition according to type,
$\Delta: \Omega^{\ell,m} \to \Omega^{\ell,m}$.
Hence, harmonic forms are also bigraded
\[
        \cH^k = \bigoplus_{\ell+m=k} \cH^{\ell,m}\ .
\]
and satisfy a K\"unneth formula $\cH^{\ell,m} (M \times N) \simeq
\bigoplus_{p+r=\ell,q+s=m} \cH^{p,q} (M) \otimes \cH^{r,s} (N)$.

\begin{theorem}\index{Hodge ! theorem}\index{theorem ! Hodge}
\textbf{(Hodge)} $\;$
Every Dolbeault cohomology class on a compact K\"ahler
manifold $(M,\omega)$ possesses a unique harmonic representative, i.e.,
there is an isomorphism
$\cH^{\ell,m} \simeq H_{\mathrm{Dolbeault}}^{\ell,m}(M)$.
\end{theorem}

Combining the two theorems of Hodge, we find the decomposition
of cohomology groups for a compact K\"ahler manifold
\[
   H_{\mathrm{deRham}}^k (M;\CC)
%\simeq \cH^\ell = \bigoplus_{\ell+m=k} \cH^{\ell,m}
   \simeq \bigoplus_{\ell+m=k} H_{\mathrm{Dolbeault}}^{\ell,m}(M) \ ,
\]
known as the \textbf{Hodge decomposition}\index{Hodge ! decomposition}.
In particular, the Dolbeault cohomology groups
$H_{\mathrm{Dolbeault}}^{\ell,m}$
are finite-dimensional and $H^{\ell,m} \simeq \overline{H^{m,\ell}}$.

%%%%%%%%%%%%%%%%%%%%%%%%%%%%%%%%%%%%%%%%%%%%%%%%%%%%%%%%%%%%%%%%%%%%%%%%%%%%%

Let $b^k(M) := \dim H_{\mathrm{deRham}}^k(M)$
be the usual \textbf{Betti numbers}\index{Betti number}\index{number !
Betti} of $M$, and let
$h^{\ell,m}(M) := \dim H_{\mathrm{Dolbeault}}^{\ell,m}(M)$ 
be the \textbf{Hodge numbers}\index{Hodge ! number}\index{number !
Hodge} of $M$.

For an arbitrary compact symplectic manifold $(M,\omega)$,
the even Betti numbers must be positive,
because $\omega^k$ is closed but not exact $(k = 0,1,\ldots,n)$.
In fact, if it were $\omega ^k = d\alpha$, by Stokes' theorem
we would have $\int_M \omega^n =
\int_M d (\alpha \wedge \omega^{n-k}) = 0$,
which contradicts $\omega^n$ being a volume form.

For a compact K\"ahler manifold $(M,\omega)$, there are finer
topological consequences coming from the Hodge theorems,\index{Hodge !
theorem}\index{theorem ! Hodge}
as we must have $b^k = \sum_{\ell+m=k} h^{\ell,m}$
and $h^{\ell,m} = h^{m,\ell}$.
The odd Betti numbers must be even because
$b^{2k+1} = \sum_{\ell+m=2k+1} h^{\ell,m} =
2 \sum_{\ell=0}^k h^{\ell,(2k+1-\ell)}$.
The number $h^{1,0} = \frac {1}{2} b^1$ must be a topological invariant.
The numbers $h^{\ell,\ell}$ are positive,
because $0 \neq [\omega^\ell] \in H^{\ell,\ell}_{\mathrm{Dolbeault}} (M)$.
First of all, $[\omega^\ell]$ defines an element of
$H_{\mathrm{Dolbeault}}^{\ell,\ell}$ as
$\omega \in \Omega^{1,1}$ implies that $\omega^\ell \in \Omega^{\ell,\ell}$,
and the closedness of $\omega^\ell$ implies that
${\bar \partial}\omega^\ell = 0$.
If it were $\omega^\ell = {\bar \partial}\beta$ for some
$\beta \in \Omega^{\ell-1,\ell}$, then
$\omega^n = \omega^\ell \wedge \omega^{n-\ell}
= {\bar \partial}(\beta \wedge \omega^{n-\ell})$
would be ${\bar \partial}$-exact.
But $[\omega^n] \ne 0$ in
$H^{2n}_{\mathrm{deRham}}(M ; \CC) \simeq
H^{n,n}_{\mathrm{Dolbeault}}(M)$ since it is a volume form.
A popular diagram to describe relations among Hodge numbers is the
\textbf{Hodge diamond}\index{Hodge ! diamond}:
\[
\begin{array}{ccccccc}
        & & & h^{n,n} \\
        & & h^{n,n-1} & & h^{n-1,n} \\
        & h^{n,n-2} & & h^{n-1,n-1} & & h^{n-2,n} \\
        \ldots & & & \vdots & & & \ldots \\
        \\
        & h^{2,0} & & h^{1,1} & & h^{0,2} \\
        & & h^{1,0} & & h^{0,1} \\
        & & & h^{0,0} \\
\end{array}
\]
Complex conjugation gives symmetry with respect to the middle vertical,
whereas the Hodge star operator
induces symmetry about the center of the diamond.
The middle vertical axis is all non-zero.

There are further symmetries and ongoing research
on how to compute $H_{\mathrm{Dolbeault}}^{\ell,m}$
for a compact K\"ahler manifold $(M,\omega)$.
In particular, the \textbf{hard Lefschetz theorem}
states isomorphisms $L^k : H^{n-k}_{\mathrm{deRham}}(M)
\stackrel{\simeq}{\longrightarrow} H^{n+k}_{\mathrm{deRham}}(M)$
given by wedging with $\omega^ k$ at the level of forms
and the \textbf{Lefschetz decompositions} $H^m_{\mathrm{deRham}}(M)
\simeq \oplus_k L^k (\ker L^{n-m+2k+1} |_{H^{m-2k}})$.
The \textbf{Hodge conjecture}\index{Hodge !
conjecture}\index{conjecture ! Hodge} claims,
for projective manifolds $M$
(i.e., complex submanifolds of complex projective space),
that the Poincar\'e duals of elements in
$H_{\mathrm{Dolbeault}}^{\ell,\ell}(M) \cap H^{2\ell}(M;\QQ)$
are rational linear combinations of classes of complex
codimension $\ell$ subvarieties of $M$.
This has been proved only for the $\ell = 1$ case
(it is the Lefschetz theorem on $(1,1)$-classes;
see for instance~\cite{gr-ha:principles}).

%%%%%%%%%%%%%%%%%%%%%%%%%%%%%%%%%%%%%%%%%%%%%%%%%%%%%%%%%%%%%%%%%%%%%%%%%%%%%
%%%%%%%%%%%%%%%%%%%%%%%%%%%%%%%%%%%%%%%%%%%%%%%%%%%%%%%%%%%%%%%%%%%%%%%%%%%%%

\ssubsection{Pseudoholomorphic Curves}
\label{sec:pseudoholomorphic}
\index{pseudoholomorphic curves}

Whereas an almost complex manifold $(M,J)$ tends to have
no $J$-holomorphic functions $M \to \CC$ at all,\footnote{However,
the study of {\em asymptotically $J$-holomorphic functions}
has been recently developed to obtain important
results~\cite{do:almost,do:pencils,au:branched};
see Section~\ref{sec:pencils}.}
it has plenty of {\em $J$-holomorphic curves}\index{J-holomorphic
curve@$J$-holomorphic curve} $\CC \to M$.
Gromov\index{Gromov !
pseudoholomorphic curve}\index{pseudoholomorphic curve}
first realized that {\em pseudoholomorphic curves}
provide a powerful tool in symplectic topology in
an extremely influential paper~\cite{gr:pseudo}.
Fix a closed Riemann surface $(\Sigma, j)$, that is,
a compact complex 1-dimensional manifold $\Sigma$ without boundary
and equipped with the canonical almost complex structure $j$.

\begin{definition}
A parametrized \textbf{pseudoholomorphic curve}
(or \textbf{ $J$-holomorphic curve})
in $(M,J)$ is a (smooth) map $u : \Sigma \to M$
whose differential intertwines $j$ and $J$,
that is, $du_p \circ j_p = J_p \circ du_p$, $\forall p \in \Sigma$.
\end{definition}

In other words, the \textbf{Cauchy-Riemann equation}
$du + J \circ du \circ j = 0$ holds.

Pseudoholomorphic curves are related to
parametrized 2-dimensional symplectic submanifolds.
If a pseudoholomorphic curve $u : \Sigma \to M$
is an embedding, then its image $S:=u(\Sigma)$
is a 2-dimensional almost complex
submanifold, hence a symplectic submanifold.
Conversely, the inclusion $i : S \hookrightarrow M$ of
a 2-dimensional symplectic submanifold can be seen
as a pseudoholomorphic curve.
An appropriate compatible almost complex structure $J$
on $(M, \omega)$ can be constructed starting from $S$,
such that $TS$ is $J$-invariant.
The restriction $j$ of $J$ to $TS$ is necessarily integrable
because $S$ is 2-dimensional.

The group $G$ of complex diffeomorphisms of $(\Sigma, j)$
acts on (parametrized) pseudoholomorphic curves by reparametrization:
$u \mapsto u \circ \gamma$, for $\gamma \in G$.
This normally means that each curve $u$ has a noncompact orbit under $G$.
The orbit space $\cM_g (A,J)$ is the set
of unparametrized pseudoholomorphic curves
in $(M,J)$ whose domain $\Sigma$ has genus $g$ and whose
image $u(\Sigma)$ has homology class $A \in H_2 (M;\ZZ)$.
The space $\cM_g (A,J)$ is called the
\textbf{moduli space of unparametrized pseudoholomorphic curves}
of genus $g$ representing the class $A$.
For generic $J$, Fredholm theory shows that pseudoholomorphic curves
occur in finite-dimensional smooth families,
so that the moduli spaces $\cM_g (A,J)$ can be manifolds,
after avoiding singularities given by {\em multiple
coverings}.\footnote{A curve $u : \Sigma \to M$ is a
\textbf{multiple covering} if
$u$ factors as $u = u' \circ \sigma$ where $\sigma : \Sigma \to \Sigma'$
is a holomorphic map of degree greater than 1.}

\begin{example}
Usually $\Sigma$ is the Riemann sphere $\CC \PP^1$,
whose complex diffeomorphisms are
those given by {\em fractional linear transformations}
(or {\em M\"obius transformations}).
So the 6-dimensional noncompact
group of projective linear transformations
$\PSL (2;\CC )$ acts on \textbf{pseudoholomorphic spheres}
by reparametrization $u \mapsto u \circ \gamma_A$, where
$A = \textstyle{\left[ \begin{array}{cc} a & b \\ c & d \end{array} \right]}
\in \PSL (2;\CC )$ acts by
$\gamma_A : \CC \PP^1 \to \CC \PP^1$,
$\gamma_A [z,1] = [ \textstyle{\frac{az+b}{cz+d}},1 ]$.
\end{example}

When $J$ is an almost complex structure {\em compatible}
with a symplectic form $\omega$, the area of the image
of a pseudoholomorphic curve $u$ (with respect to the
metric $g_{_J} ( \cdot , \cdot) = \omega (\cdot , J \cdot)$)
is determined by the class $A$ that it represents.
The number
\[
   E(u) := \omega (A) = \int_\Sigma u^* \omega
   = \mbox{ area of the image of } u \mbox{ with respect to } g_{_J}
\]
is called the \textbf{energy}\index{energy} of the curve $u$
and is a topological invariant: it only depends on $[\omega]$
and on the homotopy class of $u$.
Gromov proved that the constant energy of all the
pseudoholomorphic curves representing a homology class $A$
ensured that the space $\cM_g (A,J)$,
though not necessarily compact, had natural \textbf{compactifications}
$\overline \cM_g (A,J)$ by including what he called {\em cusp-curves}.

\begin{theorem} \textbf{(Gromov's compactness theorem)} $\;$
If $(M,\omega)$ is a compact manifold equipped with a generic
compatible almost complex structure $J$, and if $u_j$ is a
sequence of pseudoholomorphic curves in $\cM_g (A,J)$,
then there is a subsequence that weakly converges to
a cusp-curve in $\overline \cM_g (A,J)$.
\end{theorem}
%Gromov-Hausdorff topology

Hence the cobordism class of the compactified moduli space
$\overline \cM_g (A,J)$ might be a nice symplectic invariant of $(M,\omega)$,
as long as it is not empty or null-cobordant.
Actually a nontrivial regularity criterion for $J$
ensures the existence of pseudoholomorphic curves.
And even when $\overline \cM_g (A,J)$ is null-cobordant,
we can define an invariant to be the (signed) number of
pseudoholomorphic curves of genus $g$ in class $A$
that intersect a specified set of representatives of
homology classes in $M$~\cite{ru:sigma_models,ta:sw=gr,wi:sigma_models}.
For more on pseudoholomorphic curves, see for
instance~\cite{mc-sa:curves} (for a comprehensive discussion
of the genus 0 case) or~\cite{au-la:holomorphic} (for higher genus).
Here is a selection of applications of
(developments from) pseudoholomorphic curves:

\begin{itemize}
\item
   Proof of the \textbf{nonsqueezing theorem}~\cite{gr:pseudo}:
for $R > r$ there is no symplectic embedding of a ball $B^{2n}_R$
of radius $R$ into a cylinder $B^{2}_r \times \RR^{2n-2}$
of radius $r$, both in $(\RR^{2n},\omega_0)$.

\item
    Proof that there are {\em no lagrangian spheres}
in $(\CC^n,\omega_0)$, except for the circle in $\CC^2$, and more generally
{\em no compact exact lagrangian submanifolds}, in the sense that
the tautological 1-form $\alpha$ restricts to an exact form~\cite{gr:pseudo}.

\item
   Proof that if $(M,\omega)$ is a connected symplectic 4-manifold
symplectomorphic to $(\RR^4,\omega_0)$ outside a compact set
and containing no symplectic $S^2$'s,
then $(M,\omega)$ symplectomorphic to $(\RR^4,\omega_0)$~\cite{gr:pseudo}.

\item
   Study questions of \textbf{symplectic
packing}~\cite{bi:packing,mc-po:packing,tr:packing} such as:
for a given $2n$-dimensional symplectic manifold $(M,\omega)$,
what is the maximal radius $R$ for which there is a symplectic embedding
of $N$ disjoint balls $B^{2n}_R$ into $(M,\omega)$?

\item
   Study \textbf{groups of symplectomorphisms} of 4-manifolds
(for a review see~\cite{mc:groups}).
Gromov~\cite{gr:pseudo} showed that
$\mathrm{Sympl} (\CC \PP^2, \omega_{_{\mathrm{FS}}})$ and
$\mathrm{Sympl} (S^2 \times S^2, {\mathrm{pr}}_1^* \sigma
\oplus {\mathrm{pr}}_2^* \sigma)$
deformation retract onto the corresponding groups of standard isometries.

\item
   Development of \textbf{Gromov-Witten invariants} allowing to
prove, for instance, the nonexistence of symplectic forms
on $\CC\PP^2 \# \CC\PP^2 \# \CC\PP^2$ or the
classification of symplectic structures on
{\em ruled surfaces} (Section~\ref{sec:blow_up}).

\item
   Development of \textbf{Floer homology} to prove the Arnold conjecture
on the fixed points of symplectomorphisms of compact symplectic manifolds,
or on the intersection of lagrangian submanifolds
(Section~\ref{sec:arnold_floer}).

\item
   Development of \textbf{symplectic field theory}
introduced by Eliashberg, Givental and Hofer~\cite{el-gi-ho:field}
extending Gromov-Witten theory, exhibiting a rich algebraic structure
and also with applications to {\em contact geometry}.

\end{itemize}

%%%%%%%%%%%%%%%%%%%%%%%%%%%%%%%%%%%%%%%%%%%%%%%%%%%%%%%%%%%%%%%%%%%%%%%%%%%%%
%%%%%%%%%%%%%%%%%%%%%%%%%%%%%%%%%%%%%%%%%%%%%%%%%%%%%%%%%%%%%%%%%%%%%%%%%%%%%
% --> Section 4
%%%%%%%%%%%%%%%%%%%%%%%%%%%%%%%%%%%%%%%%%%%%%%%%%%%%%%%%%%%%%%%%%%%%%%%%%%%%%
%%%%%%%%%%%%%%%%%%%%%%%%%%%%%%%%%%%%%%%%%%%%%%%%%%%%%%%%%%%%%%%%%%%%%%%%%%%%%

\newpage

\ssection{Symplectic Geography}
\index{symplectic geography}
\label{section4}

%%%%%%%%%%%%%%%%%%%%%%%%%%%%%%%%%%%%%%%%%%%%%%%%%%%%%%%%%%%%%%%%%%%%%%%%%%%%%
%%%%%%%%%%%%%%%%%%%%%%%%%%%%%%%%%%%%%%%%%%%%%%%%%%%%%%%%%%%%%%%%%%%%%%%%%%%%%

\ssubsection{Existence of Symplectic Forms}
\label{sec:existence}

The utopian goal of symplectic classification
addresses the standard questions:

\begin{itemize}
\item
   {\em (Existence)}
Which manifolds carry symplectic forms?
\item
   {\em (Uniqueness)}
What are the distinct symplectic structures on a given manifold?
\end{itemize}

Existence is tackled through central examples in this section
and symplectic constructions in the next two sections.
Uniqueness is treated in the remainder of this chapter
dealing with invariants that allow to distinguish symplectic manifolds.

A K\"ahler structure naturally yields both
a symplectic form and a complex structure (compatible ones).
Either a symplectic or a complex structure on a manifold
implies the existence of an almost complex structure.
The following diagram represents the relations
among these structures.
In dimension 2, orientability trivially guarantees
the existence of all other structures, so the picture collapses.
In dimension 4, the first interesting dimension,
the picture above is faithful -- we will see that
there are {\em closed} 4-dimensional examples in each region.
\textbf{Closed} here means compact and without boundary.

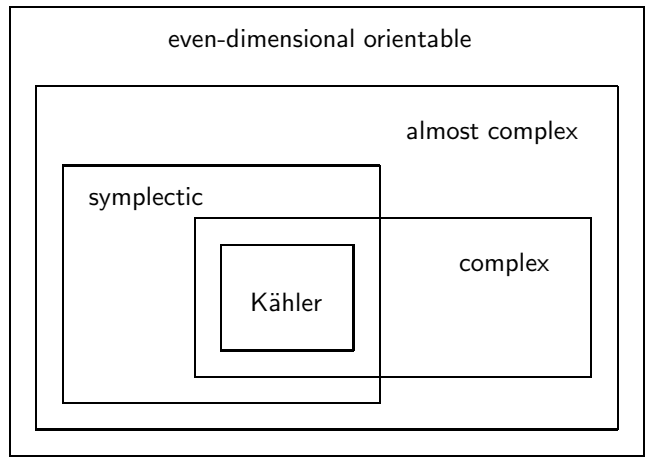
\begin{picture}(200,200)(-20,-5)
% horizontal lines
\put(10,10){\line(1,0){240}}
\put(20,20){\line(1,0){220}}
\put(30,30){\line(1,0){120}}
\put(80,40){\line(1,0){150}}
\put(90,50){\line(1,0){50}}
\put(90,90){\line(1,0){50}}
\put(80,100){\line(1,0){150}}
\put(30,120){\line(1,0){120}}
\put(20,150){\line(1,0){220}}
\put(10,180){\line(1,0){240}}
% vertical lines
\put(10,10){\line(0,1){170}}
\put(20,20){\line(0,1){130}}
\put(30,30){\line(0,1){90}}
\put(80,40){\line(0,1){60}}
\put(90,50){\line(0,1){40}}
\put(140,50){\line(0,1){40}}
\put(230,40){\line(0,1){60}}
\put(150,30){\line(0,1){90}}
\put(240,20){\line(0,1){130}}
\put(250,10){\line(0,1){170}}
% words
\put(40,105){symplectic}
\put(70,165){even-dimensional orientable}
\put(101,65){K\"ahler}
\put(160,130){almost complex}
\put(180,80){complex}
\thicklines
\end{picture}

Not all 4-dimensional manifolds are almost complex.
A result of Wu~\cite{wu:classes} gives a necessary and sufficient
condition in terms of the signature $\sigma$ and the
Euler characteristic $\chi$ of a 4-dimensional closed manifold $M$
for the existence of an almost complex structure:
{\em $3 \sigma + 2 \chi = h^2$ for some $h \in H^2 (M;\ZZ)$
congruent with the second Stiefel-Whitney class $w_2 (M)$ modulo 2}.
For example, $S^4$ and
$(S^2 \times S^2) \# (S^2 \times S^2)$ are not almost complex.
When an almost complex structure exists,
the first Chern class of the tangent bundle
(regarded as a complex vector bundle) satisfies the condition for $h$.
The sufficiency of Wu's condition is the remarkable
part.\footnote{Moreover, such solutions $h$ are in one-to-one
correspondence with
{\em isomorphism} classes of almost complex structures.}

According to Kodaira's classification of closed complex
surfaces~\cite{ko:surfacesI}\index{Kodaira ! complex surfaces},
such a surface admits a K\"ahler structure
if and only if its first Betti number $b_1$ is even.
The necessity of this condition is a Hodge relation on the Betti numbers
(Section~\ref{sec:hodge}).
The complex projective plane $\CC \PP^2$ with the Fubini-Study form
(Section~\ref{sec:kahler})
might be called the simplest example of a closed K\"ahler 4-manifold.

The \textbf{Kodaira-Thurston example}~\cite{th:examples}\index{example !
Kodaira-Thurston}\index{Kodaira ! Kodaira-Thurston
example}\index{Thurston ! Kodaira-Thurston example}
first demonstrated that a manifold
that admits both a symplectic and a complex structure
does not have to admit any K\"ahler structure.
Take $\RR^4$ with $dx_1 \wedge dy_1 + dx_2 \wedge dy_2$, and
$\Gamma$ the discrete group generated by the four symplectomorphisms:
\[
\begin{array}{rcl}
   (x_1,x_2,y_1,y_2) & \longmapsto & (x_1+1,x_2,y_1,y_2) \\
   (x_1,x_2,y_1,y_2) & \longmapsto & (x_1,x_2+1,y_1,y_2) \\
   (x_1,x_2,y_1,y_2) & \longmapsto & (x_1,x_2+y_2,y_1+1,y_2) \\
   (x_1,x_2,y_1,y_2) & \longmapsto & (x_1,x_2,y_1,y_2+1)
\end{array}
\]
Then $M = \RR^4/\Gamma$ is a symplectic manifold
that is a 2-torus bundle over a 2-torus.
Kodaira's classification~\cite{ko:surfacesI}\index{Kodaira ! complex surfaces}
shows that $M$ has a complex structure.
However, $\pi_1 (M) = \Gamma$,
hence $H_1(\RR^4/\Gamma;\ZZ) = \Gamma / [ \Gamma,\Gamma]$
has rank 3, so $b_1 = 3$ is {\em odd}.

Fern\'andez-Gotay-Gray~\cite{fe-go-gr:symplectic}\index{Fern\'andez-Gotay-Gray
example}\index{example ! Fern\'andez-Gotay-Gray}\index{Gotay !
Fern\'andez-Gotay-Gray}\index{Gray ! Fern\'andez-Gotay-Gray (A.\ Gray)}
first exhibited symplectic manifolds that do not admit any
complex structure at all.
Their examples are circle bundles over circle bundles
(i.e., a {\em tower} of circle bundles) over a 2-torus.

The \textbf{Hopf surface}\index{example !
Hopf surface}\index{Hopf ! surface}
is the complex surface diffeomorphic to $S^1 \times S^3$
obtained as the quotient $\CC^2 \backslash \{0\}/\Gamma$ where
$\Gamma = \{2^n \Id \mid n \in \ZZ\}$ is a group of {\em complex}
transformations, i.e., we factor $\CC^2 \backslash \{0\}$
by the equivalence relation $(z_1,z_2) \sim (2z_1,2z_2)$.
The Hopf surface is not symplectic because $H^2(S^1 \times S^3) = 0$.

The manifold $\CC\PP^2 \# \CC\PP^2 \# \CC\PP^2$ is almost complex
but is neither complex (since it does not fit Kodaira's
classification~\cite{ko:surfacesI})\index{Kodaira ! complex
surface}\index{complex surface}, nor symplectic
as shown by Taubes~\cite{ta:invariants}\index{example !
Taubes}\index{Taubes ! CP@$\CC\PP^2 \# \CC\PP^2 \# \CC\PP^2$ is not
complex}\index{Taubes ! CP@$\CC\PP^2 \# \CC\PP^2 \# \CC\PP^2$ is not
complex} using Seiberg-Witten invariants\index{Seiberg-Witten
invariants}\index{Witten ! Seiberg-Witten invariants}
(Section~\ref{sec:invariants}).

We could go through the previous discussion
restricting to closed 4-dimensional examples
{\em with a specific fundamental group}.
We will do this restricting to simply connected examples,
where the following picture holds.

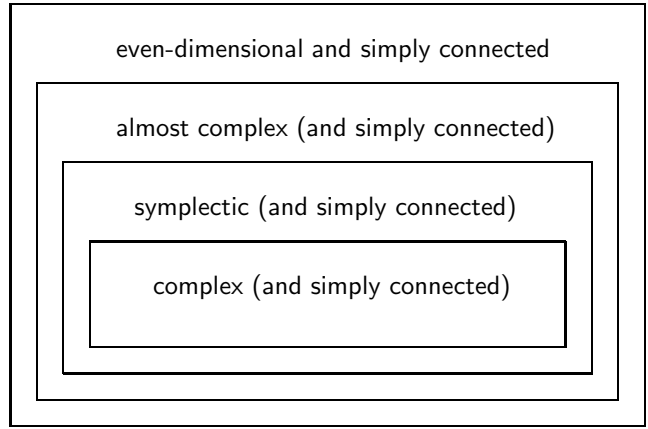
\begin{picture}(200,180)(-20,10)
% horizontal lines
\put(10,20){\line(1,0){240}}
\put(20,30){\line(1,0){220}}
\put(30,40){\line(1,0){200}}
\put(40,50){\line(1,0){180}}
\put(40,90){\line(1,0){180}}
\put(30,120){\line(1,0){200}}
\put(20,150){\line(1,0){220}}
\put(10,180){\line(1,0){240}}
% vertical lines
\put(10,20){\line(0,1){160}}
\put(20,30){\line(0,1){120}}
\put(30,40){\line(0,1){80}}
\put(40,50){\line(0,1){40}}
\put(220,50){\line(0,1){40}}
\put(230,40){\line(0,1){80}}
\put(240,30){\line(0,1){120}}
\put(250,20){\line(0,1){160}}
% words
\put(50,160){even-dimensional and simply connected}
\put(50,130){almost complex (and simply connected)}
\put(57,100){symplectic (and simply connected)}
\put(64,70){complex (and simply connected)}
\thicklines
\end{picture}

It is a consequence of Wu's result~\cite{wu:classes}
that a simply connected manifold admits an almost complex structure
if and only if $b_2^+$ is odd.\footnote{The \textbf{intersection form}
of an oriented {\em topological} closed 4-manifold $M$ is the bilinear
pairing $Q_M : H^2 (M;\ZZ) \times H^2 (M;\ZZ) \to \ZZ$,
$Q_M (\alpha,\beta) := \langle \alpha \cup \beta , [M] \rangle$,
where $\alpha \cup \beta$ is the {\em cup product}
and $[M]$ is the {\em fundamental class}.
Since $Q_M$ always vanishes on torsion elements, descending to
$H^2 (M;\ZZ) / \mathrm{torsion}$ it can be represented by a matrix.
When $M$ is smooth and simply connected, this pairing is
$Q_M (\alpha,\beta) := \int_M \alpha \wedge \beta$ since
non-torsion elements are representable by 2-forms.
As $Q_M$ is symmetric (in the smooth case,
the wedge product of 2-forms is symmetric)
and unimodular (the determinant of a matrix representing
$Q_M$ is $\pm 1$ by Poincar\'e duality), it is diagonalizable over $\RR$
with eigenvalues $\pm 1$.
We denote by $b_2^+$ (respectively, $b_2^-$) the number of
positive (resp.\ negative) eigenvalues of $Q_M$
counted with multiplicities,
i.e., the dimension of a maximal subspace where $Q_M$
is positive-definite (resp.\ negative-definite).
The \textbf{signature} of $M$ is the difference
$\sigma := b_2^+ - b_2^-$, whereas the second Betti number
is the sum $b_2 = b_2^+ + b_2^-$, i.e., the \textbf{rank} of $Q_M$.
The \textbf{type} of an intersection form is \textbf{definite}
if it is positive or negative definite (i.e., $|\sigma| = b_2$)
and \textbf{indefinite} otherwise.}
In particular, the connected sum $\#_m \CC \PP^2 \#_n \overline{\CC \PP^2}$
(of $m$ copies of $\CC \PP^2$
with $n$ copies of $\overline{\CC \PP^2}$)
has an almost complex structure if and only if $m$ is
odd.\footnote{The intersection form of a connected sum $M_0 \# M_1$
is (isomorphic to) $Q_{M_0} \oplus Q_{M_1}$.}

By Kodaira's classification~\cite{ko:surfacesI},
a simply connected complex surface always admits a
compatible symplectic form (since $b^1 = 0$ is even),
i.e., it is always K\"ahler.

Since they are simply connected,
$S^4$, $\CC\PP^2 \# \CC\PP^2 \# \CC\PP^2$ and $\CC \PP^2$
live in three of the four regions in the picture
for simply connected examples.
All of $\CC \PP^2 \#_m \overline{\CC \PP^2}$
are also simply connected K\"ahler manifolds
because they are {\em pointwise blow-ups} $\CC \PP^2$
and the {\em blow-down map} is holomorphic;
see Section~\ref{sec:blow_up}.

There is a family of manifolds obtained from
$\CC \PP^2 \#_9 \overline{\CC \PP^2} =: E(1)$
by a {\em knot surgery}~\cite{fi-st:knots}
that were shown by Fintushel and Stern to be symplectic
and confirmed not to admit a complex
structure~\cite{pa:non-complex}.
The first example
of a closed simply connected symplectic manifold
that cannot be K\"ahler, was a 10-dimensional manifold
obtained by McDuff~\cite{mc:examples} as follows.
The Kodaira-Thurston example\index{example !
Kodaira-Thurston}\index{Kodaira ! Kodaira-Thurston
example}\index{Thurston ! Kodaira-Thurston example}
$\RR^4 / \Gamma$ (not simply connected)
embeds symplectically
in $(\CC \PP^5,\omega_{FS})$~\cite{gr:partial,ti:embedding}.
McDuff's example is a {\em blow-up} of $(\CC \PP^5,\omega_{FS})$
along the image of $\RR^4 / \Gamma$.

\textbf{Geography problems} are problems on the existence
of simply connected closed oriented 4-dimensional manifolds
with some additional structure (such as,
a symplectic form or a complex structure)
for each pair of {\em topological coordinates}.
As a consequence of the work of Freedman~\cite{fr:topology}
and Donaldson~\cite{do:gauge} in the 80's,
it became known that the homeomorphism class of a connected simply connected
closed oriented {\em smooth} 4-manifold is determined
by the two integers -- the second Betti number and the signature
$(b_2, \sigma)$ -- and the {\em parity}\footnote{We say that
the \textbf{parity} of an intersection form $Q_M$ is \textbf{even}
when $Q_M (\alpha,\alpha)$ is even for all $\alpha \in H^2 (M;\ZZ)$,
and \textbf{odd} otherwise.} of the intersection form.
Forgetting about the parity, the numbers $(b_2, \sigma)$
can be treated as \textbf{topological coordinates}.
For each pair $(b_2, \sigma)$ there could well be infinite different
(i.e., nondiffeomorphic) smooth manifolds.
Using riemannian geometry,
Cheeger~\cite{ch:finiteness} showed that there are at most {\em countably many}
different smooth types for closed 4-manifolds.
There are no known finiteness results for the smooth types
of a given topological 4-manifold, in contrast to other dimensions.

Traditionally, the numbers used are
$(c_1^2, c_2):= (3 \sigma + 2 \chi, \chi) = (3 \sigma + 4 + 2b_2 , 2 + b_2)$,
and frequently just the {\em slope} $c_1^2/c_2$ is considered.
If $M$ admits an almost complex structure $J$,
then $(TM,J)$ is a complex vector bundle,
hence has Chern classes $c_1 = c_1 (M,J)$ and $c_2 = c_2 (M,J)$.
Both $c_1^2 := c_1 \cup c_1$ and  $c_2$ may be regarded as numbers
since $H^4 (M;\ZZ) \simeq \ZZ$.
They satisfy $c_1^2 = 3 \sigma + 2 \chi$ (by Hirzebruch's signature formula)
and $c_2 = \chi$ (because the top Chern class is always the Euler class),
justifying the notation for the
topological coordinates in this case.

\begin{examples}
The manifold $\CC \PP^2$
has $(b_2, \sigma)=(1,1)$, i.e., $(c_1^2, c_2) = (9,3)$.
Reversing the orientation $\overline{\CC \PP^2}$
has $(b_2, \sigma)=(1,-1)$, i.e., $(c_1^2, c_2) = (3,3)$.
Their connected sum $\CC \PP^2 \# \overline{\CC \PP^2}$
has $(b_2, \sigma)=(2,0)$, i.e., $(c_1^2, c_2) = (8,0)$.
The product $S^2 \times S^2$
also has $(b_2, \sigma)=(2,0)$ i.e., $(c_1^2, c_2) = (8,4)$.
But $\CC \PP^2 \# \overline{\CC \PP^2}$ has an {\em odd}
intersection form whereas $S^2 \times S^2$ has an {\em even}
intersection form: $\left[ \begin{array}{cc}
   1 & 0 \\ 0 & -1 \end{array} \right]$
vs.\ $\left[ \begin{array}{cc}
   0 & 1 \\ 1 & 0 \end{array} \right]$.
\end{examples}

\textbf{Symplectic geography}~\cite{go-st:calculus,st:geography}
addresses the following question:
What is the set of pairs of integers $(m,n) \in \ZZ \times \ZZ$
for which there exists a connected simply connected
closed {\em symplectic} 4-manifold $M$
having second Betti number $b_2 (M) = m$ and
signature $\sigma (M) = n$?
This problem includes the usual geography of simply connected
complex surfaces, since all such surfaces are K\"ahler
according to Kodaira's
classification~\cite{ko:surfacesI}\index{Kodaira ! complex surfaces}.
Often, instead of the numbers $(b_2,\sigma)$, the question
is equivalently phrased in terms of the Chern numbers
$(c_1^2,c_2)$ for a compatible almost complex structure,
which satisfy $c_1^2 = 3 \sigma + 2 \chi$~\cite{wu:classes}
and $c_2 = \chi$, where $\chi = b_2 + 2$ is the {\em Euler number}.
Usually only {\em minimal} (Section~\ref{sec:blow_up})
or {\em irreducible} manifolds
are considered to avoid trivial examples.
A manifold is \textbf{irreducible} when it is not a
connected sum of other manifolds, except when one of
the summands is a homotopy sphere.

It was speculated that perhaps any simply connected
closed smooth 4-manifold other than $S^4$ is diffeomorphic to
a connected sum of symplectic manifolds,
where any orientation is allowed on each summand
(the so-called {\em minimal conjecture} for smooth 4-manifolds).
Szab\'o~\cite{sz:exotic,sz:irreducible} provided counterexamples in a family
of irreducible simply connected closed non-symplectic smooth 4-manifolds.

All these problems could be posed for other fundamental groups.
Gompf~\cite{go:new} used {\em symplectic sums}
(Section~\ref{sec:constructions}) to prove the following theorem.
He also proved that his surgery construction can be adapted
to produce {\em non}K\"ahler examples.
Since finitely-presented groups are not classifiable,
this shows that compact symplectic 4-manifold are not classifiable.

\begin{theorem}
\label{thm:gompf}\index{Gompf ! theorem}\index{theorem ! Gompf}\index{example !
Gompf}\index{Gompf construction}
\textbf{(Gompf)} $\;$
Every finitely-presented group occurs as the fundamental group
$\pi_1 (M)$ of a compact symplectic 4-manifold $(M,\omega)$.
\end{theorem}

%%%%%%%%%%%%%%%%%%%%%%%%%%%%%%%%%%%%%%%%%%%%%%%%%%%%%%%%%%%%%%%%%%%%%%%%%%%%%
%%%%%%%%%%%%%%%%%%%%%%%%%%%%%%%%%%%%%%%%%%%%%%%%%%%%%%%%%%%%%%%%%%%%%%%%%%%%%

\ssubsection{Fibrations and Sums}
\label{sec:constructions}

Products of symplectic manifolds are naturally symplectic.
As we will see, special kinds of {\em twisted products}, i.e.,
fibrations,\footnote{A
\textbf{fibration} (or {\em fiber bundle}) is a manifold $M$
(called the \textbf{total space}\index{space !
total}\index{total space}) with a submersion $\pi : M \to X$
to a manifold $X$ (the \textbf{base}\index{base})
that is locally trivial in the sense that
there is an open covering of $X$, such that, to each set $\cU$
in that covering corresponds a diffeomorphism of the form
$\varphi_\cU = (\pi,s_\cU) : \pi ^{-1} (\cU) \to \cU \times F$
(a \textbf{local trivialization}\index{trivialization}) where $F$
is a fixed manifold (the \textbf{model fiber}\index{fiber}).
A collection of local trivializations such that the sets $\cU$
cover $X$ is called a \textbf{trivializing cover} for $\pi$.
Given two local trivializations, the second entry of the
composition $\varphi_\cV \circ \varphi_\cU^{-1} = (\id , \psi_{\cU\cV})$
on $(\cU \cap \cV) \times F$ gives the corresponding
\textbf{transition function} $\psi_{\cU\cV} (x): F \to F$
at each $x \in \cU \cap \cV$.}
are also symplectic.

\begin{definition}
A \textbf{symplectic fibration} is a fibration
$\pi : M \to X$ where the model fiber is a symplectic manifold
$(F,\sigma)$ and with a trivializing cover for which
all the transition functions are symplectomorphisms $F \to F$.
\end{definition}

In a symplectic fibration each \textbf{fiber} $\pi^{-1} (x)$
carries a \textbf{canonical symplectic form} $\sigma_x$
defined by the restriction of $s_\cU^* \sigma$,
for any domain $\cU$ of a trivialization
covering $x$ (i.e., $x \in \cU$).
A symplectic form $\omega$ on the total space $M$ of
a symplectic fibration is called \textbf{compatible} with the fibration
if each fiber $(\pi^{-1} (x), \sigma_x)$
is a symplectic submanifold of $(M,\omega)$,
i.e., $\sigma_x$ is the restriction of $\omega$ to $\pi^{-1} (x)$.

\begin{examples}
\begin{enumerate}
\item
{\em Every compact oriented\footnote{An \textbf{oriented fibration}
is a fibration whose model fiber is oriented
and there is a trivializing cover for which all transition functions
preserve orientation.}
fibration whose model fiber $F$ is an oriented surface
admits a structure of symplectic fibration}
for the following reason.
Let $\sigma_0$ be an area form on $F$.
Each transition function $\psi_{\cU\cV} (x): F \to F$ pulls
$\sigma_0$ back to a cohomologous area form $\sigma_1$
(depending on $\psi_{\cU\cV} (x)$).
Convex combinations $\sigma_t = (1-t) \sigma_0 + t \sigma_1$
give a path of area forms from $\sigma_0$ to $\sigma_1$
with constant class $[\sigma_t]$.
By Moser's argument (Section~\ref{sec:trick}),
there exists a diffeomorphism $\rho (x):F \to F$
isotopic to the identity, depending smoothly on $x \in \cU \cap \cV$,
such that $\psi_{\cU\cV} (x) \circ \rho (x)$
is a symplectomorphism of $(F,\sigma_0)$.
By successively adjusting local trivializations for a finite
covering of the base, we can make all transition functions
into symplectomorphisms.

\item
{\em Every fibration with connected base and compact fibers
having a symplectic form $\omega$ for which all
fibers are symplectic submanifolds admits a structure of
symplectic fibration compatible with $\omega$.}
Indeed, under trivializations, the restrictions of $\omega$ to the fibers
give cohomologous symplectic forms in the model fiber $F$.
So by Moser's Theorem~\ref{thm:moser},
all fibers are strongly isotopic to $(F,\sigma)$
where $\sigma$ is the restriction of $\omega$ to a chosen fiber.
These isotopies can be used to produce a trivializing cover
where each $s_\cU (x)$ is a symplectomorphism.
\end{enumerate}
\end{examples}

In the remainder of this section,
assume that for a fibration $\pi : M \to X$
the total space is compact and the base is connected.
For the existence of a compatible symplectic form
on a symplectic fibration, a necessary condition
is the existence of a cohomology class in $M$
that restricts to the classes of the fiber symplectic forms.
Thurston~\cite{th:examples} showed that,
when the base admits also a symplectic form, this condition is sufficient.
Yet not all symplectic fibrations with a compatible symplectic form
have a symplectic base~\cite{we:fat}.

\begin{theorem}
\label{thm:thurston}\index{Thurston ! theorem}\index{theorem ! Thurston}
\textbf{(Thurston)} $\;$
Let $\pi : M \to X$ be a compact symplectic fibration
with connected symplectic base $(X,\alpha)$
and model fiber $(F,\sigma)$.
If there is a class $[\nu] \in H^2 (M)$ pulling back to $[\sigma]$,
then, for sufficiently large $k > 0$,
there exists a symplectic form $\omega_k$ on $M$
that is compatible with the fibration and is in $[\nu + k \pi^* \alpha]$.
\end{theorem}

\vspace*{-2ex}

\begin{proof}
We first find a form $\tau$ on $M$ in the class $[\nu]$
that restricts to the canonical symplectic form on each fiber.
Pick a trivializing cover
$\{\varphi_i = (\pi,s_i) \mid i \in I\}$
with contractible domains $\cU_i$.
Let $\rho_i$, $i \in I$, be a partition of unity subordinate to this covering
and let $\widetilde \rho_i := \rho_i \circ \pi : M \to \RR$.
Since $[\nu]$ always restricts to the class
of the canonical symplectic form $[\sigma_x]$,
and the $\cU_i$'s are contractible,
on each $\pi_i^{-1} (\cU_i)$ the forms
$s_i^* \sigma - \nu$ are exact.
Choose 1-forms $\lambda_i$ such that
$s_i^* \sigma = \nu + d \lambda_i$, and set
\[
   \tau := \nu + \sum_{i \in I} d (\widetilde \rho_i \lambda_i) \ .
\]
Since $\tau$ is nondegenerate on the (vertical) subbundle
given by the kernel of $d \pi$, for $k > 0$ large enough
the form $\tau + k \pi^* \alpha$ is nondegenerate on $M$.
\end{proof}

\vspace*{-1ex}

\begin{corollary}
Let $\pi : M \to X$ be a compact oriented fibration
with connected symplectic base $(X,\alpha)$
and model fiber an oriented surface $F$ of genus $g(F) \neq 1$.
Then $\pi$ admits a compatible symplectic form.
\end{corollary}

\vspace*{-2ex}

\begin{proof}
By Example~1 above, $\pi : M \to X$ admits a structure
of symplectic fibration with model fiber $(F,\sigma)$.
Since the fiber is not a torus ($g(F) \neq 1$),
the Euler class of the tangent bundle $TF$
(which coincides with $c_1 (F,\sigma)$)
is $\lambda [\sigma]$ for some $\lambda \neq 0$.
Hence, the first Chern class $[c]$
of the {\em vertical} subbundle given by the kernel of $d\pi$
(assembling the tangent bundles to the fibers)
restricts to $\lambda [\sigma_x]$ on the fiber over $x \in X$.
We can apply Theorem~\ref{thm:thurston} using the class
$[\nu] = \lambda^{-1} [c]$.
\end{proof}

%%%%%%%%%%%%%%%%%%%%%%%%%%%%%%%%%%%%%%%%%%%%%%%%%%%%%%%%%%%%%%%%%%%%%%%%%%%%%

A pointwise connected sum $M_0 \# M_1$
of symplectic manifolds $(M_0,\omega_0)$ and $(M_1,\omega_1)$
tends to not admit a symplectic form, even if we only require the
eventual symplectic form to be isotopic to $\omega_i$ on
each $M_i$ minus a ball.
The reason~\cite{au:exemples} is that such a
symplectic form on $M_0 \# M_1$
would allow to construct an almost complex structure on
the sphere formed by the union of the two removed balls,
which is known not to exist except on $S^2$ and $S^6$.
Therefore:

\begin{proposition}
Let $(M_0,\omega_0)$ and $(M_1,\omega_1)$ be two compact
symplectic manifolds of dimension not 2 nor 6.
Then the connected sum $M_0 \# M_1$ does not admit
any symplectic structure isotopic to $\omega_i$ on
$M_i$ minus a ball, $i=1,2$.
\end{proposition}

For connected sums to work in the symplectic category,
they should be done 
along codimension 2 symplectic submanifolds.
The following construction, already mentioned in~\cite{gr:partial},
was dramatically explored and popularized by Gompf~\cite{go:new}
(he used it to prove Theorem~\ref{thm:gompf}).
Let $(M_0,\omega_0)$ and $(M_1,\omega_1)$ be two
$2n$-dimensional symplectic manifolds.
Suppose that a compact symplectic manifold $(X,\alpha)$ of
dimension $2n-2$ admits symplectic embeddings to both
$i_0 : X \hookrightarrow M_0$, $i_1 : X \hookrightarrow M_1$.
For simplicity, assume that the corresponding
normal bundles are trivial (in general, they need to have
symmetric Euler classes).
By the symplectic neighborhood theorem
(Theorem~\ref{thm:weinstein_symplectic}),
there exist symplectic embeddings
$j_0 : X \times B_\varepsilon \to M_0$ and
$j_1 : X \times B_\varepsilon \to M_1$ (called \textbf{framings})
where $B_\varepsilon$ is a ball of radius $\varepsilon$
and centered at the origin in $\RR^2$ such that
$j_k^* \omega_k = \alpha + dx \wedge dy$ and
$j_k (p,0) = i_k (p) $ $\forall p \in X$, $k=0,1$.
Chose an area- and orientation-preserving diffeomorphism $\phi$
of the annulus $B_\varepsilon \setminus B_\delta$ for
$0 < \delta < \varepsilon$ that interchanges the two
boundary components.
Let $\cU_k = j_k (X \times B_\delta) \subset M_k$, $k=0,1$.
A \textbf{symplectic sum} of $M_0$ and $M_1$ along $X$ is defined to be
\[
   M_0 \#_X M_1 := \left( M_0 \setminus \cU_0 \right)
   \cup_\phi \left( M_1 \setminus \cU_1 \right)
\]
where the symbol $\cup_\phi$ means that we identify
$j_1(p,q)$ with $j_0 (p,\phi(q))$ for all $p \in X$ and
$\delta < |q| < \varepsilon$.
As $\omega_0$ and $\omega_1$ agree on the regions under identification,
they induce a symplectic form on $M_0 \#_X M_1$.
The result depends on $j_0$, $j_1$, $\delta$ and $\phi$.

\textbf{Rational blowdown} is a surgery on 4-manifolds
that replaces a neighborhood of a chain of embedded $S^2$'s
with boundary a {\em lens space} $L(n^2,n-1)$ by a manifold with the
same rational homology as a ball.
This simplifies the homology
possibly at the expense of complicating the fundamental group.
Symington~cite{sy:blowdown} showed that rational blowdown
preserves a symplectic structure if the original spheres are
symplectic surfaces in a symplectic 4-manifold.

%%%%%%%%%%%%%%%%%%%%%%%%%%%%%%%%%%%%%%%%%%%%%%%%%%%%%%%%%%%%%%%%%%%%%%%%%%%%%
%%%%%%%%%%%%%%%%%%%%%%%%%%%%%%%%%%%%%%%%%%%%%%%%%%%%%%%%%%%%%%%%%%%%%%%%%%%%%

\ssubsection{Symplectic Blow-Up}
\label{sec:blow_up}

{\em Symplectic blow-up} is the extension to the symplectic category
of the blow-up operation in algebraic geometry.
It is due to Gromov according to the first
printed exposition of this operation in~\cite{mc:examples}.

Let $L$ be the \textbf{tautological line bundle} over $\CC \PP^{n-1}$, that is,
\[
   L = \{ ([p], z) \mid
   p \in \CC^n \setminus \{ 0 \} \ , \ z = \lambda p
   \mbox{ for some } \lambda \in \CC \}
\]
with projection to $\CC \PP^{n-1}$ given by $\pi: ([p], z) \mapsto [p]$.
The fiber of $L$ over the point $[p] \in \CC \PP^{n-1}$
is the complex line in $\CC^n$ represented by that point.
The \textbf{blow-up of $\CC^n$ at the origin}\index{blow-up ! definition}
is the total space of the bundle $L$, sometimes denoted $\widetilde \CC^n$.
The corresponding \textbf{blow-down map}\index{blow-down map} is the map
$\beta : L \to \CC^n$ defined by $\beta ([p], z) = z$.
The total space of $L$ may be decomposed as
the disjoint union of two sets: the zero section
\[
   E:= \{ ([p],0) \mid p \in \CC^n \setminus \{ 0 \} \}
\]
and
\[
   S:= \{ ([p],z) \mid p \in \CC^n \setminus \{ 0 \}  \ ,
   \ z = \lambda p \mbox{ for some } \lambda \in \CC^* \} \ .
\]
The set $E$ is called the
\textbf{exceptional divisor}\index{exceptional divisor};
it is diffeomorphic to $\CC \PP^{n-1}$ and gets mapped to the
origin by $\beta$.
On the other hand, the restriction of $\beta$ to the
complementary set $S$ is a diffeomorphism onto $\CC^n \setminus \{ 0 \}$.
Hence, we may regard $L$ as being obtained from $\CC^n$
by smoothly replacing the origin by a copy of $\CC \PP^{n-1}$.
Every biholomorphic map $f: \CC^n \to \CC^n$ with $f(0)=0$
lifts uniquely to a biholomorphic map $\widetilde f : L \to L$
with $\widetilde f (E) = E$.
The lift is given by the formula
\[
   {\widetilde f} ([p],z) = \left\{ \begin{array}{ll}
   ([f(z)],f(z)) & \mbox{ if } z \neq 0 \\
   ([p],0) & \mbox{ if } z = 0 \ . \end{array} \right.
\]
There are actions of the unitary group $\UU (n)$ on $L$, $E$ and $S$
induced by the standard linear action on $\CC^n$,
and the map $\beta$ is $\UU (n)$-equivariant.
For instance, $\beta^* \omega_0 + \pi^* \omega_{FS}$
is a $\UU (n)$-invariant K\"ahler form on $L$.

\begin{definition}
A \textbf{blow-up symplectic form}\index{blow-up ! symplectic
form}\index{symplectic ! blow-up}
on the tautological line bundle $L$ is a
$\UU (n)$-invariant symplectic form $\omega$ such that
the difference $\omega - \beta^* \omega_0$ is compactly supported,
where $\omega _0 = \frac i2 \sum_{k=1}^n dz_k \wedge d\bar z_k$
is the standard symplectic form on $\CC^n$.
\end{definition}

Two blow-up symplectic forms are \textbf{equivalent}\index{blow-up !
equivalence} if one is the pullback of the other by a
$\UU (n)$-equivariant diffeomorphism of $L$.
Guillemin and Sternberg~\cite{gu-st:birational} showed that
two blow-up symplectic forms are equivalent if and only if
they have equal restrictions to the exceptional divisor $E \subset L$.
Let $\Omega^\varepsilon$ ($\varepsilon > 0$) be the set of all
blow-up symplectic forms on $L$ whose restriction to the
exceptional divisor $E \simeq \CC \PP^{n-1}$ is
$\varepsilon \omega_{_{\mathrm{FS}}}$,
where $\omega_{_{\mathrm{FS}}}$ is the Fubini-Study form
(Section~\ref{sec:kahler}).
An \textbf{$\varepsilon$-blow-up} of $\CC^n$ at the origin
is a pair $(L,\omega)$ with $\omega \in \Omega ^\varepsilon$.

Let $(M, \omega)$ be a $2n$-dimensional symplectic manifold.
It is a consequence of Darboux's Theorem (Theorem~\ref{thm:darboux})
that, for each point $p \in M$, there exists a complex chart
$(\cU , z_1 , \ldots , z_n)$ centered at $p$ and
with image in $\CC^n$ where
$\left. \omega \right|_\cU = \frac i2 \sum_{k=1}^n dz_k \wedge d\bar z_k$.
It is shown in~\cite{gu-st:birational} that, for
$\varepsilon$ small enough, we can perform an
$\varepsilon$-blow-up of $M$ at $p$ modeled on $\CC^n$
at the origin, without changing the symplectic structure
outside of a small neighborhood of $p$.
The resulting manifold is called an
\textbf{$\varepsilon$-blow-up of $M$ at $p$}.\index{blow-up ! definition}
As a manifold, the blow-up of $M$ at a point is diffeomorphic
to the {\em connected sum}\footnote{The \textbf{connected sum}
of two oriented $m$-dimensional manifolds $M_0$ and $M_1$
is the manifold, denoted $M_0 \# M_1$, obtained from the union
of those manifolds each with a small ball removed $M_i \setminus B_i$
by identifying the boundaries via a (smooth) map
$\phi : \partial B_1 \to \partial B_2$ that extends
to an orientation-preserving diffeomorphism of neighborhoods
of $\partial B_1$ and $\partial B_2$ (interchanging the inner
and outer boundaries of the annuli).}
$M \# \overline{\CC \PP^n}$, where $\overline{\CC \PP^n}$ is
the manifold $\CC \PP^n$ equipped with the orientation
opposite to the natural complex one.

\begin{example}
Let $\PP (L \oplus \CC)$ be the $\CC \PP^1$-bundle over $\CC \PP^{n-1}$
obtained by projectivizing the direct sum of the tautological
line bundle $L$ with a trivial complex line bundle.
Consider the map
\[
\begin{array}{rrcl}
   \beta: & \CC \PP (L \oplus \CC) & \longrightarrow & \CC \PP^n \\
   & ([p],[\lambda p : w]) & \longmapsto & [\lambda p : w] \ ,
\end{array}
\]
where $[\lambda p : w]$ on the right
represents a line in $\CC^{n+1}$,
forgetting that, for each $[p] \in \CC \PP^{n-1}$,
that line sits in the 2-complex-dimensional subspace
$L_{[p]} \oplus \CC \subset \CC^n \oplus \CC$.
Notice that $\beta$ maps the {\em exceptional divisor}
\[
   E:= \{ ([p],[0:\ldots :0:1]) \mid [p] \in \CC \PP^{n-1} \}
   \simeq \CC \PP^{n-1}
\]
to the point $[0:\ldots :0:1] \in \CC \PP^n$,
and $\beta$ is a diffeomorphism on the complement
\[
   S:= \{ ([p],[\lambda p : w]) \mid [p] \in \CC \PP^{n-1} \ ,
   \ \lambda \in \CC^* \ , \ w \in \CC \} \simeq \CC \PP^n
   \setminus \{ [0:\ldots : 0:1] \} \ .
\]
Therefore, we may regard $\CC \PP (L \oplus \CC)$ as being obtained
from $\CC \PP^n$ by smoothly replacing the point $[0:\ldots :0:1]$
by a copy of $\CC \PP^{n-1}$.
The space $\CC \PP (L \oplus \CC)$ is the blow-up of $\CC \PP^n$
at the point $[0:\ldots :0:1]$, and $\beta$ is the corresponding
blow-down map.
The manifold $\CC \PP (L \oplus \CC)$ for $n=2$ is
a \textbf{Hirzebruch surface}.\index{example !
Hirzebruch surface}\index{Hirzebruch surface}
\end{example}

%%%%%%%%%%%%%%%%%%%%%%%%%%%%%%%%%%%%%%%%%%%%%%%%%%%%%%%%%%%%%%%%%%%%%%%%%%%%%

When $(\CC \PP^{n-1} , \omega_{FS})$ is symplectically
embedded in a symplectic manifold $(M,\omega)$ with image $X$
and normal bundle isomorphic to the tautological bundle $L$,
it can be subject to a {\em blow-down} operation.
By the symplectic neighborhood theorem
(Theorem~\ref{thm:weinstein_symplectic}),
some neighborhood $\cU \subset M$ of the image $X$ is
symplectomorphic to a neighborhood $\cU_0 \subset L$ of the zero section.
It turns out that some neighborhood of $\partial \cU_0$ in $L$
is symplectomorphic to a spherical shell in $(\CC^n,\omega_0)$.
The \textbf{blow-down of $M$ along $X$} is a manifold obtained
from the union of $M \setminus \cU$ with a ball in $\CC^n$.
For more details, see~\cite[\S 7.1]{mc-sa:introduction}.

Following algebraic geometry, we call \textbf{minimal} a
$2n$-dimensional symplectic manifold $(M,\omega)$ without any
symplectically embedded $(\CC \PP^{n-1} , \omega_{FS})$,
so that $(M,\omega)$ is not the blow-up at a point of another
symplectic manifold.
In dimension 4, a manifold is minimal
if it does not contain any embedded sphere $S^2$ with self-intersection $-1$.
Indeed, by the work of Taubes~\cite{ta:invariants,ta:sw=>gr},
if such a sphere $S$ exists, then either the homology class $[S]$
or its symmetric $-[S]$ can be represented by a {\em symplectically}
embedded sphere with self-intersection $-1$.

%%%%%%%%%%%%%%%%%%%%%%%%%%%%%%%%%%%%%%%%%%%%%%%%%%%%%%%%%%%%%%%%%%%%%%%%%%%%%

For a symplectic manifold $(M , \omega)$,
let $i : X \hookrightarrow M$ be the inclusion of a symplectic submanifold.
The normal bundle $NX$ to $X$ in $M$ admits a structure
of complex vector bundle (as it is a symplectic vector bundle).
Let $\PP (NX) \to X$ be the projectivization of the bundle $NX \to X$,
let $Z$ be the zero section of $NX$,
let $L (NX)$ be the corresponding \textbf{tautological line bundle}
(given by assembling the tautological line bundles over each fiber)
and let $\beta : L (NX) \to NX$ be the blow-down map.
On the {\em exceptional divisor}
\[
   E:= \{ ([p],0) \in L(NX) \mid p \in NX \setminus Z \}
   \simeq \PP (NX)
\]
the map $\beta$ is just projection to the zero section $Z$.
The restriction of $\beta$ to the complement $L (NX) \setminus E$
is a diffeomorphism to $NX \setminus Z$.
Hence, $L (NX)$ may be viewed as being obtained from $NX$
by smoothly replacing each point of the zero section by
the projectivization of its normal space.
We symplectically identify some tubular neighborhood $\cU$ of $X$ in $M$
with a tubular neighborhood $\cU_0$ of the zero section $Z$ in $NX$.
A \textbf{blow-up of the symplectic
manifold $(M , \omega)$ along the symplectic submanifold $X$}\index{blow-up !
along a submanifold} is the manifold obtained
from the union of $M \setminus \cU$ and $\beta^{-1} (\cU_0)$
by identifying neighborhoods of $\partial \cU$, and equipped with a
symplectic form that restricts to $\omega$ on
$M \setminus \cU$~\cite{mc:examples}.
When $X$ is one point, this construction reduces to
the previous symplectic blow-up at a point.

Often symplectic geography concentrates on minimal examples.
McDuff~\cite{mc:rational} showed that a minimal symplectic 4-manifold
with a symplectically embedded $S^2$ with nonnegative
self-intersection is symplectomorphic either to $\CC \PP^2$
or to an $S^2$-bundle over a surface.
Using Seiberg-Witten theory it was proved:

\begin{theorem}
Let $(M,\omega)$ be a minimal closed symplectic 4-manifold.
\begin{itemize}
\item[(a)]
\textbf{(Taubes~\cite{ta:sw=>gr})}
If $b_2^+ > 1$, then $c_1^2 \geq 0$.
\item[(b)]
\textbf{(Liu~\cite{li:wall})}
If $b_2^+ = 1$ and $c_1^2 < 0$,
then $M$ is the total space of an $S^2$-fibration
over a surface of genus $g$ where $\omega$ is nondegenerate on the fibers,
and $(c_1^2,c_2) = (8-8g, 4-4g)$,
i.e., $(M,\omega)$ is a {\em symplectic ruled surface}.
\end{itemize}
\end{theorem}

A \textbf{symplectic ruled surface}\footnote{A (rational)
\textbf{ruled surface}
is a complex (K\"ahler) surface that is the total space of a holomorphic
fibration over a Riemann surface with fiber $\CC \PP^1$.
When the base is also a sphere, these are the \textbf{Hirzebruch surfaces}
$\PP (L \oplus \CC)$ where $L$ is a holomorphic line bundle over $\CC \PP^1$.}
is a symplectic 4-manifold
$(M,\omega)$ that is the total space of an $S^2$-fibration
where $\omega$ is nondegenerate on the fibers.

A \textbf{symplectic rational surface} is a symplectic 4-manifold $(M,\omega)$
that can be obtained from the standard $(\CC \PP ^2,\omega_{FS})$
by blowing up and blowing down.

With $b_2^+ = 1$ and $c_1^2 = 0$, we have
symplectic manifolds $\CC \PP^2 \#_9 \overline{\CC \PP^2} =: E(1)$,
the {\em Dolgachev surfaces} $E(1,p,q)$,
the results $E(1)_K$ of surgery on a fibered knot $K \subset S^3$, etc.
With $b_2^+ = 1$ and $c_1^2 > 0$, we have
symplectic manifolds $\CC \PP^2$, $S^2 \times S^2$,
$\CC \PP^2 \#_n \overline{\CC \PP^2}$ for $n \leq 8$
and the {\em Barlow surface}.
For $b_2^+ = 1$ and $c_1^2 \geq 0$, Park~\cite{pa:non-complex}
gave a criterion for a symplectic 4-manifold to be rational
or ruled in terms of Seiberg-Witten theory.

%%%%%%%%%%%%%%%%%%%%%%%%%%%%%%%%%%%%%%%%%%%%%%%%%%%%%%%%%%%%%%%%%%%%%%%%%%%%
%%%%%%%%%%%%%%%%%%%%%%%%%%%%%%%%%%%%%%%%%%%%%%%%%%%%%%%%%%%%%%%%%%%%%%%%%%%%%

\ssubsection{Uniqueness of Symplectic Forms}
\label{sec:classification}

Besides the notions listed in Section~\ref{sec:trick},
the following equivalence relation for symplectic manifolds is considered.
As it allows the cleanest statements about uniqueness,
this relation is simply called {\em equivalence}.

\begin{definition}
Symplectic manifolds $(M,\omega_0)$ and $(M,\omega_1)$ are
\textbf{equivalent} if they are related by a combination
of deformation-equivalences and symplectomorphisms.
\end{definition}

Recall that $(M,\omega_0)$ and $(M,\omega_1)$ are
{\em deformation-equivalent} when there is a smooth family $\omega_t$
of symplectic forms joining $\omega_0$ to $\omega_1$
(Section~\ref{sec:trick}),
and they are {\em symplectomorphic} when there is a diffeomorphism
$\varphi : M \to M$ such that $\varphi^* \omega_1 = \omega_0$
(Section~\ref{symplectic_forms}).
Hence, equivalence is the relation
generated by deformations and diffeomorphisms.
The corresponding equivalence classes can be viewed as
the connected components of the moduli space of symplectic forms
up to diffeomorphism.
This is a useful notion when focusing on topological properties.

\begin{examples}
\begin{enumerate}
\item
The complex projective plane $\CC \PP ^2$ has a unique
symplectic structure up to symplectomorphism
and scaling.
This was shown by Taubes~\cite{ta:sw=gr}\index{Taubes !
unique symplectic structure on
$\CC \PP ^2$}\index{unique symplectic structure on $\CC \PP ^2$}
relating Seiberg-Witten invariants (Section~\ref{sec:invariants})
to pseudoholomorphic curves to prove the existence
of a pseudoholomorphic sphere.
Previous work of Gromov~\cite{gr:pseudo} and McDuff~\cite{mc:structure}
showed that the existence of
a pseudoholomorphic sphere implies that the symplectic form is standard.

Lalonde and McDuff~\cite{la-mc:classification}
concluded similar classifications for symplectic ruled surfaces
and for symplectic rational surfaces (Section~\ref{sec:blow_up}).
The symplectic form on a symplectic ruled surface
is unique up to symplectomorphism in its cohomology class,
and is isotopic to a standard K\"ahler form.
In particular, any symplectic form on $S^2 \times S^2$ is
symplectomorphic to $a \pi_1^* \sigma + b \pi_2^* \sigma$
for some $a,b >0$ where $\sigma$ is the standard area form on $S^2$.

Li-Liu~\cite{li-liu:symplectic}
showed that the symplectic structure on 
$\CC \PP^2 \#_n \overline{\CC \PP^2}$ for $2 \leq n \leq 9$ is unique
up to equivalence.

\item
McMullen and Taubes~\cite{mc-ta:inequivalent}
first exhibited simply connected closed 4-manifolds
admitting inequivalent symplectic structures.
Their examples were constructed using 3-dimensional topology,
and distinguished by analyzing the structure of Seiberg-Witten
invariants to show that the first Chern classes
(Section~\ref{sec:compatible_almost}) of the two symplectic
structures lie in disjoint orbits of the diffeomorphism group.
In higher dimensions there were previously
examples of manifolds with inequivalent symplectic forms;
see for instance~\cite{ru:algebraic}.

With symplectic techniques and avoiding gauge theory,
Smith~\cite{sm:moduli} showed that, for each $n \geq 2$,
there is a simply connected closed 4-manifold that
admits at least $n$ inequivalent symplectic forms,
also distinguished via the first Chern classes.
It is not yet known whether there exist inequivalent symplectic
forms on a 4-manifold with the same first Chern class.
\end{enumerate}
\end{examples}

%%%%%%%%%%%%%%%%%%%%%%%%%%%%%%%%%%%%%%%%%%%%%%%%%%%%%%%%%%%%%%%%%%%%%%%%%%%%%
%%%%%%%%%%%%%%%%%%%%%%%%%%%%%%%%%%%%%%%%%%%%%%%%%%%%%%%%%%%%%%%%%%%%%%%%%%%%%

\ssubsection{Invariants for 4-Manifolds}
\label{sec:invariants}

Very little was known about 4-dimensional manifolds until 1981,
when Freedman~\cite{fr:topology} provided a complete classification of
closed simply connected {\em topological} 4-manifolds, and shortly thereafter
Donaldson~\cite{do:gauge} showed that the panorama for
{\em smooth} 4-manifolds was much wilder.\footnote{It had been proved
by Rokhlin in 1952 that if such a smooth manifold $M$
has even intersection form $Q_M$ (i.e., $w_2 =0$),
then the signature of $Q_M$ must be a multiple of 16.
It had been proved by Whitehead and Milnor
that two such topological manifolds are homotopy equivalent
if and only if they have the same intersection form.}
Freedman showed that, modulo homeomorphism,
such topological manifolds are essentially classified
by their intersection forms
(for an {\em even} intersection form there is exactly one class,
whereas for an {\em odd} intersection form there are exactly two
classes distinguished by the {\em Kirby-Siebenmann invariant} $KS$,
at most one of which admits smooth representatives
-- smoothness requires $KS = 0$).
Donaldson showed that, whereas the existence of a smooth structure
imposes strong constraints on the topological type of a manifold,
for the same topological manifold there can be infinite
different smooth structures.\footnote{It
is known that in dimensions $\leq 3$, each topological manifold has exactly
one smooth structure, and in dimensions $\geq 5$
each topological manifold has at most finitely many smooth structures.
For instance, whereas each topological $\RR^n$, $n \neq 4$,
admits a unique smooth structure, the topological $\RR^4$
admits uncountably many smooth structures.}
In other words, by far not all intersection forms can occur
for smooth 4-manifolds and the same intersection form
may correspond to nondiffeomorphic manifolds.

Donaldson's key tool was a set of gauge-theoretic invariants,
defined by counting with signs the equivalence classes
(modulo gauge equivalence) of connections on $\SU (2)$- (or $\SO (3)$-)
bundles over $M$ whose curvature has vanishing self-dual part.
For a dozen years there was hard work on the invariants discovered by
Donaldson but limited advancement on the understanding of smooth 4-manifolds.

%%%%%%%%%%%%%%%%%%%%%%%%%%%%%%%%%%%%%%%%%%%%%%%%%%%%%%%%%%%%%%%%%%%%%%%%%%%%%

\begin{examples}
Finding {\em exotic}\footnote{A manifold homeomorphic but not
diffeomorphic to a smooth manifold $M$ is called an \textbf{exotic} $M$.}
smooth structures on closed simply connected manifolds
with small $b_2$ has long been an interesting problem,
especially in view of the smooth Poincar\'e conjecture for 4-manifolds.
The first exotic smooth structures on a rational surface
$\CC\PP^2 \#_n \overline{\CC \PP^2}$ were found in the late 80's
for $n=9$ by Donaldson~\cite{do:irrationality}
and for $n=8$ by Kotschick~\cite{ko:homeomorphic}.
There was no progress until the recent work of Park~\cite{pa:symplectic}
constructing a symplectic exotic $\CC \PP^2 \#_7 \overline{\CC \PP^2}$
and using this to exhibit a third distinct smooth structure
$\CC \PP^2 \#_8 \overline{\CC \PP^2}$,
thus illustrating how the existence of symplectic forms is
tied to the existence of different smooth structures.
This stimulated research by Fintushel, Ozsv\'ath, Park, Stern,
Stipsicz and Szab\'o, which together shows that there are infinitely
many exotic smooth structures on $\CC\PP^2 \#_n \overline{\CC \PP^2}$
for $n=5,6,7,8$ (the case $n=9$ had been shown in the late 80's
by Friedman-Morgan and by Okonek-Van de Ven).
\end{examples}

%%%%%%%%%%%%%%%%%%%%%%%%%%%%%%%%%%%%%%%%%%%%%%%%%%%%%%%%%%%%%%%%%%%%%%%%%%%%%

In 1994 Witten brought about a revolution in Donaldson theory
by introducing a new set of invariants
-- the \textbf{Seiberg-Witten invariants} --
which are much simpler to calculate and to apply.
This new viewpoint was inspired by developments due to Seiberg and Witten
in the understanding of {\em $N=2$ supersymmetric Yang-Mills}.

Let $M$ be a smooth oriented closed 4-dimensional manifold
with $b_2^+ (M) > 1$ (there is a version for $b_2^+ (M) = 1$).
All such 4-manifolds $M$ (with any $b_2^+ (M)$) admit a spin-c structure,
i.e., a $\mathrm{Spin}^c (4)$-bundle over $M$ with an isomorphism
of the associated $\SO (4)$-bundle to the bundle of oriented frames
on the tangent bundle for some chosen riemannian metric.
Let $\mathcal C_M = \{ a \in H^2 (M;\ZZ) \mid a \equiv w_2 (TM) (2) \}$
be the set of characteristic elements, and let
$\mathrm{Spin}^c (M)$ be the set of spin-c structures on $M$.
For simplicity, assume that $M$ is simply connected
(or at least that $H_1 (M;\ZZ)$ has no 2-torsion), so that
$\mathrm{Spin}^c (M)$ is isomorphic to $\mathcal C_M$ with isomorphism
given by the first Chern class of the {\em determinant line bundle}
(the \textbf{determinant line bundle} is the line bundle associated
by a natural group homomorphism $\mathrm{Spin}^c (4) \to \UU(1)$).
Fix an orientation of a maximal-dimensional positive-definite
subspace $H_+^2 (M ; \RR) \subset H^2 (M ; \RR)$.
The Seiberg-Witten invariant is the function
\[
   \mathrm{SW}_M : \mathcal C_M \longrightarrow \ZZ
\]
defined as follows.
Given a spin-c structure $\alpha \in \mathrm{Spin}^c (M) \simeq \mathcal C_M$,
the image
$\mathrm{SW}_M (\alpha) = [\cM] \in H_d (\cB^*;\ZZ)$
is the homology class of the moduli space $\cM$ of solutions
(called \textbf{monopoles}) of the Seiberg-Witten (SW) equations
modulo gauge equivalence.
The SW equations are non-linear differential equations on
a pair of a connection $A$ on the determinant line bundle of $\alpha$
and of a section $\varphi$ of an associated $\UU(2)$-bundle,
called the positive (half) spinor bundle:
\[
   F_A^+ = i q (\varphi) \qquad \mbox{ and } \qquad D_A \varphi = 0 \ ,
\]
where $F_A^+$ is the self-dual part of the (imaginary) curvature of $A$,
$q$ is a squaring operation taking sections of the positive
spinor bundle to self-dual 2-forms,
and $D_A$ is the corresponding Dirac operator.
For a generic perturbation of the equations
(replacing the first equation by $F_A^+ = i q (\varphi) +i\nu$,
where $\nu$ is a self-dual 2-form)
and of the riemannian metric, a transversality argument shows that
the moduli space $\cM$ is well-behaved and actually inside
the space $\cB^*$ of gauge-equivalence classes
of irreducible pairs (those $(A,\varphi)$ for which $\varphi \neq 0$),
which is homotopy-equivalent to $\CC \PP^\infty$
and hence has even-degree homology groups
$H_d (\cB^*;\ZZ) \simeq \ZZ$.
When the dimension $d$ of $\cM$ is odd or when $\cM$ is empty,
the invariant $\mathrm{SW}_M (\alpha)$ is set to be zero.
The \textbf{basic classes} are the classes $\alpha \in \mathcal C_M$
for which $\mathrm{SW}_M (\alpha) \neq 0$.
The set of basic classes is always finite,
and if $\alpha$ is a basic class then so is $-\alpha$.
The main results are that the Seiberg-Witten invariants
are invariants of the diffeomorphism type of the 4-manifold $M$
and satisfy vanishing and nonvanishing
theorems, which allowed to answer an array of
questions about specific manifolds.

%%%%%%%%%%%%%%%%%%%%%%%%%%%%%%%%%%%%%%%%%%%%%%%%%%%%%%%%%%%%%%%%%%%%%%%%%%%%%

Taubes~\cite{ta:sw=gr} discovered an equivalence between Seiberg-Witten
and Gromov invariants (using pseudoholomorphic curves) for
symplectic 4-manifolds, by proving
the existence of pseudoholomorphic curves from solutions
of the Seiberg-Witten equations and vice-versa.
As a consequence, he proved:

\begin{theorem}
\label{thm:taubes}\index{Taubes ! theorem}\index{theorem !
Taubes}
\textbf{(Taubes)} $\;$
Let $(M,\omega)$ be a compact symplectic 4-manifold.

If $b_2^+ > 1$, then $c_1 (M,\omega)$ admits a %(possibly disconnected)
smooth pseudoholomorphic representative.

If $M = M_1 \# M_2$, then one of the $M_i$'s has
negative definite intersection form.
\end{theorem}

There are results also for $b_2^+ = 1$,
and follow-ups describe the set of basic classes of
a connected sum $M \# N$ in terms of the set of
basic classes of $M$ when $N$ is a manifold with
negative definite intersection form
(starting with $\overline{\CC \PP^2}$).

%%%%%%%%%%%%%%%%%%%%%%%%%%%%%%%%%%%%%%%%%%%%%%%%%%%%%%%%%%%%%%%%%%%%%%%%%%%%%

In an attempt to understand other 4-manifolds via Seiberg-Witten
and Gromov invariants, some analysis of pseudoholomorphic curves
has been extended to nonsymplectic 4-manifolds
by equipping these with a {\em nearly nondegenerate closed 2-form}.
In particular, Taubes~\cite{ta:harmonic} has related Seiberg-Witten invariants
to pseudoholomorphic curves for compact oriented 4-manifolds with $b_2^+ > 0$.
Any compact oriented 4-manifold $M$ with $b_2^+ > 0$ admits a closed
2-form that vanishes along a union of circles and is symplectic
elsewhere~\cite{ga-ki:circles,ho:harmonic}.
In fact, for a generic metric on $M$,
there is a self-dual harmonic form $\omega$ which is transverse
to zero as a section of $\Lambda^2 T^* M$.
The vanishing locus of $\omega$ is the union of a finite number
of embedded circles, and $\omega$ is symplectic elsewhere.

The generic behavior of closed 2-forms on orientable 4-manifolds
is partially understood~\cite[pp.23-24]{ar-gi:symplectic}.
Here is a summary.
Let $\omega$ be a generic closed 2-form on a 4-manifold $M$.
At the points of some hypersurface $Z$, the form $\omega$ has rank 2.
At a generic point of $M$, $\omega$ is nondegenerate;
in particular, has the Darboux normal form
$dx_1 \wedge dy_1 + dx_2 \wedge dy_2$.
There is a codimension-1 submanifold $Z$
where $\omega$ has rank 2, and there are no
points where $\omega$ vanishes.
%$Z := \{ p \in M \mid \omega_p \wedge \omega_p = 0 \}$
At a generic point of $Z$, the kernel of $\widetilde \omega$
is transverse to $Z$; the normal form near such a point
is $x_1 dx_1 \wedge dy_1 + dx_2 \wedge dy_2$.
There is a curve $C$ where the kernel of $\widetilde \omega$
is not transverse to $Z$, hence sits in $TZ$.
At a generic point of $C$, the kernel of $\widetilde \omega$
is transverse to $C$;
there are two possible normal forms near such points,
called {\em elliptic} and {\em hyperbolic},
$d(x-\frac{z^2}{2}) \wedge dy + d(xz \pm ty-\frac{z^3}{3}) \wedge dt$.
The hyperbolic and elliptic sections of $C$ are separated
by {\em parabolic} points, where the kernel is tangent to $C$.
It is known that there exists at least one continuous family
of inequivalent degeneracies in a parabolic
neighborhood~\cite{go-ti:moduli}.

%%%%%%%%%%%%%%%%%%%%%%%%%%%%%%%%%%%%%%%%%%%%%%%%%%%%%%%%%%%%%%%%%%%%%%%%%%%%%
%%%%%%%%%%%%%%%%%%%%%%%%%%%%%%%%%%%%%%%%%%%%%%%%%%%%%%%%%%%%%%%%%%%%%%%%%%%%%

\ssubsection{Lefschetz Pencils}
\label{sec:pencils}

{\em Lefschetz pencils} in symplectic geometry imitate
linear systems in complex geometry.
Whereas holomorphic functions on a projective surface must be
constant, there are interesting functions on the complement of
a finite set, and generic such functions have only quadratic singularities.
A Lefschetz pencil can be viewed as a complex Morse function
or as a very singular fibration, in the sense that, not only
some fibers are singular (have ordinary double points)
but all {\em fibers} go through some points.

\begin{definition}
A \textbf{Lefschetz pencil} on an oriented 4-manifold $M$ is a map
$f : M \setminus \{ b_1 , \ldots , b_n \} \to \CC \PP ^1$
defined on the complement of a finite set in $M$,
called the \textbf{base locus},
that is a submersion away from a finite set $\{ p_1 , \ldots , p_{n+1} \}$,
and obeying local models $(z_1,z_2) \mapsto z_1 / z_2$
near the $b_j$'s and $(z_1,z_2) \mapsto z_1 z_2$ near the $p_j$'s,
where $(z_1,z_2)$ are oriented local complex coordinates.
\end{definition}

Usually it is also required that each fiber contains at most
one singular point.
By blowing-up $M$ at the $b_j$'s, we obtain a map to $\CC \PP^1$
on the whole manifold, called a \textbf{Lefschetz fibration}.
Lefschetz pencils and Lefschetz fibrations
can be defined on higher dimensional manifolds
where the $b_j$'s are replaced by codimension 4 submanifolds.
By working on the Lefschetz fibration,
Gompf~\cite{go:characterization,go:Lefschetz}
proved that a structure of Lefschetz pencil
(with a nontrivial base locus)
gives rise to a symplectic form, canonical up to isotopy,
such that the fibers are symplectic.

Using asymptotically holomorphic techniques~\cite{au:asymptotically,do:almost},
Donaldson~\cite{do:pencils} proved that symplectic 4-manifolds
admit Lefschetz pencils.
More precisely:

\begin{theorem}
\label{thm:donaldson_lefschetz}
\textbf{(Donaldson)} $\;$
Let $J$ be a compatible almost complex structure
on a compact symplectic 4-manifold $(M,\omega)$
where the class $[\omega]/2\pi$ is integral.
Then $J$ can be deformed through almost complex structures
to an almost complex structure $J'$ such that $M$ admits
a Lefschetz pencil with $J'$-holomorphic fibers.
\end{theorem}

The closure of a smooth fiber of the Lefschetz pencil is a
symplectic submanifold Poincar\'e dual to $k [\omega]/2\pi$;
cf.\ Theorem~\ref{thm:donaldson_submanifolds}.
Other perspectives on Lefschetz pencils have been explored,
including in terms of representations of the free group
$\pi_1 (\CC \PP ^1 \setminus \{ p_1 , \ldots , p_{n+1} \})$
in the mapping class group $\Gamma_g$ of the generic fiber
surface~\cite{sm:monodromy}.

Similar techniques were used
by Auroux~\cite{au:branched} to realize symplectic 4-manifolds
as {\em branched covers} of $\CC \PP ^2$,
and thus reduce the classification of symplectic 4-manifolds
to a (hard) algebraic question about factorization in the braid group.
Let $M$ and $N$ be compact oriented 4-manifolds,
and let $\nu$ be a symplectic form on $N$.

\begin{definition}
A map $f:M \to N$ is a \textbf{symplectic branched cover}
if for any $p \in M$ there are complex charts centered at $p$
and $f(p)$ %with domains $\cU$ and $\cV$, respectively, and 
such that $\nu$ is positive on each complex line
and where $f$ is given by: %one of the following three models:
a local diffeomorphism $(x,y) \to (x,y)$, or
a simple branching $(x,y) \to (x^2,y)$,
or an ordinary cusp $(x,y) \to (x^3 -xy,y)$.
\end{definition}

\vspace*{-2ex}

\begin{theorem}
\label{thm:auroux_branched}
\textbf{(Auroux)} $\;$
Let $(M,\omega)$ be a compact symplectic 4-manifold
where the class $[\omega]$ is integral,
and let $k$ be a sufficiently large integer.
Then there is a symplectic branched cover $f_k :(M,k\omega) \to \CC\PP^2$,
that is canonical up to isotopy for $k$ large enough.
Conversely, given a symplectic branched cover $f:M \to N$,
the domain $M$ inherits a symplectic form
canonical up to isotopy in the class $f^*[\nu]$.
\end{theorem}

%%%%%%%%%%%%%%%%%%%%%%%%%%%%%%%%%%%%%%%%%%%%%%%%%%%%%%%%%%%%%%%%%%%%%%%%%%%%%
%%%%%%%%%%%%%%%%%%%%%%%%%%%%%%%%%%%%%%%%%%%%%%%%%%%%%%%%%%%%%%%%%%%%%%%%%%%%%
% --> Section 5
%%%%%%%%%%%%%%%%%%%%%%%%%%%%%%%%%%%%%%%%%%%%%%%%%%%%%%%%%%%%%%%%%%%%%%%%%%%%%%
%%%%%%%%%%%%%%%%%%%%%%%%%%%%%%%%%%%%%%%%%%%%%%%%%%%%%%%%%%%%%%%%%%%%%%%%%%%%%

\newpage

\ssection{Hamiltonian Geometry}
\index{moment map}
\label{section5}

%%%%%%%%%%%%%%%%%%%%%%%%%%%%%%%%%%%%%%%%%%%%%%%%%%%%%%%%%%%%%%%%%%%%%%%%%%%%%
%%%%%%%%%%%%%%%%%%%%%%%%%%%%%%%%%%%%%%%%%%%%%%%%%%%%%%%%%%%%%%%%%%%%%%%%%%%%%

\ssubsection{Symplectic and Hamiltonian Vector Fields}
\label{sec:symplectic_hamiltonian_fields}
\index{hamiltonian ! vector field}\index{symplectic ! vector field}
\index{vector field ! hamiltonian}\index{vector field ! symplectic}

Let $(M,\omega)$ be a symplectic manifold and let
$H : M \to \RR$ be a smooth function.
By nondegeneracy, there is
a unique vector field $X_{_H}$ on $M$ such that
$\imath_{X_{H}}\omega=dH$.
Supposing that $X_{_H}$ is complete
(this is always the case when $M$ is compact),
let $\rho_{t} : M \to M$, $t \in \RR$,
be its flow (cf.\ Section~\ref{cotangent_bundles}).
Each diffeomorphism $\rho_{t}$ preserves $\omega$,
i.e., $\rho_{t}^* \omega = \omega$, because
$\frac{d}{dt}\rho_{t}^{*}\omega = \rho_{t}^{*}\cL_{X_{_H}}\omega
= \rho_{t}^{*}(d \imath_{X_{_H}}\omega + \imath_{X_{_H}} d\omega ) = 0$.
Therefore, every function on $(M,\omega)$ produces a
family of symplectomorphisms.
Notice how this feature involves
both the {\em nondegeneracy} and the {\em closedness} of $\omega$.

\begin{definition}
A vector field $X_{_H}$ such that
$\imath_{X_{H}}\omega=dH$ for some $H \in C^\infty (M)$
is a \textbf{hamiltonian vector field}\index{hamiltonian !
vector field}\index{vector field ! hamiltonian}
with \textbf{hamiltonian function}\index{hamiltonian ! function} $H$.
\end{definition}

Hamiltonian vector fields preserve their hamiltonian functions
($\cL_{X_{H}}H=\imath_{X_{H}}dH$ $=\imath_{X_{H}}\,\imath_{X_{H}}\omega=0$),
so each integral curve\index{integral ! curve}
$\{ \rho_t (x) \mid t \in \RR \}$
of a hamiltonian vector field $X_{_H}$ must be contained in a level set
of the hamiltonian function $H$.
%, i.e., $H(x) = (\rho_t^* H)(x) = H (\rho_t (x))$, for all $t$.
In $(\RR^{2n}, \omega_0 = \sum dx_j \wedge dy_j)$,
the {\em symplectic gradient}
$X_{_H} = \sum \left( \frac{\partial H}{\partial
y_j} \frac{\partial} {\partial x_j}-\frac{\partial H}{\partial x_j}
\frac{\partial} {\partial y_j} \right)$ and the usual (euclidean) gradient
$\nabla H = \sum_j \left (\frac{\partial H}{\partial x_j}
\frac{\partial}{\partial x_j}+\frac{\partial H}{\partial y_j}
\frac{\partial}{\partial y_j}\right)$
of a function $H$ are\index{gradient vector
field}\index{vector field ! gradient}
related by $JX_{_H}=\nabla H$,
where $J$ is the standard almost complex structure.

\begin{examples}
\begin{enumerate}

\item
For the height function $H(\theta,h)=h$ on the sphere
$(M,\omega) = (S^{2},d\theta \wedge dh)$, from
$\imath_{X_H}(d\theta \wedge dh)=dh$ we get
$X_{_H}=\frac{\partial}{\partial\theta}$.
Thus, $\rho_{t}(\theta,h)=(\theta+t,h)$,
which is rotation about the vertical axis,
preserving the height $H$.

\item
Let $X$ be any vector field on a manifold $W$.
There is a unique vector field $X_{\sharp}$ on the
cotangent bundle $T^*W$ whose flow is the lift\index{lift !
of a vector field} of the flow of $X$.
Let $\alpha$ be the tautological form
and $\omega = -d \alpha$ the canonical symplectic
form on $T^*W$.
The vector field $X_{\sharp}$ is hamiltonian
with hamiltonian function $H := \imath_{X_{\sharp}} \alpha$.

\item
Consider euclidean space $\RR^{2n}$ with coordinates
$(q_{1},\ldots,q_{n},p_{1},\ldots,p_{n})$
and $\omega_{0}=\sum dq_j\wedge dp_j$.
The curve $\rho_{t}=(q(t),p(t))$ is an integral curve for
a hamiltonian vector field $X_{_H}$ exactly
when it satisfies the \textbf{Hamilton equations}:\index{Hamilton
equations}\index{classical mechanics}\index{mechanics ! classical}
\[
        \left \{ \begin{array}{l}
        \frac{dq_i}{dt}(t) =
        \phantom{-}\frac{\partial H}{\partial p_i} \\
        \frac{dp_i}{dt}(t) =
        -\frac{\partial H}{\partial q_i}
        \end{array} \right.
\]

\item
\textbf{Newton's second law}\index{Newton ! second law}
states that a particle of mass $m$ moving
in \textbf{configuration space}\index{space !
configuration}\index{configuration space}
$\RR^3$ with coordinates $q=(q_1,q_2,q_3)$ under a potential $V(q)$
moves along a curve $q(t)$ such that 
\[ 
        m\frac{d^2 q}{d t^2} = -\nabla V(q)\ .
\]
Introduce the \textbf{momenta}\index{momentum}
$p_i=m\frac{dq_i}{dt}$ for $i=1,2,3$, and
\textbf{energy}\index{energy ! classical mechanics}
function $H(q,p)=\frac{1}{2m}|p|^2+V(q)$ on the
\textbf{phase space}\footnote{The \textbf{phase space} of a system
of $n$ particles is the space parametrizing the position and momenta
of the particles.
The mathematical model for a phase space
is a symplectic manifold.}\index{space ! phase}\index{phase space}
$\RR^6 = T^* \RR^3$
with coordinates $(q_{1},q_{2},q_{3},p_{1},p_{2},p_{3})$.
The energy $H$ is conserved by the motion
and Newton's second law\index{Newton ! second law}
in $\RR^3$ is then equivalent to the Hamilton
equations\index{Hamilton equations} in $\RR^6$:
\[
\left\{ \begin{array}{l}
        \frac{dq_i}{dt}=\frac{1}{m}p_i=\frac{\partial H}{\partial p_i} \\
        \frac{dp_i}{dt}=m \frac{d^2q_i}{dt^2}=
        -\frac{\partial V}{\partial q_i}=-\frac{\partial H}{\partial q_i}
\end{array} \right.
\]
\end{enumerate}
\end{examples}

\vspace*{-2ex}

\begin{definition}
A vector field $X$ on $M$ preserving $\omega$
(i.e., such that $\cL_{X}\omega=0$) is a
\textbf{symplectic vector field}.\index{symplectic !
vector field}\index{vector field ! symplectic}
\end{definition}

Hence, a vector field $X$ on $(M,\omega)$
is called \textbf{symplectic} when $\imath_X\omega$ is closed,
and \textbf{hamiltonian}\index{vector field !
hamiltonian}\index{hamiltonian ! vector field}
when $\imath_X \omega$ is exact.
In the latter case, a {\em primitive} $H$ of $\imath_X\omega$
is called a \textbf{hamiltonian function}\index{hamiltonian !
function}\index{function ! hamiltonian} of $X$.
On a contractible open set
every symplectic vector field is hamiltonian.
Globally, the group $H_{\rm de Rham}^{1}(M)$ measures the obstruction
for symplectic vector fields to be hamiltonian.
For instance, the vector field $X_{1}=\frac{\partial}{\partial\theta_{1}}$
on the 2-torus $(M,\omega) = (\TT^{2},d\theta_{1}\wedge d\theta_{2})$
is symplectic but not hamiltonian.

%%%%%%%%%%%%%%%%%%%%%%%%%%%%%%%%%%%%%%%%%%%%%%%%%%%%%%%%%%%%%%%%%%%%%%%%%%%%%

A vector field $X$ is a differential operator on functions:
$X \cdot f := \cL_{_X} f = df (X)$ for $f \in C^\infty (M)$.
As such, the bracket $W = [X,Y]$ is the commutator:
$\cL_W = [ \cL_X , \cL_Y ] = \cL_X \cL_Y - \cL_Y \cL_X$
(cf.\ Section~\ref{sec:integrability}).
This endows the set $\chi(M)$ of vector fields on a manifold $M$ with
a structure of {\em Lie algebra}.\footnote{A (real)
\textbf{Lie algebra}\index{Lie ! algebra}
is a (real) vector space $\fg$ together with a
\textbf{Lie bracket}\index{Lie ! bracket} $[ \cdot, \cdot ]$,
i.e., a bilinear map $[ \cdot, \cdot ] : \fg \times \fg \to \fg$
satisfying \textbf{antisymmetry},
$[x,y] = -[y,x]$, $\forall x,y \in \fg$,\index{antisymmetry}
and the \textbf{Jacobi identity},
$[x,[y,z]] + [y,[z,x]] + [z,[x,y]] = 0$,
$\forall x,y,z \in \fg$.\index{Jacobi ! identity}}
For a symplectic manifold $(M,\omega)$,
using $\imath_{[X,Y]} = [ \cL_X , \imath_Y ]$
and Cartan's magic formula, we find that
$\imath_{[X,Y]}\omega =
%\cL_{X}\imath_{Y}\omega-\imath_{Y}\cL_{X}\omega =
d\imath_{X}\imath_{Y}\omega + \imath_{X} d\imath_{Y}\omega -
\imath_{Y} d\imath_{X}\omega - \imath_{Y}\imath_{X} d\omega
= d(\omega(Y,X))$.
Therefore:

\begin{proposition}
\label{prop:brackets}
If $X$ and $Y$ are symplectic vector fields on
a symplectic manifold $(M,\omega)$,
then $[X,Y]$ is hamiltonian with hamiltonian function $\omega(Y,X)$.
\end{proposition}

Hence, hamiltonian vector fields and symplectic vector fields
form Lie subalgebras for the Lie bracket $[\cdot,\cdot]$.

\begin{definition}
The \textbf{Poisson bracket}\index{bracket !
Poisson}\index{Poisson ! bracket} of two functions $f,g\in C^{\infty}(M)$
is the function $\{f,g\} := \omega(X_f,X_g) = \cL_{X_g} f$.
\end{definition}

By Proposition~\ref{prop:brackets} we have $X_{\{f,g\}}=-[X_{f},X_{g}]$.
Moreover, the bracket $\{\cdot,\cdot\}$
satisfies the \textbf{Jacobi identity}\index{Jacobi identity},
$\{f,\{g,h\}\}+\{g,\{h,f\}\}+\{h,\{f,g\}\}=0$,
and the \textbf{Leibniz rule}\index{Leibniz rule},
$\{ f, gh \} = \{ f,g \} h + g \{ f,h \}$.

\begin{definition}
A \textbf{Poisson algebra}\index{Poisson ! algebra}
$(\cP , \{\cdot,\cdot\})$
is a commutative associative algebra $\cP$
with a Lie bracket $\{\cdot,\cdot\}$ satisfying the
Leibniz rule\index{Leibniz rule}.
\end{definition}

When $(M,\omega)$ is a symplectic manifold,
$(C^{\infty}(M), \{\cdot,\cdot\})$ is a Poisson algebra,
and the map $C^{\infty}(M) \to \chi(M)$,
$H \mapsto X_{_H}$ is a Lie algebra anti-homomorphism.

\begin{examples}
\begin{enumerate}
\item
For the prototype $(\RR^{2n}, \sum dx_i \wedge dy_i)$,
we have $X_{x_i} = - \frac{\partial}{\partial y_i}$ and
$X_{y_i} = \frac{\partial}{\partial x_i}$, so that
$\{ x_i , x_j \} = \{ y_i , y_j \} = 0$
and $\{ x_i , y_j \} = \delta_{ij}$ for all $i,j$.
Arbitrary functions $f,g \in C^\infty (\RR^{2n})$ have
the \textbf{classical Poisson bracket}
\[
   \{ f, g \} = \sum \limits_{i=1}^{n}
   \left( \frac{\partial f}{\partial x_i}\frac{\partial g}{\partial y_i}
   - \frac{\partial f}{\partial y_i}\frac{\partial g}{\partial x_i} \right) \ .
\]

\item
Let $G$ be a Lie group,\footnote{A \textbf{Lie group} is
a manifold $G$ equipped with a group structure where the
group operation $G \times G \to G$ and inversion $G \to G$
are smooth maps.
An \textbf{action} of a Lie group $G$ on a manifold $M$ is a
group homomorphism $G \to \mathrm{Diff}(M)$, $g \mapsto \psi_g$,
where the \textbf{evaluation map} $M \times G \to M$,
$(p,g) \mapsto \psi_g (p)$ is a smooth map.
The \textbf{orbit}\index{orbit ! definition}
of $G$ through $p \in M$ is $\{\psi_g(p) \mid g \in G\}$.
The \textbf{stabilizer}\index{stabilizer}
(or {\em isotropy}\index{isotropy}) of $p \in M$ is
$G_p := \{g \in G \mid \psi_g(p) = p\}$.}
${\fg}$ its Lie algebra and
${\fg}^*$ the dual vector space of ${\fg}$.
The vector field $^\fg X^\#$ generated by $X \in \fg$
for the adjoint action\footnote{Any Lie group $G$ acts on itself by
\textbf{conjugation}:\index{conjugation}
$g \in G \mapsto \psi_g \in \mathrm{Diff}(G)$,
$\psi_g(a) = g \cdot a \cdot g^{-1}$.
Let $\mathrm{Ad}_g: \fg \to \fg$ be the derivative
at the identity of $\psi_g: G \to G$.
We identify the Lie algebra $\fg$ with the tangent space $T_eG$.
For matrix groups, $\mathrm{Ad}_g X = g X g^{-1}$.
Letting $g$ vary, we obtain the
\textbf{adjoint action}\index{adjoint ! action}\index{action ! adjoint}
of $G$ on its Lie algebra
$\mathrm{Ad}: G \to \GL(\fg)$.
%, $g \mapsto \mathrm{Ad}_g$.
Let $\langle \cdot,\cdot \rangle : \fg^* \times \fg \to \RR$
be the natural pairing $\langle \xi,X\rangle = \xi(X)$.
Given $\xi \in \fg^*$, we define $\mathrm{Ad}_g^*\xi$ by
$\langle\mathrm{Ad}_g^*\xi,X\rangle = 
\langle \xi,\mathrm{Ad}_{g^{-1}}X\rangle$, for any $X \in \fg$.
The collection of maps $\mathrm{Ad}_g^*$ forms the \textbf{coadjoint
action}\index{coadjoint ! action}\index{action !
coadjoint} of $G$ on the dual of its Lie algebra
$\mathrm{Ad}^*: G \to \GL(\fg^*)$.
%, $g \mapsto \mathrm{Ad}_g^*$.
These satisfy $\mathrm{Ad}_g \circ \mathrm{Ad}_h = \mathrm{Ad}_{gh}$
and $\mathrm{Ad}_g^* \circ \mathrm{Ad}_h^* = \mathrm{Ad}_{gh}^*$.}
of $G$ on $\fg$ has value $[X,Y]$ at $Y \in \fg$.
The vector field $X^\#$ generated by $X \in \fg$
for the coadjoint action of $G$ on $\fg^*$ is
$\langle X^\#_{_\xi} , Y \rangle = \langle \xi , [Y,X] \rangle$,
$\forall \ \xi \in \fg^* , Y \in \fg$.
The skew-symmetric pairing $\omega$ on $\fg$ defined at $\xi \in \fg^*$ by
\[
   \omega_{_\xi} (X,Y) := \langle \xi , [X,Y] \rangle
\]
has kernel at $\xi$ the Lie algebra $\fg_{_\xi}$ of the
stabilizer of $\xi$ for the coadjoint action.
Therefore, $\omega$ restricts to a nondegenerate 2-form
on the tangent spaces to the orbits of the coadjoint action.
As the tangent spaces to an orbit are generated by the vector
fields $X^\#$, the Jacobi identity in $\fg$ implies
that this form is closed.
It is called the \textbf{canonical symplectic form}\index{canonical !
symplectic form on a coadjoint orbit}\index{symplectic !
canonical symplectic form on a coadjoint orbit} (or the
\textbf{Lie-Poisson}\index{Lie-Poisson symplectic form}\index{Poisson !
Lie-Poisson symplectic form} or
\textbf{Kirillov-Kostant-Souriau symplectic structure}\index{Kostant-Kirillov
symplectic form}\index{Kirillov|see{Kostant-Kirillov}})
on the \textbf{coadjoint orbits}.
The corresponding Poisson structure on
$\fg^*$\index{Poisson ! structure on $\fg^*$}
is the canonical one induced by the Lie bracket:
\[
        \{ f , g \} (\xi) = \langle \xi, [df_{_\xi}, dg_{_\xi}] \rangle
\]
for $f,g \in C^\infty (\fg^*)$ and $\xi \in {\fg}^*$.
The differential
$df_{_\xi} : T_{_\xi} {\fg}^* \simeq {\fg}^* \to \RR$ is identified with
an element of $\fg \simeq \fg^{**}$.

\end{enumerate}
\end{examples}

%%%%%%%%%%%%%%%%%%%%%%%%%%%%%%%%%%%%%%%%%%%%%%%%%%%%%%%%%%%%%%%%%%%%%%%%%%%%%
%%%%%%%%%%%%%%%%%%%%%%%%%%%%%%%%%%%%%%%%%%%%%%%%%%%%%%%%%%%%%%%%%%%%%%%%%%%%%

\ssubsection{Arnold Conjecture and Floer Homology}
\label{sec:arnold_floer}

There is an important generalization of Poincar\'e's last
geometric theorem (Theorem~\ref{thm:poincare_birkhoff})
conjectured by Arnold starting around~1966.\index{symplectomorphism !
Arnold conjecture}\index{Arnold !
conjecture}\index{conjecture ! Arnold}
Let $(M,\omega)$ be a compact symplectic manifold,
and $h_t: M \to \RR$ a 1-periodic (i.e., $h_t = h_{t+1}$)
smooth family of functions.
Let $\rho: M \times \RR \to M$ be the isotopy generated by
the time-dependent hamiltonian vector field $v_t$
defined by the equation $\omega(v_t,\cdot) = dh_t$.
The symplectomorphism $\varphi = \rho_1$ is then said to be
\textbf{exactly homotopic to the identity}.\index{symplectomorphism !
Arnold conjecture}\index{Arnold ! conjecture}\index{conjecture !
Arnold}\index{symplectomorphism !
exactly homotopic to the identity}\index{exactly homotopic to the identity}
In other words, a symplectomorphism exactly homotopic to the identity
is the time-1 map of the isotopy generated by some
time-dependent 1-periodic hamiltonian function.
There is a one-to-one correspondence between the
fixed points of $\varphi$ and the period-1 orbits of $\rho$.
When all the fixed points of such $\varphi$ are nondegenerate (generic case),
we call $\varphi$ \textbf{nondegenerate}.
The \textbf{Arnold conjecture}~\cite[Appendix~9]{ar:mathematical}
predicted that
\[
        \# \{ \mbox{fixed points of a nondegenerate } \varphi\} \geq
        \displaystyle{\sum_{i=0}^{2n}} \dim H^i(M;\RR)
\]
(or even that the number of fixed points of a nondegenerate $\varphi$
is at least the minimal number of critical points of
a Morse function\footnote{A
\textbf{Morse function} is a smooth function
$f:M \to \RR$ all of whose critical points are nondegenerate,
i.e., at any critical point the hessian matrix is nondegenerate.}).
When the hamiltonian $h: M \to \RR$ is independent of $t$,
this relation is trivial:
a point $p$ is critical for $h$ if and only if $dh_p = 0$,
if and only if $v_p = 0$, if and only if
$\rho (t,p) = p$, $\forall t \in \RR$, which implies
that $p$ is a fixed point of $\rho_1=\varphi$,
so the Arnold conjecture reduces to a Morse inequality.
Notice that, according to the Lefschetz fixed point theorem,\index{theorem !
Lefschetz fixed point}\index{Lefschetz fixed point theorem}
the Euler characteristic of $M$, i.e., the {\em alternating}
sum of the Betti numbers, $\sum (-1)^ i \dim H^i(M;\RR)$,
is a (weaker) lower bound for the number of fixed points of $\varphi$.
%The Arnold conjecture can be viewed as an instance of
%symplectic rigidity related to some Morse theory,
%as the sum of the Betti numbers $\sum \dim H^i(M;\RR)$ is a lower bound
%for the number of critical points of a Morse function on $M$.

The Arnold conjecture\index{Arnold ! conjecture}\index{conjecture !
Arnold} was gradually proved from the late~70's to the late~90's
by Eliashberg~\cite{el:fixed_points}, Conley-Zehnder~\cite{co-ze:arnold},
Floer~\cite{fl:holomorphic}, Sikorav~\cite{si:points},
Weinstein~\cite{we:conley_zehnder}, Hofer-Salamon~\cite{ho-sa:floer},
Ono~\cite{on:arnold}, culminating with independent proofs by
Fukaya-Ono~\cite{fu-on:arnold} and Liu-Tian~\cite{li-ti:arnold}.
There are open conjectures for sharper bounds on the number of fixed points.
The breakthrough tool for establishing the Arnold conjecture was
\textbf{Floer homology} -- an $\infty$-dimensional analogue of Morse theory.
Floer homology was defined by
Floer~\cite{fl:relative,fl:gradient,fl:lagrangian,fl:holomorphic,fl:witten}
and developed through
the work of numerous people after Floer's death.
It combines the variational approach of Conley and Zehnder~\cite{co-ze:morse},
with Witten's Morse-Smale complex~\cite{wi:morse_theory},
and with Gromov's compactness theorem for
pseudo-holomorphic curves~\cite{gr:pseudo}.

%Instead of a real function $f$,
Floer theory starts from a symplectic action functional
on the space of loops $\cL M$ of a symplectic manifold $(M,\omega)$
whose zeros of the differential $dF : T (\cL M) \to \RR$
are the period-1 orbits of the isotopy $\rho$ above.
The tangent bundle $T (\cL M)$ is the space of loops with vector
fields over them: pairs $(\ell,v)$, where $\ell : S^1 \to M$
and $v : S^1 \to \ell^* (TM)$ is a section.
Then $df (\ell,v) = \int_0^1 \omega
( \dot \ell (t) - X_{h_t} (\ell (t) , v(t) ) \, dt$.
The {\em Floer complex}\footnote{The \textbf{Morse complex}
for a Morse function on a compact manifold, $f:M \to \RR$,
is the chain complex freely generated by the critical points of $f$,
graded by the {\em Morse index} $\imath$
and with differential given by counting the number $n(x,y)$ of flow lines
of the negative gradient $- \nabla f$ (for a metric on $X$)
from the point $x$ to the point $y$ whose indices differ by 1:
\[
   C_* = \oplus_{x \in \mathrm{Crit} (f)} \ZZ \langle x \rangle
   \quad \mbox{ and } \quad
   \partial \langle x \rangle = \sum \limits_{y \in \mathrm{Crit} (f),
   \imath (y) = \imath (x) -1} n(x,y) \langle y \rangle \ .
\]
The coefficient $n(x,y)$ is thus the number of solutions
(modulo $\RR$-reparametrization) $u : \RR \to X$
of the ordinary differential equation
$\frac{d}{dt} u(t) = - \nabla f (u(t))$ with conditions
$\lim_{t \to -\infty} u(t) = x$, $\lim_{t \to +\infty} u(t) = y$.
The \textbf{Morse index} of a critical point
of $f$ is the dimension of its unstable manifold, i.e., the number of
negative eigenvalues of the hessian of $f$ at that point.
For a generic metric, the unstable manifold of a critical point
$W^u (x)$ intersects transversally with the stable manifold of
another critical point $W^s (y)$.
When $\imath (x) - \imath (y) = 1$,
the intersection $W^u (x) \cap W^s (y)$ has dimension 1,
so when we quotient out by the $\RR$-reparametrization
(to count actual image curves) we get
a discrete set, which is finite by compactness.
That $(C_*,\partial)$ is indeed a complex, i.e., $\partial ^2 =0$,
follows from counting broken flow lines between points
whose indices differ by 2.
Morse's theorem states that the homology of the Morse complex
coincides with the ordinary homology of $M$.
In particular, the sum of all the Betti numbers
$\sum \dim H^i(M;\RR)$ is a lower bound for the number of
critical points of a Morse function.}
is the chain complex freely generated by the critical points of $F$
(corresponding to the fixed points of $\varphi$),
with {\em relative grading} $\mathrm{index} (x,y)$
given by the difference in the number
of positive eigenvalues from the spectral flow.
The Floer differential is given by counting
the number $n(x,y)$ of pseudo-holomorphic surfaces
(the {\em gradient flow lines} joining two fixed points):
\[
   C_* = \oplus_{x \in \mathrm{Crit} (F)} \ZZ \langle x \rangle
   \quad \mbox{ and } \quad
   \partial \langle x \rangle = \sum \limits_{\footnotesize{\begin{array}{l}
   y \in \mathrm{Crit} (F) \\ \mathrm{index} (x,y) = 1
   \end{array}}} n(x,y) \langle y \rangle \ .
\]
Pondering transversality, compactness and orientation,
Floer's theorem states that the homology of $(C_* , \partial)$
is isomorphic to the ordinary homology of $M$.
In particular, the sum of the Betti numbers
is a lower bound for the number of fixed points of $\varphi$.

From the above {\em symplectic Floer homology},
Floer theory has branched out to tackle other differential
geometric problems in symplectic geometry and 3- and 4-dimensional topology.
It provides a rigorous definition of invariants
viewed as homology groups of infinite-dimensional Morse-type theories,
with relations to gauge theory and quantum field theory.
There is {\em lagrangian Floer homology} (for
%a version of the formulation before for --- Seidel
the case of lagrangian intersections,
i.e., intersection of a lagrangian submanifold
with a hamiltonian deformation of itself),
{\em instanton Floer homology} (for invariants of 3-manifolds),
{\em Seiberg-Witten Floer homology}, {\em Heegaard Floer homology}
and {\em knot Floer homology}.
For more on Floer homology, see for instance~\cite{do:floer,sa:pcmi}.

%%%%%%%%%%%%%%%%%%%%%%%%%%%%%%%%%%%%%%%%%%%%%%%%%%%%%%%%%%%%%%%%%%%%%%%%%%%%%
%%%%%%%%%%%%%%%%%%%%%%%%%%%%%%%%%%%%%%%%%%%%%%%%%%%%%%%%%%%%%%%%%%%%%%%%%%%%%

\ssubsection{Euler-Lagrange Equations}
\label{sec:euler_lagrange}

The equations of motion in classical mechanics arise
from \textbf{variational principles}.\index{principle !
of least action}\index{equations !
of motion}\index{motion ! equations}\index{variational ! problem}
%as solutions of {\em variational problems}.
The physical path of a general mechanical
system\index{mechanical system}\index{system ! mechanical}
of $n$ particles is the path that {\em minimizes}
a quantity called the {\em action}.
When dealing with systems with constraints,
such as the simple pendulum,\index{pendulum ! simple}
or two point masses attached by a rigid rod, or a rigid body,
the language of variational principles becomes
more appropriate than the explicit analogues of
Newton's second laws.\index{Newton ! second law}
Variational principles\index{variational ! principle}
are due mostly to D'Alembert,\index{D'Alembert !
variational principle}\index{principle !
variational} Maupertius\index{Maupertius !
variational principle}, Euler\index{Euler !
variational principle} and Lagrange\index{Lagrange ! variational principle}.

%%%%%%%%%%%%%%%%%%%%%%%%%%%%%%%%%%%%%%%%%%%%%%%%%%%%%%%%%%%%%%%%%%%%%%%%%%%%%

Let $M$ be an $n$-dimensional manifold,
and let $F: TM \to \RR$ be a function on its tangent bundle.
If $\gamma: [a,b] \to M$ is a curve on $M$,
the \textbf{lift of $\gamma$ to $TM$}\index{lift ! of a path}
is the curve on $TM$ given by
${\widetilde \gamma}: [a,b] \to TM$,
$t \mapsto \left( \gamma(t), \frac{d\gamma}{dt} (t) \right)$.
The \textbf{action}\index{action ! of a path} of $\gamma$ is
\[
        \cA_{\gamma} := 
        \displaystyle{\int_a^b ({\widetilde \gamma}^*F)(t)dt} =
        \displaystyle{\int_a^b
        F\left( \gamma(t), \frac{d\gamma}{dt} (t) \right) dt}\ .
\]
For fixed $p,q$, let
$\cP(a,b,p,q) = \{\gamma: [a,b] \to M \mbox{ smooth} \mid
\gamma(a) = p, \gamma(b) = q\}$.
The goal is to find, among all $\gamma \in \cP(a,b,p,q)$,
the curve that {\em locally minimizes} $\cA_{\gamma}$.
(Minimizing curves are always locally
minimizing.)\index{action ! minimizing}\index{minimizing !
action}\index{minimizing ! locally}
Assume that $p$, $q$ and the image of $\gamma$
lie in a coordinate neighborhood $(\cU,x_1,\dots,x_n)$.
On $T\cU$ we have coordinates
$(x_1,\dots,x_n,v_1,\dots,v_n)$ associated with a trivialization of $T\cU$
by $\frac {\partial}{\partial x_1},\dots,\frac {\partial}{\partial x_n}$.
Using this trivialization, a curve $\gamma: [a,b] \to \cU$,
$\gamma(t) = (\gamma_1(t),\dots,\gamma_n(t))$ lifts to
\[
        {\widetilde \gamma}: [a,b] \longrightarrow T\cU\ , \qquad
        {\widetilde \gamma}(t) = \left(\gamma_1(t),\dots,\gamma_n(t),
        \frac {d\gamma_1}{dt} (t),\dots,\frac {d\gamma_n}{dt} (t)\right)\ .
\]
Consider infinitesimal variations of $\gamma$.
Let $c_1,\dots,c_n \in C^{\infty}([a,b])$ be such that $c_k(a) = c_k(b) = 0$.
For $\varepsilon$ small, let $\gamma_{\varepsilon}: [a,b] \to \cU$
be the curve $\gamma_{\varepsilon}(t) =
(\gamma_1(t) + \varepsilon c_1(t),\dots,\gamma_n(t) + \varepsilon c_n(t))$.
Let $\cA_{\varepsilon} := \cA_{\gamma_{\varepsilon}}$.
A necessary condition for $\gamma = \gamma_0 \in \cP(a,b,p,q)$
to minimize the action is that $\varepsilon = 0$
be a critical point of $\cA_{\varepsilon}$.
By the Leibniz rule and integration by parts, we have that
\[
\begin{array}{rl}
        \displaystyle{\frac {d\cA_{\varepsilon}}{d\varepsilon} (0)} & =
        \displaystyle{\int_a^b \sum_k \left[ \frac {\partial F}{\partial x_k}
        \left( \gamma_0(t), \frac {d\gamma_0}{dt} (t)\right)
        c_k(t) + \frac {\partial F}{\partial v_k}
        \left( \gamma_0, \frac {d\gamma_0}{dt} \right) \frac
        {dc_k}{dt} (t)\right] dt} \\
        & =
        \displaystyle{\int_a^b \sum_k \left[ \frac {\partial F}{\partial x_k}
        (\dots) - \frac {d}{dt} \frac {\partial F}{\partial v_k}
        (\dots)\right] c_k(t)\, dt \ .}
\end{array}
\]
For $\frac {d\cA_{\varepsilon}}{d\varepsilon} (0)$ to vanish
for all $c_k$'s satisfying boundary conditions $c_k(a) = c_k(b) = 0$,
the path $\gamma_0$
must satisfy the \textbf{Euler-Lagrange equations}\index{Euler !
Euler-Lagrange equations}:
\[
        \frac {\partial F}{\partial x_k}
        \left( \gamma_0(t), \frac {d\gamma_0}{dt} (t)\right) =
        \frac {d}{dt} \frac {\partial F}{\partial v_k}
        \left( \gamma_0(t), \frac{d\gamma_0}{dt} (t)\right) \ ,
        \quad k=1, \ldots , n \ .
\]

%%%%%%%%%%%%%%%%%%%%%%%%%%%%%%%%%%%%%%%%%%%%%%%%%%%%%%%%%%%%%%%%%%%%%%%%%%%%%

\begin{examples}
\begin{enumerate}
\item
Let $(M,g)$ be a riemannian manifold\index{riemannian !
manifold}\index{manifold ! riemannian}.
Let $F: TM \to \RR$ be the function
whose restriction to each tangent space is the
quadratic form defined by the riemannian metric.\index{riemannian !
metric}\index{metric}
On a coordinate chart $F(x,v) = |v|^2 = \sum g_{ij}(x) v^i v^j$.
Let $p,q \in M$ and $\gamma : [a,b] \to M$ a curve joining $p$ to $q$.
The {\em action}\index{action ! of a path} of $\gamma$ is
\[
   \cA_\gamma = \displaystyle{\int_a^b
   \left| {d \gamma \over dt} \right|^2 dt} \ .
\]
The Euler-Lagrange equations\index{Euler ! Euler-Lagrange
equations}\index{equations ! Euler-Lagrange}\index{Lagrange !
Euler-Lagrange equations} become
the \textbf{Christoffel equations}\index{Christoffel !
equations}\index{equations ! Christoffel} for a geodesic
\[
        {d^2 \gamma^k \over dt^2} + \sum (\Gamma_{ij}^k \circ \gamma)
        {d \gamma^i \over dt} {d \gamma^j \over dt} = 0\ ,
\]
where the \textbf{Christoffel symbols}\index{Christoffel ! symbols}
$\Gamma_{ij}^k$'s
are defined in terms of the coefficients of the riemannian metric
($g^{ij}$ is the matrix inverse to $g_{ij}$) by
\[
   \Gamma_{ij}^k = \frac{1}{2} \sum \limits_\ell
   g^{\ell k} \left( \frac{\partial g_{\ell i}}{\partial x_j} +
   \frac{\partial g_{\ell j}}{\partial x_i} -
   \frac{\partial g_{ij}}{\partial x_\ell } \right) \ .
\]

%%%%%%%%%%%%%%%%%%%%%%%%%%%%%%%%%%%%%%%%%%%%%%%%%%%%%%%%%%%%%%%%%%%%%%%%%%%%%

\item
Consider\index{example ! of mechanical system}
a point-particle of mass $m$ moving in $\RR^3$
under a \textbf{force field} $G$.
The \textbf{work}\index{work} of $G$
on a path $\gamma: [a,b] \to \RR^3$ is
$W_{\gamma} := \int_a^b G(\gamma(t)) \cdot \frac {d\gamma}{dt} (t) \, dt$.
Suppose that $G$ is \textbf{conservative},\index{system !
conservative}\index{conservative system} i.e.,
$W_{\gamma}$ depends only on the initial and final points,
$p= \gamma(a)$ and $q = \gamma(b)$.
We can define the \textbf{potential energy}\index{energy !
potential}\index{potential ! energy} as
$V: \RR^3 \to \RR$, $V(q) := W_{\gamma}$,
where $\gamma$ is a path joining a fixed base point $p_0 \in \RR^3$ to $q$.
%(the {\em origin})
Let $\cP$
be the set of all paths going from $p$ to $q$ over time $t \in [a,b]$.
By the \textbf{principle of least action},\index{action !
principle of least action}\index{principle ! of least action}
the physical path is the path $\gamma \in \cP$ that minimizes a kind of
mean value of kinetic minus potential energy\index{energy !
kinetic}\index{energy ! potential}\index{potential !
energy}\index{kinetic energy},
known as the \textbf{action}\index{action ! of a path}:
\[
        \cA_{\gamma} :=
        \int_a^b \left( \frac {m}{2} \left| \frac{d\gamma}{dt} (t)
        \right|^2 - V(\gamma(t))\right) dt\ .
\]
The Euler-Lagrange equations are then equivalent to
\textbf{Newton's second law}:\index{Newton ! second law}
\[
   m \frac {d^2x}{dt^2} (t) -\frac {\partial V}{\partial x} (x(t)) = 0
   \quad \iff \quad m \frac {d^2x}{dt^2} (t) = G(x(t)) \ .
\]
In the case of the earth moving about the sun,
both regarded as point-masses and assuming that the sun
to be stationary at the origin,
the \textbf{gravitational potential}\index{gravitational
potential}\index{potential ! gravitational}
$V(x) = \frac {\mbox{const.}}{|x|}$ yields the
\textbf{inverse square law} for the motion.\index{inverse square law}

%%%%%%%%%%%%%%%%%%%%%%%%%%%%%%%%%%%%%%%%%%%%%%%%%%%%%%%%%%%%%%%%%%%%%%%%%%%%%

\item
Consider now $n$ point-particles of masses $m_1,\dots,m_n$
moving in $\RR^3$ under a conservative force corresponding
to a potential energy $V \in C^{\infty}(\RR^{3n})$.
At any instant $t$, the configuration of this system
is described by a vector $x = (x_1,\dots,x_n)$
in configuration space $\RR^{3n}$,
where $x_k \in \RR^3$ is the position of the $k$th particle.
For fixed $p,q \in \RR^{3n}$, let $\cP$ be the set of all paths
$\gamma = (\gamma_1 , \ldots , \gamma_n) : [a,b] \to \RR^{3n}$
from $p$ to $q$.
The \textbf{action}\index{action ! of a path} of a path $\gamma \in \cP$ is
\[
        \cA_{\gamma} :=
        \int_a^b \left( \sum \limits_{k=1}^n
        \frac {m_k}{2} \left| \frac{d\gamma_k}{dt} (t)
        \right|^2 - V(\gamma(t))\right) dt\ .
\]
The Euler-Lagrange equations
reduce to Newton's law for each particle.
%\[
%        m_k \frac {d^2x_{k,i}}{dt^2} (t) =
%        -\frac {\partial V}{\partial x_{k,i}}
%        (x_1(t),\dots,x_n(t))\ , \quad k=1, \ldots , n \ , i=1,2,3 \ .
%\]
Suppose that the particles are restricted to move on a
submanifold $M$ of $\RR ^{3n}$ called the
\textbf{constraint set}.\index{constraint set}\index{Newton !
second law}\index{system ! constrained}\index{constrained system}
By the \textbf{principle of least action for a constrained system},
the physical path has
minimal action among all paths satisfying the rigid constraints.
I.e., we single out the actual physical path as the one
that minimizes $\cA_{\gamma}$ among all
$\gamma: [a,b] \to M$ with $\gamma(a) = p$ and $\gamma(b) = q$.
\end{enumerate}
\end{examples}

%%%%%%%%%%%%%%%%%%%%%%%%%%%%%%%%%%%%%%%%%%%%%%%%%%%%%%%%%%%%%%%%%%%%%%%%%%%%%

In the case where $F(x,v)$ does not depend on $v$,
the Euler-Lagrange equations\index{Euler ! Euler-Lagrange equations}
are simply $\frac {\partial F}{\partial x_i} \left( \gamma_0(t),
\frac {d\gamma_0}{dt} (t)\right) = 0$.
These are satisfied if and only if
the curve $\gamma_0$ sits on the critical set of $F$.
For generic $F$, the critical points are isolated,
hence $\gamma_0(t)$ must be a constant curve.
In the case where $F(x,v)$ depends affinely on $v$,
$F(x,v) = F_0(x) + \sum_{j=1}^n F_j(x)v_j$,
the Euler-Lagrange equations become
\[
        \frac {\partial F_0}{\partial x_i} (\gamma(t)) = \sum_{j=1}^n
        \left( \frac {\partial F_i}{\partial x_j}
        - \frac {\partial F_j}{\partial x_i} \right)
        (\gamma(t)) \frac {d\gamma_j}{dt} (t)\ .
\]
If the $n \times n$ matrix $\left( \frac {\partial F_i}{\partial x_j}
- \frac {\partial F_j}{\partial x_i}\right)$ has an inverse $G_{ij}(x)$,
we obtain the system of first order ordinary differential equations
$\frac {d\gamma_j}{dt} (t) = \sum G_{ji}(\gamma(t))
\frac {\partial F_0}{\partial x_i} (\gamma(t))$.
Locally it has a unique solution through each point $p$.
If $q$ is not on this curve, there is no solution at all to the
Euler-Lagrange equations belonging to $\cP(a,b,p,q)$.

Therefore, we need non-linear dependence of $F$ on the $v$
variables in order to have appropriate solutions.
From now on, assume the
\textbf{Legendre condition}:\index{Legendre ! condition}
\[
   \displaystyle{ \det
   \left( \frac {\partial^2F}{\partial v_i \partial v_j} \right)} \ne 0\ .
\]
Letting $G_{ij}(x,v) = \left( \frac {\partial^2F}{\partial v_i
\partial v_j} (x,v) \right)^{-1}$, the Euler-Lagrange equations become
\[
        \frac {d^2\gamma_j}{dt^2} =
        \sum_i G_{ji} \frac {\partial F}{\partial x_i}
        \left( \gamma,\frac {d\gamma}{dt} \right) -
        \sum_{i,k} G_{ji} \frac {\partial^2F}{\partial
        v_i \partial x_k} \left( \gamma, \frac {d\gamma}{dt} \right)
        \frac {d\gamma_k}{dt}\ .
\]
This second order ordinary differential equation
has a unique solution given initial conditions
$\gamma(a) = p$ and $\frac {d\gamma}{dt} (a) = v$.
Assume that
$\left( \frac {\partial^2F}{\partial v_i \partial v_j} (x,v) \right) \gg 0$,
$\forall (x,v)$, i.e.,
with the $x$ variable frozen, the function $v \mapsto F(x,v)$ is
\textbf{strictly convex}\index{strictly convex function}.
Then the path $\gamma_0 \in \cP(a,b,p,q)$ satisfying the above
Euler-Lagrange equations does indeed locally
minimize\index{minimizing ! locally} $\cA_{\gamma}$
(globally it is only critical):

\begin{proposition}
For every sufficiently small subinterval $[a_1,b_1]$ of $[a,b]$,
$\gamma_0|_{[a_1,b_1]}$ is locally minimizing in $\cP(a_1,b_1,p_1,q_1)$
where $p_1 =
\gamma_0(a_1)$, $q_1 = \gamma_0(b_1)$.
\end{proposition}

\vspace*{-2ex}

\begin{proof}
Take $c = (c_1,\dots,c_n)$ with $c_i \in C^{\infty}([a,b])$,
$c_i(a) = c_i(b) = 0$.
Let $\gamma_{\varepsilon} = \gamma_0 + \varepsilon c \in \cP(a,b,p,q)$,
and let $\cA_{\varepsilon} = \cA_{\gamma_{\varepsilon}}$.
Suppose that $\gamma_0: [a,b] \to \cU$ satisfies the
Euler-Lagrange equations, i.e.,
$\frac {d\cA_{\varepsilon}}{d\varepsilon} (0) = 0$.
Then
\[
\begin{array}{llll}
        \displaystyle{\frac {d^2\cA_{\varepsilon}}{d\varepsilon^2} (0)} &
        = & \displaystyle{\int_a^b
        \sum_{i,j} \frac {\partial^2F}{\partial x_i\partial x_j}
        \left( \gamma_0, \frac {d\gamma_0}{dt} \right) \ c_i \ c_j \ dt}
        & \quad \mbox{(A)} \\
        & + & \displaystyle{2 \int_a^b \sum_{i,j} \frac {\partial^2F}{\partial
        x_i\partial v_j} \left( \gamma_0, \frac {d\gamma_0}{dt} \right)
        \ c_i \ \frac {dc_j}{dt} \ dt}
        & \quad \mbox{(B)} \\
        & + & \displaystyle{\int_a^b \sum_{i,j} \frac
        {\partial^2F}{\partial v_i\partial v_j}
        \left( \gamma_0, \frac {d\gamma_0}{dt} \right)
        \ \frac {dc_i}{dt} \ \frac {dc_j}{dt} \ dt}
        & \quad \mbox{(C)} \ .
\end{array}
\]
Since $\left( \frac {\partial^2F}{\partial v_i\partial v_j} (x,v)\right)
\gg 0$ at all $x,v$, we have
\[
        |\mbox{(A)}| \le
        \displaystyle{K_{_{\mathrm{A}}} |c|_{L^2[a,b]}^2} \ ,
        \quad
        |\mbox{(B)}| \le 
        \displaystyle{K_{_{\mathrm{B}}}
        |c|_{L^2[a,b]} \left| \frac {dc}{dt}
        \right|_{L^2[a,b]}}
        \; \mbox{ and } \;
        \mbox{(C)} \ge
        \displaystyle{K_{_{\mathrm{C}}}
        \left| \frac {dc}{dt} \right|^2_{L^2[a,b]}} \ .
\]
where $K_{_{\mathrm{A}}}, K_{_{\mathrm{B}}}, K_{_{\mathrm{C}}}$
are positive constants.
By the Wirtinger inequality\footnote{The
\textbf{Wirtinger inequality}\index{Wirtinger inequality}
states that, for $f \in C^1([a,b])$ with $f(a) = f(b) = 0$, we have
\[
        \int_a^b \left| \frac {df}{dt} \right|^2 dt
        \ge \frac {\pi^2}{(b-a)^2} \int_a^b |f|^2dt\ .
\]
This can be proved with Fourier series.}\index{Wirtinger inequality},
if $b-a$ is very small, then $\mbox{(C)} > |\mbox{(A)}| + |\mbox{(B)}|$
when $c \not\equiv 0$.
Hence, $\gamma_0$ is a local minimum.
\end{proof}

In Section~\ref{sec:symplectic_hamiltonian_fields} we saw that
solving Newton's second law\index{Newton ! second law}
in {\em configuration space}\index{configuration space}\index{space !
configuration} $\RR^3$ is
equivalent to solving in {\em phase space}\index{phase space}\index{space !
phase} for the integral curve\index{integral ! curve} in
$T^*\RR^3 = \RR^6$ of the hamiltonian vector field with
hamiltonian function $H$.
In the next section we will see how
this correspondence extends to more general Euler-Lagrange equations.

%%%%%%%%%%%%%%%%%%%%%%%%%%%%%%%%%%%%%%%%%%%%%%%%%%%%%%%%%%%%%%%%%%%%%%%%%%%%%
%%%%%%%%%%%%%%%%%%%%%%%%%%%%%%%%%%%%%%%%%%%%%%%%%%%%%%%%%%%%%%%%%%%%%%%%%%%%%

\ssubsection{Legendre Transform}
\label{sec:legendre}
\index{Legendre ! transform}

The Legendre transform gives the relation
between the variational (Euler-Lagrange)\index{equations !
Euler-Lagrange}\index{Euler ! Euler-Lagrange equations}
and the symplectic (Hamilton-Jacobi)\index{equations !
Hamilton-Jacobi}\index{Hamilton-Jacobi equations}\index{Jacobi !
Hamilton-Jacobi equations}
formulations of the equations of motion.

Let $V$ be an $n$-dimensional vector space,
with $e_1,\dots,e_n$ a basis of $V$ and $v_1,\dots,v_n$
the associated coordinates.
Let $F: V \to \RR$, $F = F(v_1,\dots,v_n)$, be a smooth function.
The function $F$ is \textbf{strictly convex}\index{function !
strictly convex}\index{strictly convex function} if and only if for
every pair of elements $p,v \in V$, $v \neq 0$, the restriction of
$F$ to the line $\{ p + xv \, | \, x \in \RR \}$ is strictly
convex.\footnote{A function $F:V \to \RR$ is \textbf{strictly convex} if
at every $p \in V$ the {\em hessian} $d^2 F_p$ is positive definite.
Let $u = \sum_{i=1}^n u_ie_i \in V$.
The \textbf{hessian}\index{hessian} of $F$ at $p$
is the quadratic function on $V$
\[
        (d^2F)_p(u) :=
        \sum_{i,j} \frac {\partial^2F}{\partial v_i\partial v_j} (p) u_iu_j
        = \left. \frac {d^2}{dt^2} F(p+tu) \right|_{t=0} \ .
\]}
It follows from the case of real functions on $\RR$ that,
for a strictly convex function $F$ on $V$, the following are
equivalent:
\footnote{A smooth function $f: \RR \to \RR$ is
\textbf{strictly convex}\index{function !
strictly convex}\index{strictly convex function}
if $f''(x) >0$ for all $x \in \RR$.
Assuming that $f$ is strictly convex,
the following four conditions are equivalent:
$f' (x) = 0$ at some point,
$f$ has a local minimum,
$f$ has a unique (global) minimum, and
$f(x) \to +\infty$ as $x \to \pm \infty$.
The function $f$ is \textbf{stable}\index{function !
stable}\index{stable ! function} if it satisfies one
(and hence all) of these conditions.
For instance, $e^x + ax$ is strictly convex for any $a \in \RR$,
but it is stable only for $a < 0$.
The function $x^2 + ax$ is strictly convex and stable for any $a \in \RR$.}
\begin{itemize}
\item[(a)]
$F$ has a critical point, i.e., a point where $dF_p = 0$;
\item[(b)]
$F$ has a local minimum at some point;
\item[(c)]
$F$ has a unique critical point (global minimum); and
\item[(d)]
$F$ is proper\index{proper function}, that is,
$F(p) \to +\infty$ as $p \to \infty$ in $V$.
\end{itemize}
A strictly convex function $F$ is \textbf{stable}\index{function !
stable}\index{stability ! definition}
when it satisfies conditions (a)-(d) above.

\begin{definition}\index{Legendre ! transform}
The \textbf{Legendre transform} associated to
$F \in C^{\infty}(V)$ is the map
\[
\begin{array}{rrcl}
        L_F: & V & \longrightarrow & V^* \\
        & p & \longmapsto & dF_p \in T_p^*V \simeq V^*\ ,
\end{array}
\]
where $T_p^*V \simeq V^*$ is the canonical
identification for a vector space $V$.
\end{definition}

From now on, assume that $F$ is a strictly convex function on $V$.
Then, for every point $p \in V$, $L_{_F}$ maps a neighborhood of $p$
diffeomorphically onto a neighborhood of $L_{_F}(p)$.
Given $\ell \in V^*$, let
\[
        F_\ell: V \longrightarrow \RR\ , \qquad F_\ell(v) = F(v) - \ell(v)\ .
\]
Since $(d^2F)_p = (d^2F_\ell)_p$,
$F$ is strictly convex if and only if $F_\ell$ is strictly convex.
The \textbf{stability set}\index{stability ! set}
of $F$ is
\[
        S_F = \{\ell \in V^* \mid F_\ell \mbox{ is stable}\}\ .
\]
The set $S_{_F}$ is open and convex,
and $L_{_F}$ maps $V$ diffeomorphically onto $S_{_F}$.
(A way to ensure that $S_{_F} = V^*$ and hence that
$L_{_F}$ maps $V$ diffeomorphically onto $V^*$,
is to assume that a strictly convex function $F$ has
\textbf{quadratic growth at infinity}\index{quadratic growth at infinity},
i.e., there exists a positive-definite quadratic form $Q$ on $V$
and a constant $K$ such that $F(p) \geq Q(p) - K$, for all $p$.)
The inverse to $L_F$ is the map $L_F^{-1}: S_F \to V$ described
as follows: for $\ell \in S_F$, the value $L_F^{-1}(\ell)$ is the
unique minimum point $p_\ell \in V$ of $F_\ell$.
Indeed $p$ is the minimum of $F(v) - dF_p(v)$.

\begin{definition}
The \textbf{dual function}\index{function !
dual}\index{dual function} $F^*$ to $F$ is
\[
        F^*: S_F \longrightarrow \RR \ ,
        \quad F^*(\ell) = -\min_{p \in V} F_\ell(p)\ .
\]
\end{definition}

The dual function $F^*$ is smooth
and, for all $p \in V$ and all $\ell \in S_{_F}$,
satisfies the \textbf{Young inequality}\index{Young inequality}
$F(p) + F^*(\ell) \geq \ell(p)$.

On one hand we have $V \times V^* \simeq T^*V$, and on the other hand,
since $V = V^{**}$, we have $V \times V^* \simeq V^* \times V \simeq T^*V^*$.
Let $\alpha_1$ be the tautological 1-form on $T^*V$ and
$\alpha_2$ be the tautological 1-form on $T^*V^*$.
Via the identifications above,
we can think of both of these forms as living on $V \times V^*$.
Since $\alpha_1 = d\beta - \alpha_2$, where
$\beta: V \times V^* \to \RR$ is the function $\beta (p,\ell) = \ell(p)$,
we conclude that the forms $\omega_1 = - d \alpha_1$
and $\omega_2 = - d \alpha_2$ satisfy $\omega_1 = - \omega_2$.

\begin{theorem}
For a strictly convex function $F$ we have that $L_F^{-1} = L_{F^*}$.
\end{theorem}

\vspace*{-2ex}

\begin{proof}
The graph $\Lambda_{_F}$ of the Legendre transform $L_{_F}$
is a lagrangian submanifold of
$V \times V^*$ with respect to the symplectic form $\omega_1$.
Hence, $\Lambda_{_F}$ is also lagrangian for $\omega_2$.
Let ${\mathrm{pr}}_1 : \Lambda_{_F} \to V$ and
${\mathrm{pr}}_2 : \Lambda_{_F} \to V^*$ be the restrictions of the
projection maps $V \times V^* \to V$ and $V \times V^* \to V^*$,
and let $i : \Lambda_{_F} \hookrightarrow V \times V^*$ be the inclusion map.
Then $i^* \alpha_1 = d ({\mathrm{pr}}_1)^* F$
as both sides have value $dF_p$ at $(p,dF_p) \in \Lambda_{_F}$.
It follows that $i^* \alpha_2 = d (i^* \beta - ({\mathrm{pr}}_1)^* F)
= d ({\mathrm{pr}}_2)^* F^*$,
which shows that $\Lambda_{_F}$ is the graph of the inverse of $L_{F^*}$.
From this we conclude that the inverse of the Legendre transform
associated with $F$ is the Legendre transform\index{Legendre !
transform} associated with $F^*$.
\end{proof}

%%%%%%%%%%%%%%%%%%%%%%%%%%%%%%%%%%%%%%%%%%%%%%%%%%%%%%%%%%%%%%%%%%%%%%%%%%%%%

Let $M$ be a manifold and $F: TM \to \RR$.
We return to the Euler-Lagrange equations for minimizing the action
$\cA_{\gamma} = \int {\widetilde \gamma}^*F$.\index{variational ! problem}
At $p \in M$, let $F_p := F|_{T_pM}: T_p M \to \RR$.
Assume that $F_p$ is strictly convex for all $p \in M$.  To simplify
notation, assume also that $S_{F_p} = T_p^*M$.
The Legendre transform on each tangent space
$L_{F_p}: T_pM \stackrel{\simeq}{\longrightarrow} T_p^*M$
is essentially given by the first derivatives of $F$ in the $v$ directions.
Collect these and the dual functions $F_p^*: T_p^*M \to \RR$ into maps
\[
        \cL: TM \longrightarrow T^*M \ , \ \cL|_{T_pM} = L_{F_p}
        \quad \mbox{ and } \quad
        H: T ^*M \longrightarrow \RR \ , \ H|_{T_p^*M} = F_p^*\ .
\]
The maps $H$ and $\cL$ are smooth, and $\cL$ is a diffeomorphism.

\begin{theorem}
\index{theorem ! Euler-Lagrange equations}\index{equations !
Euler-Lagrange}\index{Euler ! Euler-Lagrange equations}
Let $\gamma: [a,b] \to M$ be a curve, and
${\widetilde \gamma}: [a,b] \to TM$ its lift.
Then $\gamma$ satisfies the Euler-Lagrange equations
on every coordinate chart
if and only if $\cL \circ {\widetilde \gamma}: [a,b] \to T^*M$
is an integral curve of the hamiltonian vector field $X_H$.
\end{theorem}

\vspace*{-2ex}

\begin{proof}
Let $(\cU,x_1,\dots,x_n)$ be a coordinate chart in $M$,
with associated tangent $(T\cU,x_1,\dots,x_n,v_1,\dots,v_n)$
and cotangent $(T^*\cU,x_1,\dots,x_n,\xi_1,\dots,\xi_n)$ coordinates.
On $T\cU$ we have $F = F(x,v)$, on $T^*\cU$ we have $H = H(x,\xi)$, and
\[
\begin{array}{rrclcrrcl}
        \cL: & T\cU & \longrightarrow & T^*\cU & \qquad &
        H: & T^*\cU & \longrightarrow & \RR \\
        & (x,v) & \longmapsto & (x,\xi) & &
        & (x,\xi) & \longmapsto & F_x^*(\xi) = \xi \cdot v - F(x,v)
\end{array}
\]
where $\xi := L_{F_x}(v) = \frac {\partial F}{\partial v} (x,v)$
is called the \textbf{momentum}\index{momentum}.
Integral curves $(x(t),\xi(t))$ of $X_H$ satisfy the
Hamilton equations\index{Hamilton equations}\index{equations ! Hamilton}:
\[
   \mbox{(H)} \qquad \qquad \qquad
   \left\{ \begin{array}{rll}
   \frac {dx}{dt} & = & \phantom{-} \frac {\partial H}{\partial \xi} (x,\xi) \\
   \frac {d\xi}{dt} & = & - \frac {\partial H}{\partial x} (x,\xi)
\end{array} \right.
\]
whereas the physical path $x(t)$ satisfies the Euler-Lagrange
equations:
\[
        \mbox{(E-L)} \qquad \qquad
        \frac {\partial F}{\partial x} \left( x,\frac {dx}{dt} \right) =
        \frac {d}{dt} \frac {\partial F}{\partial v}
        \left( x,\frac {dx}{dt} \right)\ .
\]
Let $(x(t),\xi(t)) = \cL\left(x(t),\frac {dx}{dt} (t)\right)$.
For an arbitrary curve $x(t)$, we want to prove that $t \mapsto (x(t),\xi(t))$
satisfies~(H) if and only if $t \mapsto \left(x(t), \frac {dx}{dt} (t)\right)$
satisfies~(E-L).
The first line of~(H) comes automatically from the definition of $\xi$:
\[
        \xi = L_{F_x}\left( \frac {dx}{dt} \right)
        \qquad \iff \qquad
        \frac {dx}{dt} = L_{F_x}^{-1}(\xi) =
        L_{F_x^*}(\xi) = \frac {\partial H}{\partial \xi} (x,\xi) \ .
\]
If $(x,\xi) = \cL(x,v)$, by differentiating both sides of
$H(x,\xi) = \xi \cdot v - F(x,v)$ with respect to $x$,
where $\xi = L_{F_x}(v) = \xi(x,v)$
and $v = \frac {\partial H}{\partial \xi}$, we obtain
\[
        \frac {\partial H}{\partial x} +
        \frac {\partial H}{\partial \xi}
        \frac {\partial \xi}{\partial x}
        = \frac {\partial \xi}{\partial x}
        \cdot v - \frac {\partial F}{\partial x}
        \qquad \iff \qquad
        \frac {\partial F}{\partial x} (x,v) =
        -\frac {\partial H}{\partial x} (x,\xi) \ .
\]
Using the last equation and the definition of $\xi$,
the second line of~(H) becomes~(E-L):
\[
        \frac {d\xi}{dt}
        = -\frac {\partial H}{\partial x} (x,\xi)
        \qquad \iff \qquad 
        \frac {d}{dt} \frac {\partial F}{\partial v} (x,v)
        = \frac {\partial F}{\partial x} (x,v) \ .
\]
\end{proof}

%%%%%%%%%%%%%%%%%%%%%%%%%%%%%%%%%%%%%%%%%%%%%%%%%%%%%%%%%%%%%%%%%%%%%%%%%%%%
%%%%%%%%%%%%%%%%%%%%%%%%%%%%%%%%%%%%%%%%%%%%%%%%%%%%%%%%%%%%%%%%%%%%%%%%%%%%

\ssubsection{Integrable Systems}
\label{sec:integrable}
\index{integrable ! system}

\begin{definition}
A \textbf{hamiltonian system}\index{hamiltonian ! system}
is a triple $(M,\omega,H)$,
where $(M,\omega)$ is a symplectic manifold and $H \in C^\infty (M)$
is the \textbf{hamiltonian function}.\index{hamiltonian ! function}
\end{definition}

\vspace*{-2ex}

\begin{proposition}\label{prop:integral}
For a function $f$ on a symplectic manifold $(M,\omega)$
we have that $\{f,H\}=0$ if and only if
$f$ is constant along integral curves of $X_{_H}$.
\end{proposition}

\vspace*{-2ex}

\begin{proof}
Let $\rho_{t}$ be the flow of $X_{_H}$.  Then
\[
   \frac{d}{dt}(f\circ \rho_{t})
   = \rho_{t}^{*}\cL_{X_{H}}f
   = \rho_{t}^{*}\imath_{X_{H}}df
   = \rho_{t}^{*}\imath_{X_{H}}\imath_{X_{f}}\omega \\
   = \rho_{t}^{*}\omega(X_{f},X_{_H})
   = \rho_{t}^{*} \{f,H\} \ .
\]
\end{proof}

A function $f$ as in Proposition~\ref{prop:integral} is called an
\textbf{integral of motion}\index{integral ! of motion}\index{motion !
integral of motion} (or a \textbf{first integral}\index{first
integral}\index{integral ! first}
or a \textbf{constant of motion}).\index{motion ! constant of motion}
In general, hamiltonian systems do not admit
integrals of motion that are {\em independent}
of the hamiltonian function.
Functions $f_1, \ldots, f_n$ are said to be \textbf{independent}
if their differentials $(df_1)_p, \ldots, (df_n)_p$
are linearly independent at all points $p$ in some dense subset of $M$.
Loosely speaking, a hamiltonian system is
{\em (completely) integrable}\index{integrable ! system} if it has
as many {\em commuting} integrals of motion as possible.
\textbf{Commutativity} is with respect to the Poisson
bracket.\index{Poisson ! bracket}\index{bracket ! Poisson}
If $f_1, \ldots, f_n$ are commuting integrals of motion
for a hamiltonian system $(M,\omega,H)$, then
$\omega (X_{f_i}, X_{f_j}) = \{ f_i, f_j \} = 0$,
so at each $p \in M$ the hamiltonian vector fields
generate an isotropic subspace of $T_pM$.
When $f_1, \ldots, f_n$ are independent,
by symplectic linear algebra
$n$ can be at most half the dimension of $M$.

\begin{definition}
A hamiltonian system $(M,\omega,H)$ where $M$ is a
$2n$-dimensional manifold
is \textbf{(completely) integrable}\index{completely
integrable system}\index{integrable ! system} if it possesses
$n$ independent commuting integrals of motion,
$f_1=H, f_2,\ldots, f_n$.% pairwise in involution with respect to
%the Poisson bracket, i.e., $\{ f_i, f_j \} = 0$, for all $i,j$.
\end{definition}

Any 2-dimensional hamiltonian system (where the set
of non-fixed points is dense) is trivially integrable.
Basic examples are the simple pendulum and the harmonic oscillator.
A hamiltonian system $(M,\omega,H)$ where $M$ is 4-dimensional
is integrable if there is an integral of motion independent of $H$
(the commutativity condition is automatically satisfied).
A basic example is the spherical pendulum.
Sophisticated examples of integrable systems
can be found in~\cite{au:tops,hi-se-wa:integrable}.

\begin{examples}
\begin{enumerate}
\item
The \textbf{simple pendulum}\index{example !
simple pendulum}\index{pendulum ! simple}\index{simple
pendulum} is a mechanical system
consisting of a massless rigid rod
of length $\ell$, fixed at one end, whereas the other end has
a bob of mass $m$, which may oscillate in the vertical plane.
We assume that the force of gravity\index{gravity} is constant pointing
vertically downwards and the only external force acting on this system.
Let $\theta$ be the oriented angle between the rod
and the vertical direction.
Let $\xi$ be the coordinate along the fibers of
$T^* S^1$ induced by the standard angle coordinate on $S^1$.
The energy function $H: T^* S^1 \to \RR$,
$H(\theta, \xi) = \frac{\xi^2}{2m\ell^2} + m\ell (1-\cos \theta)$,
is an appropriate hamiltonian
function to describe the simple pendulum.
Gravity is responsible for the potential
energy\index{energy ! potential}\index{potential ! energy}
$V(\theta) = m\ell (1-\cos \theta)$,
and the kinetic energy\index{energy !
kinetic}\index{kinetic energy} is given by
$K(\theta,\xi) = \frac{1}{2m\ell^2} \xi^2$.

\item
The \textbf{spherical pendulum}\index{example !
spherical pendulum}\index{pendulum !
spherical}\index{spherical pendulum}
consists of a massless rigid rod
of length $\ell$, fixed at one end, whereas the other end has
a bob of mass $m$, which may oscillate {\em freely in all directions}.
For simplicity let $m=\ell=1$.
Again assume that gravity\index{gravity} is
the only external force.
Let $\varphi, \theta$ ($0 < \varphi < \pi$, $0 < \theta < 2\pi$)
be spherical coordinates for the bob,
inducing coordinates $\eta, \xi$ along the fibers of $T^* S^2$.
An appropriate hamiltonian function for this system
is the energy function $H: T^* S^2 \to \RR$, $H(\varphi, \theta, \eta, \xi) =
\frac{1}{2} \left( \eta^2 + \frac{\xi^2}{(\sin \varphi)^2} \right)
+ \cos \varphi$.
The function $J(\varphi, \theta, \eta, \xi) = \xi$
is an independent integral of motion
corresponding to the group of symmetries
given by rotations about the vertical axis (Section~\ref{sec:actions}).
The points $p \in T^*S^2$ where $dH_p$ and $dJ_p$
are linearly dependent are:
\begin{itemize}
\item
the two critical points of $H$
(where both $dH$ and $dJ$ vanish);
\item
if $x \in S^2$ is in the southern hemisphere
($x_3 < 0$), then there exist exactly two points,
$p_+ = (x,\eta,\xi)$ and $p_- = (x,-\eta,-\xi)$,
in the cotangent fiber above $x$ where $dH_p$ and $dJ_p$
are linearly dependent;
\item
since $dH_p$ and $dJ_p$ are linearly dependent
along the trajectory of the hamiltonian vector field of $H$
through $p_+$,
this trajectory is also a trajectory of
the hamiltonian vector field of $J$
and hence its projection onto $S^2$ is a latitudinal
(or horizontal) circle.
The projection of the trajectory through $p_-$
is the same latitudinal circle traced in the opposite direction.
\end{itemize}
\end{enumerate}
\end{examples}

Let $(M,\omega,H)$ be an integrable system of dimension $2n$
with integrals of motion $f_1=H, f_2,\ldots, f_n$.
Let $c \in \RR^n$ be a regular value of $f:= (f_1, \ldots, f_n)$.
The corresponding level set $f^{-1} (c)$ is a lagrangian submanifold,
as it is $n$-dimensional and its tangent bundle is isotropic.
If the flows are complete on $f^{-1} (c)$,
by following them we obtain global coordinates.
Any compact component of $f^{-1} (c)$ must hence be a torus.
These components, when they exist, are called
\textbf{Liouville tori}\index{Liouville ! torus}.
A way to ensure that compact components exist
is to have one of the $f_i$'s proper.

\begin{theorem}
\label{thm:arnold}\index{theorem ! Arnold-Liouville}\index{Liouville !
Arnold-Liouville theorem}\index{Arnold ! Arnold-Liouville theorem}
\textbf{(Arnold-Liouville~\cite{ar:mathematical})} $\;$
Let $(M,\omega,H)$ be an integrable system of dimension $2n$
with integrals of motion $f_1=H, f_2,\ldots, f_n$.
Let $c \in \RR^n$ be a regular value of $f:= (f_1, \ldots, f_n)$.
The level $f^{-1} (c)$ is a lagrangian
submanifold of $M$.
\begin{itemize}
\item[(a)]
If the flows of the hamiltonian vector fields
$X_{f_1}, \ldots, X_{f_n}$ starting at
a point $p \in f^{-1} (c)$ are complete,
then the connected component of $f^{-1} (c)$ containing $p$
is a homogeneous space for $\RR^n$,
i.e., is of the form $\RR^{n-k} \times \TT^k$
for some $k$, $0 \leq k \leq n$, where $\TT^k$ is
a $k$-dimensional torus..
With respect to this affine structure, that component has
coordinates $\varphi_1, \ldots, \varphi_n$, known as
\textbf{angle coordinates}\index{angle coordinates},
in which the flows of $X_{f_1}, \ldots, X_{f_n}$ are linear.
\item[(b)]
There are coordinates $\psi_1, \ldots, \psi_n$, known as
\textbf{action coordinates}\index{action ! coordinates},
complementary to the angle coordinates, such that
the $\psi_i$'s are integrals of motion and
$\varphi_1, \ldots, \varphi_n,\psi_1, \ldots, \psi_n$
form a Darboux chart.
\end{itemize}
\end{theorem}

Therefore, the dynamics of an integrable system
has a simple explicit solution
in action-angle coordinates\index{action-angle coordinates}.
The proof of part~(a) -- the easy part of the theorem -- is sketched above.
For the proof of part~(b), see for instance~\cite{ar:mathematical,du:global}.
Geometrically, regular levels being lagrangian submanifolds
implies that, in a neighborhood of a regular value,
the map $f: M \to \RR^n$
collecting the given integrals of motion is a
\textbf{lagrangian fibration}\index{lagrangian fibration}, i.e.,
it is locally trivial and its fibers are lagrangian submanifolds.
Part~(a) states that there are coordinates along the fibers,
the angle coordinates,\footnote{The name {\em angle coordinates}
is used even if the fibers are not tori.}
in which the flows of $X_{f_1}, \ldots, X_{f_n}$ are linear.
Part (b) guarantees the existence of coordinates
on $\RR^n$, the action coordinates, $\psi_1, \ldots, \psi_n$,
complementary to the angle coordinates,
that (Poisson) commute among themselves and satisfy
$\{ \varphi _i , \psi _j \} = \delta_{ij}$.
%with respect to the angle coordinates.
The action coordinates are generally not
the given integrals of motion because
$\varphi_1, \ldots, \varphi_n,f_1, \ldots, f_n$
do not form a Darboux chart.

%%%%%%%%%%%%%%%%%%%%%%%%%%%%%%%%%%%%%%%%%%%%%%%%%%%%%%%%%%%%%%%%%%%%%%%%%%%%%
%%%%%%%%%%%%%%%%%%%%%%%%%%%%%%%%%%%%%%%%%%%%%%%%%%%%%%%%%%%%%%%%%%%%%%%%%%%%%

\ssubsection{Symplectic and Hamiltonian Actions}
\index{action ! symplectic}\index{action ! hamiltonian}\index{hamiltonian !
action}\index{symplectic ! action}
\label{sec:actions}

Let $(M,\omega)$ be a symplectic manifold, and $G$ a Lie group.

\begin{definition}\index{action ! symplectic}\index{symplectic ! action}
An action\footnote{A (smooth) \textbf{action} of $G$ on $M$ is a group
homomorphism $G \to \mathrm{Diff}(M)$, $g \mapsto \psi_g$,
whose evaluation map $M \times G \to M$, $(p,g) \mapsto \psi_g (p)$,
is smooth.} $\psi: G \to \mathrm{Diff}(M)$, $g \mapsto \psi_g$,
is a \textbf{symplectic action} if each $\psi_g$ is a symplectomorphism,
i.e., $\psi: G \to \mathrm{Sympl}(M,\omega) \subset \mathrm{Diff}(M)$.
\end{definition}

In particular, symplectic actions of $\RR$ on $(M,\omega)$\index{symplectic !
action}\index{action ! symplectic} are in one-to-one correspondence
with complete symplectic vector fields on $M$:\index{complete vector
field}\index{vector field ! symplectic}\index{symplectic ! vector field}
\[
        \psi = \exp tX \quad \longleftrightarrow \quad
        X_p = \left. \frac {d\psi_t(p)}{dt} \right|_{t=0}\ , \ p \in M \ .
\]
We may define a symplectic action $\psi$ of $S^1$ or $\RR$
on $(M,\omega)$ to be \textbf{hamiltonian}\index{action !
hamiltonian}\index{hamiltonian ! action}\index{action ! hamiltonian}
if the vector field $X$ generated by $\psi$ is hamiltonian,
that is, when there is $H:M \to \RR$ with
$d H = \imath_X \omega$.
An action of $S^1$ may be viewed as a periodic action of $\RR$.

\begin{examples}
\begin{enumerate}
\item
On $(\RR^{2n}, \omega_0)$, the orbits of the action
generated by $X = -\frac {\partial}{\partial y_1}$ are lines
parallel to the $y_1$-axis,
$\{(x_1,y_1-t,x_2,y_2,\dots,x_n,y_n) \mid t \in \RR\}$.
Since $X$ is hamiltonian with hamiltonian function $x_1$,
this is a
{\em hamiltonian action}\index{action ! hamiltonian}\index{hamiltonian !
action}\index{example ! of hamiltonian actions} of $\RR$.

\item
On the 2-sphere $(S^{2},d\theta \wedge dh)$
in cylindrical coordinates, the one-parameter group of
diffeomorphisms given by rotation around the vertical axis,
$\psi_{t}(\theta,h)=(\theta+t,h)$ ($t \in \RR$)
is a symplectic action of the group
$S^1 \simeq \RR / \langle 2\pi \rangle$,
as it preserves the area form $d\theta \wedge dh$.
Since the vector field corresponding to $\psi$
is hamiltonian with hamiltonian function $h$,
this is a {\em hamiltonian action}\index{action !
hamiltonian}\index{hamiltonian ! action}\index{example !
of hamiltonian actions} of $S^1$.
\end{enumerate}
\end{examples}

When $G$ is a product of $S^1$'s or $\RR$'s,
an action $\psi: G \to \mathrm{Sympl}(M,\omega)$
is called {\em hamiltonian}\index{hamiltonian !
action}\index{action ! hamiltonian} when the restriction
to each 1-dimensional factor
is hamiltonian in the previous sense
{\em with hamiltonian function preserved by
the action of the rest of $G$}.

For an arbitrary Lie group $G$,
we use an upgraded hamiltonian function $\mu$,
known as a {\em moment map}\index{moment map !
upgraded hamiltonian function},
determined up to an additive local constant
by coordinate functions $\mu_i$ indexed by a basis
%satisfying $d \mu_i = \imath_{X_i^\#} \omega$ for a basis $X_i$
of the Lie algebra of $G$.
We require that the constant be such that $\mu$ is {\em equivariant},
i.e., $\mu$ intertwines the action of $G$ on $M$ and the
coadjoint action of $G$ on the dual of its Lie algebra.
(If $M$ is compact, equivariance can be achieved
by adjusting the constant so that $\int_M \mu \omega^n = 0$.
Similarly when there is a fixed point $p$
(on each component of $M$) by imposing $\mu (p) = 0$.)

%%%%%%%%%%%%%%%%%%%%%%%%%%%%%%%%%%%%%%%%%%%%%%%%%%%%%%%%%%%%%%%%%%%%%%%%%%%%%

Let $G$ be a Lie group, $\fg$ the Lie algebra of $G$,
and $\fg^*$ the dual vector space of $\fg$.

\begin{definition}
An action $\psi: G \to \mathrm{Diff}(M)$ on a
symplectic manifold $(M,\omega)$
is a \textbf{hamiltonian action}\index{hamiltonian !
action}\index{action ! hamiltonian} if there exists a map
$\mu: M \to \fg^*$ satisfying:
\begin{itemize}
\item
For each $X \in \fg$, we have $d\mu^X = \imath_{X^{\#}}\omega$,
i.e., $\mu^X$ is a hamiltonian function\index{hamiltonian !
function}\index{function ! hamiltonian}\index{hamiltonian !
moment map} for the vector field $X^{\#}$, where
\begin{itemize}
\item
$\mu^X: M \to \RR$, $\mu^X(p) := \langle
\mu(p),X\rangle$, is the component of $\mu$ along $X$,
\item
$X^{\#}$ is the vector field on $M$ generated by the one-parameter
subgroup $\{\exp tX \mid t \in \RR\} \subseteq G$.
\end{itemize}
\item
The map $\mu$ is {\em equivariant}\index{equivariant !
moment map}\index{moment map ! equivariance}
with respect to the given action
$\psi$ on $M$ and the
coadjoint action:
$\mu \circ \psi_g = \mathrm{Ad}_g^* \circ \mu$, for all $g \in G$.
\end{itemize}
Then $(M,\omega,G,\mu)$ is a
\textbf{hamiltonian $G$-space}\index{hamiltonian !
G-space@$G$-space}\index{G-space@$G$-space} and $\mu$ is
a \textbf{moment map}\index{moment map ! definition}\index{moment map !
hamiltonian G-space@hamiltonian $G$-space}.
\end{definition}

This definition matches the previous one when
$G$ is an abelian group $\RR$, $S^1$ or $\TT^n$,
for which equivariance becomes
invariance since the coadjoint action is trivial.

%%%%%%%%%%%%%%%%%%%%%%%%%%%%%%%%%%%%%%%%%%%%%%%%%%%%%%%%%%%%%%%%%%%%%%%%%%%%%

\begin{examples}\index{moment map ! example}
\begin{enumerate}
\item
Let $\TT ^n = \{ (t_1, \ldots, t_n) \in \CC^n\, : \,
|t_j| =1, \mbox{ for all }j \, \}$ be a torus acting on $\CC^n$ by
$(t_1, \ldots, t_n) \cdot (z_1, \ldots, z_n) =
(t_1^{k_1} z_1, \ldots, t_n^{k_n} z_n)$,
where $k_1, \ldots, k_n \in \ZZ$ are fixed.
This action is hamiltonian with a moment map
$\mu : \CC^n \to (\ft ^n)^* \simeq \RR^n$,
$\mu (z_1, \ldots, z_n) = - \textstyle{\frac12}
(k_1 |z_1|^2, \ldots, k_n |z_n|^2)$.

\item
When a Lie group $G$ acts
on two symplectic manifolds $(M_j, \omega_j)$, $j = 1,2$,
with moment maps $\mu_j : M_j \to \fg^*$,
the diagonal action of $G$ on $M_1 \times M_2$
% is hamiltonian
has moment map $\mu : M_1 \times M_2 \to \fg^*$,
$\mu (p_1,p_2) = \mu_1 (p_1) + \mu_2 (p_2)$.

\item
Equip the coadjoint orbits\index{coadjoint ! orbit} of a Lie group $G$
with the canonical symplectic form\index{canonical !
symplectic form on a coadjoint orbit}\index{symplectic !
canonical symplectic form on a coadjoint orbit}
(Section~\ref{sec:symplectic_hamiltonian_fields}).
Then, for each $\xi \in \fg^*$, the coadjoint
action\index{coadjoint ! action}
on the orbit $G \cdot \xi$ is hamiltonian with
moment map simply the inclusion map
$\mu : G \cdot \xi \hookrightarrow \fg^*$.

\item
Identify the Lie algebra of the unitary group $\UU (n)$ with its dual via
the inner product $\langle A, B \rangle = \mathrm{trace} (A ^* B)$.
The natural action of $\UU (n)$ on $(\CC^n, \omega_0)$
is hamiltonian\index{moment map ! example} with moment map
$\mu : \CC^n \to \fu (n)$ given by $\mu (z) = \textstyle{i \over {2}} z z^*$.
Similarly, a moment map for the natural action of $\UU (k)$ on the space
$(\CC^{k\times n}, \omega_0)$ of complex $(k\times n)$-matrices is given by
$\mu (A) = \textstyle{{i} \over {2}} A A^*$ for $A \in \CC^{k\times n}$.
Thus the $\UU (n)$-action by conjugation on the space
$(\CC^{n^2},\omega_0)$ of complex $(n\times n)$-matrices
is hamiltonian, with moment map given by
$\mu (A) = \textstyle{i \over 2} [A, A^*]$.

\item
For the spherical pendulum (Section~\ref{sec:integrable}),
the {\em energy-momentum map}\index{energy !
energy-momentum map} $(H,J) : T^* S^2 \to \RR^2$
is a moment map for the $\RR \times S^1$ action
given by time flow and rotation about the vertical axis.

\item
Suppose that a compact Lie group acts on a symplectic
manifold $(M,\omega)$ in a hamiltonian way, and that $q \in M$
is a fixed point for the $G$-action.
Then, by an equivariant version of Darboux's
theorem,\footnote{\textbf{Equivariant Darboux
Theorem~\cite{we:lagrangian}:}\index{theorem ! equivariant
Darboux}\index{Darboux ! equivariant Darboux
theorem}\index{equivariant ! Darboux theorem}
{\em Let $(M, \omega)$ be a $2n$-dimensional symplectic
manifold equipped with a symplectic action of a
compact Lie group $G$, and let $q$ be a fixed point.
Then there exists a $G$-invariant chart
$(\cU,x_1,\dots,x_n,y_1,\dots,y_n)$ centered at $q$
and $G$-equivariant with respect to a linear action
of $G$ on $\RR^{2n}$ such that
\[
   \left. \omega \right|_\cU =
   \sum\limits_{k=1}^n dx_k \wedge dy_k \ .
\]}A suitable linear action on $\RR^{2n}$ is
equivalent to the induced action of $G$ on $T_qM$.
The proof relies on an equivariant version of the
Moser trick and may be found in~\cite{gu-st:techniques}.}
there exists a Darboux chart $(\cU , z_1 , \ldots , z_n)$
centered at $q$ that is $G$-equivariant with respect to
a linear action of $G$ on $\CC^n$.
Consider an $\varepsilon$-blow-up of $M$ relative to this chart,
for $\varepsilon$ sufficiently small.
Then $G$ acts on the blow-up in a hamiltonian way.
\end{enumerate}
\end{examples}

%%%%%%%%%%%%%%%%%%%%%%%%%%%%%%%%%%%%%%%%%%%%%%%%%%%%%%%%%%%%%%%%%%%%%%%%%%%%%

The concept of a moment map\index{moment map ! origin}
was introduced by Souriau~\cite{so:structure}
under the french name {\em application moment};
besides the more standard english translation
to {\em moment map}, the alternative {\em momentum map} is also used,
and recently James Stasheff has proposed the short
unifying new word \textbf{momap}.
The name comes from being the generalization of
{\em linear and angular momenta}\index{momentum}\index{angular
momentum}\index{linear momentum} in classical mechanics.

Let $\RR^3$ act on
$(\RR^6 \simeq T^*\RR^3, \omega_0 = \sum dx_i \wedge dy_i)$ by
\textbf{translations}:
\[
        a \in \RR^3 \; \longmapsto \;
        \psi_a \in \mathrm{Sympl}(\RR^6,\omega_0)\ , \; 
        \psi_a (x,y) = (x + a,y)\ .
\]
The vector field generated by $X = a = (a_1,a_2,a_3)$ is
$X^{\#} = a_1 \frac {\partial}{\partial x_1} + a_2 \frac
{\partial}{\partial x_2} + a_3 \frac {\partial}{\partial x_3}$,
and the \textbf{linear momentum}\index{linear momentum} map
\[
        \mu: \RR^6 \longrightarrow \RR^3\ , \quad
        \mu(x,y) = y
\]
is a moment map, with
$\mu^{a}(x,y) = \langle \mu(x,y),a \rangle = y \cdot a$.
Classically, $y$ is
called the \textbf{momentum vector}\index{momentum vector}
corresponding to the \textbf{position vector}\index{momentum vector} $x$.

The $\SO (3)$-action on $\RR^3$ by \textbf{rotations} lifts
to a symplectic action $\psi$ on the cotangent bundle $\RR^6$.
The infinitesimal version of this action is\footnote{The Lie group
$\SO (3) = \{A \in \GL (3;\RR) \mid A^t A = \mathrm{Id} \mbox{ and }
\mathrm{det} A = 1\}$,\index{example ! coadjoint orbits}
has Lie algebra, $\fg = \{A \in \mathfrak{gl}(3;\RR) \mid A + A^t = 0\}$,
the space of $3 \times 3$ skew-symmetric matrices.
The standard identification of $\fg$ with $\RR^3$
carries the Lie bracket to the exterior product:
\[
\begin{array}{rcl}
        A = \left[ \begin{array}{ccc}
        0 & -a_3 & a_2 \\
        a_3 & 0 & -a_1 \\
        -a_2 & a_1 & 0
        \end{array} \right] & \longmapsto &
        a = (a_1,a_2,a_3) \\
        \\
        {[ A,B ]} = AB - BA & \longmapsto & a \times b \ .
\end{array}
\]}
\[
        a \in \RR^3 \; \longmapsto \;
        d \psi (a) \in \chi^{\mathrm{sympl}}(\RR^6) \ , \;
        d\psi (a) (x,y) = (a \times x, a \times y)\ .
\]
Then the \textbf{angular momentum}\index{angular
momentum} map
\[
        \mu: \RR^6 \longrightarrow \RR^3 \ , \quad
        \mu(x,y) = x \times y
\]
is a moment map, with
$\mu^{a}(x,y) = \langle \mu(x,y),a \rangle = (x \times y) \cdot a$.

%%%%%%%%%%%%%%%%%%%%%%%%%%%%%%%%%%%%%%%%%%%%%%%%%%%%%%%%%%%%%%%%%%%%%%%%%%%%%

The notion of a moment map associated to a group action on a
symplectic manifold formalizes the \textbf{Noether principle}\index{Noether !
principle}\index{principle ! Noether},
which asserts that there is a one-to-one
correspondence between {\em symmetries} (or one-parameter group actions)
and {\em integrals of motion}
(or conserved quantities) for a mechanical system.

\begin{definition}
An \textbf{integral of motion}\index{motion !
integral of motion}\index{integral ! of motion} of a
hamiltonian $G$-space $(M,\omega,G,\mu)$ is
a $G$-invariant function $f: M \to \RR$.
When $\mu$ is constant on the trajectories of a hamiltonian
vector field $X_f$,
the corresponding flow $\{\exp tX_f \mid t \in \RR\}$
(regarded as an $\RR$-action) is a \textbf{symmetry}
of the hamiltonian $G$-space $(M,\omega,G,\mu)$.
\end{definition}

\vspace*{-2ex}

\begin{theorem}\index{theorem ! Noether}\index{Noether !
theorem} \textbf{(Noether)} $\;$
Let $(M,\omega,G,\mu)$ be a hamiltonian $G$-space
where $G$ is connected.
If $f$ is an integral of motion,
the flow of its hamiltonian vector field $X_f$ is a symmetry.
If the flow of some hamiltonian vector field $X_f$ is a symmetry,
then a corresponding hamiltonian function $f$ is an integral of motion.
\end{theorem}

\vspace*{-2ex}

\begin{proof}
Let $\mu^X = \langle \mu,X\rangle: M \to \RR$ for $X \in \fg$.
We have
$\cL_{X_f} \mu^X = \imath_{X_f} d\mu^X = \imath_{X_f}\imath_{X^{\#}}\omega
= -\imath_{X^{\#}} \imath_{X_f}\omega = -\imath_{X^{\#}} df = -\cL_{X^{\#}} f$.
So $\mu$ is invariant over the flow of $X_f$
if and only if $f$ is invariant under the infinitesimal $G$-action.
\end{proof}

%%%%%%%%%%%%%%%%%%%%%%%%%%%%%%%%%%%%%%%%%%%%%%%%%%%%%%%%%%%%%%%%%%%%%%%%%%%%%
%%%%%%%%%%%%%%%%%%%%%%%%%%%%%%%%%%%%%%%%%%%%%%%%%%%%%%%%%%%%%%%%%%%%%%%%%%%%%

%\ssubsection{Existence and Uniqueness of Moment Maps}
%\index{moment map ! properties}\index{moment map ! existence}

We now turn to the questions of existence and uniqueness of moment maps.

Let $\fg$ be a Lie algebra, and let
$C^k := \Lambda^k\fg^*$ be the set of \textbf{$k$-cochains} on $\fg$,
that is, of alternating $k$-linear maps $\fg \times \dots \times \fg \to \RR$.
The linear operator $\delta: C^k \to C^{k+1}$ defined by
$\delta c(X_0,\dots,X_k) = \sum_{i < j} (-1)^{i+j} c([X_i,X_j],X_0,\dots,
{\widehat X}_i,\dots,{\widehat X}_j,\dots,X_k)$ satisfies $\delta^2 = 0$.
The \textbf{Lie algebra cohomology groups}\index{cohomology !
Lie algebra}\index{Lie ! algebra cohomology}
(or \textbf{Chevalley cohomology groups}\index{cohomology !
Chevalley}\index{Chevalley cohomology})
of $\fg$ are the cohomology groups of the complex
$0 \stackrel{\delta}{\longrightarrow} C^0 \stackrel{\delta}{\longrightarrow}
C^1 \stackrel{\delta}{\longrightarrow} \dots$:
\[
        H^k(\fg;\RR) :=
        \frac {\ker \delta: C^k \to C^{k+1}}{\mathrm{im}\
        \delta: C^{k-1} \to C^k}\ .
\]
It is always $H^0(\fg;\RR) = \RR$.
If $c \in C^1 = \fg^*$, then $\delta c(X,Y) = -c([X,Y])$.
The \textbf{commutator ideal}\index{commutator ideal} $[\fg,\fg]$
is the subspace of $\fg$ spanned by $\{ [X,Y] \mid X,Y \in \fg \}$.
Since $\delta c = 0$ if and only if $c$ vanishes on $[\fg,\fg]$,
we conclude that $H^1(\fg;\RR) = [\fg,\fg]^0$,
where $[\fg,\fg]^0 \subseteq \fg^*$ is the \textbf{annihilator}
of $[\fg,\fg]$.
An element of $C^2$ is an alternating bilinear map
$c: \fg \times \fg \to \RR$, and $\delta c(X,Y,Z) =
-c([X,Y],Z) + c([X,Z],Y) - c([Y,Z],X)$.
If $c = \delta b$ for some $b \in C^1$, then
$c(X,Y) = (\delta b) (X,Y) = -b([X,Y])$.

If $\fg$ is the Lie algebra of a compact connected Lie group $G$,
then by averaging one can show that the de Rham cohomology
may be computed from the subcomplex of $G$-invariant forms,
and hence $H^k(\fg;\RR) = H_{\mathrm{de Rham}}^k(G)$.

\begin{proposition}
If $H^1(\fg;\RR) = H^2(\fg,\RR) = 0$,
then any symplectic $G$-action is hamiltonian.
\end{proposition}

\vspace*{-2ex}

\begin{proof}
Let $\psi: G \to \mathrm{Sympl}(M,\omega)$
be a symplectic action of $G$ on a symplectic manifold $(M,\omega)$.
Since $H^1(\fg;\RR) = 0$ means that $[\fg,\fg] = \fg$, and
since commutators of symplectic vector fields are hamiltonian, we have
$d\psi: \fg = [\fg,\fg] \to \chi^{\mathrm{ham}}(M)$.
The action $\psi$ is hamiltonian if and only if there is a Lie algebra
homomorphism $\mu^*: \fg \to C^{\infty}(M)$
such that the hamiltonian vector field of $\mu^* (\xi)$
is $d\psi (\xi)$.
We first take an arbitrary vector space lift
$\tau: \fg \to C^{\infty}(M)$ with this property,
i.e., for each basis vector $X \in \fg$, we choose
$\tau(X) = \tau^X \in C^{\infty}(M)$ such that $v_{(\tau^X)} = d\psi(X)$.
The map $X \mapsto \tau^X$ may not be a Lie algebra homomorphism.
By construction, $\tau^{[X,Y]}$ is a hamiltonian function for
$[X,Y]^{\#}$, and (as computed in Section~\ref{sec:integrable})
$\{\tau^X,\tau^Y\}$ is a hamiltonian function for $-[X^{\#},Y^{\#}]$.
Since $[X,Y]^{\#} = -[X^{\#},Y^{\#}]$, the
corresponding hamiltonian functions must differ by a constant:
\[
        \tau^{[X,Y]} - \{\tau^X,\tau^Y\} = c(X,Y) \in \RR \ .
\]
By the Jacobi identity, $\delta c = 0$.
Since $H^2(\fg;\RR) = 0$, there is $b \in \fg^*$
satisfying $c = \delta b$, $c(X,Y) = -b([X,Y])$.  We define
\[
\begin{array}{rrcl}
        \mu^*: & \fg & \longrightarrow & C^{\infty}(M) \\
        & X & \longmapsto & \mu^*(X) = \tau^X + b(X) = \mu^X \ .
\end{array}
\]
Now $\mu^*$ is a Lie algebra homomorphism:
$\mu^*([X,Y]) =
%\tau^{[X,Y]} + b([X,Y]) =
\{\tau^X,\tau^Y\} = \{\mu^X,\mu^Y\}$.
\end{proof}

By the Whitehead lemmas\index{Whitehead lemmas}
(see for instance~\cite[pages 93-95]{ja:lie})
a semisimple Lie group $G$ has $H^1(\fg;\RR) = H^2(\fg;\RR) = 0$.
As a corollary, {\em when $G$ is semisimple,
any symplectic $G$-action is hamiltonian.}\footnote{A
compact Lie group $G$ has $H^1(\fg;\RR) = H^2(\fg;\RR) = 0$
if and only if it is semisimple.
In fact, a compact Lie group $G$ is
semisimple\index{semisimple} when $\fg = [\fg,\fg]$.
The unitary group $\UU(n)$ is not semisimple because
the multiples of the identity, $S^1 \cdot \mathrm{Id}$,
form a nontrivial center; at the level of the Lie algebra, this
corresponds to the subspace $\RR \cdot \mathrm{Id}$ of
scalar matrices, which are not commutators since they are not traceless.
Any abelian Lie group is {\em not} semisimple.
Any direct product of the other compact classical groups
$\SU (n)$, $\SO (n)$ and $\Sp (n)$ is semisimple.
An arbitrary compact Lie group admits a finite cover by a direct
product of tori and semisimple Lie groups.}

\begin{proposition}
For a connected Lie group $G$,\index{moment map ! uniqueness}
if $H^1(\fg;\RR) = 0$, then moment maps for
hamiltonian $G$-actions are unique.
\end{proposition}

\vspace*{-2ex}

\begin{proof}
Suppose that $\mu_1$ and $\mu_2$ are two moment maps
for an action $\psi$.
For each $X \in \fg$, $\mu_1^X$ and $\mu_2^X$ are both
hamiltonian functions for $X^{\#}$, thus
$\mu_1^X - \mu_2^X = c(X)$ is locally constant.
This defines $c \in \fg^*$, $X \mapsto c(X)$.
Since the corresponding $\mu_i^* : \fg \to C^\infty (M)$
are Lie algebra homomorphisms,
we have $c([X,Y]) = 0$, $\forall X,Y
\in \fg$, i.e., $c \in [\fg,\fg]^0 = \{0\}$.
Hence, $\mu_1 = \mu_2$.
\end{proof}

In general, if $\mu: M \to \fg^*$ is a
moment map, then given any $c \in [\fg,\fg]^0$, $\mu_1 = \mu + c$ is
another moment map.
In other words, moment maps are unique up to elements of the
dual of the Lie algebra that annihilate the commutator ideal.

The two extreme cases are when
\[
\begin{array}{ll}
\mbox{$\bullet$ $G$ is semisimple:} &
\mbox{any symplectic action is hamiltonian} \ ,\\
& \mbox{moment maps are unique} \ ; \\
\mbox{$\bullet$ $G$ is abelian:} &
\mbox{symplectic actions may not be hamiltonian} \ ,\\
& \mbox{moment maps are unique up to a constant $c \in \fg^*$}\ .
\end{array}
\]

%%%%%%%%%%%%%%%%%%%%%%%%%%%%%%%%%%%%%%%%%%%%%%%%%%%%%%%%%%%%%%%%%%%%%%%%%%%%%
%%%%%%%%%%%%%%%%%%%%%%%%%%%%%%%%%%%%%%%%%%%%%%%%%%%%%%%%%%%%%%%%%%%%%%%%%%%%%

\ssubsection{Convexity}
\index{moment map ! properties}\index{moment map ! convexity}

Atiyah, Guillemin and Sternberg~\cite{at:convexity, gu-st:convexity}
showed that the image of the moment map for a
hamiltonian torus action on a compact connected symplectic
manifold is always a polytope.\footnote{A \textbf{polytope}
in $\RR^n$ is the convex hull of a finite number of points in $\RR^n$.
A \textbf{convex polyhedron} is a subset of $\RR^n$ that
is the intersection of a finite number of affine half-spaces.
Hence, polytopes coincide with bounded convex polyhedra.}
A proof of this theorem can also be found
in~\cite{mc-sa:introduction}.

\begin{theorem}\label{thm:convexity}\index{theorem !
Atiyah-Guillemin-Sternberg}\index{Atiyah-Guillemin-Sternberg
theorem}\index{Guillemin|see{Atiyah-Guillemin-Sternberg}}\index{theorem !
convexity}\index{Sternberg|see{Atiyah-Guillemin-Sternberg}}\index{convexity}
\textbf{(Atiyah, Guillemin-Sternberg)} $\;$
Let $(M,\omega)$ be a compact connected symplectic manifold.
Suppose that $\psi: \TT^m \to \mathrm{Sympl}(M,\omega)$ is a
hamiltonian action of an $m$-torus with moment map $\mu: M \to \RR^m$.
Then:
\begin{itemize}
\item[(a)]
the levels $\mu^{-1} (c)$ are connected ($c \in \RR^m$);\index{connectedness}
\item[(b)]
the image $\mu (M)$ is convex;
\item[(c)]
$\mu (M)$ is the convex hull
of the images of the fixed points of the action.
\end{itemize}
\end{theorem}

The image $\mu(M)$ of the moment map is called the
\textbf{moment polytope}\index{moment polytope}\index{polytope ! moment}.

\begin{examples}
\begin{enumerate}

\item
Suppose that $\TT^m$ acts linearly on $(\CC^n, \omega_0)$.
Let $\lambda^{(1)}, \ldots, \lambda^{(n)} \in \ZZ^m$ be the
{\em weights} appearing in the corresponding weight space decomposition,
that is,
\[
   \CC^n \simeq \displaystyle{ \bigoplus _{k=1}^n V_{\lambda^{(k)}}} \ ,
\]
where, for $\lambda^{(k)} = (\lambda^{(k)}_1,\ldots,\lambda^{(k)}_m)$,
the torus $\TT^m$ acts on the complex line $V_{\lambda^{(k)}}$ by
$( e^{i t_1} , \ldots, e^{i t_m} ) \cdot v =
e ^{i \sum_j \lambda^{(k)}_j t_j} v$.
If the action is effective\footnote{An action of a group $G$
on a manifold $M$ is called
\textbf{effective}\index{action ! effective}\index{effective ! action}
if each group element $g \ne e$ moves at least one point $p \in M$,
that is, $\cap_{p \in M} G_p = \{e\}$,
where $G_p = \{g \in G \mid g \cdot p = p\}$ is the stabilizer of $p$.},
then $m \leq n$ and the weights $\lambda^{(1)} , \ldots, \lambda^{(n)}$
are part of a $\ZZ$-basis of $\ZZ^m$.
If the action is symplectic (hence hamiltonian in this case),
then the weight spaces $V_{\lambda^{(k)}}$ are symplectic subspaces.
In this case, a moment map is given by
\[
   \mu (v) = - \textstyle{\frac12} \sum \limits_{k=1}^{n}
   \lambda^{(k)} | v_{\lambda^{(k)}} |^2 \ ,
\]
where $| \cdot |$ is the standard norm\footnote{The
standard inner product satisfies
$\langle v,w \rangle = \omega_0 (v, Jv)$ where
$J \frac{\partial}{\partial z} = i \frac{\partial}{\partial z}$
and $J \frac{\partial}{\partial \bar z} =
-i \frac{\partial}{\partial \bar z}$.
In particular, the standard norm is invariant for a symplectic
complex-linear action.}
and $v = v_{\lambda^{(1)}} + \ldots + v_{\lambda^{(n)}}$
is the weight space decomposition of $v$.
We conclude that, if $\TT^n$ acts on $\CC^n$ in a linear,
effective and hamiltonian way, then any moment map $\mu$
is a submersion, i.e., each differential
$d \mu_v : \CC^n \to \RR^n$ ($v \in \CC^n$) is surjective.

\item
Consider a coadjoint orbit $\cO_{\lambda}$ for the unitary group $\UU (n)$.
Multiplying by $i$, the orbit $\cO_{\lambda}$ can be viewed as
the set of hermitian matrices with a given eigenvalue
spectrum $\lambda = (\lambda_1 \geq \ldots \geq \lambda_n)$.
The restriction of the coadjoint action to the maximal torus
$\TT^n$ of diagonal unitary matrices is hamiltonian
with moment map $\mu : \cO_\lambda \to \RR^n$ taking
a matrix to the vector of its diagonal entries.
Then the moment polytope $\mu (\cO_\lambda)$ is the convex hull $C$
of the points given by all the permutations of
$(\lambda_1, \ldots , \lambda_n)$.
This is a rephrasing of the classical theorem of Schur
($\mu (\cO_\lambda) \subseteq C$) and Horn ($C \subseteq \mu (\cO_\lambda)$).
\end{enumerate}
\end{examples}

%%%%%%%%%%%%%%%%%%%%%%%%%%%%%%%%%%%%%%%%%%%%%%%%%%%%%%%%%%%%%%%%%%%%%%%%%%%%%

Example~1 is related to the
universal local picture for a moment map
near a fixed point of a hamiltonian torus action:

\begin{theorem}
\label{thm:moment_darboux}
Let $(M^{2n}, \omega, \TT^m , \mu)$ be a hamiltonian
$\TT^m$-space, where $q$ is a fixed point.
Then there exists a chart $(\cU,x_1,\dots,x_n,y_1,\dots,y_n)$
centered at $q$ and weights
$\lambda^{(1)}, \ldots, \lambda^{(n)} \in \ZZ^m$ such that
\[
   \left. \omega \right|_\cU =
   \sum\limits_{k=1}^n dx_k \wedge dy_k \quad \mbox{ and } \quad
   \left. \mu \right|_\cU =
   \mu (q) - \frac 12 \sum\limits_{k=1}^n
   \lambda^{(k)} (x_k^2 + y_k^2) \ .
\]
\end{theorem}

%%%%%%%%%%%%%%%%%%%%%%%%%%%%%%%%%%%%%%%%%%%%%%%%%%%%%%%%%%%%%%%%%%%%%%%%%%%%%

The following two results use the crucial fact that
any effective action of an $m$-torus on a manifold
has orbits of dimension $m$;
a proof may be found in~\cite{br:groups}.

\begin{corollary}
Under the conditions of the convexity theorem,
if the $\TT^m$-action is effective, then there must be at least
$m+1$ fixed points.
\end{corollary}

\vspace*{-2ex}

\begin{proof}
At a point $p$ of an $m$-dimensional orbit the moment map is
a submersion, i.e., $(d\mu_1)_p,\dots,(d\mu_m)_p$ are linearly independent.
Hence, $\mu(p)$ is an interior point of $\mu(M)$, and $\mu(M)$ is a
nondegenerate polytope.
A nondegenerate polytope in
$\RR^m$ has at least $m+1$ vertices.
The vertices of $\mu(M)$ are images of fixed points.
\end{proof}

\vspace*{-1ex}

\begin{proposition}
\label{prop:dimension}
Let $(M,\omega,\TT^m,\mu)$ be a hamiltonian $\TT^m$-space.
If the $\TT^m$-action is effective, then $\dim M \ge 2m$.
\end{proposition}

\vspace*{-2ex}

\begin{proof}
Since the moment map is constant on an orbit $\cO$,
for $p \in \cO$ the differential $d\mu_p: T_pM \to \fg^*$
maps $T_p\cO$ to $0$.
Thus $T_p\cO \subseteq \ker d\mu_p = (T_p\cO)^{\omega}$,
where $(T_p\cO)^{\omega}$ is the symplectic orthogonal
of $T_p\cO$.
This shows that orbits $\cO$ of a hamiltonian torus
action are isotropic submanifolds of $M$.
In particular, by symplectic linear algebra we have
that $\dim \cO \le \frac {1}{2} \dim M$.
Now consider an $m$-dimensional orbit.
\end{proof}

%%%%%%%%%%%%%%%%%%%%%%%%%%%%%%%%%%%%%%%%%%%%%%%%%%%%%%%%%%%%%%%%%%%%%%%%%%%%%

For a hamiltonian action of an arbitrary compact Lie group $G$
on a compact symplectic manifold $(M,\omega)$,
the following {\em nonabelian} convexity theorem
was proved by Kirwan~\cite{ki:convexity_III}:
if $\mu: M \to \fg^*$ is a moment map,
then the intersection $\mu(M) \cap \ft_+^*$ of the image of $\mu$
with a Weyl chamber for a Cartan subalgebra $\ft \subseteq \fg$
is a convex polytope.
This had been conjectured by Guillemin and Sternberg
and proved by them in particular cases.

%%%%%%%%%%%%%%%%%%%%%%%%%%%%%%%%%%%%%%%%%%%%%%%%%%%%%%%%%%%%%%%%%%%%%%%%%%%%%
%%%%%%%%%%%%%%%%%%%%%%%%%%%%%%%%%%%%%%%%%%%%%%%%%%%%%%%%%%%%%%%%%%%%%%%%%%%%%
% --> Section 6
%%%%%%%%%%%%%%%%%%%%%%%%%%%%%%%%%%%%%%%%%%%%%%%%%%%%%%%%%%%%%%%%%%%%%%%%%%%%%
%%%%%%%%%%%%%%%%%%%%%%%%%%%%%%%%%%%%%%%%%%%%%%%%%%%%%%%%%%%%%%%%%%%%%%%%%%%%%

\newpage

\ssection{Symplectic Reduction}
\index{symplectic reduction}
\label{section6}

%%%%%%%%%%%%%%%%%%%%%%%%%%%%%%%%%%%%%%%%%%%%%%%%%%%%%%%%%%%%%%%%%%%%%%%%%%%%%
%%%%%%%%%%%%%%%%%%%%%%%%%%%%%%%%%%%%%%%%%%%%%%%%%%%%%%%%%%%%%%%%%%%%%%%%%%%%%

\ssubsection{Marsden-Weinstein-Meyer Theorem}
\index{theorem ! Marsden-Weinstein-Meyer}\index{Weinstein !
Marsden-Weinstein-Meyer theorem}\index{Marsden-Weinstein-Meyer !
theorem}\index{Meyer|see{Marsden-Weinstein-Meyer}}\index{symplectic !
reduction|see{reduction}}
\label{sec:m-w-m-thm}

Classical physicists realized that,
whenever there is a symmetry group of dimension $k$
acting on a mechanical system,
the number of degrees of freedom for the
position and momenta of the particles
may be reduced by $2k$.
Symplectic reduction formulates this process mathematically.

%   Symplectic reduction is useful for
%building new symplectic manifolds, producing interesting
%invariants, simplifying computations, providing
%a beautiful geometry, etc.

%%%%%%%%%%%%%%%%%%%%%%%%%%%%%%%%%%%%%%%%%%%%%%%%%%%%%%%%%%%%%%%%%%%%%%%%%%%%%

\begin{theorem}\index{theorem ! Marsden-Weinstein-Meyer}\index{Weinstein !
Marsden-Weinstein-Meyer theorem}\index{Marsden-Weinstein-Meyer !
theorem}\label{thm:reduction}
\textbf{(Marsden-Weinstein, Meyer~\cite{ma-we:reduction,me:symmetries})} $\;$
Let $(M,\omega,G,\mu)$ be a hamiltonian $G$-space
(Section~\ref{sec:actions}) for a compact Lie group $G$.
Let $i: \mu^{-1}(0) \hookrightarrow M$ be the inclusion map.
Assume that $G$ acts freely on $\mu^{-1}(0)$.  Then
\begin{itemize}
\item[(a)]
the orbit space $M_{\mathrm{red}} = \mu^{-1}(0)/G$ is a manifold,
\item[(b)]
$\pi: \mu^{-1}(0) \rightarrow M_{\mathrm{red}}$ is a principal
$G$-bundle, and
\item[(c)]
there is a symplectic form $\omega_{\mathrm{red}}$ on
$M_{\mathrm{red}}$ satisfying $i^*\omega = \pi^*\omega_{\mathrm{red}}$.
\end{itemize}
\end{theorem}

\vspace*{-1ex}

%%%%%%%%%%%%%%%%%%%%%%%%%%%%%%%%%%%%%%%%%%%%%%%%%%%%%%%%%%%%%%%%%%%%%%%%%%%%%

\begin{definition}\index{quotient ! Marsden-Weinstein-Meyer}\index{Weinstein !
Marsden-Weinstein-Meyer quotient}\index{Marsden-Weinstein-Meyer !
quotient}\index{quotient ! symplectic}\index{symplectic !
quotient}\index{reduced ! space}
The symplectic manifold $(M_{\mathrm{red}},\omega_{\mathrm{red}})$
is the \textbf{reduction} (or \textbf{reduced space},
or \textbf{symplectic quotient})
of $(M,\omega)$ with respect to $G,\mu$.
\end{definition}

%%%%%%%%%%%%%%%%%%%%%%%%%%%%%%%%%%%%%%%%%%%%%%%%%%%%%%%%%%%%%%%%%%%%%%%%%%%%%

When $M$ is K\"ahler\index{K\"ahler structure}
and the action of $G$ preserves the complex structure,
we can show that the symplectic reduction has a natural K\"ahler structure.

%%%%%%%%%%%%%%%%%%%%%%%%%%%%%%%%%%%%%%%%%%%%%%%%%%%%%%%%%%%%%%%%%%%%%%%%%%%%%

Let $(M,\omega,G,\mu)$ be a hamiltonian $G$-space
for a compact Lie group $G$.
To reduce at a level $\xi \in \fg^*$ of $\mu$,\index{reduction ! other levels}
we need $\mu^{-1}(\xi)$ to be preserved by $G$,
or else take the $G$-orbit of $\mu^{-1}(\xi)$,
or else take the quotient by the maximal subgroup of $G$
that preserves $\mu^{-1}(\xi)$.
Since $\mu$ is equivariant,
$G$ preserves $\mu^{-1}(\xi)$ if and only if
$\mathrm{Ad}_g^*\xi = \xi$, $\forall g \in G$.
Of course, the level $0$ is always preserved.
Also, when $G$ is a torus, any level is
preserved and \textbf{reduction at $\xi$}
for the moment map $\mu$, is equivalent to
reduction at $0$ for a shifted moment map
$\phi: M \rightarrow \fg^*$, $\phi(p) := \mu(p) - \xi$.
In general, let $\cO$ be a coadjoint orbit in $\fg^*$ equipped with the
\textbf{canonical symplectic form}\index{canonical !
symplectic form on a coadjoint orbit}\index{symplectic !
canonical symplectic form on a coadjoint orbit}
$\omega_{\cO}$ (defined in Section~\ref{sec:symplectic_hamiltonian_fields}).
Let $\cO^-$ be the orbit $\cO$ equipped with $-\omega_{\cO}$.
The natural product action of $G$ on $M \times \cO^-$ is hamiltonian
with moment map $\mu_{\cO}(p,\xi) = \mu(p) - \xi$.
If the hypothesis of Theorem~\ref{thm:reduction} is satisfied
for $M \times \cO^-$, then
one obtains a \textbf{reduced space with respect to the coadjoint orbit
$\cO$}.\index{reduced ! space}

%%%%%%%%%%%%%%%%%%%%%%%%%%%%%%%%%%%%%%%%%%%%%%%%%%%%%%%%%%%%%%%%%%%%%%%%%%%%%

\begin{examples}
\index{reduction ! examples}\index{example ! reduction}
\begin{enumerate}
\item
The standard symplectic form on $\CC^n$ is
$\omega_0 = \frac {i}{2} \sum dz_i \wedge d{\bar z}_i =
\sum dx_i \wedge dy_i = \sum r_idr_i \wedge d\theta_i$
in polar coordinates.
The $S^1$-action on $(\CC^n,\omega_0)$ where $e^{it} \in S^1$
acts as multiplication by $e^{it}$ has vector field
$X^{\#} = \frac {\partial}{\partial\theta_1} + \frac
{\partial}{\partial\theta_2} + \dots + \frac {\partial}{\partial\theta_n}$.
This action is hamiltonian with moment map
$\mu: \CC^n \to \RR$, $\mu(z) = -\frac {|z|^2}{2}$,
since $\imath_{X^\#}\omega = \sum r_idr_i
= -\frac {1}{2} \sum dr_i^2 = d\mu$.
The level $\mu^{-1}(- \frac 12)$ is the unit sphere $S^{2n-1}$,
whose orbit space is the projective space,\index{complex !
projective space}\index{reduction ! reduced space}\index{reduced ! space}
\[
   \mu^{-1}(\textstyle{- \frac 12})/S^1 = S^{2n-1}/S^1 = \CC\PP^{n-1} \ .
\]
%which is a $(2n-2)$-dimensional symplectic manifold.
The reduced symplectic form at level $- \frac 12$ is 
$\omega_{_{\mathrm{red}}} = \omega_{_{\mathrm{FS}}}$
the Fubini-Study symplectic form.\index{Fubini-Study
form}\index{Study|see{Fubini-Study}}\index{symplectic
! Fubini-Study form}\index{form ! Fubini-Study}
Indeed, if $\mathrm{pr} : \CC^{n+1} \setminus \{0\} \to \CC \PP ^n$
is the standard projection, the forms
$\mathrm{pr}^* \omega_{_{\mathrm{FS}}} = \textstyle{\frac{i}{2}}
\partial \bar{\partial} \log (|z|^2)$ and $\omega_0$
have the same restriction to $S^{2n+1}$.

\item
Consider the natural action of $\UU (k)$ on $\CC^{k\times n}$
with moment map $\mu (A) = \textstyle{{i} \over {2}} A A^*
+ {\mathrm{Id} \over {2i}}$ for $A \in \CC^{k\times n}$
(Section~\ref{sec:actions}).
Since
$\mu^{-1} (0) = \{ A \in \CC^{k\times n} \, | \, A A^* = \mathrm{Id} \}$,
the reduced manifold is the grassmannian of $k$-planes in $\CC^n$:
\[
        \mu^{-1} (0) / \UU (k) = \GG (k,n) \ .
\]
\end{enumerate}
\end{examples}

%The natural actions of $\TT ^{n+1}$ and $\UU (n+1)$ on
%$(\CC \PP ^n, \omega_{_{\mathrm{FS}}})$ are hamiltonian,
%and have moment maps...

%%%%%%%%%%%%%%%%%%%%%%%%%%%%%%%%%%%%%%%%%%%%%%%%%%%%%%%%%%%%%%%%%%%%%%%%%%%%%

For the case where $G = S^1$ and $\dim M = 4$,
here is a glimpse of reduction.\index{reduction ! low-brow proof}
Let $\mu: M \rightarrow \RR$ be the moment map and $p \in \mu^{-1}(0)$.
Choose local coordinates near $p$:
$\theta$ along the orbit through $p$,
$\mu$ given by the moment map, and $\eta_1,\eta_2$ the pullback
of coordinates on $M_{\mathrm{red}} = \mu^{-1}(0)/S^1$.
Then the symplectic form can be written
\[
        \omega = A \ d\theta \wedge d\mu + \sum B_j \ d\theta \wedge d\eta_j
        + \sum C_j \ d\mu \wedge d\eta_j + D \ d\eta_1 \wedge d\eta_2 \ .
\]
As $d\mu = \imath \left( \frac {\partial}{\partial\theta}
\right)\omega$, we must have $A = 1$, $B_j = 0$.
Since $\omega$ is symplectic, it must be $D \ne 0$.
Hence, $i^*\omega = D \ d\eta_1 \wedge d\eta_2$ is the pullback of a
symplectic form on $M_{\mathrm{red}}$.

%%%%%%%%%%%%%%%%%%%%%%%%%%%%%%%%%%%%%%%%%%%%%%%%%%%%%%%%%%%%%%%%%%%%%%%%%%%%%

The actual proof of Theorem~\ref{thm:reduction}
requires some preliminary ingredients.

Let $\mu : M \to \fg^*$ be the moment map for an (hamiltonian) action
of a Lie group $G$ on a symplectic manifold $(M,\omega)$.
Let $\fg_p$ be the Lie algebra of the stabilizer of a point $p \in M$,
let $\fg_p^0 = \{\xi \in \fg^* \mid \langle \xi,X \rangle =
0,\ \forall X \in \fg_p\}$ be the annihilator of $\fg_p$,
and let $\cO _p$ be the $G$-orbit through $p$.
Since $\omega_p(X_p^{\#},v) = \langle d\mu_p(v),X \rangle$,
for all $v \in T_pM$ and all $X \in \fg$,
%and count dimensions.
the differential $d\mu_p: T_pM \rightarrow \fg^*$ has
\[
   {\mathrm{ker}} \ d\mu_p = (T_p \cO _p)^{\omega_p}
   \quad \mbox{ and } \quad
   {\mathrm{im}} \ d\mu_p = \fg_p^0 \ .
\]
Consequently,
the action is locally free\footnote{The action is
\textbf{locally free} at $p$ when $\fg_p = \{ 0 \}$, i.e.,
the stabilizer of $p$ is a discrete group.
The action is \textbf{free} at $p$ when the stabilizer of $p$
is trivial, i.e., $G_p = \{ e \}$.} at $p$
if and only if $p$ is a regular point of $\mu$
(i.e., $d\mu_p$ is surjective), and we obtain:

\begin{lemma}
\label{lem:free}
If $G$ acts freely on $\mu^{-1}(0)$,
then $0$ is a regular value of $\mu$,
the level $\mu^{-1}(0)$ is a submanifold of $M$
of codimension $\dim G$, and,
for $p \in \mu^{-1}(0)$, the tangent space
$T_p\mu^{-1}(0) = \ker d\mu_p$ is the symplectic
orthogonal to $T_p{\cal O}_p$ in $T_pM$.
\end{lemma}

In particular, {\em orbits in $\mu^{-1}(0)$ are isotropic}.
Since any tangent vector to the orbit is the value of a vector field
generated by the group, we can show this directly by computing,
for any $X,Y \in \fg$ and $p \in \mu^{-1}(0)$,
the hamiltonian function for $[Y^\#,X^\#] = [Y,X]^\#$
at that point: $\omega_p(X_p^\#,Y_p^\#) = \mu^{[Y,X]}(p) = 0$.

\begin{lemma}
\label{lem:isotropic}
Let $(V,\Omega)$ be a symplectic vector space, and $I$ an isotropic subspace.
Then $\Omega$ induces a canonical symplectic structure
$\Omega_{\mathrm{red}}$ on $I^{\Omega}/I$.
\end{lemma}

\vspace*{-2ex}

\begin{proof}
Let $[u],[v]$ be the classes in $I^{\Omega}/I$ of $u,v \in I^{\Omega}$.
We have $\Omega(u+i,v+j) = \Omega(u,v)$, $\forall i,j \in I$, because
$\Omega(u,j) = \Omega(i,v) = \Omega(i,j) = 0$.
Hence, we can define $\Omega_{\mathrm{red}} ([u],[v]) := \Omega(u,v)$.
This is nondegenerate: if $u \in I^{\Omega}$ has $\Omega(u,v) = 0$,
for all $v \in I^{\Omega}$, then
$u \in (I^{\Omega})^{\Omega} = I$, i.e., $[u] = 0$.
\end{proof}

\vspace*{-1ex}

\begin{proposition}
\label{prop:quotient}
If a compact Lie group $G$ acts freely on a manifold $M$,
then $M/G$ is a manifold and the map $\pi: M \rightarrow M/G$
is a principal $G$-bundle.
\end{proposition}

\vspace*{-2ex}

\begin{proof}
We first show that, for any $p \in M$,
the $G$-orbit through $p$ is a compact
submanifold of $M$ diffeomorphic to $G$.\footnote{Even
if the action is not free, the orbit through $p$
is a compact submanifold of $M$.
In that case, the orbit of a point $p$ is diffeomorphic
to the quotient $G / G_p$ of $G$ by the stabilizer of $p$.}
The $G$-orbit through $p$ is the
image of the smooth injective map $\mathrm{ev}_p: G \to M$,
$\mathrm{ev}_p(g) = g \cdot p$.
%Injectivity of $\mathrm{ev}_p$ follows from the action being free.
The map $\mathrm{ev}_p$ is proper because, if $A$ is
a compact, hence closed, subset of $M$, then its inverse image
$(\mathrm{ev}_p)^{-1} (A)$, being a closed subset of
the compact Lie group $G$, is also compact.
The differential $d (\mathrm{ev}_p )_e$ is injective
because $d (\mathrm{ev}_p )_e (X) = 0 \Leftrightarrow
X_p^{\#} = 0 \Leftrightarrow X = 0$, $\forall X \in T_e G$,
as the action is free.
At any other point $g \in G$, for $X \in T_g G$ we have
$d (\mathrm{ev}_p )_g (X) = 0 \Leftrightarrow
d (\mathrm{ev}_p \circ R_g )_e  \circ (d R_{g^{-1}} )_g (X) = 0$,
where $R_g : G \to G$, $h \mapsto hg$, is right multiplication by $g$.
But $\mathrm{ev}_p \circ R_g = \mathrm{ev}_{g \cdot p}$
has an injective differential at $e$,
and $(d R_{g^{-1}} )_g$ is an isomorphism.
It follows that $d (\mathrm{ev}_p )_g$ is always injective,
so $\mathrm{ev}_p$ is an immersion.
We conclude that $\mathrm{ev}_p$ is a closed embedding.

We now apply the slice
theorem\footnote{\textbf{Slice Theorem: }\index{theorem !
slice}\index{slice theorem}\label{thm:slice}
{\em Let $G$ be a compact Lie group acting on a manifold $M$
such that $G$ acts freely at $p \in M$.
Let $S$ be a transverse section to $\cO_p$ at $p$
(this is called a \textbf{slice}).
Choose a coordinate chart $x_1,\dots,x_n$ centered at $p$ such that
$\cO_p \simeq G$ is given by $x_1 = \dots = x_k = 0$ and
$S$ by $x_{k+1} = \dots = x_n = 0$.
Let $S_{\varepsilon} = S \cap B_{\varepsilon}$ where
$B_{\varepsilon}$ is the ball of radius $\varepsilon$
centered at $0$ with respect to these coordinates.
Let $\eta: G \times S \rightarrow M$, $\eta(g,s) = g\cdot s$.
Then, for sufficiently small $\varepsilon$,
the map $\eta: G \times S_{\varepsilon} \rightarrow M$ takes
$G \times S_{\varepsilon}$ diffeomorphically onto a $G$-invariant
neighborhood $\cU$ of the $G$-orbit through $p$.}
In particular, if the action of $G$ is free at $p$,
then the action is free on $\cU$,
so the set of points where $G$ acts freely is open.}
which is an equivariant tubular neighborhood
theorem.\index{equivariant ! tubular neighborhood theorem}\index{tubular
neighborhood ! equivariant}
For $p \in M$, let $q = \pi(p) \in M/G$.
Choose a $G$-invariant neighborhood $\cU$ of $p$ as in the
slice theorem, so that $\cU \simeq G \times S$
where $S$ is an appropriate slice.
Then $\pi(\cU) = \cU/G =: \cV$ is a neighborhood of $q$ in $M/G$
homeomorphic\footnote{We equip the orbit space $M/G$ with the
\textbf{quotient topology}\index{quotient ! topology}, i.e.,
$\cV \subseteq M/G$ is open if and only if $\pi^{-1}(\cV)$
is open in $M$.} to $S$.
Such neighborhoods $\cV$ are used as charts on $M/G$.
To show that
the associated transition maps are smooth, consider two $G$-invariant
open sets $\cU_1,\cU_2$ in $M$ and corresponding slices $S_1,S_2$.
Then $S_{12} = S_1 \cap \cU_2$, $S_{21} = S_2 \cap \cU_1$ are
both slices for the $G$-action on $\cU_1 \cap \cU_2$.
To compute the transition map $S_{12} \rightarrow S_{21}$, consider
the sequence
$S_{12} \stackrel{\simeq}{\longrightarrow} \{ e \} \times S_{12}
\hookrightarrow G \times S_{12} \stackrel{\simeq}{\longrightarrow}
\cU_1 \cap \cU_2$ and similarly for $S_{21}$.
The composition $S_{12} \hookrightarrow \cU_1 \cap \cU_2
\stackrel{\simeq}{\longrightarrow}
G \times S_{21} \stackrel{pr}{\longrightarrow} S_{21}$ is smooth.

Finally, we show that $\pi: M \rightarrow M/G$
is a principal $G$-bundle.
For $p \in M$, $q = \pi(p)$, choose a $G$-invariant neighborhood $\cU$
of $p$ of the form $\eta: G \times S \stackrel{\simeq}{\longrightarrow} \cU$.
Then $\cV = \cU/G \simeq S$ is the corresponding neighborhood of $q$ in $M/G$:
\[
\begin{array}{rccc}
        M \supseteq & \cU &
        \stackrel{\eta}{\simeq} \; \; G \times S \; \; \simeq &
        G \times \cV \\
        & \phantom{\pi} \downarrow \pi & & \downarrow \\
        M/G \supseteq & \cV & = & \cV
\end{array}
\]
Since the projection on the right is smooth, $\pi$ is smooth.
By considering the overlap of two trivializations
$\phi_1: \cU_1 \to G \times \cV_1$ and $\phi_2: \cU_2 \to G \times \cV_2$,
we check that the transition map
$\phi_2 \circ \phi_1^{-1} = (\sigma_{12}, \id) :
G \times (\cV_1 \cap \cV_2 ) \to G \times (\cV_1 \cap \cV_2 )$ is smooth.
\end{proof}

%%%%%%%%%%%%%%%%%%%%%%%%%%%%%%%%%%%%%%%%%%%%%%%%%%%%%%%%%%%%%%%%%%%%%%%%%%%%%

\noindent
\textbf{Proof of Theorem~\ref{thm:reduction}.}
Since $G$ acts freely on $\mu^{-1}(0)$,
by Lemma~\ref{lem:free} the level $\mu^{-1}(0)$
is a submanifold.
Applying Proposition~\ref{prop:quotient} to the free action
of $G$ on the manifold $\mu^{-1}(0)$, we conclude the
assertions~(a) and~(b).

At $p \in \mu^{-1}(0)$ the tangent space to the orbit
$T_p\cO_p$ is an isotropic subspace of the
symplectic vector space $(T_pM,\omega_p)$.
By Lemma~\ref{lem:isotropic} there is a canonical symplectic structure
on the quotient $T_p\mu^{-1}(0)/T_p\cO_p$.
The point $[p] \in M_{\mathrm{red}} = \mu^{-1}(0)/G$
has tangent space $T_{[p]}
M_{\mathrm{red}} \simeq T_p\mu^{-1}(0)/T_p\cO_p$.
This gives a well-defined nondegenerate 2-form
$\omega_{\mathrm{red}}$ on $M_{\mathrm{red}}$
because $\omega$ is $G$-invariant.
By construction $i^*\omega = \pi^*\omega_{\mathrm{red}}$ where
\[
\begin{array}{cll}
        \mu^{-1}(0) & \stackrel{i}{\hookrightarrow} & M \\
        \downarrow \pi \\
        M_{\mathrm{red}}
\end{array}
\]
The injectivity of $\pi^*$ yields closedness:
% of $\omega_{\mathrm{red}}$:
$\pi^* d \omega_{\mathrm{red}} = d \pi^* \omega_{\mathrm{red}}
= d \imath ^* \omega = \imath ^* d \omega = 0$.
\hfill $\Box$

%%%%%%%%%%%%%%%%%%%%%%%%%%%%%%%%%%%%%%%%%%%%%%%%%%%%%%%%%%%%%%%%%%%%%%%%%%%%%
%%%%%%%%%%%%%%%%%%%%%%%%%%%%%%%%%%%%%%%%%%%%%%%%%%%%%%%%%%%%%%%%%%%%%%%%%%%%%

\ssubsection{Applications and Generalizations}
\label{sec:generalizations}

Let $(M,\omega,G,\mu)$ be a hamiltonian $G$-space
for a compact Lie group $G$.
Suppose that another Lie group $H$ acts on $(M,\omega)$
in a hamiltonian way with moment map $\phi: M \rightarrow \fh^*$.
Suppose that the $H$-action commutes with the $G$-action,
that $\phi$ is $G$-invariant and that $\mu$ is $H$-invariant.
Assuming that $G$ acts freely on $\mu^{-1}(0)$,
let $(M_{\mathrm{red}}, \omega_{\mathrm{red}})$ be the
corresponding reduced space.
Since the action of $H$ preserves $\mu^{-1}(0)$ and $\omega$
and commutes with the $G$-action,
the reduced space $(M_{\mathrm{red}}, \omega_{\mathrm{red}})$
inherits a symplectic action of $H$.
Since $\phi$ is preserved by the $G$-action,
the restriction of this moment map to $\mu^{-1}(0)$ descends
to a moment map $\phi_{\mathrm{red}}: M_{\mathrm{red}} \rightarrow \fh^*$
satisfying $\phi_{\mathrm{red}} \circ \pi = \phi \circ i$,
where $\pi : \mu^{-1}(0) \to M_{\mathrm{red}}$
and $i : \mu^{-1}(0) \hookrightarrow M$.
Therefore, $(M_{\mathrm{red}}, \omega_{\mathrm{red}}, H, \phi_{\mathrm{red}})$
is a hamiltonian $H$-space.
%Notice that $\phi$ is $G$-invariant, since the $H$-action
%commutes with the $G$-action; this can be already
%seen at the level of hamiltonian functions as $\{\phi^Y,\mu^X\} = 0$.

%%%%%%%%%%%%%%%%%%%%%%%%%%%%%%%%%%%%%%%%%%%%%%%%%%%%%%%%%%%%%%%%%%%%%%%%%%%%%

Consider now the action of a \textbf{product group} $G = G_1 \times G_2$,
where $G_1$ and $G_2$ are compact connected Lie groups.\index{reduction !
for product groups}\index{product group}\index{group ! product}
We have $\fg = \fg_1 \oplus \fg_2$ and $\fg^* = \fg_1^* \oplus \fg_2^*$.
Suppose that $(M,\omega,G,\psi)$ is a hamiltonian $G$-space with moment map
\[
        \psi = (\psi_1,\psi_2) : M \longrightarrow \fg_1^* \oplus \fg_2^*\ ,
\]
where $\psi_i : M \rightarrow \fg_i^*$ for $i=1,2$.
The fact that $\psi$ is equivariant
implies that $\psi_1$ is invariant under $G_2$
and $\psi_2$ is invariant under $G_1$.
Assume that $G_1$ acts freely on $Z_1 := \psi_1^{-1}(0)$.
Let $(M_1 = Z_1/G_1,\omega_1)$ be the
reduction of $(M,\omega)$ with respect to $G_1,\psi_1$.
From the observation above, $(M_1,\omega_1)$
inherits a hamiltonian $G_2$-action with moment map
$\mu_2 : M_1 \rightarrow \fg_2^*$
such that $\mu_2 \circ \pi = \psi_2 \circ i$,
where $\pi : Z_1 \to M_1$ and $i : Z_1 \hookrightarrow M$.
If $G$ acts freely on $\psi ^{-1} (0,0)$,
then $G_2$ acts freely on $\mu_2 ^{-1} (0)$,
and there is a natural symplectomorphism
\[
   \mu_2 ^{-1} (0) / G_2 \; \simeq \; \psi ^{-1} (0,0) / G \ .
\]
This technique of performing reduction with respect to
one factor of a product group at a time is called
\textbf{reduction in stages}\index{reduction ! in stages}.
It may be extended to reduction by a normal subgroup $H \subset G$
and by the corresponding quotient group $G / H$.

%%%%%%%%%%%%%%%%%%%%%%%%%%%%%%%%%%%%%%%%%%%%%%%%%%%%%%%%%%%%%%%%%%%%%%%%%%%%%

\begin{example}\index{reduction ! symmetry}
Finding symmetries for a mechanical problem
may reduce degrees of freedom by two at a time:
an integral of motion $f$ for a $2n$-dimensional hamiltonian system
$(M,\omega,H)$ may allow to understand the trajectories of this
system in terms of the trajectories of a $(2n-2)$-dimensional
hamiltonian system $(M_{\mathrm{red}},\omega_{\mathrm{red}},H_{\mathrm{red}})$.
Locally this process goes as follows.
Let $(\cU,x_1,\dots,x_n,\xi_1,\dots,\xi_n)$ be a Darboux chart for $M$
such that $f = \xi_n$.\footnote{To obtain such a chart,
in the proof of Darboux's Theorem~\ref{thm:darboux}
start with coordinates $(x'_1,\ldots ,x'_n,$ $y'_1,\ldots y'_n)$
such that $y'_n = f$ and $\frac{\partial}{\partial x'_n} = X_f$.}
Since $\xi_n$ is an integral of motion,
$0 = \{\xi_n,H\} = -\frac {\partial H}{\partial x_n}$,
the trajectories of the hamiltonian vector field $X_H$
lie on a constant level $\xi_n = c$ (Proposition~\ref{prop:integral}),
and $H$ does not depend on $x_n$.
The \textbf{reduced space}\index{reduced !
phase space}\index{phase space} is
$\cU_{\mathrm{red}} =
\{(x_1,\dots,x_{n-1},\xi_1,\dots,\xi_{n-1}) \mid
\exists a: (x_1,\dots,x_{n-1},a,\xi_1,\dots,\xi_{n-1},c) \in \cU \}$
and the \textbf{reduced hamiltonian}\index{reduced !
hamiltonian}\index{hamiltonian ! reduced} is
$H_{\mathrm{red}}: \cU_{\mathrm{red}} \to \RR$,
$H_{\mathrm{red}}(x_1,\dots,x_{n-1},\xi_1,\dots,\xi_{n-1}) =$
$H(x_1,\dots,x_{n-1},a,\xi_1,\dots,\xi_{n-1},c)$ for some $a$.
In order to find the trajectories of the original system on the
hypersurface $\xi_n = c$, we look for the trajectories
$(x_1(t),\dots,x_{n-1}(t),\xi_1(t),\dots,\xi_{n-1}(t))$
of the reduced system on $\cU_{\mathrm{red}}$,
and integrate the equation
$\frac {dx_n}{dt} (t) = \frac {\partial H}{\partial \xi_n}$
to obtain the original trajectories where
\[
\left\{ \begin{array}{rcl}
        x_n(t) & = & x_n(0) +
        \displaystyle{\int_0^t \frac {\partial H}{\partial \xi_n}
        (x_1(t),\dots,x_{n-1}(t),\xi_1(t),\dots,\xi_{n-1}(t),c)\ dt} \\
        \xi_n(t) & = & c \ .
\end{array} \right.
\]
\end{example}

%%%%%%%%%%%%%%%%%%%%%%%%%%%%%%%%%%%%%%%%%%%%%%%%%%%%%%%%%%%%%%%%%%%%%%%%%%%%%

By Sard's theorem, the singular values
of a moment map $\mu : M \to \fg^*$ form a set of measure zero.
So, perturbing if necessary,
we may assume that a level of $\mu$ is regular
hence, when $G$ is compact, that any point $p$ of that level
has finite stabilizer $G_p$.
Let $\cO_p$ be the orbit of $p$.
By the slice theorem
for the case of orbifolds, near $\cO_p$
the orbit space of the level is modeled by $S/G_p$,
where $S$ is a $G_p$-invariant disk in the level
and transverse to $\cO_p$ (a {\em slice}).
Thus, the orbit space is an {\em orbifold}.\index{orbifold !
reduced space}\footnote{Let $|M|$ be a Hausdorff
topological space satisfying the second axiom of countability.
An \textbf{orbifold chart}\index{orbifold ! chart}
on $|M|$ is a triple $(\cV, \Gamma, \varphi)$,
where $\cV$ is a connected open subset of some euclidean
space $\RR^m$, $\Gamma$ is a finite group that acts
linearly on $\cV$ so that the set of points where
the action is not free has codimension at least two,
and $\varphi: \cV \to |M|$ is a $\Gamma$-invariant
map inducing a homeomorphism from
$\cV / \Gamma$ onto its image $\cU \subset |M|$.
An \textbf{orbifold atlas}\index{orbifold ! atlas}
$\cA$ for $|M|$ is a collection of
orbifold charts on $|M|$ such that:
the collection of images $\cU$ forms a basis of open sets
in $|M|$, and the charts are compatible in the sense that,
whenever two charts $(\cV_1, \Gamma_1, \varphi_1)$ and
$(\cV_2, \Gamma_2, \varphi_2)$ satisfy $\cU_1 \subseteq \cU_2$,
there exists an injective homomorphism
$\lambda :\Gamma_1 \to \Gamma_2$ and a $\lambda$-equivariant
open embedding $\psi : \cV_1 \to \cV_2$ such that
$\varphi_2 \circ \psi = \varphi_1$.
Two orbifold atlases are \textbf{equivalent}\index{orbifold ! equivalence}
if their union is still an atlas.
An $m$-dimensional \textbf{orbifold}\index{orbifold ! definition}
$M$ is a Hausdorff topological
space $|M|$ satisfying the second axiom of countability,
plus an equivalence class of orbifold atlases on $|M|$.
We do not require the action of
each group $\Gamma$ to be effective.
Given a point $p$ on an orbifold $M$,
let $(\cV, \Gamma, \varphi)$ be an orbifold chart
for a neighborhood $\cU$ of $p$.
The \textbf{orbifold structure group}\index{orbifold ! structure group}
of $p$, $\Gamma_p$, is (the isomorphism class of) the stabilizer
of a pre-image of $p$ under $\phi$.
Orbifolds were introduced by Satake in~\cite{sa:orbifolds}.}
This implies that, when $G = \TT^n$ is an $n$-torus,
for most levels reduction goes through,
however the quotient space is not necessarily a manifold
but an orbifold.
Roughly speaking, orbifolds are singular manifolds where each
singularity is locally modeled on $\RR^m / \Gamma$,
for some finite group $\Gamma \subset \GL (m;\RR)$.
The differential-geometric notions of vector fields, differential
forms, exterior differentiation, group actions, etc., extend naturally
to orbifolds by gluing corresponding local $\Gamma$-invariant or
$\Gamma$-equivariant objects.
In particular, a \textbf{symplectic orbifold}\index{orbifold ! symplectic}
is a pair $(M,\omega)$
where $M$ is an orbifold and $\omega$ is a closed 2-form on $M$
that is nondegenerate at every point.

\begin{examples}\index{orbifold ! examples}
The $S^1$-action on $\CC^2$ given by
$e^{i\theta} \cdot (z_1,z_2) =
(e^{ik\theta}z_1,e^{i\ell\theta}z_2)$, for some integers $k$ and $\ell$,
has moment map $\mu: \CC^2 \to \RR$,
$(z_1,z_2) \mapsto -\frac {1}{2} (k|z_1|^2 + \ell |z_2|^2)$.
Any $\xi < 0$ is a regular value and $\mu^{-1}(\xi)$ is a
3-dimensional ellipsoid.

When $\ell = 1$ and $k \ge 2$, the stabilizer of $(z_1,z_2)$
is $\{1\}$ if $z_2 \ne 0$ and is $\ZZ_k =
\left\{ e^{i \frac {2\pi m}{k}} \mid m = 0,1,\dots,k-1 \right\}$ if $z_2 = 0$.
The reduced space $\mu^{-1}(\xi)/S^1$ is then called a
\textbf{teardrop}\index{orbifold ! teardrop}\index{teardrop
orbifold} orbifold or {\em conehead}\index{orbifold !
conehead}\index{conehead orbifold};
it has one \textbf{cone} (or {\em dunce cap}\index{orbifold !
dunce cap}\index{dunce cap orbifold})
singularity with cone angle $\frac {2\pi}{k}$,
that is, a point with orbifold structure group $\ZZ_k$.

When $k, \ell \ge 2$ are relatively prime,
for $z_1,z_2 \ne 0$ the stabilizer of $(z_1,0)$ is $\ZZ_k$,
of $(0,z_2)$ is $\ZZ_{\ell}$ and of $(z_1,z_2)$ is $\{1\}$.
The quotient $\mu^{-1}(\xi)/S^1$ is called a \textbf{football} orbifold:
it has two cone singularities,
with angles $\frac {2\pi}{k}$ and $\frac {2\pi}{\ell}$.

For $S^1$ acting on $\CC^n$ by $e^{i\theta} \cdot (z_1,\dots,z_n) =
(e^{ik_1\theta}z_1,\dots,e^{ik_n\theta}z_n)$
%with weights $k_j$ relatively prime in pairs,
the reduced spaces are orbifolds called \textbf{weighted}
(or \textbf{twisted}) \textbf{projective spaces}.\index{example !
weighted projective
space}\index{weighted projective space}\index{twisted projective space}
\end{examples}

%%%%%%%%%%%%%%%%%%%%%%%%%%%%%%%%%%%%%%%%%%%%%%%%%%%%%%%%%%%%%%%%%%%%%%%%%%%%%

Let $(M,\omega)$ be a symplectic manifold where $S^1$ acts
in a hamiltonian way, $\rho: S^1 \to {\rm Diff} (M)$,
with moment map $\mu : M \to \RR$.
Suppose that:
\begin{itemize}
\item
$M$ has a unique nondegenerate minimum at $q$
where $\mu (q) = 0$, and
\item
for $\varepsilon$ sufficiently small,
$S^1$ acts freely on the level set $\mu^{-1} (\varepsilon)$.
\end{itemize}
Let $\CC$ be equipped with the symplectic form $-i dz \wedge d \bar z$.
Then the action of $S^1$ on the product
$\psi: S^1 \to {\rm Diff} (M \times \CC)$,
$\psi_t (p,z) = (\rho_t (p) , t \cdot z)$,
is hamiltonian with moment map
\[
   \phi : M \times \CC \longrightarrow \RR \ , \qquad
   \phi (p,z) = \mu (p) - |z|^2 \ .
\]
Observe that $S^1$ acts freely on the $\varepsilon$-level
of $\phi$ for $\varepsilon$ small enough:
\[
\begin{array}{rcl}
   \phi^{-1} (\varepsilon) & = &
   \{ (p,z) \in M \times \CC \mid \mu (p) - |z|^2 = \varepsilon \} \\
   & = & \{ (p,0) \in M \times \CC \mid \mu (p) = \varepsilon \} \\
   & & \quad \cup \quad
   \{ (p,z) \in M \times \CC \mid |z|^2 = \mu (p) -\varepsilon >0\} \ .
\end{array}
\]
The reduced space is hence
\[
   \phi^{-1} (\varepsilon) / S^1 \simeq
   \mu^{-1} (\varepsilon) / S^1 \cup
   \{ p \in M \mid \mu (p) > \varepsilon \} \ .
\]
The open submanifold of $M$ given by
$\{ p \in M \mid \mu (p) > \varepsilon \}$ embeds
as an open dense submanifold into $\phi^{-1} (\varepsilon) / S^1$.
The reduced space $\phi^{-1} (\varepsilon) / S^1$
is the $\varepsilon$-blow-up of $M$ at $q$ (Section~\ref{sec:actions}).
This global description of blow-up for hamiltonian
$S^1$-spaces is due to Lerman~\cite{le:cuts},
as a particular instance of his {\em cutting} technique.
\textbf{Symplectic cutting}\index{symplectic ! cutting}
is the application of symplectic
reduction to the product of a hamiltonian $S^1$-space
with the standard $\CC$ as above, in a way that the
reduced space for the original hamiltonian $S^1$-space
embeds symplectically as a codimension 2 submanifold
in a symplectic manifold.
As it is a local construction, the cutting operation may
be more generally performed at a local minimum (or maximum)
of the moment map $\mu$.
There is a remaining $S^1$-action on the cut space
$M_{\mathrm{cut}}^{\geq \varepsilon} := \phi^{-1} (\varepsilon) / S^1$
induced by
\[
   \tau: S^1 \longrightarrow {\rm Diff} (M \times \CC) \ , \qquad
   \tau_t (p,z) = (\rho_t (p) , z) \ .
\]
In fact, $\tau$ is a hamiltonian $S^1$-action on $M \times \CC$
that commutes with $\psi$, thus descends to an action
$\widetilde \tau : S^1 \to {\rm Diff} (M_{\mathrm{cut}}^{\geq \varepsilon})$,
which is also hamiltonian.

Loosely speaking, the cutting technique provides a hamiltonian
way to close the open manifold $\{ p \in M \mid \mu (p) > \varepsilon \}$,
by using the reduced space at level $\varepsilon$,
$\mu^{-1} (\varepsilon) / S^1$.
We may similarly close $\{ p \in M \mid \mu (p) < \varepsilon \}$.
The resulting hamiltonian $S^1$-spaces are called
\textbf{cut spaces}\index{cut spaces},
and denoted $M_{\mathrm{cut}}^{\geq \varepsilon}$
and $M_{\mathrm{cut}}^{\leq \varepsilon}$.
If another group $G$ acts on $M$ in a hamiltonian way
that commutes with the $S^1$-action, then the cut spaces are
also hamiltonian $G$-spaces.

%%%%%%%%%%%%%%%%%%%%%%%%%%%%%%%%%%%%%%%%%%%%%%%%%%%%%%%%%%%%%%%%%%%%%%%%%%%%%
%%%%%%%%%%%%%%%%%%%%%%%%%%%%%%%%%%%%%%%%%%%%%%%%%%%%%%%%%%%%%%%%%%%%%%%%%%%%%

\ssubsection{Moment Map in Gauge Theory}
\index{moment map ! in gauge theory}\index{gauge ! theory}

Let $G$ be a Lie group and $P$ a principal $G$-bundle
over $B$.\footnote{Let $G$
be a Lie group and $B$ a manifold.
A \textbf{principal $G$-bundle over $B$} is a
fibration $\pi : P \to B$ (Section~\ref{sec:constructions})
with a free action of $G$ (the \textbf{structure group}\index{group !
structure}) on the total space $P$,
such that the base $B$ is the orbit space,
the map $\pi$ is the point-orbit projection and
the local trivializations are of the form 
$\varphi_\cU = (\pi,s_\cU) : \pi ^{-1} (\cU) \to \cU \times G$ with
$s_\cU (g \cdot p) = g \cdot s_\cU (p)$
for all $g \in G$ and all $p \in \pi ^{-1} (\cU)$.
A principal $G$-bundle is represented by a diagram
\[
\begin{array}{cll}
        G & \hookrightarrow & P \\
        & & \downarrow \pi \\
        & & B
\end{array}
\]
For instance, the \textbf{Hopf fibration}\index{Hopf ! fibration}
is a principal $S^1$-bundle over $S^2$($=\CC \PP^1$) with total space $S^3$
regarded as unit vectors in $\CC ^2$ where circle elements
act by complex multiplication.}
%$e^{i\theta} \cdot (z_1,z_2) = (e^{i\theta} z_1, e^{i\theta} z_2)$.
If $A$ is a connection (form)\footnote{An action
$\psi: G \to \mathrm{Diff}(P)$
induces an infinitesimal action\index{action !
infinitesimal}\index{infinitesimal action}
$d\psi: \fg \to \chi(P)$ mapping $X \in \fg$ to the
vector field $X^{\#}$ generated by the
one-parameter group $\{\exp tX(e) \mid t \in \RR\} \subseteq G$.
Fix a basis $X_1, \ldots, X_k$ of $\fg$.
Let $P$ be a principal $G$-bundle over $B$.
Since the $G$-action is free, the vector fields
$X_1^{\#}, \ldots, X_k^{\#}$ are linearly independent
at each $p \in P$.
The \textbf{vertical bundle} $V$ is the rank $k$ subbundle
of $TP$ generated by $X_1^{\#}, \ldots, X_k^{\#}$.
Alternatively, $V$ is the set of vectors
tangent to $P$ that lie in the kernel of the derivative of
the bundle projection $\pi$, so $V$ is indeed independent
of the choice of basis for $\fg$.
An \textbf{(Ehresmann) connection}\index{Ehresmann !
connection}\index{connection ! on a principal bundle}\index{principal bundle !
connection} on $P$ is a choice of a splitting $TP = V \oplus H$,
where $H$ (called the \textbf{horizontal bundle})
is a $G$-invariant subbundle of $TP$
complementary to the vertical bundle $V$.
A \textbf{connection form}\index{connection ! form}\index{form !
connection} on $P$ is a Lie-algebra-valued 1-form
$A = \sum_{i=1}^k A_i \otimes X_i \in \Omega^1 (P) \otimes \fg$
such that $A$ is $G$-invariant, with respect to
the product action of $G$ on $\Omega^1 (P)$ (induced by the
action on $P$) and on $\fg$ (the adjoint action), and
$A$ is vertical, in the sense that
$\imath_{X^{\#}} A = X$ for any $X \in \fg$.
A connection $TP = V \oplus H$ determines
a connection (form) $A$ and vice-versa by the formula
$H = \ker A = \{ v \in TP \mid \imath_v A = 0 \}$.
Given a connection on $P$, the splitting $TP = V \oplus H$
induces splittings for bundles
$T^*P = V^* \oplus H^*$, $\wedge ^2 T^*P =
(\wedge ^2 V^*) \oplus (V^* \wedge H^*) \oplus (\wedge ^2 H^*)$, etc.,
and for their sections:
$\Omega^1 (P) = \Omega^1_{\mathrm{vert}} \oplus \Omega^1_{\mathrm{horiz}}$,
$\Omega^2 (P) = \Omega^2_{\mathrm{vert}} \oplus \Omega^2_{\mathrm{mix}}
\oplus \Omega^2_{\mathrm{horiz}}$, etc.
The corresponding connection form $A$ is in
$\Omega^1_{\mathrm{vert}} \otimes \fg$.}
on $P$, and if $a \in \Omega^1_{\mathrm{horiz}} \otimes \fg$
is $G$-invariant for the product action,
then $A+a$ is also a connection on $P$.
Reciprocally, any two connections on $P$
differ by an $a \in (\Omega^1_{\mathrm{horiz}} \otimes \fg)^G$.
We conclude that the \textbf{set $\cA$ of all connections} on the
principal $G$-bundle $P$ is an affine space\index{space ! affine}
modeled on the linear space
$\fa = (\Omega^1_{\mathrm{horiz}} \otimes \fg)^G$.
\index{connection ! space}\index{symplectic ! structure on the
space of connections}\index{example !
of infinite-dimensional symplectic manifold}\index{space ! of connections}

Now let $P$ be a principal $G$-bundle over a compact
Riemann surface.\index{Riemann ! surface}
Suppose that the group $G$ is compact or semisimple\index{group !
semisimple}\index{semisimple}.
Atiyah and Bott~\cite{at-bo:surfaces}\index{Atiyah !
moduli space}\index{Bott ! moduli space} noticed that
the corresponding space $\cA$ of all connections
may be treated as an {\em infinite-dimensional
symplectic manifold}.\index{manifold !
infinite-dimensional}
This requires choosing a $G$-invariant
inner product $\langle \cdot,\cdot \rangle$ on $\fg$,
which always exists, either by averaging any
inner product when $G$ is compact, or by using the
{\em Killing form}\index{form ! Killing}\index{Killing form}
on semisimple groups.

Since $\cA$ is an affine space, its tangent space at any point $A$
is identified with the model linear space $\fa$.
With respect to a basis $X_1, \ldots, X_k$
for the Lie algebra $\fg$, elements $a,b \in \fa$ are written
\[
   a = \sum a_i \otimes X_i \quad \mbox{ and } \quad
   b = \sum b_i \otimes X_i \ .
\]
If we wedge $a$ and $b$, and then integrate over $B$, we obtain a real number:
\[
\begin{array}{rrclcl}
   \omega : & \fa \times \fa & \longrightarrow &
   \left( \Omega ^2_{\mathrm{horiz}} (P) \right) ^G \simeq \Omega ^2 (B)
   & \longrightarrow & \RR \\
   & (a,b) & \longmapsto & \sum \limits_{i,j} a_i \wedge b_j
   \langle X_i , X_j \rangle &
   \longmapsto & \int \limits_B \sum \limits_{i,j} a_i \wedge b_j
   \langle X_i , X_j \rangle \ .
\end{array}
\]
We used that the pullback $\pi^* : \Omega ^2 (B) \to \Omega ^2 (P)$
is an isomorphism onto its image
$\left( \Omega ^2_{\mathrm{horiz}} (P) \right) ^G$.
When $\omega(a,b) =0$ for all $b \in \fa$,
then $a$ must be zero.
The map $\omega$ is nondegenerate,
skew-symmetric, bilinear and constant in the sense that
it does not depend on the base point $A$.
Therefore, it has the right to be called a
symplectic form on $\cA$, so the pair $(\cA, \omega)$
is an \textbf{infinite-dimensional symplectic manifold}.

%%%%%%%%%%%%%%%%%%%%%%%%%%%%%%%%%%%%%%%%%%%%%%%%%%%%%%%%%%%%%%%%%%%%%%%%%%%%%

A diffeomorphism $f:P \to P$ commuting with the $G$-action
determines a diffeomorphism $f_{\mathrm{basic}} : B \to B$
by projection.
Such a diffeomorphism $f$ is called a
\textbf{gauge transformation}\index{gauge ! transformation}
if the induced $f_{\mathrm{basic}}$ is the identity.
The \textbf{gauge group}\index{group ! gauge}\index{gauge ! group}
of $P$ is the group $\cG$ of all gauge transformations of $P$.

The derivative of an $f \in \cG$ takes an {\em Ehresmann connection}
$TP = V \oplus H$ to another connection $TP = V \oplus H_f$,
and thus induces an action of $\cG$ in the space $\cA$ of
all connections.
Atiyah and Bott~\cite{at-bo:surfaces} noticed that
the action of $\cG$\index{action ! gauge group}\index{principal bundle !
gauge group} on $(\cA , \omega)$ is hamiltonian,
where the moment map (appropriately interpreted) is
\[
\begin{array}{rrcl}
   \mu: & \cA & \longrightarrow &
   \left( \Omega ^2 (P) \otimes \fg \right) ^G \\
   & A & \longmapsto & \curv \ A \ ,
\end{array}
\]
i.e., the moment map {\em is} the curvature.\footnote{The
exterior derivative of a connection $A$ decomposes into three components,
\[
   dA = (dA)_{\mathrm{vert}} +
   (dA)_{\mathrm{mix}} + (dA)_{\mathrm{horiz}}
   \in \left( \Omega^2_{\mathrm{vert}}
   \oplus \Omega^2_{\mathrm{mix}}
   \oplus \Omega^2_{\mathrm{horiz}} \right) \otimes \fg
\]
satisfying $(dA)_{\mathrm{mix}} = 0$ and
$(dA)_{\mathrm{vert}} (X,Y) = [X,Y]$, i.e.,
$(dA)_{\mathrm{vert}} = \frac12 \sum_{i,\ell,m}
c_{\ell m}^i A_\ell \wedge A_m \otimes X_i$,
where the $c_{\ell m}^i$'s are the \textbf{structure constants}
of the Lie algebra with respect to the chosen basis, and defined by
$[X_\ell, X_m] = \sum_{i,\ell,m} c_{\ell m}^i X_i$.
So the relevance of $dA$ may come only from its horizontal component,
called the \textbf{curvature form}\index{curvature form}\index{form !
curvature} of the connection $A$, and denoted
$\curv \ A = (dA)_{\mathrm{horiz}} \in
\Omega^2_{\mathrm{horiz}} \otimes \fg$.
A connection is called \textbf{flat} if its curvature is zero.}
The reduced space $\cM = \mu ^{-1} (0) / \cG$
is the space of {\em flat connections} modulo gauge equivalence,
known as the \textbf{moduli space of flat
connections},\index{connection ! flat}\index{flat
connection}\index{connection !
moduli space}\index{moduli space}\index{space ! moduli}
which is a finite-dimensional symplectic orbifold.

%%%%%%%%%%%%%%%%%%%%%%%%%%%%%%%%%%%%%%%%%%%%%%%%%%%%%%%%%%%%%%%%%%%%%%%%%%%%%

\begin{example}
We describe the Atiyah-Bott construction for the
case of a circle bundle\index{circle bundle}
\[
\begin{array}{cll}
        S^1 & \hookrightarrow & P \\
        & & \downarrow \pi \\
        & & B
\end{array}
\]
Let $v$ be the generator of the $S^1$-action on $P$,
corresponding to the basis $1$ of $\fg \simeq \RR$.
A connection form on $P$ is an ordinary 1-form
$A \in \Omega ^1 (P)$ such that
$\cL _v A = 0$ and $\imath_v A = 1$.
If we fix one particular connection $A_0$,
then any other connection is of the form $A = A_0 + a$
for some $a \in \fa = \left( \Omega^1_{\mathrm{horiz}} (P) \right)^G
= \Omega^1 (B)$.
The symplectic form on $\fa = \Omega^1 (B)$ is simply
\[
\begin{array}{rrcccl}
   \omega : & \fa \times \fa & \longrightarrow &
   \Omega ^2 (B) & \longrightarrow & \RR \\
   & (a,b) & \longmapsto & a \wedge b & \longmapsto & \int_B a \wedge b \ .
\end{array}
\]
The gauge group is $\cG = \mathrm{Maps}(B,S^1)$,
because a gauge transformation is multiplication
by some element of $S^1$ over each point in $B$
encoded in a map $h : B \to S^1$.
The action $\phi : \cG \to \mathrm{Diff} (P)$
takes $h \in \cG$ to the diffeomorphism
\[
\begin{array}{rrcl}
   \phi_h : & p & \longmapsto & h(\pi(p)) \cdot p \ .
\end{array}
\]
The Lie algebra of $\cG$ is
$\mathrm{Lie} \ \cG = \mathrm{Maps}(B,\RR) = C^\infty (B)$
with dual $\left( \mathrm{Lie} \ \cG \right) ^* = \Omega ^2 (B)$,
where the (smooth) duality is provided by integration
$C^\infty (B) \times \Omega ^2 (B) \to \RR$,
$(h,\beta) \mapsto \int_B h \beta$.
The gauge group acts on the space of all connections by
\[
\begin{array}{rrcl}
   \psi & \cG & \longrightarrow & \mathrm{Diff} (\cA) \\
   & (h: x \mapsto e^{i\theta(x)}) & \longmapsto &
   ( \psi_h: A \mapsto A - \pi^* d\theta )
\end{array}
\]
(In the case where $P = S^1 \times B$
is a trivial bundle, every connection can be written
$A = dt + \beta$, with $\beta \in \Omega ^1 (B)$.
A gauge transformation $h \in \cG$ acts on $P$ by
$\phi_h : (t,x) \mapsto (t + \theta (x) , x)$
and on $\cA$ by $A \mapsto \phi^*_{h^{-1}} (A)$.)
The infinitesimal action is
\[
\begin{array}{rrcl}
   d\psi: & \mathrm{Lie} \ \cG & \longrightarrow & \chi(\cA) \\
   & X & \longmapsto & X^{\#} = \mbox{ vector field described by }
   ( A \mapsto A -dX )
\end{array}
\]
so that $X^{\#} = -dX$.
It remains to check that
\[
\begin{array}{rrcl}
   \mu: & \cA & \longrightarrow &
   \left( \mathrm{Lie} \ \cG \right) ^* = \Omega ^2 (B) \\
   & A & \longmapsto & \curv \ A
\end{array}
\]
is indeed a moment map for the action of the gauge group on $\cA$.
Since in this case $\curv \ A = dA \in
\left( \Omega_{\mathrm{horiz}} ^2 (P) \right) ^G = \Omega ^2 (B)$,
the action of $\cG$ on $\Omega ^2 (B)$ is trivial
and $\mu$ is $\cG$-invariant,
the equivariance condition is satisfied.
Take any $X \in \mathrm{Lie} \ \cG = C^\infty (B)$.
Since the map $\mu ^X: A \mapsto \langle X , dA \rangle = \int_B X \cdot dA$
is linear in $A$, its differential is
\[
\begin{array}{rrcl}
   d\mu ^X: & \fa & \longrightarrow & \RR \\
   & a & \longmapsto & \int_B X da \ .
\end{array}
\]
By definition of $\omega$ and the Stokes theorem,\index{Stokes
theorem}\index{theorem ! Stokes} we have that
\[
   \omega (X^{\#} , a) = \displaystyle{\int_B} X^{\#} \cdot a
   = - \displaystyle{\int_B} dX \cdot a
   = \displaystyle{\int_B} X \cdot da
   = d \mu ^X (a) \ , \qquad \forall a \in \Omega ^1 (B) \ .
\]
so we are done in proving that $\mu$ is the moment map.
\end{example}

The function $||\mu||^2 : \cA \to \RR$ giving the square
of the $L^2$ norm of the curvature
is the \textbf{Yang-Mills functional}, whose
Euler-Lagrange equations are the {\em Yang-Mills equations}.
Atiyah and Bott~\cite{at-bo:surfaces} studied the topology of $\cA$
by regarding $||\mu||^2$ as an equivariant Morse function.
In general, it is a good idea to apply Morse theory
to the norm square of a moment map~\cite{ki:quotients}.

%%%%%%%%%%%%%%%%%%%%%%%%%%%%%%%%%%%%%%%%%%%%%%%%%%%%%%%%%%%%%%%%%%%%%%%%%%%%%
%%%%%%%%%%%%%%%%%%%%%%%%%%%%%%%%%%%%%%%%%%%%%%%%%%%%%%%%%%%%%%%%%%%%%%%%%%%%%

\ssubsection{Symplectic Toric Manifolds}
\label{sec:stm}
\index{symplectic ! toric manifolds}\index{toric manifolds}

Toric manifolds are smooth {\em toric varieties}.\footnote{Toric
varieties were introduced by Demazure in~\cite{de:subgroups}.
There are many nice surveys of the theory of
toric varieties in algebraic geometry;
see, for instance,~\cite{da:toric,fu:toric,ke-kn-mu-sa:toroidal,od:toric}.
Toric geometry has recently become an
important tool in physics in connection with
mirror symmetry~\cite{co:recent}.}
When studying the symplectic features of these spaces,
we refer to them as {\em symplectic toric manifolds}.
Relations between the algebraic and symplectic viewpoints
on toric manifolds are discussed in~\cite{ca:toric}.

\begin{definition}
A \textbf{symplectic toric manifold} is a compact
connected symplectic manifold $(M,\omega)$ equipped
with an effective hamiltonian action of a torus $\TT$
of dimension equal to half the dimension of the manifold,
$\dim \TT = \frac {1}{2} \dim M$,
and with a choice of a corresponding moment map $\mu$.
Two symplectic toric manifolds,
$(M_i,\omega_i,\TT_i,\mu_i)$, $i=1,2$,
are \textbf{equivalent} if there exists an isomorphism
$\lambda : \TT_1 \to \TT_2$ and a $\lambda$-equivariant
symplectomorphism $\varphi : M_1 \to M_2$ such that
$\mu_1 = \mu_2 \circ \varphi$.
\end{definition}

\vspace*{-1ex}

%Equivalent symplectic toric manifolds are often not distinguished.

\begin{examples}
\begin{enumerate}

\item
The circle $S^1$ acts on the 2-sphere
$(S^2,\omega_{\mathrm{standard}} = d\theta \wedge dh)$ by rotations,
$e^{i \nu} \cdot (\theta, h) = (\theta + \nu ,h)$.
with moment map $\mu = h$ equal to the height function
and moment polytope $[-1,1]$.

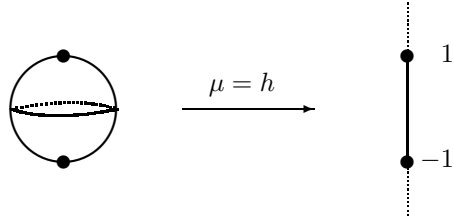
\begin{picture}(200,100)(-40,0)
%\put(180,10){\line(0,1){80}}
\qbezier[50](180,10)(180,50)(180,90)
\put(95,50){\vector(1,0){50}}
\put(105,57){$\mu=h$}
\put(185,28){$-1$}
\put(185,68){$\phantom{-}1$}
\thicklines
\put(180,30){\line(0,1){40}}
\put(50,50){\circle{50}}
\qbezier(30,50)(45,45)(70,50)
\qbezier[25](30,50)(55,55)(70,50)
\put(180,30){\circle*{5}}
\put(180,70){\circle*{5}}
\put(50,30){\circle*{5}}
\put(50,70){\circle*{5}}
\end{picture}

Analogously, $S^1$ acts on
the Riemann sphere $\CC \PP^1$ with the Fubini-Study form
$\omega_{\mathrm{FS}} = \frac {1}{4} \omega_{\mathrm{standard}}$, by
$e ^{i\theta} \cdot [z_0,z_1] = [z_0,e ^{i\theta}z_1]$.
This is hamiltonian with moment map
$\mu[z_0,z_1] = -\frac {1}{2} \cdot \frac {|z_1|^2}{|z_0|^2 + |z_1|^2}$,
and moment polytope $\left[ -\frac {1}{2} ,0\right]$.

\item
For the $\TT^n$-action on the product of $n$ Riemann spheres
$\CC \PP^1 \times \ldots \times \CC \PP^1$ by
\[
   (e ^{i\theta_1}, \ldots ,e ^{i\theta_n})
   \cdot ([z_1,w_1],\ldots , [z_n,w_n]) =
   ([z_1,e ^{i\theta_1}w_1],\ldots , [w_0,e ^{i\theta_n}w_1]) \ ,
\]
the moment polytope is an $n$-dimensional cube.

\item
   Let $(\CC \PP^2, \omega_{\mathrm{FS}})$ be 2-(complex-)dimensional
complex projective space equipped with the Fubini-Study form
defined in Section~\ref{sec:kahler}.
The $\TT^2$-action on $\CC \PP^2$ by
$(e ^{i\theta_1},e ^{i\theta_2}) \cdot [z_0,z_1,z_2] =
[z_0,e ^{i\theta_1}z_1,e ^{i\theta_2}z_2]$
has moment map
\[
   \mu[z_0,z_1,z_2] = -\frac {1}{2} \left( \frac
   {|z_1|^2}{|z_0|^2+|z_1|^2+|z_2|^2} ,
   \frac {|z_2|^2}{|z_0|^2+|z_1|^2+|z_2|^2} \right) \ .
\]
The image is the isosceles triangle with vertices
$(0,0)$, $( -\frac {1}{2} ,0)$ and $(0, -\frac {1}{2})$.

\item
For the $\TT^n$-action on $(\CC \PP^n,\omega_{\mathrm{FS}})$ by
\[
   (e ^{i\theta_1}, \ldots ,e ^{i\theta_n})
   \cdot [z_0,z_1,\ldots,z_n] = [z_0,e ^{i\theta_1}z_1,
   \ldots , e ^{i\theta_n}z_n]
\]
the moment polytope is an $n$-dimensional simplex.
\end{enumerate}
\end{examples}

Since the coordinates of the moment map are
commuting integrals of motion,
a symplectic toric manifold gives rise to
a completely integrable system\index{integrable ! system}.
By Proposition~\ref{prop:dimension}, symplectic
toric manifolds are optimal hamiltonian torus-spaces.
By Theorem~\ref{thm:convexity}, they have an associated polytope.
It turns out that the moment polytope contains enough
information to sort all symplectic toric manifolds.
We now define the class of polytopes
that arise in the classification.
For a symplectic toric manifold
the weights $\lambda^{(1)}, \ldots, \lambda^{(n)}$
in Theorem~\ref{thm:moment_darboux} form a $\ZZ$-basis of $\ZZ^m$,
hence the moment polytope is a {\em Delzant polytope}:

\begin{definition}
A \textbf{Delzant polytope}\index{polytope !
Delzant}\index{Delzant ! polytope}
in $\RR^n$ is a polytope satisfying:
\begin{itemize}
\item
\textbf{simplicity}\index{simple polytope}\index{polytope !
simple}, i.e., there are $n$ edges meeting at each vertex;
\item
\textbf{rationality}\index{rational polytope}\index{polytope !
rational}, i.e., the edges meeting at the vertex $p$
are rational in the sense that each edge is of the form $p + tu_i$,
$t \ge 0$, where $u_i \in \ZZ^n$;
\item
\textbf{smoothness}\index{smooth polytope}\index{polytope !
smooth}, i.e., for each vertex, the corresponding $u_1,\dots,u_n$
can be chosen to be a $\ZZ$-basis of $\ZZ^n$.
\end{itemize}
\end{definition}

In $\RR^2$ the simplicity condition is always satisfied
(by nondegenerate polytopes).
In $\RR^3$, for instance a square pyramid fails the simplicity condition.

\begin{examples}
The pictures below represent Delzant polytopes in $\RR^2$.\index{example !
of Delzant polytope}\index{Delzant ! example
of Delzant polytope}\index{polytope ! example
of Delzant polytope}

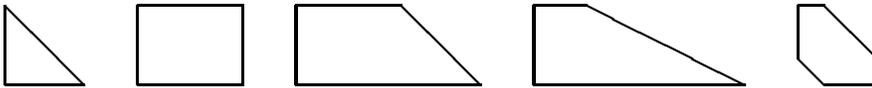
\begin{picture}(250,60)(5,-20)
\thicklines
% triangle
\put(0,0){\line(1,0){30}}
\put(0,0){\line(0,1){30}}
\put(0,30){\line(1,-1){30}}
% rectangle
\put(50,0){\line(1,0){40}}
\put(50,0){\line(0,1){30}}
\put(90,30){\line(-1,0){40}}
\put(90,30){\line(0,-1){30}}
% hirzebruch
\put(110,0){\line(1,0){70}}
\put(110,0){\line(0,1){30}}
\put(150,30){\line(1,-1){30}}
\put(110,30){\line(1,0){40}}
% hirzebruch
\put(200,0){\line(1,0){80}}
\put(200,0){\line(0,1){30}}
\put(200,30){\line(1,0){20}}
\put(220,30){\line(2,-1){60}}
% blow-up of $\CC \PP^2$ at the three fixed points
\put(310,0){\line(1,0){20}}
\put(310,0){\line(-1,1){10}}
\put(300,30){\line(1,0){10}}
\put(300,10){\line(0,1){20}}
\put(330,10){\line(-1,1){20}}
\put(330,0){\line(0,1){10}}
\end{picture}
\end{examples}

The following theorem
classifies (equivalence classes of) symplectic toric manifolds
in terms of the combinatorial data encoded by a Delzant polytope.

\begin{theorem}\index{Delzant ! theorem}\index{theorem ! Delzant}
\textbf{(Delzant~\cite{de:hamiltoniens})} $\;$
Toric manifolds are classified by Delzant polytopes,
and their bijective correspondence is given by the moment map:
\[
\begin{array}{rcl}
        \{\mbox{toric manifolds}\} & \longleftrightarrow
        &\{\mbox{Delzant polytopes}\}
        \\
        (M^{2n},\omega,\TT^n,\mu) &\longmapsto &\mu(M) \ .
\end{array}
\]
\end{theorem}

Delzant's construction (Section~\ref{sec:construction})
shows that for a toric manifold the moment map
takes the fixed points bijectively to the vertices of the moment polytope
and takes points with a $k$-dimensional stabilizer to the
codimension $k$ faces of the polytope.
The moment polytope is exactly the orbit space,
i.e., the preimage under $\mu$ of each point in the polytope
is exactly one orbit.
For instance, consider $(S^2,\omega=d\theta \wedge dh, S^1, \mu = h)$,
where $S^1$ acts by rotation.
The image of $\mu$ is the line segment $I = [-1,1]$.
The product $S^1 \times I$ is an open-ended cylinder.
We can recover the 2-sphere
by collapsing each end of the cylinder to a point.
Similarly, we can build $\CC \PP^2$ from 
$\TT^2 \times \Delta$
where $\Delta$ is a rectangular isosceles triangle, and so on.

\begin{examples}
\begin{enumerate}
\item
By a linear transformation in $\SL (2;\ZZ)$,
we can make one of the angles in a Delzant triangle into a right angle.
Out of the rectangular triangles, only the isosceles one
satisfies the smoothness condition.
Therefore, up to translation, change of scale and the action of $\SL (2;\ZZ)$,
there is just one 2-dimensional Delzant polytope\index{polytope !
Delzant}\index{Delzant ! polytope} with three vertices,
namely an {\em isosceles triangle}.
We conclude that the projective space $\CC \PP^2$
is the only 4-dimensional toric manifold with three fixed points,
up to choices of a constant in the moment map,
of a multiple of $\omega_{_{\mathrm{FS}}}$
and of a lattice basis in the Lie algebra of $\TT^2$.

\item
Up to translation, change of scale and the action of $\SL (n;\ZZ)$,
the {\em standard $n$-simplex} $\Delta$ in $\RR^n$ (spanned by the origin
and the standard basis vectors $(1,0,\ldots,0),\ldots,(0,\ldots,0,1)$)
is the only $n$-dimensional Delzant polytope\index{polytope !
Delzant}\index{Delzant ! polytope} with $n+1$ vertices.
Hence, $M_\Delta = \CC \PP^n$ is the only $2n$-dimensional toric manifold
with $n+1$ fixed points,
up to choices of a constant in the moment map,
of a multiple of $\omega_{_{\mathrm{FS}}}$
and of a lattice basis in the Lie algebra of $\TT^N$.

\item
A transformation in $\SL (2;\ZZ)$
makes one of the angles in a Delzant quadrilateral into a right angle.
Automatically an adjacent angle also becomes $90^o$.
Smoothness imposes that the slope of the skew side be integral.
Thus, up to translation, change of scale and $\SL (2;\ZZ)$-action,
the 2-dimensional Delzant polytopes\index{polytope !
Delzant}\index{Delzant ! polytope} with four vertices are
trapezoids with vertices
$(0,0)$, $(0,1)$, $(\ell,1)$ and $(\ell +n,0)$,
for $n$ a nonnegative integer and $\ell > 0$.
Under Delzant's construction (that is, under
symplectic reduction of $\CC^4$ with respect to an action of $(S^1)^2$),
these correspond to
the so-called \textbf{Hirzebruch surfaces}\index{example !
Hirzebruch surfaces}\index{Hirzebruch surface}
-- the only 4-dimensional symplectic toric
manifolds\index{toric manifold ! four-dimensional@4-dimensional}
that have four fixed points up to equivalence as before.
Topologically, they are $S^2$-bundles over $S^2$,
either the trivial bundle $S^2 \times S^2$
when $n$ is even or the nontrivial bundle
(given by the blow-up of $\CC \PP^2$ at a point;
see Section~\ref{sec:blow_up}) when $n$ is odd.
\end{enumerate}
\end{examples}

Let $\Delta$ be an $n$-dimensional Delzant polytope,
and let $(M_\Delta,\omega_\Delta, \TT^n, \mu_\Delta)$ be the
associated symplectic toric manifold.
The $\varepsilon$-blow-up of $(M_\Delta,\omega_\Delta)$
at a fixed point of the $\TT^n$-action is a new symplectic
toric manifold (Sections~\ref{sec:blow_up} and~\ref{sec:actions}).
Let $q$ be a fixed point of the $\TT^n$-action
on $(M_\Delta,\omega_\Delta)$, and let $p = \mu_\Delta (q)$
be the corresponding vertex of $\Delta$.
Let $u_1, \ldots ,u_n$ be the primitive (inward-pointing)
edge vectors at $p$, so that the rays $p + t u_i$, $t \geq 0$,
form the edges of $\Delta$ at $p$.

\begin{proposition}
\index{blow-up ! of toric manifold}
The $\varepsilon$-blow-up of $(M_\Delta,\omega_\Delta)$
at a fixed point $q$ is the symplectic toric manifold
associated to the polytope $\Delta_\varepsilon$
obtained from $\Delta$ by replacing the vertex $p$
by the $n$ vertices $p + \varepsilon u_i$, $i=1,\ldots , n$.
\end{proposition}

In other words, the moment polytope for the blow-up of
$(M_\Delta,\omega_\Delta)$ at $q$ is obtained from $\Delta$
by chopping off the corner corresponding to $q$, thus substituting
the original set of vertices by the same set with the vertex
corresponding to $q$ replaced by exactly $n$ new vertices.
The truncated polytope is Delzant.
We may view the $\varepsilon$-blow-up of $(M_\Delta, \omega_\Delta)$
as being obtained from $M_\Delta$ by smoothly replacing $q$ by
$(\CC \PP^{n-1}, \varepsilon \omega_{_{\mathrm{FS}}})$
(whose moment polytope is an $(n-1)$-dimensional simplex).

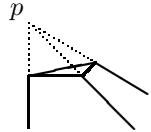
\begin{picture}(200,70)(-95,-10)
\put(43,43){$p$}
\qbezier[10](50,20)(50,30)(50,40)
\qbezier[15](70,20)(60,30)(50,40)
\qbezier[15](75,25)(65,31)(50,40)
\thicklines
\put(50,20){\line(5,1){25}}
\put(50,20){\line(1,0){20}}
%\put(70,20){\line(1,1){5}}
\qbezier[15](70,20)(72,22)(75,25)
\put(50,0){\line(0,1){20}}
\put(90,0){\line(-1,1){20}}
\put(75,25){\line(5,-3){20}}
\end{picture}

\begin{example}
The moment polytope for the standard $\TT^2$-action
on $(\CC \PP^2, \omega_{_{\mathrm{FS}}})$ is a right
isosceles triangle $\Delta$.
If we blow-up $\CC \PP^2$ at $[0:0:1]$ we obtain
a symplectic toric manifold associated to the trapezoid below:
a {\em Hirzebruch surface}.

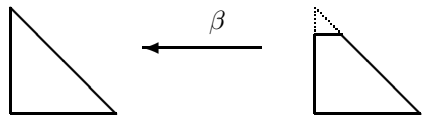
\begin{picture}(200,60)(-30,-5)
\qbezier[7](165,30)(165,35)(165,40)
\qbezier[10](175,30)(170,35)(165,40)
\thicklines
\put(50,0){\line(1,0){40}}
\put(50,0){\line(0,1){40}}
\put(90,0){\line(-1,1){40}}
\put(145,25){\vector(-1,0){45}}
\put(125,32){$\beta$}
\put(165,0){\line(1,0){40}}
\put(165,30){\line(1,0){10}}
\put(165,0){\line(0,1){30}}
\put(205,0){\line(-1,1){30}}
\end{picture}

\end{example}

Let $(M, \omega, \TT^n , \mu)$ be a $2n$-dimensional symplectic
toric manifold.
Choose a suitably generic direction in $\RR^n$ by picking a
vector $X$ whose components are independent over $\QQ$.
This condition ensures that:
\begin{itemize}
\item
the one-dimensional subgroup $\TT^X$ generated
by the vector $X$ is dense in $\TT^n$,
\item
$X$ is not parallel to the facets of the moment polytope
$\Delta := \mu (M)$, and
\item
the vertices of $\Delta$ have different projections along $X$.
\end{itemize}

Then the fixed points for the $\TT^n$-action are exactly
the fixed points of the action restricted to $\TT^X$, that is,
are the zeros of the vector field, $X^\#$ on $M$ generated by $X$.
The projection of $\mu$ along $X$,
$\mu ^X := \langle \mu , X \rangle : M \to \RR$,
is a hamiltonian function for the vector field $X^\#$ generated by $X$.
We conclude that the critical points
of $\mu ^X$ are precisely the fixed points of the $\TT^n$-action.

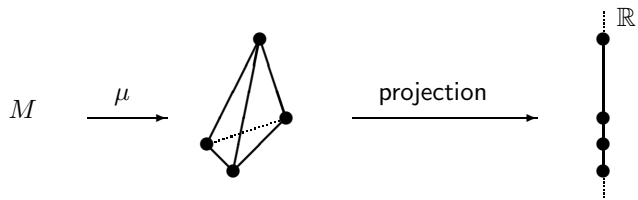
\begin{picture}(200,90)(-5,10)
\put(25,50){$M$}
\put(55,50){\vector(1,0){30}}
\put(65,57){$\mu$}
\put(155,50){\vector(1,0){70}}
\put(165,57){projection}
\qbezier[20](100,40)(115,45)(130,50)
\qbezier[50](250,20)(250,55)(250,90)
\put(255,85){$\RR$}
\put(100,40){\circle*{5}}
\put(110,30){\circle*{5}}
\put(120,80){\circle*{5}}
\put(130,50){\circle*{5}}
\put(250,30){\circle*{5}}
\put(250,40){\circle*{5}}
\put(250,50){\circle*{5}}
\put(250,80){\circle*{5}}
\thicklines
\put(100,40){\line(1,-1){10}}
\put(100,40){\line(1,2){20}}
\put(110,30){\line(1,5){10}}
\put(110,30){\line(1,1){20}}
\put(120,80){\line(1,-3){10}}
\put(250,30){\line(0,1){50}}
\end{picture}

By Theorem~\ref{thm:moment_darboux},
if $q$ is a fixed point for the $\TT^n$-action,
then there exists a chart $(\cU,x_1,\dots,x_n,y_1,\dots,y_n)$
centered at $q$ and weights
$\lambda^{(1)}, \ldots, \lambda^{(n)} \in \ZZ^n$ such that
\[
   \left. \mu ^X \right|_\cU =
   \left. \langle \mu , X \rangle \right|_\cU =
   \mu ^X (q) - \frac 12 \sum\limits_{k=1}^n
   \langle \lambda^{(k)} , X \rangle (x_k^2 + y_k^2) \ .
\]
Since the components of $X$ are independent over $\QQ$,
all coefficients $\langle \lambda^{(k)} , X \rangle$
are nonzero, so $q$ is a {\em nondegenerate}\index{nondegenerate
critical point} critical point of $\mu^X$.
Moreover, the {\em index}\footnote{A
\textbf{Morse function}\index{Morse ! function}
on an $m$-dimensional manifold $M$ is a smooth function $f: M \to \RR$
all of whose critical points (where $df$ vanishes)
are nondegenerate\index{nondegenerate critical point}
(i.e., the {\em hessian matrix}\index{hessian} is nonsingular).
Let $q$ be a nondegenerate critical point for $f : M \to \RR$.
The \textbf{index of $f$ at $q$} is the index of
the hessian $H_q : \RR^m \times \RR^m \to \RR$
regarded as a symmetric bilinear function,
that is, the the maximal dimension of a subspace of $\RR$ where
$H$ is negative definite.}
of $q$ is twice the number of labels $k$ such that
$- \langle \lambda^{(k)} , X \rangle < 0$.
But the $-\lambda^{(k)}$'s are precisely the
edge vectors $u_i$ which satisfy Delzant's conditions.
Therefore, geometrically, the index of $q$
can be read from the moment polytope $\Delta$,
by taking twice the number of edges whose
inward-pointing edge vectors at $\mu (q)$
{\em point up relative to $X$}, that is,
whose inner product with $X$ is positive.
In particular, $\mu ^X$ is a {\em perfect Morse
function}\footnote{A \textbf{perfect Morse function}\index{perfect
Morse function} is a Morse function $f$ for which the
{\em Morse inequalities}~\cite{mi:morse,mo:calculus}
are equalities, i.e., $b_\lambda (M) = C_\lambda$ and
$b_\lambda (M) - b_{\lambda -1} (M) + \ldots \pm b_0 (M) =
C_\lambda - C_{\lambda -1} + \ldots \pm C_0$
where $b_\lambda (M) = \dim H_\lambda (M)$
and $C_\lambda$ be the number of critical points
of $f$ with index $\lambda$.
If all critical points of a Morse function $f$
have even index, then $f$ is a perfect Morse function.}
and we have:

\begin{proposition}
Let $X\in \RR^n$ have components independent over $\QQ$.
The degree-$2k$ homology group of the symplectic toric
manifold $(M, \omega, \TT , \mu)$ has dimension equal to the number
of vertices of the moment polytope where there are exactly
$k$ (primitive inward-pointing) edge vectors that point up
relative to the projection along the $X$.
All odd-degree homology groups of $M$ are zero.
\end{proposition}

By Poincar\'e duality (or by taking $-X$ instead of $X$),
the words {\em point up} may be replaced by {\em point down}.
The Euler characteristic of a symplectic toric manifold
is simply the number of vertices of the corresponding polytope.
There is a combinatorial way of understanding the
cohomology ring~\cite{fu:toric}.

A \textbf{symplectic toric orbifold}\index{symplectic !
toric orbifold}\index{orbifold ! toric}
is a compact connected symplectic orbifold $(M,\omega)$ equipped
with an effective hamiltonian action of a torus
of dimension equal to half the dimension of the orbifold,
and with a choice of a corresponding moment map.
Symplectic toric orbifolds were classified by Lerman and
Tolman~\cite{le-to} in a theorem that generalizes Delzant's:
a symplectic toric orbifold is determined by its moment polytope
plus a positive integer label attached to each of the polytope facets.
The polytopes that occur are more general than the Delzant polytopes
in the sense that only simplicity and rationality
are required; the edge vectors $u_1,\dots,u_n$
need only form a rational basis of $\ZZ^n$.
When the integer labels are all equal to 1,
the failure of the polytope smoothness accounts
for all orbifold singularities.

%%%%%%%%%%%%%%%%%%%%%%%%%%%%%%%%%%%%%%%%%%%%%%%%%%%%%%%%%%%%%%%%%%%%%%%%%%%%%
%%%%%%%%%%%%%%%%%%%%%%%%%%%%%%%%%%%%%%%%%%%%%%%%%%%%%%%%%%%%%%%%%%%%%%%%%%%%%

\ssubsection{Delzant's Construction}
\index{Delzant ! construction}
\label{sec:construction}

Following~\cite{de:hamiltoniens,gu:moment},
we prove the existence part
(or surjectivity) in Delzant's theorem,\index{Delzant ! construction}
by using symplectic reduction to associate to an $n$-dimensional
Delzant polytope $\Delta$ a symplectic toric manifold
$(M_\Delta,\omega_\Delta,\TT^n,\mu_\Delta)$.

Let $\Delta$ be a Delzant polytope in
$(\RR^n)^*$\footnote{Although we identify $\RR^n$ with its dual
via the euclidean inner product, it may be more clear
to see $\Delta$ in $(\RR^n)^*$ for Delzant's construction.}
and with $d$ facets.\footnote{A \textbf{face} of a polytope
$\Delta$ is a set of the form
$F = P \cap \{ x \in \RR^n \mid f(x) = c \}$
where $c \in \RR$ and $f \in (\RR^n)^*$
satisfies $f(x) \geq c$, $\forall x \in P$.
A \textbf{facet}\index{facet}\index{polytope !
facet} of an $n$-dimensional polytope is an $(n-1)$-dimensional face.}
We can algebraically describe $\Delta$
as an intersection of $d$ halfspaces.
Let $v_i \in \ZZ^n$, $i = 1,\dots,d$, be the
primitive\footnote{A lattice vector $v \in \ZZ^n$ is
\textbf{primitive}\index{primitive
vector} if it cannot be written as $v = ku$
with $u \in \ZZ^n$, $k \in \ZZ$ and $|k| > 1$; for instance, $(1,1)$, $(4,3)$,
$(1,0)$ are primitive, but $(2,2)$, $(3,6)$ are not.}
outward-pointing normal vectors to the facets of $\Delta$.
Then, for some $\lambda_i \in \RR$, we can write
$\Delta = \{x \in (\RR^n)^* \mid \langle x,v_i\rangle
\le \lambda_i,\ i = 1,\dots,d\}$.

\begin{example}
When $\Delta$ is the triangle below, we have
\[
   \Delta = \{x \in (\RR^2)^* \mid \langle x,(-1,0)\rangle \le 0 \ ,
   \ \langle x,(0,-1) \rangle \le 0 \ , \ \langle x,(1,1)\rangle \le 1\} \ .
\]

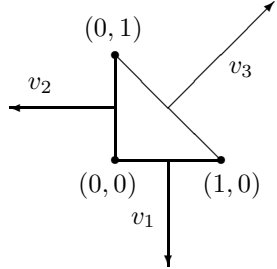
\begin{picture}(100,80)(-140,-25)
\put(0,0){\line(1,0){40}}
\put(0,0){\line(0,1){40}}
\put(0,40){\line(1,-1){40}}
% labels
\put(-14,-12){$(0,0)$}
\put(33,-12){$(1,0)$}
\put(-12,47){$(0,1)$}
%\put(80,67){$\overrightarrow{(1,1)}$}
%\put(43,-60){$\overrightarrow{(0,-1)}$}
%\put(-67,27){$\overrightarrow{(-1,0)}$}
\put(43,33){$v_3$}
\put(6,-24){$v_1$}
\put(-33,27){$v_2$}
% dots
\put(0,0){\circle*{3}}
\put(0,40){\circle*{3}}
\put(40,0){\circle*{3}}
% arrows
\put(0,20){\vector(-1,0){40}}
\put(20,0){\vector(0,-1){40}}
\put(20,20){\vector(1,1){40}}
\end{picture}

\end{example}

For the  standard basis
$e_1 = (1,0,\dots,0),\dots,e_d = (0,\dots,0,1)$ of $\RR^d$, consider
\[
\begin{array}{rrcl}
        \pi: &\RR^d &\longrightarrow &\RR^n \\
        &e_i &\longmapsto &v_i \ .
\end{array}
\]

\begin{lemma}
\label{le:onto}
The map $\pi$ is onto and maps $\ZZ^d$ onto $\ZZ^n$.
\end{lemma} 

\vspace*{-2ex}

\begin{proof}
We need to show that the set $\{v_1,\dots,v_d\}$ spans $\ZZ^n$.
At a vertex $p$, the edge vectors $u_1,\dots,u_n \in (\RR^n)^*$
form a basis for $(\ZZ^n)^*$ which, by a change of basis
if necessary, we may assume is the standard basis.
Then the corresponding primitive normal vectors to the facets
meeting at $p$ are $-u_1,\dots,-u_n$.
\end{proof}

We still call $\pi$ the induced surjective map
$\TT^d = \RR^d/(2\pi \ZZ^d) \stackrel{\pi}{\rightarrow}
\TT^n = \RR^n/(2\pi \ZZ^n)$.
The kernel $N$ of $\pi$ is a $(d-n)$-dimensional Lie
subgroup of $\TT^d$ with inclusion $i : N \hookrightarrow \TT^d$.
Let $\fn$ be the Lie algebra of $N$.
The exact sequence of tori
\[
   1 \longrightarrow N \stackrel{i}{\longrightarrow} \TT^d
   \stackrel{\pi}{\longrightarrow} \TT^n \longrightarrow 1
\]
induces an exact sequence of Lie algebras
\[
        0 \longrightarrow \fn \stackrel{i}{\longrightarrow} \RR^d
        \stackrel{\pi}{\longrightarrow} \RR^n \longrightarrow 0
\]
with dual exact sequence
\[
        0 \longrightarrow (\RR^n)^* \stackrel{\pi^*}{\longrightarrow}
        (\RR^d)^*
        \stackrel{i^*}{\longrightarrow} \fn^* \longrightarrow 0 \ .
\]
Consider $\CC^d$ with symplectic form
$\omega_0 = \frac {i}{2} \sum dz_k \wedge
d{\bar z}_k$, and standard hamiltonian action of $\TT^d$ given by
$(e^{i t_1},\dots,e^{i t_d}) \cdot (z_1,\dots,z_d) =
(e^{i t_1}z_1,\dots,e^{i t_d}z_d)$.
A moment map is $\phi: \CC^d \to (\RR^d)^*$ defined by
\[
        \phi(z_1,\dots,z_d) =
        - \frac 12 (|z_1|^2,\dots,|z_d|^2) + (\lambda_1,\dots,\lambda_d) \ ,
\]
where the constant is chosen for later convenience.
The subtorus $N$ acts on $\CC^d$ in a hamiltonian way with moment map
$i^* \circ \phi: \CC^d \to \fn^*$.
Let $Z = (i^* \circ \phi)^{-1}(0)$.

In order to show that $Z$ (a closed set) is compact
it suffices (by the Heine-Borel theorem) to show that $Z$ is bounded.
Let $\Delta'$ be the image of $\Delta$ by $\pi^*$.
First we show that $\phi(Z) = \Delta'$.
A value $y\in(\RR^d)^*$ is in the image of $Z$ by $\phi$
if and only if
\[
   \mbox{(a) $y$ is in the image of $\phi$}
   \qquad \mbox{ and } \qquad
   \mbox{(b) $i^* y = 0$}
\]
if and only if (using the expression for $\phi$
and the third exact sequence)
\[
   \mbox{(a) $\ip{y}{e_i} \le \lambda_i$ for $i=1,\ldots,d$}
   \quad \mbox{ and } \quad
   \mbox{(b) $y = \pi^*(x)$ for some $x\in(\RR^n)^*$} \ .
\]
Suppose that $y = \pi^*(x)$.
Then
\[
   \ip{y}{e_i} \le \lambda_i, \forall i
   \iff \ip{x}{\pi(e_i)} \le \lambda_i, \forall i
   \iff \ip{x}{v_i} \le \lambda_i, \forall i
   \iff x \in \Delta \ .
\]
Thus, $y \in \phi(Z) \Leftrightarrow y\in\pi^*(\Delta) = \Delta'$.
Since $\Delta '$ is compact,
$\phi$ is proper and $\phi (Z) = \Delta'$,
we conclude that $Z$ must be bounded, and hence compact.

In order to show that $N$ acts freely on $Z$,
pick a vertex $p$ of $\Delta$, and let
$I= \{ i_1 , \ldots , i_n \}$ be the set
of indices for the $n$ facets meeting at $p$.
%so that $v_{i_1},\ldots,v_{i_n}$ are the
%normal primitive vectors pointing outwards
%with respect to those facets.
Pick $z \in Z$ such that $\phi (z) = \pi ^* (p)$.
Then $p$ is characterized by $n$ equations
$\langle p , v_i \rangle = \lambda_i$ where $i \in I$:
\[
\begin{array}{rcl}
   \langle p , v_i \rangle = \lambda_i
   & \iff & \langle p , \pi (e_i) \rangle = \lambda_i \\
   & \iff & \langle \pi ^* (p) ,e_i \rangle = \lambda_i \\
   & \iff & \langle \phi (z) ,e_i \rangle = \lambda_i \\
   & \iff & \mbox{$i$-th coordinate of
   $\phi (z)$ is equal to $\lambda_i$} \\
   & \iff & -\frac 12 |z_i|^2 + \lambda_i = \lambda_i \\
   & \iff & z_i = 0 \ .
\end{array}
\]
Hence, those $z$'s are points whose coordinates in the
set $I$ are zero, and whose other coordinates are nonzero.
Without loss of generality, we may assume that
$I = \{ 1,\ldots ,n\}$.
The stabilizer of $z$ is
\[
   (\TT^d)_z = \{ (t_1,\ldots , t_n,1,\ldots ,1) \in \TT^d \} \ .
\]
As the restriction $\pi : (\RR^d)_z \to \RR^n$
maps the vectors $e_1, \ldots , e_n$ to a $\ZZ$-basis
$v_1, \ldots , v_n$ of $\ZZ^n$ (respectively), at the
level of groups $\pi : (\TT^d)_z \to \TT^n$ must be bijective.
Since $N = \ker (\pi : \TT^d \to \TT^n)$,
we conclude that $N \cap (\TT^d)_z = \{ e \}$,
i.e., $N_z = \{ e \}$.
Hence all $N$-stabilizers at points mapping to
vertices are trivial.
But this was the worst case, since other stabilizers
$N_{z'}$ ($z' \in Z$) are contained in stabilizers
for points $z$ that map to vertices.
We conclude that $N$ acts freely on $Z$.

We now apply reduction.
Since $i^*$ is surjective, $0 \in \fn^*$ is a
regular value of $i^* \circ \phi$.
Hence, $Z$ is a compact submanifold of $\CC^d$ of
(real) dimension $2d - (d-n) = d+n$.
The orbit space $M_{\Delta} = Z/N$ is a compact manifold
of (real) dimension $\dim Z - \dim N = (d+n) - (d-n) = 2n$.
The point-orbit map $p: Z \to M_{\Delta}$
is a principal $N$-bundle over $M_{\Delta}$.
Consider the diagram
\[
\begin{array}{ccc}
        Z &\stackrel{j}{\hookrightarrow} &\CC^d \\
        {\scriptstyle{p}} \downarrow \phantom{\scriptstyle{p}} \\
        M_{\Delta}
\end{array}
\]
where $j: Z \hookrightarrow \CC^d$ is inclusion.
The Marsden-Weinstein-Meyer theorem (Theorem~\ref{thm:reduction})
guarantees the existence of a
symplectic form $\omega_{\Delta}$ on $M_{\Delta}$ satisfying
\[
        p^*\omega_{\Delta} = j^*\omega_0 \ .
\]
Since $Z$ is connected, the  symplectic
manifold $(M_{\Delta},\omega_{\Delta})$ is also connected.

It remains to show that $(M_\Delta,\omega_\Delta)$
is a hamiltonian $\TT^n$-space
with a moment map $\mu_\Delta$ having image $\mu_\Delta(M_\Delta) = \Delta$.
Let $z$ be such that $\phi (z) = \pi ^* (p)$
where $p$ is a vertex of $\Delta$.
Let $\sigma : \TT^n \to (\TT^d)_z$ be the
inverse for the earlier bijection
$\pi : (\TT^d)_z \to \TT^n$.
This is a {\em section}, i.e.,
a right inverse for $\pi$, in the sequence
\[
\begin{array}{ccccccccc}
   1 & \longrightarrow & N & \stackrel{i}{\longrightarrow}
   & \TT^d & \stackrel{\pi}{\longrightarrow}
   & \TT^n & \longrightarrow & 1 \ , \\
   & & & & & \stackrel{\sigma}{\longleftarrow}
\end{array}
\]
so it {\em splits}, i.e., becomes like a sequence for a product,
as we obtain an isomorphism
$(i , \sigma) : N \times \TT^n \stackrel{\simeq}{\longrightarrow} \TT^d$.
The action of the $\TT^n$ factor (or, more
rigorously, $\sigma (\TT^n) \subset \TT^d$)
descends to the quotient $M_\Delta = Z / N$.
Consider the diagram
\[
\begin{array}{rl}
   Z & \stackrel{j}{\hookrightarrow} \CC^d
   \stackrel{\phi}{\longrightarrow} (\RR^d)^*
   \simeq \eta^* \oplus (\RR^n)^*
   \stackrel{\sigma^*}{\longrightarrow} (\RR^n)^* \\
   p \downarrow \\
   M_\Delta
\end{array}
\]
where the last horizontal map is projection onto the second factor.
Since the composition of the horizontal maps
is constant along $N$-orbits, it descends to a map
\[
   \mu_\Delta : M_\Delta \longrightarrow (\RR^n)^*
\]
which satisfies $\mu_\Delta \circ p = \sigma^* \circ \phi \circ j$.
By reduction for product groups (Section~\ref{sec:generalizations}),
this is a moment map for the action of $\TT^n$
on $(M_\Delta, \omega_\Delta)$.
The image of $\mu_\Delta$ is
\[
   \mu_\Delta (M_\Delta) = (\mu_\Delta \circ p) (Z)
   = (\sigma ^* \circ \phi \circ j ) (Z)
   = (\sigma ^* \circ \pi ^* ) (\Delta) = \Delta \ ,
\]
because $\phi (Z) = \pi^* (\Delta)$ and $\pi \circ \sigma = \mbox{id}$.
We conclude that $(M_\Delta,\omega_\Delta, \TT^n , \mu_\Delta)$ is 
the required toric manifold corresponding to $\Delta$.
This construction via reduction also shows that
symplectic toric manifolds are in fact K\"ahler.

\begin{example}\index{example ! complex projective space}\index{complex !
projective space}\index{example ! Delzant construction}
Here are the details of Delzant's construction
for the case of a segment $\Delta = [0,a] \subset \RR^*\ (n = 1,d = 2)$.
Let $v(=1)$ be the standard basis vector in $\RR$.
Then $\Delta$ is described by
$\langle x,-v \rangle \le 0$ and $\langle x, v \rangle \le a$,
where $v_1 = -v$, $v_2 = v$, $\lambda_1 =0$ and $\lambda_2 =a$.
The projection $\RR^2 \stackrel{\pi}{\longrightarrow} \RR$,
$e_1 \mapsto -v$, $e_2 \mapsto v$,
has kernel equal to the span of $(e_1 + e_2)$,
so that $N$ is the diagonal subgroup of $\TT^2 = S^1 \times S^1$.
The exact sequences become
\[
\begin{array}{ccccccccc}
   1 & \longrightarrow & N & \stackrel{i}{\longrightarrow} & \TT^2
   & \stackrel{\pi}{\longrightarrow} & S^1 & \longrightarrow & 1 \\
   & & t & \longmapsto & (t,t) \\
   & & & & (t_1,t_2) & \longmapsto & t_1^{-1}t_2 \\
   \\
   0 & \longrightarrow & \fn & \stackrel{i}{\longrightarrow}
   & \RR^2 & \stackrel{\pi}{\longrightarrow} & \RR
   & \longrightarrow & 0 \\
   & & x & \longmapsto & (x,x) \\
   & & & & (x_1,x_2) & \longmapsto & x_2 - x_1 \\
   \\
   0 & \longrightarrow & \RR^* & \stackrel{\pi^*}{\longrightarrow}
   & (\RR^2)^* & \stackrel{i^*}{\longrightarrow} & \fn^*
   & \longrightarrow & 0 \\
   & & x & \longmapsto & (-x,x) \\
   & & & & (x_1,x_2) & \longmapsto & x_1 + x_2 \ .
\end{array}
\]
The action of the diagonal subgroup
$N = \{(e^{i t},e^{i t}) \in S^1 \times S^1\}$ on $\CC^2$ by
\[
   (e^{i t},e^{i t}) \cdot (z_1,z_2) = (e^{i t}z_1,e^{i t}z_2)
\]
has moment map $(i^* \circ\phi)(z_1,z_2) =
\textstyle{-\frac 12} (|z_1|^2 + |z_2|^2) + a$,
with zero-level set
\[
   (i^* \circ \phi)^{-1}(0) =
   \{(z_1,z_2) \in \CC^2 :|z_1|^2 + |z_2|^2 = 2a \} \ .
\]
Hence, the reduced space is a projective space,
$(i^* \circ \phi)^{-1}(0)/N = \CC \PP^1$.
\end{example}

%%%%%%%%%%%%%%%%%%%%%%%%%%%%%%%%%%%%%%%%%%%%%%%%%%%%%%%%%%%%%%%%%%%%%%%%%%%%%
%%%%%%%%%%%%%%%%%%%%%%%%%%%%%%%%%%%%%%%%%%%%%%%%%%%%%%%%%%%%%%%%%%%%%%%%%%%%%

\ssubsection{Duistermaat-Heckman Theorems}%\label{}
\index{Duistermaat-Heckman !
polynomial}\index{Heckman|see{Duistermaat-Heckman}}

Throughout this section, let $(M, \omega, G, \mu)$ be a hamiltonian $G$-space,
where $G$ is an $n$-torus\footnote{The discussion in this section
may be extended to hamiltonian actions of other compact Lie
groups, not necessarily tori; see~\cite[Exercises 2.1-2.10]{gu:moment}.}
and the moment map $\mu$ is proper.%case of noncompact groups?

If $G$ acts freely on $\mu^{-1} (0)$, it also acts freely
on nearby levels $\mu^{-1} (t)$, $t \in \fg^*$ and $t \approx 0$.
(Otherwise, assume only that $0$ is a regular value of $\mu$
and work with orbifolds.)
We study the variation of the reduced spaces by
relating\index{reduction ! local form}\index{local form}
\[
   (M_{\red} = \mu^{-1} (0) / G, \omega_{\red}) \qquad \mbox{ and } \qquad
   (M_t = \mu^{-1} (t) / G , \omega_t) \ .
\]
For simplicity, assume $G$ to be the circle $S^1$.
Let $Z = \mu^{-1} (0)$ and let $i : Z \hookrightarrow M$
be the inclusion map.
Fix a connection form $\alpha \in \Omega^1 (Z)$
for the principal bundle
\[
\begin{array}{cll}
        S^1 & \hookrightarrow & Z \\
        & & \downarrow \pi \\
        & & M_{\mathrm{red}}
\end{array}
\]
that is, $\cL _{X^\#} \alpha = 0$ and $\imath _{X^\#} \alpha = 1$,
where $X^\#$ is the infinitesimal generator for the $S^1$-action.
Construct a 2-form on the product
manifold $Z \times (-\varepsilon, \varepsilon)$ by the recipe
\[
   \sigma = \pi^* \omega_{\red} - d(x \alpha) \ ,
\]
where $x$ is a linear coordinate on the interval
$(-\varepsilon, \varepsilon) \subset \RR \simeq \fg^*$.
(By abuse of notation, we shorten the symbols for
forms on $Z \times (-\varepsilon, \varepsilon)$
that arise by pullback via projection onto each factor.)

\begin{lemma}
The 2-form $\sigma$ is symplectic for $\varepsilon$ small enough.
\end{lemma}

\vspace*{-2ex}

\begin{proof}
At points where $x=0$, the form
$\sigma |_{x=0} = \pi^* \omega_{\red} + \alpha \wedge dx$ satisfies
$\sigma |_{x=0} \left( X^{\#}, \frac{\partial}{\partial x} \right) = 1$,
so $\sigma$ is nondegenerate along $Z \times \{ 0 \}$.
Since nondegeneracy is an open condition,
we conclude that $\sigma$ is nondegenerate for $x$
in a sufficiently small neighborhood of $0$.
Closedness is clear.
%The form $\sigma$ is clearly closed.
\end{proof}

Notice that $\sigma$ is invariant with respect to
the $S^1$-action on the first factor of
$Z \times (-\varepsilon, \varepsilon)$.
This action is hamiltonian with moment map
$x : Z \times (-\varepsilon, \varepsilon)
\to (-\varepsilon, \varepsilon)$
given by projection onto the second factor
(since $\cL _{X^\#} \alpha =0$ and $\imath _{X^\#} \alpha =1$):
\[
   \imath _{X^\#} \sigma = - \imath _{X^\#} d (x \alpha)
   = - \cL _{X^\#} (x \alpha) + d {\imath _{X^\#} (x \alpha)} = dx \ .
\]

\begin{lemma}
\label{lem:model}
There exists an equivariant symplectomorphism between
a neighborhood of $Z$ in $M$ and a neighborhood of
$Z \times \{ 0 \}$ in $Z \times (-\varepsilon, \varepsilon)$,
intertwining the two moment maps, for $\varepsilon$ small enough.
\end{lemma}

\vspace*{-2ex}

\begin{proof}
The inclusion
$i_0 : Z \hookrightarrow Z \times (-\varepsilon, \varepsilon)$
as $Z \times \{ 0 \}$ and the natural inclusion
$i : Z \hookrightarrow M$ are $S^1$-equivariant
coisotropic embeddings.
Moreover, they satisfy $i_0^* \sigma = i^* \omega$
since both sides are equal to $\pi^* \omega_{\red}$,
and the moment maps coincide on $Z$ because $i_0^* x = 0 = i^* \mu$.
Replacing $\varepsilon$ by a smaller positive number if necessary,
the result follows from the equivariant version
of the coisotropic embedding theorem
(Theorem~\ref{thm:coisotropic}).\footnote{\textbf{Equivariant
Coisotropic Embedding Theorem:}\index{theorem !
equivariant coisotropic embedding}\index{equivariant !
coisotropic embedding} {\em Let $(M_0, \omega_0)$, $(M_1, \omega_1)$
be symplectic manifolds of dimension $2n$,
$G$ a compact Lie group acting on $(M_i, \omega_i)$, $i=0,1$,
in a hamiltonian way with moment maps $\mu_0$ and $\mu_1$,
respectively, $Z$ a manifold of dimension $k \geq n$
with a $G$-action, and $\iota_i : Z \hookrightarrow M_i$, $i=0,1$,
$G$-equivariant coisotropic embeddings.
Suppose that $\iota_0^* \omega_0 = \iota_1^* \omega_1$
and $\iota_0^* \mu_0 = \iota_1^* \mu_1$.
Then there exist $G$-invariant neighborhoods
$\cU_0$ and $\cU_1$ of $\iota_0 (Z)$ and $\iota_1 (Z)$
in $M_0$ and $M_1$, respectively,
and a $G$-equivariant symplectomorphism
$\varphi : \cU_0 \rightarrow \cU_1$ such that
$\varphi \circ \iota_0 = \iota_1$ and $\mu_0 = \varphi ^* \mu_1$.}}
\end{proof}

Therefore, in order to compare the reduced spaces
$M_t = \mu ^{-1} (t) / S^1$ for $t \approx 0$,
we can work in $Z \times (-\varepsilon, \varepsilon)$
and compare instead the reduced spaces $x ^{-1} (t) / S^1$.

\begin{proposition}
\label{prop:curvature}
The space $(M_t,\omega_t)$ is symplectomorphic to
$(M_{\red}, \omega_{\red} -t \beta)$
where $\beta$ is the curvature form of the connection $\alpha$.
\end{proposition}

\vspace*{-2ex}

\begin{proof}
By Lemma~\ref{lem:model}, $(M_t,\omega_t)$ is symplectomorphic
to the reduced space at level $t$ for the hamiltonian space
$(Z \times (-\varepsilon, \varepsilon) , \sigma , S^1, x)$.
Since $x ^{-1} (t) = Z \times \{ t \}$,
where $S^1$ acts on the first factor, all the manifolds
$x ^{-1} (t) / S^1$ are diffeomorphic to $Z / S^1 = M_{\red}$.
As for the symplectic forms,
let $\iota_t : Z \times \{ t \} \hookrightarrow
Z \times (-\varepsilon, \varepsilon)$ be the inclusion map.
The restriction of $\sigma$ to $Z \times \{ t \}$ is
\[
   \iota_t ^* \sigma = \pi^* \omega_{\red} -t d \alpha \ .
\]
By definition of curvature, $d \alpha = \pi ^* \beta$.
Hence, the reduced symplectic form on $x ^{-1} (t) / S^1$ is
$\omega_{\red} - t \beta$.
\end{proof}

In loose terms, Proposition~\ref{prop:curvature} says that
the reduced forms $\omega_t$ vary linearly in $t$,
for $t$ close enough to $0$.
However, the identification of $M_t$ with $M_{\red}$ as abstract
manifolds is not natural.
Nonetheless, any two such identifications are isotopic.
By the homotopy invariance of de Rham classes, we obtain:

\begin{theorem}\index{theorem !
Duistermaat-Heckman}\index{Duistermaat-Heckman ! theorem}\label{thm:dh2}
\textbf{(Duistermaat-Heckman~\cite{du-he:variation})} $\;$
Under the hypotheses and notation before,
the cohomology class of the reduced symplectic form
$[ \omega_t ]$ varies linearly in $t$.
More specifically,
if $c = [- \beta] \in H_{\mathrm{deRham}}^2 (M_{\red})$
is the first Chern class\footnote{Often
the Lie algebra of $S^1$ is identified with $2\pi i \RR$
under the exponential map $\exp : \fg \simeq 2\pi i \RR \rightarrow S^1$,
$\xi \mapsto e^{\xi}$.
Given a principal $S^1$-bundle, by this identification
the infinitesimal action maps the generator $2\pi i$ of $2\pi i \RR$
to the generating vector field $X^\#$.
A connection form $A$ is then an imaginary-valued 1-form
on the total space satisfying $\cL _{X^\#} A =0$
and $\imath_ {X^\#} A = 2\pi i$.
Its curvature form $B$ is an imaginary-valued 2-form
on the base satisfying $\pi^* B = dA$.
By the Chern-Weil isomorphism, the \textbf{first Chern class}\index{Chern !
first Chern class}\index{first Chern class} of the principal
$S^1$-bundle is $c= [\frac{i}{2\pi} B]$.

Here we identify the Lie algebra of $S^1$ with $\RR$
and implicitly use the exponential map
$\exp : \fg \simeq \RR \rightarrow S^1$, $t \mapsto e^{2\pi i t}$.
Hence, given a principal $S^1$-bundle,
the infinitesimal action maps the generator 1 of $\RR$ to $X^\#$,
and here a connection form $\alpha$ is an ordinary 1-form
on the total space satisfying $\cL _{X^\#} \alpha =0$
and $\imath_ {X^\#} \alpha = 1$.
The curvature form $\beta$ is an ordinary 2-form
on the base satisfying $\pi^* \beta = d\alpha$.
Consequently, we have $A=2\pi i \alpha$, $B=2\pi i \beta$
and the first Chern class is given by $c= [-\beta]$.}\index{Chern !
first Chern class}\index{first Chern class}
of the $S^1$-bundle $Z \rightarrow M_{\red}$,
we have
\[
   [ \omega_t ] = [ \omega_{\red} ] + t c \ .
\]
\end{theorem}

\vspace*{-1ex}

%%%%%%%%%%%%%%%%%%%%%%%%%%%%%%%%%%%%%%%%%%%%%%%%%%%%%%%%%%%%%%%%%%%%%%%%%%%%%

\begin{definition}
The \textbf{Duistermaat-Heckman measure}\index{measure !
Duistermaat-Heckman}\index{Duistermaat-Heckman ! measure},
$m_{DH}$, on $\fg^*$ is the push-forward of the Liouville
measure\footnote{On an arbitrary symplectic manifold $(M^{2n},\omega)$,
with symplectic volume\index{symplectic !
volume}\index{volume} $\frac{\omega^n}{n!}$,
the \textbf{Liouville measure}\index{Liouville !
measure}\index{measure ! Liouville}
(or \textbf{symplectic measure})\index{symplectic ! measure}\index{measure !
symplectic}
of a Borel subset\index{Borel subset} $\cU$ of $M$ is
\[
        m_\omega(\cU) = \int_\cU \frac{\omega^n}{n!} \ .
\]
The set $\cB$
of \textbf{Borel subsets}\index{Borel subset} is the {\em $\sigma$-ring}
generated by the set of compact subsets, i.e., if $A,B \in \cB$,
then $A \setminus B \in \cB$, and if $A_i \in \cB$, $i=1,2,\ldots$,
then $\cup_{i=1}^{\infty} A_i \in \cB$.}
by $\mu:M \rightarrow \fg^*$,
that is, for any Borel subset $\cU$ of $\fg^*$, we have
\[
        m_{DH}(\cU) =
%(\mu_* m_\omega)(\cU) =
        \int_{\mu^{-1}(\cU)} \frac{\omega^n}{n!} \ .
\]
\end{definition}

The integral with respect to the Duistermaat-Heckman measure
of a compactly-supported function $h \in C^\infty (\fg ^*)$ is
\[
   \int_{\fg ^*} h \ dm_{DH} :=
   \int_M (h \circ \mu) \frac{\omega^n}{n!} \ .
\]
On $\fg ^*$ regarded as a vector space, say $\RR ^n$,
there is also the Lebesgue (or euclidean) measure\index{euclidean !
measure}\index{Lebesgue ! measure}, $m_0$.
The relation between $m_{DH}$ and $m_0$ is governed
by the {\em Radon-Nikodym derivative}\index{Radon-Nikodym derivative},
denoted by $\frac{dm_{DH}}{dm_0}$, which is a {\em generalized
function} satisfying
\[
   \int_{\fg ^*} h \ dm_{DH} =
   \int_{\fg ^*} h \ \frac{dm_{DH}}{dm_0} \ dm_0 \ .
\]

\begin{theorem}\index{theorem !
Duistermaat-Heckman}\index{Duistermaat-Heckman ! theorem}\label{thm:dh}
\textbf{(Duistermaat-Heckman~\cite{du-he:variation})} $\;$
Under the hypotheses and notation before,
the Duistermaat-Heckman measure is
a piecewise polynomial multiple of Lebesgue
measure\index{measure ! Lebesgue}\index{Lebesgue ! measure}
on $\fg^* \simeq \RR^n$,
that is, the Radon-Nikodym derivative\index{Nikodym !
Radon-Nikodym derivative}\index{Radon-Nikodym derivative}
$f = \frac{dm_{DH}}{dm_0}$ is piecewise polynomial.
More specifically, for any Borel subset $\cU$ of $\fg^*$,
we have $m_{DH}(\cU) = \int_\cU f(x)\, dx$,
where $dx = dm_0$ is the Lebesgue volume\index{Lebesgue ! volume}
form on $\cU$ and $f: \fg^* \simeq \RR^n \to \RR$ is polynomial on
any region consisting of regular values of $\mu$.
\end{theorem}

This Radon-Nikodym derivative $f$ is called the
\textbf{Duistermaat-Heckman polynomial}.\index{Duistermaat-Heckman !
polynomial}
In the case of a toric manifold,
the Duistermaat-Heckman polynomial is a universal constant
equal to $(2\pi)^n$ when $\Delta$ is $n$-dimensional.
Thus the symplectic volume of $(M_\Delta,\omega_\Delta)$
is $(2\pi)^n$ times the euclidean volume of $\Delta$.

\begin{example}
For the standard spinning of a sphere,
$(S^2,\omega=d\theta\wedge dh, S^1, \mu=h)$,
the image of $\mu$ is the interval $[-1,1]$.
The Lebesgue measure of $[a,b] \subseteq [-1,1]$ is
$m_0([a,b]) = b-a$.
The Duistermaat-Heckman measure of $[a,b]$ is
\[
        m_{DH}([a,b]) = \int_{\{(\theta,h) \in S^2 \mid a\le h \le b \}}
        d\theta \ dh = 2\pi(b-a) \ ,
\]
i.e., $m_{DH} = 2\pi \ m_0$.
Consequently,
{\em the area of the spherical region between two parallel planes
depends only on the distance between the planes},
a result that was known to Archimedes\index{theorem !
Archimedes}\index{Archimedes} around 230 BC.
\end{example}

\vspace*{-1ex}

%%%%%%%%%%%%%%%%%%%%%%%%%%%%%%%%%%%%%%%%%%%%%%%%%%%%%%%%%%%%%%%%%%%%%%%%%%%%%

\begin{proof}
We sketch the proof of Theorem~\ref{thm:dh} for the case $G=S^1$.
The proof for the general case, which follows along
similar lines, can be found in, for instance, \cite{gu:moment},
besides the original articles.

Let $(M, \omega, S^1, \mu)$ be a hamiltonian $S^1$-space
of dimension $2n$ and let $(M_x, \omega_x)$ be its reduced
space at level $x$.
Proposition~\ref{prop:curvature} or Theorem~\ref{thm:dh2}
imply that, for $x$ in a sufficiently narrow neighborhood of $0$,
the symplectic volume\index{symplectic ! volume}\index{volume} of $M_x$,
\[
   \mathrm{vol} (M_x) =
   \int_{M_x} \frac{\omega _x ^{n-1}}{(n-1)!} =
   \int_{M_{\red}} \frac{(\omega _{\red} - x \beta) ^{n-1}}{(n-1)!} \ ,
\]
is a polynomial in $x$ of degree $n-1$.
This volume can be also expressed as
\[
   \mathrm{vol} (M_x) = \int_{Z}
   \frac{\pi^* (\omega _{\red} - x \beta) ^{n-1}}{(n-1)!}
   \wedge \alpha \ ,
\]
where $\alpha$ is a connection form for the
$S^1$-bundle $Z \rightarrow M_{\red}$ and $\beta$ is its
curvature form.
Now we go back to the computation of the
Duistermaat-Heckman measure.
For a Borel subset $\cU$ of $(-\varepsilon, \varepsilon)$,
the Duistermaat-Heckman measure is, by definition,
\[
   m_{DH}(\cU) = \int_{\mu^{-1}(\cU)} \frac{\omega^n}{n!} \ .
\]
Using the fact that $(\mu^{-1}(-\varepsilon, \varepsilon), \omega)$
is symplectomorphic to $(Z \times (-\varepsilon, \varepsilon), \sigma)$
and, moreover, they are isomorphic as hamiltonian $S^1$-spaces,
we obtain
\[
   m_{DH}(\cU) = \int_{Z \times \cU} \frac{\sigma^n}{n!} \ .
\]
Since $\sigma = \pi^* \omega_{\red} - d(x \alpha)$, its power is
$\sigma^n = n (\pi^* \omega_{\red} - x d \alpha)^{n-1}
\wedge \alpha \wedge dx$.
By the Fubini theorem\index{Fubini theorem}\index{theorem ! Fubini},
we then have
\[
   m_{DH}(\cU) = \int_{\cU} \left[ \int_{Z}
   \frac{\pi^* (\omega _{\red} - x \beta) ^{n-1}}{(n-1)!}
   \wedge \alpha \right] \wedge dx \ .
\]
Therefore, the Radon-Nikodym derivative of $m_{DH}$
with respect to the Lebesgue measure, $dx$, is
\[
   f(x) = \int_{Z}
   \frac{\pi^* (\omega _{\red} - x \beta) ^{n-1}}{(n-1)!}
   \wedge \alpha = \mathrm{vol} (M_x) \ .
\]

The previous discussion proves that, for $x \approx 0$,
$f(x)$ is a polynomial in $x$.
The same holds for a neighborhood of any other regular value of $\mu$,
because we may change the moment map $\mu$
by an arbitrary additive constant.
\end{proof}

%%%%%%%%%%%%%%%%%%%%%%%%%%%%%%%%%%%%%%%%%%%%%%%%%%%%%%%%%%%%%%%%%%%%%%%%%%%%%
%%%%%%%%%%%%%%%%%%%%%%%%%%%%%%%%%%%%%%%%%%%%%%%%%%%%%%%%%%%%%%%%%%%%%%%%%%%%%

Duistermaat and Heckman~\cite{du-he:variation} also applied these results
when $M$ is compact to provide a formula for the oscillatory integral
$\int_M e^{i\mu^X} \frac{\omega^n}{n!}$
for $X \in \fg$ as a sum of contributions of the fixed points
of the action of the one-parameter subgroup generated by $X$.
They hence showed that the {\em stationary phase
approximation}\footnote{The \textbf{stationary phase lemma}
gives the asymptotic behavior (for large $N$) of integrals
$\left( \frac{N}{2\pi} \right)^n
\int_M f e^{ig} vol$, where $f$ and $g$ are real functions
and $vol$ is a volume form on a $2n$-dimensional manifold $M$.}
is exact in the case of the moment map.
When $G$ is a maximal torus of a compact connected simple Lie group
acting on a coadjoint orbit, the Duistermaat-Heckman formula
reduces to the Harish-Chandra formula.
%reference?
It was observed by Berline and Vergne~\cite{be-ve:zeros}
and by Atiyah and Bott~\cite{at-bo:moment}
that the Duistermaat-Heckman formula can be derived
by {\em localization in equivariant cohomology}.
This is an instance of \textbf{abelian localization},
i.e., a formula for an integral (in equivariant cohomology)
in terms of data at the fixed points of the action,
and typically is used for the case of abelian groups
(or of maximal tori).
Later \textbf{non-abelian localization} formulas were found,
where integrals (in equivariant cohomology)
are expressed in terms of data at the zeros of the moment map,
normally used for the case of non-abelian groups.
Both localizations gave rise to computations of the cohomology
ring structure of reduced spaces~\cite{ki:quotients}.

%%%%%%%%%%%%%%%%%%%%%%%%%%%%%%%%%%%%%%%%%%%%%%%%%%%%%%%%%%%%%%%%%%%%%%%%%%%%%
%%%%%%%%%%%%%%%%%%%%%%%%%%%%%%%%%%%%%%%%%%%%%%%%%%%%%%%%%%%%%%%%%%%%%%%%%%%%%
% REFERENCES
%%%%%%%%%%%%%%%%%%%%%%%%%%%%%%%%%%%%%%%%%%%%%%%%%%%%%%%%%%%%%%%%%%%%%%%%%%%%%
%%%%%%%%%%%%%%%%%%%%%%%%%%%%%%%%%%%%%%%%%%%%%%%%%%%%%%%%%%%%%%%%%%%%%%%%%%%%%

%%%%%%%%%%%%%%%%%%%%%%%%%%%%%%%%%%%%%%%%%%%%%%%%%%%%%%%%%%%%%%%%%%%%%%%%%%%%%
%%%%%%%%%%%%%%%%%%%%%%%%%%%%%%%%%%%%%%%%%%%%%%%%%%%%%%%%%%%%%%%%%%%%%%%%%%%%%

}

\end{document}